\documentclass[11pt, twoside, a4paper]{amsart}
\usepackage{macro}

\title{Lagrangian Skeleta of Very Affine Complete intersections}
\author[Danil Koževnikov]{Danil Koževnikov}
\address{University of Edinburgh, Edinburgh, UK}
\email{d.kozevnikov@sms.ed.ac.uk}
\usepackage{xpatch}
\makeatletter
\xpatchcmd{\@tocline}
  {\hfil}
  {\dotfill}
  {}{}
\renewcommand{\l@subsection}{\@tocline{2}{0pt}{1.75pc}{5pc}{}}
\renewcommand{\l@subsubsection}{\@tocline{3}{0pt}{3pc}{5pc}{}}
\makeatother
\begin{document}
\begin{abstract}
    Let $Z^\circ$ be a complete intersection inside $(\Cs)^n$ that compactifies to a smooth Calabi--Yau subvariety $Z$ of a Fano toric variety. We compute the Lagrangian skeleton of $Z^\circ$ and describe its decomposition into standard pieces that are mirror to toric varieties. 
    
    This set-up was first considered by Batyrev and Borisov, who used combinatorial techniques to construct a mirror pair $(Z,\check{Z})$ of Calabi--Yau complete intersections in Fano toric varieties. We apply our main result to establish homological mirror symmetry for Batyrev–Borisov pairs in the large-volume limit. We also prove that the equivalence is compatible with toric HMS, along with further functoriality properties with respect to certain natural inclusions of very affine complete intersections. 
\end{abstract}
\maketitle
\tableofcontents

\section{Introduction} \label{section:introduction}
A \emph{Liouville domain} is a compact manifold with boundary $X$ equipped with an exact symplectic form $\omega=d\lambda$ such that its \emph{Liouville vector field} $Z$, defined implicitly through the equation $\iota_Z\omega=\lambda$, is outward pointing along $\partial X$. For a Liouville domain $(X,\lambda)$, define its \emph{skeleton} as the set of points that do not escape into the boundary under the Liouville flow, which we denote $\skel(X,\lambda)$. Equivalently, one can describe the skeleton as $\skel(X,\lambda)= X \backslash \bigcup_{t \geq 0} \Phi^{-t}_Z(\partial X)$. When $(X,\lambda)$ is \emph{Weinstein}, it is known that the skeleton is stratified by isotropic submanifolds (the reader is referred to \cite{CE} for a more detailed exposition). Skeleta are of interest to symplectic topology due to the fact that the negative Liouville flow can be used to contract $X$ onto an arbitrarily small neighbourhood of the skeleton while merely rescaling the symplectic form, so many properties of $X$ can be read off purely from $\skel(X,\lambda)$. 

In particular, it was conjectured by Kontsevich in \cite{kontsevich09} that the \emph{wrapped Fukaya category} $\mW(X,\lambda)$ can be computed in terms of global sections of a certain cosheaf of categories on $\skel(X,\lambda)$. A definition of the sheaves that should correspond to the wrapped Fukaya category appears in the work of Nadler (see, e.g., \cite{nadler16}), and the seminal work of Ganatra, Pardon and Shende on \emph{partially wrapped Fukaya categories} establishes a general connection between microlocal sheaves on skeleta and Fukaya categories in \cite{gps3}. The biggest advantage of such an approach is that the cosheaf categories have local-to-global properties and, in favourable circumstances, one can pick the Liouville form $\lambda$ so that the skeleton is particularly simple, making it possible to compute Fukaya categories from purely combinatorial and topological data. For example, homological mirror symmetry for punctured Riemann surfaces was studied through this lens in works such as \cite{stz} and \cite{pascaleff-sibilla} as a proof of concept. 

A particularly nice class of Weinstein manifolds are smooth \emph{very affine varieties}, which is just a term for subvarieties of $(\Cs)^n$ from tropical geometry that we shall adopt for brevity. Thanks to the presence of a special Lagrangian torus fibration $\Log \colon (\Cs)^n \rightarrow \R^n$, these fall within the reach of tropical geometry, which means that questions about their topology can often be reduced to studying \emph{tropical subvarieties} of $\R^n$ through the complex version of Viro patchworking developed by Mikhalkin in \cite{Mikh}. Very affine varieties also play an important role in homological mirror symmetry, since they can be obtained from subvarieties of toric varieties by removing the toric boundary. Numerous constructions of mirror pairs operate within the framework of Calabi--Yau subvarieties of toric varieties, such as Greene-Plesser mirror pairs in \cite{gp90} (complete intersections in weighted projective spaces) and Batyrev mirror pairs in \cite{Batyrev94} (hypersurfaces in Fano toric varieties). 

The first step towards understanding Lagrangian skeleta of very affine varieties was undertaken in \cite{rstz}, where the authors identified topological skeleta of a wide class of very affine hypersurfaces, leaving open the question of how these relate to the symplectic topology of such hypersurfaces. The first computation of an actual Lagrangian skeleton was performed by Nadler in \cite{nadler16} for the \emph{generalised pairs of pants $\mP_n\coloneqq\{z \in (\Cs)^{n+1} \colon 1+z_1+\dots+z_{n+1}=0\}$}. Later, that computation was extended by Gammage and Shende to a broader class of very affine hypersurfaces in \cite{GS22} (under a hypothesis that the associated fan is \emph{perfectly centred}; this assumption was removed by Zhou in \cite{Zhou}), using tropicalisation and decompositions into pairs of pants (developed in the context of symplectic geometry in \cite{Mikh} and \cite{Abou}) to bootstrap Nadler's result. The skeleta of these very affine hypersurfaces admit covers by \emph{FLTZ Lagrangians} that were introduced by Fang, Liu, Treumann and Zaslow in \cite{fltz} as mirrors to toric varieties. In \cite{GS22}, this structure is exploited to prove an HMS result, where the B-side mirror is the toric boundary of a certain toric stack that is a resolution of a Fano toric variety.

In general, computing Fukaya categories of very affine varieties can be used as the first step in proving HMS for Calabi--Yau subvarieties of Fano toric varieties following the strategy envisioned by Seidel (first implemented in the case of quartic K3 surfaces in \cite{seidel-quartic}). One begins by studying the \emph{large volume limit} on the A-side, an exact symplectic manifold obtained by removing an appropriate divisor. The B-side mirror will usually be a singular algebraic variety, the \emph{large complex structure limit}, obtained by degenerating a family of smooth varieties. Seidel's insight is that once we have an HMS statement at the limit point, it is possible to prove it over its neighbourhood by deforming both categories and matching the deformations. This leads to the notion of the \emph{relative Fukaya category}, which was introduced in the aforementioned paper \cite{seidel-quartic}. In the work of Sheridan, the framework of relative Fukaya categories was developed further and successfully applied in proofs of homological symmetry in cases such as Calabi--Yau hypersurfaces in projective spaces (\cite{nick-hyp}), Greene-Plesser mirrors (\cite{greene}) and, most recently, Batyrev mirror pairs (\cite{ghhps}). Notably, the last mentioned work indeed uses the main result from \cite{GS22} as its starting point. 

In this paper, we will mostly be interested in the construction of mirror pairs due to Batyrev and Borisov that generalises both the aforementioned ones. In their work \cite{BB96}, they provide a combinatorial construction of codimension $r$ Calabi--Yau complete intersections (CYCI's) inside Fano toric varieties $X_{\Delta}$ through \emph{nef partitions}, which are particularly nice decompositions of the reflexive polytope $\Delta$ as a Minkowski sum $\Delta_1+\dots+\Delta_r$. They also explain how to construct a \emph{dual nef partition} of another reflexive polytope $\nabla=\nabla_1+\dots+\nabla_r$, and check that the corresponding complete intersections satisfy the equalities of Hodge numbers predicted by mirror symmetry. The pairs $(Z,\check{Z})$ obtained this way will be called \emph{Batyrev--Borisov mirror pairs}\footnote{Not all CYCI's inside a Gorenstein Fano toric variety can be obtained this way, see Remark \ref{remark:gorenstein-cones}.}. For a \emph{Batyrev--Borisov complete intersection (BBCI)} $Z \subset X_\Delta$, we call the associated very affine variety $Z^\circ \subset (\Cs)^n \cong X_\Delta \backslash \partial X_\Delta$ an \emph{open Batyrev--Borisov complete intersection}. 

\subsection{Main results}
Our goal is providing an explicit description of Lagrangian skeleta of open Batyrev--Borisov complete intersections, in the same spirit as \cite{GS22} and \cite{Zhou}. Somewhat informally, our main results can be summarised as follows: the skeleton of open BBCI's still admits a description in terms of mirrors of lower dimensional toric stacks. 

To state the result precisely, we need to introduce some notation (see Sections \ref{section:nef-partitions} and \ref{section:tropical} for more details): let $M$, $N$ be a pair of dual $n$-dimensional lattices with $M_G\coloneqq M \otimes G$ and $N_G=\hom(M,G)$ for any abelian group $G$. We start from the data of an irreducible nef partition $\nabla=\nabla_1+\dots+\nabla_r$ of a reflexive polytope $\nabla \subset N$ (dual to a nef partition $\Delta=\Delta_1+\dots+\Delta_r$) and functions $h \in \R^{\Delta^\vee \cap N}$, $\check{h} \in \R^{\nabla^\vee \cap M}$ that induce refined star-shaped triangulations $\mT$ and $\check{\mT}$ of $\Delta^\vee$ and $\nabla^\vee$, respectively, and taking cones over the faces in $\mT$, $\check{\mT}$ yields a pair of fans $\Sigma\subset N_\R$ and $\check{\Sigma} \subset M_\R$. By treating $\nabla_j$'s as Newton polytopes of the defining equations, these data determine a one-parameter family of open BBCI's $Z_\beta \subset M_{\Cs}$ that admit natural toric compactifications $\overline{Z}_\beta \subset \mX_{\check{\Sigma}}$, where $\mX_{\check{\Sigma}}$ is the smooth toric Deligne--Mumford stack associated to $\check{\Sigma}$ with primitive ray generators. Our main goal will be understanding $\Lambda\coloneqq\skel(Z_\beta,\lambda)$ for large $\beta>0$ and an appropriately chosen Liouville form $\lambda$. 

In order to provide a concise description of $\Lambda$, recall that to a fan $\Sigma \subset N_\R$, one can associate the \emph{FLTZ skeleton} $\eL_\Sigma\coloneqq\bigcup_{\sigma \in \Sigma} \sigma^\perp \times \sigma \subset T^*M_{S^1}$. This is a singular Lagrangian whose open strata are strongly exact. Its usefulness for us is twofold: the skeleton will live inside $\eL_\Sigma$, while certain lower dimensional FLTZ skeleta will provide natural local models for the singularities of $\Lambda$ as a singular Lagrangian. With that in mind, we are ready to state our main result:

\begin{theorem}[Main theorem]\label{theorem:main-thm}
    Given the choice of data outlined above, there exists a Liouville form $\lambda$ on $M_{\Cs}$ and a compact smooth domain $K \subset M_{\Cs}$ such that for all $\beta>0$ large enough
    \begin{enumerate}[(i)]
        \item\label{main-thm-1} The open Batyrev--Borisov complete intersection $Z\coloneqq Z_\beta$ is a codimension $2r$ exact symplectic submanifold of $M_{\Cs}$ when equipped with the restriction of $\lambda$.
        \item\label{main-thm-2} The intersection $A\coloneqq Z \cap K$ is a Liouville domain that is isomorphic to the Liouville domain associated to the affine variety $Z$ (in the sense of \cite{based-seidel}).
        \item\label{main-thm-3} The skeleton $\Lambda=\skel(A,\lambda)$ is a singular Lagrangian ambient diffeomorphic to $(M_{S^1} \times i(S^{n-r})) \cap \eL_\Sigma \subset T^*M_{S^1}$ for a certain embedding $i \colon S^{n-r} \hookrightarrow N_\R \backslash\{0\}$ of the $(n-r)$-sphere $S^{n-r}$.
        \item\label{main-thm-4} There exists an open cover of $\Lambda$ by open sets $\Lambda(\sigma)$ anti-indexed by the poset $\Sigma_\trans$ of transversal cones in $\Sigma$, where $\Lambda(\sigma)$ is ambient diffeomorphic to $\eL_{\Sigma/\sigma}\times C_\sigma$ for some contractible set $C_\sigma$. Moreover, the identifications can be chosen so that the inclusions $\Lambda(\tau) \hookrightarrow \Lambda(\sigma)$ for $\sigma \subset \tau$ correspond to the standard\footnote{See the discussion at the start of Section \ref{section:std-covers-skeleton} or in Section 4.3 of \cite{GS22} for what precisely is meant by \enquote{standard}.} inclusions of FLTZ skeleta. 
    \end{enumerate}
\end{theorem}

\begin{remark}\label{remark:generality}
    Note that the previous results for hypersurfaces are slightly more general than the set-up of Theorem \ref{theorem:main-thm} with $r=1$. Namely, they do not require the polytope $\Delta^\vee$ to be reflexive, instead, it can be any lattice polytope containing the origin with a choice of a refined star-shaped triangulation. It is therefore a natural question whether our results could be generalised to some broader class of complete intersections (e.g., allowing $\nabla_j$ to be any lattice polytopes containing the origin rather than summands of a nef partition). The main bottleneck is that some of our methods heavily rely on the results about combinatorics of tropical complete intersections associated to nef partitions from \cite{HZ} and Section \ref{section:tropical-bbci}. Therefore, we mostly restrict our attention to the more specific case of open Batyrev--Borisov complete intersections, which also arise naturally in the setting of HMS outlined earlier in the introduction. For a number of intermediate results, the discussion of the general case is also included (e.g., Sections \ref{section:tropical-ci} and \ref{section:tailor-ci}) to present them in a more natural setting. 
\end{remark}

As advertised, this result is a natural extension of the description of skeleta of hypersurfaces provided in \cite{GS22} and \cite{Zhou} to open Batyrev--Borisov complete intersections. That means that we also recover a number of interesting corollaries.

Firstly, we exploit the recursive structure of the skeleton to prove an HMS equivalence: for the toric stack $\mX_\Sigma$, define its \emph{transversal boundary} as $\check{Z}\coloneqq\bigcup_{\sigma \in \Sigma_\trans} \overline{O(\sigma)}$, the subset of its toric boundary $\partial \mX_\Sigma$ consisting of closed orbits indexed by \emph{transversal cones} of $\Sigma$ (see Section \ref{section:tropical-bbci}) that shall also sometimes be denoted as $\partial_\trans \mX_{\Sigma}$. This will be a closed singular substack of $\mX_\Sigma$ of codimension $r$, which we prove is a B-side mirror\footnote{The opposite in the statement of Theorem \ref{theorem:hms-for-bbci} can be removed by negating the Liouville form of $Z$, we leave the statement in this form to stick to the natural symplectic form induced from $(\Cs)^n$ and the definition of the FLTZ skeleton given above.} for the open Batyrev--Borisov complete intersection $Z$.

\begin{theorem}\label{theorem:hms-for-bbci}
    There exists a quasi-equivalence of $\Z$-graded, $\C$-linear $A_\infty$-categories $\perf \mW(Z)^{op} \cong \coh(\check{Z})$. 
\end{theorem}

See Section \ref{section:hms} for more details on what exactly should be understood by the categories on both sides. We also construct a thickening of the Liouville domain $A$ from Theorem \ref{theorem:main-thm} to get a Liouville hypersurface embedding $A \times \C^{r-1} \hookrightarrow S^*M_{S^1}$ into the cosphere bundle of a real $n$-torus that interacts well with the skeleta (Theorem \ref{theorem:emb}), so we can apply \cite[Theorem 1.4]{gps3} to relate the equivalence from Theorem \ref{theorem:hms-for-bbci} with toric homological mirror symmetry. 

\begin{theorem}\label{theorem:comm-diag}
    There exists a commutative diagram
    \[\begin{tikzcd}
    	{\perf \mW(Z)^{op}} & {\perf \mW(T^*M_{S^1},\partial\mathbb{L}_\Sigma)^{op}} \\
        \coh(\check{Z}) & \coh(\mX_{\Sigma})
    	\arrow[from=1-1, to=1-2]
    	\arrow[sloped, "\sim", from=1-1, to=2-1]
    	\arrow[sloped, "\sim", from=1-2, to=2-2]
    	\arrow[from=2-1, to=2-2]
    \end{tikzcd}\]
    The vertical equivalences are Theorem \ref{theorem:hms-for-bbci} and \cite[Theorem 1.2]{coco}, while the bottom arrow is derived pushforward under the closed embedding $\check{Z} \hookrightarrow \mX_{\Sigma}$. 
\end{theorem}

It is noteworthy that since our proof of Theorems \ref{theorem:hms-for-bbci} and \ref{theorem:comm-diag} passes through microlocal sheaves via \cite{gps3}, there is no obvious description of the top functor in terms of Floer theoretically defined Fukaya categories; see Section \ref{section:functoriality} and, in particular, Remark \ref{remark:arrow-interpretation} for a more detailed discussion. 

\subsection{Relations to existing work}

Finally, we explain the relationship between the thesis and previous work applicable to the setting of open Batyrev--Borisov complete intersections. 

Most notably, the mirror construction of Abouzaid--Auroux--Katzarkov from \cite{aak} gives a Landau--Ginzburg model $(Y_{AAK},W_{AAK})$ that is SYZ-mirror to a general very affine complete intersection $Z$, which is different from the one that we consider. In the special case when $Z$ is an open BBCI, our mirror turns out to be derived equivalent to theirs by Knörrer periodicity (see Remark \ref{remark:aak-mirror}). The other direction of HMS for this construction, i.e. identifying $\coh(Z)$ with the Fukaya--Seidel category of $(Y_{AAK},W_{AAK})$, has been explored in \cite{aa} by introducing the \emph{fibrewise wrapped Fukaya category}. The direction of HMS that is of interest to us has been discussed in recent work \cite{msz}. Their methods are, however, quite different from ours, since the aforementioned paper heavily relies on a gluing construction of Weinstein manifolds to avoid computing the global skeleton. Another distinction is that we restrict our attention to a particular class of complete intersections (analogously to how \cite{GS22} only treat hypersurfaces whose tropicalisation has a single compact complementary region), but obtain more results in that specific setting. 

Another related idea is generalising the strategy of gluing together pieces that look like mirrors to smooth toric stacks, which underlies the construction of \emph{fanifolds} in \cite{GS23}. Even though the skeleta we determine here admit stratifications into the same pieces, it is not clear how our setting relates to the framework of fanifolds. The key difference is that fanifolds are constructed through an inductive procedure of generalised Weinstein handle attachments, whereas the spaces we consider naturally arise from algebraic varieties, so it is unclear how to compare the two without a neighbourhood theorem for singular Lagrangian skeleta. 

\subsection{Summary of the paper}\label{section:summary} 

We start out by giving an overview of the relevant facts from toric geometry and introducing \emph{nef partitions} (a notion originally defined in \cite{Borisov93}) to set the stage for the main act by introducing the Batyrev--Borisov construction in more detail. Despite the fact that we are exclusively concerned with the very affine varieties rather than their toric compactifications, the combinatorics of the set-up still plays a key role in our approach. 

With preliminaries out of the way, we carry on to extending the results about tropicalisation of very affine hypersurfaces from \cite{Mikh} to the setting of complete intersections in Section \ref{section:tropical}. We explain the connection between the general combinatorial model (appearing, for example, in \cite{MS}) and the computationally convenient framework of \cite{HZ} that is specialised to the Batyrev--Borisov case, and use these to describe the combinatorics of very affine complete intersections. We also review the theory of toric compactifications inside smooth toric stacks and extend some of our results (such as the tropical smoothness criterion for complete intersections) to that setting. 

Then, in Section \ref{section:potentials}, we discuss a generalisation of \emph{adapted potentials} introduced in \cite{Zhou} to get rid of the restrictive technical condition of \emph{perfectly centredness} present in \cite{GS22}. The key difference between our setting and the case of hypersurfaces is that simple restrictions on positions of minima and homogeneity of the potential no longer seem sufficient to pin down the topology of the Lagrangian skeleton. To fix this issue, we introduce a more restrictive class of \emph{strongly adapted potentials}, which are assumed to be isotopic to certain standard embeddings. 

In Section \ref{section:tailoring}, we expand the framework of \emph{tailoring} to the setting of complete intersections, which allows us to localise our calculations in the subsequent sections. The idea was originally introduced in \cite{Mikh}, but we follow the more quantitative approach from \cite{Abou} via explicit constructions of cut-off functions. As an application of the earlier discussion of toric compactifications, we also explain how our ideas can be utilised for tailoring complete intersections in general smooth toric stacks. This concludes the set-up for the problem, since we end the section by explaining how to associate a Liouville domain amenable to localised calculations to an algebraic complete intersection in $(\Cs)^n$. 

The key calculation of the paper, generalising the ones from \cite{GS22} and \cite{Zhou}, appears in Section \ref{section:skeleta}. We start by presenting a high-level outline of the argument in Section \ref{section:skeleta-outline}, emphasising how different features explored in the previous section come into play, and then fill in the computational details in Section \ref{section:skeleta-calculations}.

Section \ref{section:combinatorics} is devoted to understanding the topology of the skeleton. As in the hypersurface case, it admits a natural description as a codimension $r$ sphere embedded in $\R^n$ with tori attached to it based on its intersection pattern with a simplicial fan $\Sigma$. However, studying the embeddings of $S^{n-r}$ into the $n$-dimensional space relative to $\Sigma$ becomes substantially more delicate for $r>1$. These difficulties are addressed through the framework of adapted potentials developed in the earlier sections, along with some technical results on smoothing topological submanifolds of $\R^n$ relative to a stratification (Appendix \ref{section:appendix-smoothing}) following the methods of \cite{Whitehead1961}. 

In Section \ref{section:stabilisation}, we describe a procedure for thickening the codimension $r$ complete intersection to a Liouville hypersurface in the cosphere bundle of an $n$-torus, which produces a Liouville sector associated to an open Batyrev--Borisov complete intersection, placing us within the framework of \cite{gps1} and \cite{gps2}. This construction serves as the first step in proving Theorem \ref{theorem:comm-diag} and relating open BBCI's with toric mirror symmetry. The intermediate results about positive codimension embeddings of Liouville domains might be of independent interest. 

Finally, Section \ref{section:hms} is devoted to the applications of our results from previous sections to homological mirror symmetry. We discuss the B-side mirror in more detail and prove Theorem \ref{theorem:hms-for-bbci} by adapting the techniques from \cite{GS23} to the present setting. As an intermediate step to proving Theorem \ref{theorem:comm-diag}, we also prove that natural inclusions of open BBCI's (obtained through the simple observation that by adding up some of the defining equations, one gets a lower codimension complete intersection) correspond to inclusions of the mirror substacks inside $\partial \mX_\Sigma$.

\ack I would like to thank my advisor, Nick Sheridan, for suggesting this project, for many fruitful discussions and for helpful comments on the first draft. I would also like to thank Sheel Ganatra for explaining his work, and to Andrew Hanlon, Jeff Hicks, Kai Hugtenburg and Sukjoo Lee for useful conversations, and I am grateful to Arend Bayer and Georgios Dimitroglou Rizell for their comments on an earlier version of the manuscript. This work was supported by the ERC Starting Grant 850713 – HMS. 
\section{Batyrev--Borisov mirror construction} \label{section:nef-partitions}
In this section, we go over the combinatorics of nef partitions and their relevance to geometry through the Batyrev--Borisov mirror construction. We also review some facts from toric geometry that will be important in the rest of this paper. 

\subsection{Toric background}

We shall always take $M \cong \Z^n$ to be a lattice of rank $n$ with a dual lattice $N\coloneqq \hom_{\Z}(M,\Z)$. For an abelian group $G$, we then denote $M_G\coloneqq M \otimes_\Z G$ and $N_G\coloneqq \hom_{\Z}(M,G)$. For example, $M_\R$ and $N_\R$ are the pair of dual vector spaces in which the lattices can be embedded, $M_{S^1}\cong M_\R/M$ is an $n$-dimensional real torus with a Pontryagin dual $N_{S^1}$, and $M_{\Cs}$ is an algebraic $n$-dimensional torus with a dual group $N_{\Cs}$. Unless otherwise specified, $\langle\cdot,\cdot\rangle \colon N_{\R}\times M_{\R} \rightarrow \R$ stands for the dual pairing. 

By a \emph{polyhedron}, we shall mean an intersection of finitely many closed half-spaces in a vector space (in particular, all polyhedra are convex). The term \emph{polytope} will be synonymous with \enquote{compact polyhedron}, which can also be presented as the convex hull of a finite set of points. A polytope in $M_\R$ is called a \emph{lattice polytope} if it can be written as a convex hull of a finite set of points in the lattice $M$. The \emph{relative interior} of a polytope $P$ is the interior of $P$ inside the smallest affine subspace containing $P$. We say that a polytope $F \subset \partial P$ is a \emph{face of $P$} if there exists a hyperplane $H$ such that $F=P \cap H$. In particular, the $0$-dimensional faces of $P$ are called \emph{vertices}, $1$-dimensional faces are called \emph{edges}. 

For a full-dimensional polyhedron $P \subset M_{\R}$ and a face $F \subset \partial P$, we define\footnote{By the general discussion of embedded manifolds with corners from Appendix \ref{section:appendix-smoothing}, this should actually be called the \emph{conormal cone}, but we stick to the standard naming.} the \emph{normal cone to $P$ at $F$} as
\begin{equation*}
    \nc_P(F) \coloneqq \{ \lambda \in N_\R \colon F \textnormal{ is the set of maximisers of } \lambda \textnormal{ over } P\},
\end{equation*}
with the extra convention that $\nc_P(P)\coloneqq 0$. It is clear that $\nc_P(F)$ is a polyhedral cone and that for an inclusion $F \hookrightarrow F'$, we get an inclusion $\nc_P(F') \hookrightarrow \nc_P(F)$. Therefore, the normal cones assemble into the \emph{normal fan of $P$}, denoted as $N(P)$. This fan is going to be complete if and only if $P$ is a polytope; it will be rational if and only if the polytope $P$ is rational. If that is the case, we define the \emph{toric variety associated to $P$} as the toric variety coming from its normal fan. 

\begin{remark}\label{remark:weird-convention}
    Our convention for normal cones is opposite to the usual convention in toric geometry, where one considers the inward normals instead. The two choices clearly differ by a single sign; the main reason why we choose to diverge from the standard convention is to make everything compatible with the terminology from convex geometry (so that the terms \emph{convex function}, \emph{Legendre transform} and \emph{support function} maintain their usual meanings). The advantage becomes apparent in Section \ref{section:potentials}, where this choice allows us to phrase everything in terms of sublevel sets of convex functions, while avoiding ample minus signs. 
\end{remark}

We now explain our conventions for toric stacks following the foundational paper \cite{BCS}: given a stacky fan $(\Sigma,\beta)$, one can define the associated toric stack $\mX_{\Sigma,\beta}$, which will be a smooth Deligne--Mumford stack with coarse moduli space $X_\Sigma$ (Propositions 3.2 and 3.7 from loc. cit.). The standard facts about line bundles on toric varieties extend to that setting, provided these are interpreted as \emph{equivariant} line bundles and sections, see \cite[Section 2.1]{HHL} for a more detailed overview. Any simplicial fan $\Sigma$ has a canonical stacky fan structure by picking the primitive ray generator of each ray, so when we refer to a toric stack $\mX_\Sigma$ without the data of $\beta$, it is understood that we mean this particular `stacky enhancement'. Following \cite{torbifolds}, any toric stack over $\C$ also carries a natural Kähler form coming from a potential that is a generalisation of Guillemin's potential for toric manifolds, giving it the structure of a \emph{symplectic toric orbifold}.  

Finally, we introduce the object defined by Fang--Liu--Treumann--Zaslow in \cite{fltz} as a candidate mirror to the toric variety $X_\Sigma$, which will be of utmost importance in the rest of this work.

\begin{definition}\label{definition:fltz-lagrangian}
The \emph{FLTZ Lagrangian} of a fan of $\Sigma\subset N_\R$ (also referred to as the \emph{FLTZ skeleton}) is given by
    $$
    \eL_\Sigma\coloneqq\bigcup_{\sigma \in \Sigma} \sigma^\perp \times \sigma \subset M_{S^1} \times N_\R \cong T^*M_{S^1}.
    $$
\end{definition}

Note that we are using the identification $M_{S^1}\cong M_\R/M$ to get a sub-torus of $M_{S^1}$ from the (rational) subspace $\sigma^\perp\subset M_\R$, along with the trivialisation of the cotangent bundle of $M_{S^1}$. Each piece $\sigma^\perp \times \relint(\sigma)$ is a cylindrical strongly exact Lagrangian, and its closure $\sigma^\perp \times \sigma$ is a Lagrangian with corners (see Section \ref{section:skeleta-outline} for more details). Therefore, $\eL_\Sigma$ is a singular cylindrical Lagrangian with singular Legendrian boundary $\partial \eL_\Sigma \subset S^*M_{S^1}$. 

\begin{remark}\label{remark:stacky-fltz}
    It was shown in \cite{coco} that, unlike in the torus-equivariant setting originally considered in \cite{fltz}, for non-smooth simplicial fans $\Sigma$, the associated Lagrangian $\eL_\Sigma$ gives a natural mirror for the smooth toric stack $\mX_\Sigma$ rather than the singular variety $X_\Sigma$, which is the main reason why our results are naturally framed in the setting of toric stacks rather than varieties. Note that we are only using the `canonical enhancement' of $\Sigma$ to a stacky fan, but Kuwagaki's equivalence holds for general stacky fans $(\Sigma,\beta)$ if we modify the FLTZ Lagrangian to $\eL_{\Sigma,\beta}$ by replacing $\sigma^\perp \subset M_{S^1}$ with bigger sets $(\bigoplus_{\rho \in \sigma(1)} \Z b_\rho)^\perp_{S^1} \subset M_{S^1}$. 
\end{remark}

\subsection{Combinatorics}

We begin the section by introducing some definitions from combinatorial geometry, building up towards duality of nef partitions, which will be an important actor in the rest of this work. 

\begin{definition}\label{definition:reflexive-polytope}
    We call a lattice polytope $\Delta \subset M_{\R}$ containing the origin in its interior \emph{reflexive} if its \emph{dual polytope} $\Delta^{\vee}\coloneqq\{ y \in N_{\R}\colon \langle y, x \rangle \leq 1 \textnormal{ for all } x \in \Delta\}$ is also a lattice polytope. 
\end{definition}

We now review some standard facts about nef partitions. This notion was originally introduced in \cite{Borisov93}, but we shall follow slightly different conventions from later work (such as \cite{BB96} and \cite{BN}) instead. 

\begin{definition}\label{definition:nef-partition}
    Let $\Delta \subset M_{\R}$ be a reflexive polytope. A \emph{nef partition of $\Delta$ of length $r$} is a Minkowski sum decomposition $\Delta=\Delta_1+\dots+\Delta_r$, where $\Delta_1$, \dots, $\Delta_r$ are lattice polytopes containing the origin. 
\end{definition}
\begin{remark}\label{remark:nef_defn}
    Note that the original definition of nef partitions also demands that there exist PL convex functions $\varphi_j\colon M_{\R}\rightarrow \R$ for $j=1,\dots,r$ such that $\varphi_j(v)=1$ for $v$ a non-zero vertex of $\Delta_j$ and $\varphi_j$ identically vanishes on the other polytopes $\Delta_k$, $k\neq j$. It follows from the results in \cite[Chapter 3]{BN} that, if we restrict our attention to Minkowski decompositions of reflexive polytopes, this condition is unnecessary. 
\end{remark}

\begin{lemma}[{\cite[Corollary 3.17]{BN}}]\label{lemma:nef-summands-disjoint}
    For a nef partition $\Delta=\Delta_1+\dots+\Delta_r$, we have $\R_{\geq 0} \Delta_j \cap \R_{\geq 0} \Delta_k=\{0\}$ for all $j \neq k$. 
\end{lemma}

Due to its conciseness, we adopt the next definition from \cite{HZ} instead of following the original reference \cite{BB96}: 

\begin{definition}\label{definition:dual-nef-partition}
    Let $\Delta=\Delta_1+\dots+\Delta_r$ be a nef partition, define the convex functions $\psi_j\colon N_{\R}\rightarrow \R$ by 
    \begin{equation*}
        \psi_j(y)\coloneqq\sup\left\{\langle y, x \rangle \colon x \in \Delta_j\right\}.
    \end{equation*}
    Let $\nabla_j\subset N_{\R}$ be a polytope defined as the convex hull of $0$ and all $y \in \Delta^\vee$ such that $\psi_j(y)=1$. Then the \emph{dual nef partition} is given by $\nabla=\nabla_1+\dots+\nabla_r$.
\end{definition}
\begin{theorem}[{\cite{Borisov93}}]\label{theorem:dual-nef-partitions}
    Let $\Delta=\Delta_1+\dots+\Delta_r$ be a nef partition of a reflexive polytope $\Delta$. Then $\nabla=\nabla_1+\dots+\nabla_r$ is a reflexive polytope, so the dual nef partition is indeed a nef partition. Moreover, one also has 
    \begin{equation*}
        \begin{split}
            \nabla^{\vee}&=\conv\left\{\Delta_1,\dots,\Delta_r\right\},\\
            \Delta^{\vee}&=\conv\left\{\nabla_1,\dots,\nabla_r\right\}.
        \end{split}
    \end{equation*}
\end{theorem}

\begin{remark}\label{remark:nef_dual}
    To draw a complete picture on both sides of the duality, we can also consider the functions $\varphi_j\colon M_{\R}\rightarrow \R$ given as 
    \begin{equation*}
        \varphi_j(x)\coloneqq\sup\left\{\langle y, x \rangle : y \in \nabla_i\right\},
    \end{equation*}
    which recover the functions from Remark \ref{remark:nef_defn}, so $\Delta=\Delta_1+\dots+\Delta_r$ is the dual nef partition to $\nabla=\nabla_1+\dots+\nabla_r$ and dualising twice gets us back where we started. Moreover, we can verify that $\varphi=\varphi_1+\dots+\varphi_r$ is the support function of $\nabla^\vee$ and $\psi=\psi_1+\dots+\psi_r$ is the support function of $\Delta^\vee$. 
\end{remark}
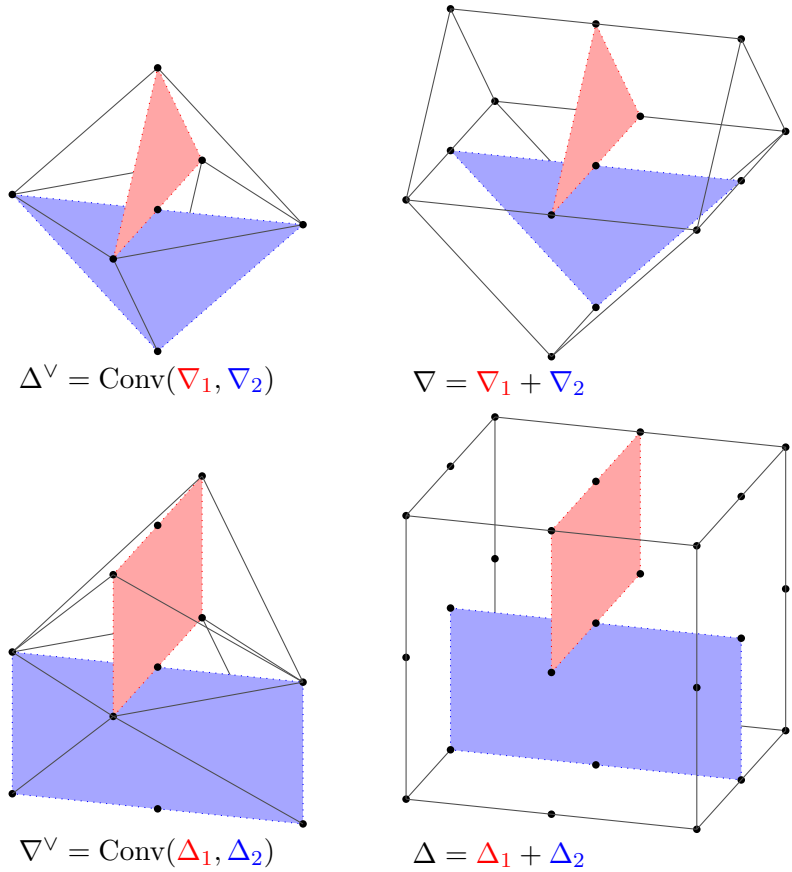
\begin{figure}[ht!]
\tdplotsetmaincoords{70}{107} 
\begin{center}
\begin{tikzpicture}[tdplot_main_coords, scale=2]
    
    \coordinate (A) at (1,0,0);
    \coordinate (B) at (0,0,1); 
    \coordinate (C) at (-1,0,0); 
    
    \coordinate (D) at (0,1,0); 
    \coordinate (E) at (0,0,-1); 
    \coordinate (F) at (0,-1,0);
    \coordinate (O) at (0,0,0);
    
    \draw[black!70] (C) -- (E); 
    \draw[black!70] (C) -- (D); 
    \draw[black!70] (B) -- (D); 
    \draw[black!70] (B) -- (F); 
    \draw[black!70] (C) -- (F);
    \filldraw[blue, fill=blue!35, dotted] (D) -- (E) -- (F) -- cycle;
    
    \filldraw[red, fill=red!35, dotted] (A) -- (B) -- (C) -- cycle;
    \draw[black!70] (A) -- (E); 
    \draw[black!70] (A) -- (D); 
    \draw[black!70] (A) -- (F);
    
    \node[fill,circle,inner sep=1pt] at (A) {};
    \node[fill,circle,inner sep=1pt] at (B) {};
    \node[fill,circle,inner sep=1pt] at (C) {};
    
    \node[fill,circle,inner sep=1pt] at (D) {};
    \node[fill,circle,inner sep=1pt] at (E) {};
    \node[fill,circle,inner sep=1pt] at (F) {};
    \node[fill,circle,inner sep=1pt] at (O) {};
    \node[below right] at (current bounding box.south west) {$\Delta^\vee=\textnormal{Conv}(\textcolor{red}{\nabla_1},\textcolor{blue}{\nabla_2})$};
\end{tikzpicture}
\hspace{1cm}
\begin{tikzpicture}[tdplot_main_coords, scale=2]

    \coordinate (A) at (1,0,0);
    \coordinate (B) at (0,0,1); 
    \coordinate (C) at (-1,0,0); 

    \coordinate (D) at (0,1,0); 
    \coordinate (E) at (0,0,-1); 
    \coordinate (F) at (0,-1,0);
    \coordinate (O) at (0,0,0);

    \draw[black!70] ($(C)+(D)$) -- ($(C)+(E)$);
    \draw[black!70] ($(C)+(D)$) -- ($(C)+(F)$);
    \draw[black!70] ($(C)+(E)$) -- ($(C)+(F)$);
    \draw[black!70] ($(A)+(E)$) -- ($(C)+(E)$);
    
    \filldraw[blue, fill=blue!35, dotted] (D) -- (E) -- (F) -- cycle;
    \filldraw[red, fill=red!35, dotted] (A) -- (B) -- (C) -- cycle;
    
    \node[fill,circle,inner sep=1pt] at (A) {};
    \node[fill,circle,inner sep=1pt] at (B) {};
    \node[fill,circle,inner sep=1pt] at (C) {};
    \node[fill,circle,inner sep=1pt] at ($(C)+(F)$) {};
    \node[fill,circle,inner sep=1pt] at ($(C)+(D)$) {};
    \node[fill,circle,inner sep=1pt] at ($(A)+(F)$) {};
    \node[fill,circle,inner sep=1pt] at ($(A)+(D)$) {};
    \node[fill,circle,inner sep=1pt] at ($(A)+(E)$) {};
    \node[fill,circle,inner sep=1pt] at ($(B)+(F)$) {};
    \node[fill,circle,inner sep=1pt] at ($(B)+(D)$) {};

    \node[fill,circle,inner sep=1pt] at (D) {};
    \node[fill,circle,inner sep=1pt] at (E) {};
    \node[fill,circle,inner sep=1pt] at (F) {};
    \node[fill,circle,inner sep=1pt] at (O) {};
    
    \draw[black!70] ($(A)+(D)$) -- ($(A)+(F)$);
    \draw[black!70] ($(A)+(E)$) -- ($(A)+(F)$);
    \draw[black!70] ($(A)+(D)$) -- ($(A)+(E)$);
  
    \draw[black!70] ($(D)+(A)$) -- ($(D)+(B)$);
    \draw[black!70] ($(D)+(B)$) -- ($(D)+(C)$);
    \draw[black!70] ($(D)+(C)$) -- ($(D)+(A)$);

    \draw[black!70] ($(F)+(A)$) -- ($(F)+(B)$);
    \draw[black!70] ($(F)+(B)$) -- ($(F)+(C)$);
    \draw[black!70] ($(F)+(C)$) -- ($(F)+(A)$);

    \draw[black!70] ($(D)+(B)$) -- ($(F)+(B)$);
    \node[below right] at (current bounding box.south west) {$\nabla=\textcolor{red}{\nabla_1}+\textcolor{blue}{\nabla_2}$};
\end{tikzpicture}
\end{center}
\begin{center}
\begin{tikzpicture}[tdplot_main_coords, scale=2]
    
    \coordinate (A) at (1,0,0);
    \coordinate (B) at (1,0,1); 
    \coordinate (C) at (-1,0,1); 
    \coordinate (D) at (-1,0,0);
    
    \coordinate (E) at (0,1,0); 
    \coordinate (F) at (0,1,-1); 
    \coordinate (G) at (0,-1,-1);
    \coordinate (H) at (0,-1,0);
    \coordinate (O) at (0,0,0);

    \draw[black!70] (D) -- (E); 
    \draw[black!70] (D) -- (F);
    \draw[black!70] (D) -- (G); 
    \draw[black!70] (D) -- (H);

    \filldraw[blue, fill=blue!35, dotted] (E) -- (F) -- (G) -- (H) -- cycle;
    
    \filldraw[red, fill=red!35, dotted] (A) -- (B) -- (C) -- (D) -- cycle;
    
    \node[fill,circle,inner sep=1pt] at (A) {};
    \node[fill,circle,inner sep=1pt] at (B) {};
    \node[fill,circle,inner sep=1pt] at (C) {};
    \node[fill,circle,inner sep=1pt] at (D) {};
    \node[fill,circle,inner sep=1pt] at (E) {};
    \node[fill,circle,inner sep=1pt] at (F) {};
    \node[fill,circle,inner sep=1pt] at (G) {};
    \node[fill,circle,inner sep=1pt] at (H) {};
    \node[fill,circle,inner sep=1pt] at (O) {};
    \node[fill,circle,inner sep=1pt] at ($(A)+(C)$) {};
    \node[fill,circle,inner sep=1pt] at ($(E)+(G)$) {};
    
    \draw[black!70] (A) -- (E); 
    \draw[black!70] (A) -- (F);
    \draw[black!70] (A) -- (G); 
    \draw[black!70] (A) -- (H);

    \draw[black!70] (B) -- (E); 
    \draw[black!70] (B) -- (H);
    \draw[black!70] (C) -- (E); 
    \draw[black!70] (C) -- (H);
    \node[below right] at (current bounding box.south west) {$\nabla^\vee=\textnormal{Conv}(\textcolor{red}{\Delta_1},\textcolor{blue}{\Delta_2})$};

\end{tikzpicture}
\hspace{1cm}
\begin{tikzpicture}[tdplot_main_coords, scale=2]
    
    \coordinate (A) at (1,0,0);
    \coordinate (B) at (1,0,1); 
    \coordinate (C) at (-1,0,1); 
    \coordinate (D) at (-1,0,0);
    
    \coordinate (E) at (0,1,0); 
    \coordinate (F) at (0,1,-1); 
    \coordinate (G) at (0,-1,-1);
    \coordinate (H) at (0,-1,0);
    \coordinate (O) at (0,0,0);

    \node[fill,circle,inner sep=1pt] at ($(D)+(G)$) {};
    \node[fill,circle,inner sep=1pt] at ($(D)+(E)+(G)$) {};
    
    \draw[black!70] ($(D)+(F)$) -- ($(D)+(G)$);
    \draw[black!70] ($(D)+(G)$) -- ($(A)+(G)$);
    
    \draw[black!70] ($(C)+(H)$) -- ($(G)+(D)$);
    
    \filldraw[blue, fill=blue!35, dotted] (E) -- (F) -- (G) -- (H) -- cycle;
    
    \filldraw[red, fill=red!35, dotted] (A) -- (B) -- (C) -- (D) -- cycle;
    
    \node[fill,circle,inner sep=1pt] at (A) {};
    \node[fill,circle,inner sep=1pt] at (B) {};
    \node[fill,circle,inner sep=1pt] at (C) {};
    \node[fill,circle,inner sep=1pt] at (D) {};
    \node[fill,circle,inner sep=1pt] at (E) {};
    \node[fill,circle,inner sep=1pt] at (F) {};
    \node[fill,circle,inner sep=1pt] at (G) {};
    \node[fill,circle,inner sep=1pt] at (H) {};
    \node[fill,circle,inner sep=1pt] at (O) {};
    \node[fill,circle,inner sep=1pt] at ($(A)+(C)$) {};
    \node[fill,circle,inner sep=1pt] at ($(E)+(G)$) {};
    
    \node[fill,circle,inner sep=1pt] at ($(A)+(E)$) {};
    \node[fill,circle,inner sep=1pt] at ($(A)+(F)$) {};
    \node[fill,circle,inner sep=1pt] at ($(A)+(G)$) {};
    \node[fill,circle,inner sep=1pt] at ($(A)+(H)$) {};
    \node[fill,circle,inner sep=1pt] at ($(A)+(E)+(G)$) {};

    \node[fill,circle,inner sep=1pt] at ($(D)+(E)$) {};
    \node[fill,circle,inner sep=1pt] at ($(D)+(F)$) {};

    \node[fill,circle,inner sep=1pt] at ($(D)+(H)$) {};

    \node[fill,circle,inner sep=1pt] at ($(E)+(B)$) {};
    \node[fill,circle,inner sep=1pt] at ($(E)+(C)$) {};
    \node[fill,circle,inner sep=1pt] at ($(H)+(B)$) {};
    \node[fill,circle,inner sep=1pt] at ($(H)+(C)$) {};
    \node[fill,circle,inner sep=1pt] at ($(E)+(A)+(C)$) {};
    \node[fill,circle,inner sep=1pt] at ($(H)+(A)+(C)$) {};

    \draw[black!70] ($(B)+(E)$) -- ($(C)+(E)$);
    \draw[black!70] ($(C)+(E)$) -- ($(C)+(H)$);
    \draw[black!70] ($(C)+(H)$) -- ($(B)+(H)$);
    \draw[black!70] ($(B)+(H)$) -- ($(B)+(E)$);
    
    \draw[black!70] ($(A)+(F)$) -- ($(D)+(F)$);

    \draw[black!70] ($(A)+(G)$) -- ($(A)+(F)$);

    \draw[black!70] ($(B)+(E)$) -- ($(F)+(A)$);
    \draw[black!70] ($(C)+(E)$) -- ($(F)+(D)$);
    \draw[black!70] ($(B)+(H)$) -- ($(G)+(A)$);
    \node[below right] at (current bounding box.south west) {$\Delta=\textcolor{red}{\Delta_1}+\textcolor{blue}{\Delta_2}$};
\end{tikzpicture}
\end{center}
    \caption{Dual nef partitions $(\Delta=\Delta_1+\Delta_2, \nabla=\nabla_1+\nabla_2)$ from Example \ref{example:running-nef}}
    \label{fig:dual-nef}
\end{figure}
\begin{example}\label{example:running-nef}
    Our running example of a nef partition (depicted on Figure \ref{fig:dual-nef}) is going to be a length two nef partition of the reflexive cube $\Delta=\conv(\pm e_1\pm e_2 \pm e_3) \subset M_\R$ into $\Delta_1=\conv(e_1,e_1+e_3,-e_1+e_3,-e_1)$ and $\Delta_2=\conv(e_2,e_2-e_3,-e_2-e_3,-e_2)$ (same as in \cite{HZ}). The dual nef partition is then given by $\nabla_1=\conv(\eta_1,\eta_3,-\eta_1)$ and $\nabla_2=\conv(\eta_2,-\eta_3,-\eta_2)$ in the dual basis of $N_\R$. In particular, we check that $\Delta^\vee=\conv(\pm \eta_1,\pm\eta_2,\pm\eta_3)=\conv(\nabla_1,\nabla_2)$. The associated toric Fano variety is $X_\Delta \cong(\Pp^1)^3$ with moment polytope $\Delta$, while the toric variety $X_\nabla$ is singular and needs to be resolved (which corresponds to subdividing the fan $\check{\Sigma}_0$ on the faces of $\nabla^\vee$ by a function $\check{h}$). 
\end{example}

Sometimes, nef partitions can be decomposed into smaller pieces, which leads to the following notion:

\begin{definition}\label{definition:irred_nef}
    We say that a nef partition $\Delta=\Delta_1+\dots+\Delta_r$ is \emph{irreducible} if there is no proper subset $\{i_1,\dots, i_s \} \subsetneq \{ 1, \dots, r\}$ such that the polytope $\Delta_{i_1}+\dots+\Delta_{i_s}$ contains $0$ in its relative interior. Otherwise, the nef partition is called \emph{reducible}.
\end{definition}

Note that if we have two nef partitions $\Delta=\Delta_1+\dots+\Delta_r$ inside $M_{\R}$ and $\Delta'=\Delta_1'+\dots+\Delta_s'$ inside $M'_{\R}$, we can consider $(\Delta_1\oplus0+\dots+\Delta_r\oplus0)+(0\oplus \Delta_1'+\dots+0\oplus \Delta_s')$ as a nef-partition of $\Delta \oplus \Delta' \subset M_{\R} \oplus M'_{\R}$. Note that the new nef partition will clearly be reducible, and in fact, any reducible nef partition can be decomposed in a similar fashion:

\begin{theorem}[{\cite[Theorem 5.8]{BB96}}]\label{theorem:irred_nef}
    Let $\{\Delta_j\}_{j=1}^r$ be a nef partition of a reflexive lattice polytope $\Delta$ in $M$. Then there exists a refinement $M \subset M'$ of finite index such that $\{\Delta_j\}_{j=1}^r$ decomposes as a direct sum of $l \geq 1$ irreducible nef partitions inside $M_1' \oplus \dots \oplus M'_l=M'$.
\end{theorem}

\begin{definition}\label{definition:induced_tri}
    Let $P$ be a lattice polytope inside $M_\R$, $A\subset P \cap M$ a set of lattice points in $P$ satisfying $P=\conv(A)$ and $h \in \R^A$ a function. Then the \emph{regular subdivision induced by $h$} is defined by projecting the lower convex hull of the graph $\Gamma(h) \subset P \times \R$ onto $P$. A subdivision of $P$ into lattice polyhedra is called \emph{regular} if it is induced by some choice of $A$ and $h$.
\end{definition}

Equivalently, we can consider the largest piecewise affine linear convex function $\psi_h \colon P \rightarrow \R $ satisfying $\psi_h(\alpha) \leq h(\alpha)$ for all $\alpha \in A$, and consider the subdivision of $P$ into its domains of linearity. We will often make use of the following result about the existence of a particular kind of regular subdivisions.

\begin{lemma}\label{lemma:triang_fns}
    Let $A$ be a finite set of lattice points containing the origin and denote $P=\conv(A)$. Then there exists a function $h\in \R^A$ such that:
    \begin{enumerate}[(1)]
        \item $h(0)=0$ and $h(\alpha)>0$ for all $\alpha \in A \backslash\{0\}$;
        \item the regular subdivision induced by $h$ is a triangulation;
        \item all the maximal simplices of the induced triangulation contain $0$, i.e. it is a \emph{star-shaped triangulation} centred at the origin. 
    \end{enumerate}
\end{lemma}
\begin{proof}
    The argument from \cite[Lemma 2.5.3]{GS22} gives the desired statement, but we shall repeat the argument for completeness: let $h_0$ be the piecewise linear function that is identically equal to $1$ on the boundary $\partial P$ and, for each facet $F$ of $P$, pick a function $h_F \in \R_{<0}^{F \cap A}$ that induces a regular triangulation of $F$ and so that $h_F|_{F \cap F' \cap A}=h_{F'}|_{F \cap F' \cap A}$ for all facets $F,F'\subset \partial P$, so they glue to a function $h_{\partial P} \in \R_{<0}^{\partial P \cap A}$ (by standard properties of secondary fans, this can be done by starting from $h_{\partial P} \equiv -1$ and perturbing generically). The desired function can then be defined through $h(\alpha)=h_0(\alpha)$ if $\alpha \notin \partial P$ and $h(\alpha)=1+\varepsilon h_{\partial P}(\alpha)$ if $\alpha \in \partial P$ for any sufficiently small $\varepsilon>0$. 
\end{proof}

We call a function $h$ (with an implicit choice of the set $A$) that induces a regular triangulation of $P$ a \emph{triangulating function on $P$}, while a function $h$ satisfying all the requirements of the Lemma \ref{lemma:triang_fns} is called a \emph{centred triangulating function}. We shall also need the following generalisation of unimodular triangulations:

\begin{definition}\label{definition:refined-triangulation}
    We say that a triangulating function on $P$ is \emph{refined} if, for any simplex $T$ in the induced regular triangulation $\mT$, the set of points $\alpha \in T \cap A$ for which $(\alpha,h(\alpha))$ lies on the lower convex hull of $\Gamma(h)$ is equal to the vertex set of $T$. 
\end{definition}

Using the equivalent definition of regular subdivisions, one can instead require that for every $T \in \mT$ that is a domain of linearity, the set of $\alpha \in T \cap A$ for which $\psi_h(\alpha)=h(\alpha)$ is equal to the vertex set of $T$. 

\begin{example}
    Inducing a regular unimodular triangulation is sufficient, but not necessary for ensuring that the corresponding $h$ is refined (if we take a $P=\conv(\pm2e_1,\pm2e_2)$ and $A=P\cap M$ with $h(0)=0$, $h(\pm2e_j)=2$ and $h(\alpha)>2$ for all other points $\alpha \in A$, this will give a refined triangulating function for which the simplices are non-unimodular). More generally, it is enough to ask that every simplex in $\mT$ contains no lattice points in $A$ other than its vertices; conversely, given any refined triangulating function, we can take $A'\subset A$ to be the set of points $\alpha \in A$ that are vertices in $\mT$, then $h|_{A'}$ will also be a refined triangulating function where all the simplices are `minimal' in this sense. 
\end{example}

Note that if we take our polytope $P$ in Lemma \ref{lemma:triang_fns} to be reflexive, we can ensure that the centred triangulating function that we construct is also refined (since the origin is the only interior lattice point of $P$ in this case, it suffices to start from a refined triangulation of $\partial P$). Note that we can not, in general, force the triangulation to be unimodular, since the faces of $\partial P$ need not admit a regular unimodular triangulation. 

In what follows, we will often need to choose a refined triangulating function as an auxiliary piece of data. In particular, for $\Delta=\Delta_1+\dots+\Delta_r$ and $\nabla=\nabla_1+\dots+\nabla_r$ dual nef partitions, pick two centred refined triangulating functions $(h,\check{h})$ on $\Delta^\vee$, $\nabla^\vee$ and denote the resulting triangulations of the boundaries $\partial \Delta^{\vee}$ and $\partial \nabla^{\vee}$ by $\mT$ and $\check{\mT}$, respectively. Note that since $\Delta^{\vee}=\conv\left\{\nabla_1,\dots,\nabla_r\right\}$, we also obtain induced triangulations $\mT_j$ of all the faces of $\partial\nabla_j$ not containing the origin for $j=1,\dots,r$ by intersecting $\mT$ with $\partial \nabla_j$ and, analogously, triangulations $\check{\mT}_j$ of the boundary faces of $\partial \Delta_j$ not containing the origin. 

\subsection{Geometry}

We shall now briefly explain the relevance of these combinatorial constructions to mirror symmetry. Recall that to a full-dimensional lattice polytope $\Delta \subset M$, we can associate a complete toric variety $X_{\Delta}$. The resulting variety turns out to be a Gorenstein Fano variety if and only if $\Delta$ is a reflexive polytope. Recall that the combinatorics of the same toric variety can also be encoded in the \emph{normal fan} of $\Delta$ $\Sigma_0 \subset N_{\R}$, whose positive-dimensional cones are $\cone(F)\coloneqq\bigcup_{t\geq0}t\cdot F$ for faces $F \subset \partial\Delta^{\vee}$, so we can also denote the associated toric variety as $X_{\Sigma_0}$. When $\Delta$ is reflexive with a nef partition, we can also repeat this story for $\nabla$ to produce a fan $\check{\Sigma}_0$ and a corresponding toric variety\footnote{Note that the notations for character and cocharacter lattices for this variety will then be reversed compared to the standard notation, so that $N$ is its character lattice and $M$ is the cocharacter lattice.} $X_\nabla=X_{\check{\Sigma}_0}$. By taking cones on simplices of $\mT$, $\check{\mT}$ instead of the entire faces, we get simplicial refinements $\Sigma$, $\check{\Sigma}$ of $\Sigma_0$, $\check{\Sigma}_0$. The nef partition naturally induces a decomposition of the ray generators for rays in $\Sigma(1)$ into sets $\Sigma_j(1)$ of rays coming from $\mT_j$ (these are disjoint by Lemma \ref{lemma:nef-summands-disjoint}).

\begin{remark}\label{remark:resolutions}
    The refinement $\Sigma \rightarrow \Sigma_0$ corresponds to a crepant toric partial resolution of singularities $X_\Sigma \rightarrow X_{\Sigma_0}$, where $X_\Sigma$ is an orbifold. Analogously, $\check{\mT}$ corresponds to a partial resolution of singularities $X_{\check{\Sigma}}\rightarrow X_{\check{\Sigma}_0}$. We can also view $\Sigma$, $\check{\Sigma}$ as stacky fans by picking primitive ray generators, which gives us smooth toric stacks $\mX_{\Sigma}$, $\mX_{\check{\Sigma}}$ with coarse moduli spaces $X_\Sigma$, $X_{\check{\Sigma}}$, so they also come with morphisms $\mX_\Sigma \rightarrow X_{\Sigma_0}$ and $\mX_{\check{\Sigma}} \rightarrow X_{\check{\Sigma}_0}$.   
\end{remark}

\begin{remark}\label{remark:irred-nef-wlog}
    Note that if we have complete intersections $Z \subset X_{\Delta}$, $Z' \subset X_{\Delta'}$ corresponding to nef partitions of $\Delta$, $\Delta'$, then the direct sum of nef partitions will yield $Z\times Z' \subset X_{\nabla}\times X_{\nabla'}=X_{\nabla \oplus \nabla'}$. The process of refining the lattice corresponds to quotienting our varieties by an action of a finite abelian group. Moreover, the cover is going to be unbranched over the algebraic torus, so we also get a corresponding unbranched cover of open BBCI's. There are general results that allow one to compute wrapped Fukaya categories of products (Künneth theorem from \cite{gps2}) and finite unbranched covers (see, e.g., \cite{seidel-quartic} and \cite{seidel-genus2}) and analogous results on the B-side, so we shall restrict our attention to the irreducible case for simplicity. 
\end{remark}

The geometric motivation for nef partitions is provided by the following result due to Batyrev and Borisov: 

\begin{proposition}[{\cite[Proposition 4.5]{BB96}}]
\label{proposition:nef_cy}
    Assume that an $n$-dimensional Gorenstein toric Fano variety $X_{\Delta}$ contains a codimension $r$ complete intersection $Z$ ($r<n$) of $r$ semi-ample Cartier divisors $Y_1$, \dots, $Y_r$ such that the canonical sheaf of $Z$ is trivial. Then there exist lattice polyhedra $\Delta_1$, \dots, $\Delta_r$ such that 
    \begin{equation*}
        \Delta=\Delta_1+\dots+\Delta_r.
    \end{equation*}
\end{proposition}

\begin{definition}\label{definition:bbci}
    We call the variety $Z \subset X_{\Delta}$ constructed above a \emph{Batyrev--Borisov complete intersection} (BBCI) corresponding to a nef-partition $\Delta=\Delta_1+\dots+\Delta_r$. We call its interior $Z^\circ=Z\backslash \partial X_{\Delta}$ an \emph{open Batyrev--Borisov complete intersection}.
\end{definition}

\begin{remark}\label{remark:coord_translation}
    One can also provide an explicit description for $Z^\circ$ in terms of coordinates on the dense torus $N_{\Cs}$ inside $X_{\Delta}$: the interior of a Cartier divisor $Y_i$ is the vanishing locus of a Laurent polynomial $f_i$ whose Newton polytope is equal to $\Delta_i$. Therefore, $Z^\circ$ is a very affine complete intersection
    \begin{equation*}
        Z^{\circ}=\V\left(f_1,\dots, f_r \right) \subset N_{\Cs},
    \end{equation*}
    and the variety $Z$ can then be recovered as the toric compactification of $Z^{\circ}$ inside $X_{\Delta}$. Note that we do not quite get a nef-partition, but a slightly more general Minkowski decomposition. By restricting our attention to the nef partition case, we only consider the cases where $0 \in \Delta_j$ for all $j=1$, \dots, $r$. This corresponds to constant terms of all the Laurent polynomials $f_j$ being non-zero, so both $Z^\circ=f^{-1}(0)$ (for $f=(f_1,\dots,f_r)\colon (\Cs)^n \rightarrow \C^r$) and all the neighbouring fibres of $f$ compactify to CYCI's inside $X_\Delta$.
\end{remark}

\begin{remark}\label{remark:gorenstein-cones}
    In the spirit of Remark \ref{remark:coord_translation}, the case when $0 \notin \Delta_j$ for some $j$ corresponds to $Z$ being \emph{rigid}, in the sense that changing $f_j$ by an arbitrarily small constant term $c_j$ will change the Newton polytope $\Delta_j$, hence the fibres of $f$ over such points no longer compactify to CYCI's inside $X_\Delta$. There is a more general framework of \emph{dual Gorenstein cones} that allows one to construct a mirror to such complete intersections (introduced in \cite{bb-cones}), with the caveat that the mirror might no longer be geometric, i.e. the mirror is not a Calabi--Yau variety, but a certain Calabi--Yau subcategory of the derived category of a higher-dimensional variety. Notably, this is the case for the \emph{generalised Greene-Plesser mirrors}, where HMS was proved in \cite{greene}. We shall not pursue this generalisation in the present work. 
\end{remark}

The key idea of the combinatorial construction from \cite{BB96} is to combine Proposition \ref{proposition:nef_cy} with duality of nef partitions: the dual nef partition $\nabla=\nabla_1+\dots+\nabla_r$ gives us a toric variety $X_\nabla$ and, inside it, another BBCI $\check{Z}$. We shall use the term \emph{Batyrev--Borisov mirror pairs} for varieties $(Z,\check{Z})$ arising this way. We are using the term liberally, since in general, one needs to consider crepant partial resolutions $X_\Sigma \rightarrow X_{\Sigma_0}$ and $X_{\check{\Sigma}} \rightarrow X_{\check{\Sigma}_0}$ as in Remark \ref{remark:resolutions} and look at the proper transforms of $Z$, $\check{Z}$ instead (or, for more generality, work with the associated smooth stacky objects). Note that we are also being intentionally vague about the relationship between coefficients of $Z$, $\check{Z}$ and the choice of Kähler forms $\omega$, $\check{\omega}$, i.e. the corresponding \emph{mirror map}, since it does not play a significant role in the bulk of this paper (see \cite[Section 6.3]{ck99} for a discussion of mirror maps for Batyrev--Borisov pairs). 

The main result of \cite{BB96} is that, under further assumptions on singularities in the set-up, the equalities of Hodge numbers $h^{p,q}(Z)=h^{n-r-p,q}(\check{Z})$ predicted by mirror symmetry hold for $p=0,1$ and $0 \leq q \leq n-r$. This provides evidence towards $(Z,\check{Z})$ being a mirror pair and, in particular, serves as an invitation to study homological mirror symmetry for it. The ultimate goal would be to prove some version of the following (intentionally somewhat imprecise) statement:

\begin{conjecture}[Batyrev--Borisov HMS]\label{conjecture:bb_hms}
    Let $(Z,\check{Z})$ be a Batyrev--Borisov mirror pair, where $Z$ and $\check{Z}$ are compact smooth subvarieties of toric varieties $X_{\Sigma}$, $X_{\check{\Sigma}}$ associated to the dual nef partitions $\Delta=\Delta_1+\dots+\Delta_r$, $\nabla=\nabla_1+\dots+\nabla_r$ and centred refined triangulating functions $h$, $\check{h}$. Then there exist quasi-equivalences:
    \begin{equation*}
        \begin{split}
            \perf \mF(Z) &\cong \coh(\check{Z}), \\
            \perf \mF(\check{Z}) &\cong \coh(Z),
        \end{split}
    \end{equation*}
    where $\mF$ stands for the compact Fukaya category of a closed symplectic manifold. 
\end{conjecture}

A special case of this statement for \emph{Batyrev mirror pairs} (corresponding to the hypersurface case with $r=1$) was recently proved in \cite{ghhps}. We make no attempts to address Conjecture \ref{conjecture:bb_hms} directly in this work, but we provide a brief explanation of its relationship with our results. The original formulation of mirror symmetry can be seen as a relationship between two moduli spaces of Calabi--Yau varieties. In particular, it is natural to ask what happens when one of the varieties degenerates. When we consider a degeneration $\{\check{X}_t\}$ on the B-side of homological mirror symmetry (also known as the \emph{large complex structure limit}), the category $\coh(\check{X}_t)$ is no longer homologically proper for the singular fibre $t=0$. This corresponds to the \emph{large volume limit} on the A-side, where the mirror $X_0$ of $\check{X}_0$ is expected to be a non-compact exact symplectic manifold, and the compact Fukaya category needs to be replaced by the \emph{wrapped Fukaya category} $\mW(X_0)$. Through this prism, the statement of Theorem \ref{theorem:hms-for-bbci} can be viewed as a limit case of Conjecture \ref{conjecture:bb_hms}.   
\section{Tropical geometry} \label{section:tropical}
In this following few sections, we summarise some existing results concerning tropicalisations of hypersurfaces in $(\Cs)^n$ and also provide their generalisations to complete intersections that will be necessary for accomplishing our goal, but might also be of independent interest. 
\subsection{Hypersurfaces}\label{section:tropical-hypersurfaces}
In this section, we review some standard results about tropicalisation of hypersurfaces inside $(\Cs)^n$.

Let $A \subset N$ be a finite set of lattice points and let $h\colon A \rightarrow \R$ be any function. Recall that the \emph{Legendre transform of $h$} is defined to be a piecewise linear function $L_h\colon M_{\R}\rightarrow \R$ given by
\begin{equation*}
    L_h(x)\coloneqq\sup_{\alpha \in A}\{\langle \alpha, x \rangle-h(\alpha)\}.
\end{equation*}
For brevity, when the choice of a function $h$ is clear from the context, we shall sometimes denote the affine function $\langle \alpha, x\rangle - h(\alpha)$ simply by $l_{\alpha}(x)$, so that we have $L_h(x)=\sup_{\alpha \in A}{l_{\alpha}(x)}$.

\begin{definition}
    The \emph{tropical hypersurface $\mA_{h,\trop}$ associated to $h$} is defined as the non-smooth locus of its Legendre transform inside $M_{\R}$. 
\end{definition}

\begin{remark}\label{remark:h_specifications}
    In what follows, we shall usually suppress the subscript $h$ to simplify the notation and write $\mA_{\trop}$ instead of $\mA_{h,\trop}$. This should not lead to any ambiguity, since the function $h$ is usually clear from the context. 
\end{remark}

The structure of the tropical hypersurface associated to $h$ is closely related to the subdivision of $\conv(A)$ induced by $h$:

\begin{proposition}[\cite{Mikh}]\label{proposition:amoeba_strata}
    Let $A$ be a finite set of lattice points in $N$ with $P=\conv(A)$ and a function $h\in \R^A$. Then we can decompose the associated tropical hypersurface $\mA_{\trop}$ into locally closed polyhedral cells 
        \begin{equation*}
            C_{S}\coloneqq\left\{ x \in M_{\R}\colon  L_h(x)=l_{\alpha}(x) \textnormal{ if and only if } \alpha \in S \right\},
        \end{equation*}
    where $S$ runs over all the positive-dimensional polytopes in the subdivision $\mS$ of $P$ induced by $h$. The poset of closures of such cells $\{\overline{C}_{S} \colon S \in \mS,  \dim(S)>0\}$ is anti-equivalent to its indexing poset $\{S \in \mS \colon  \dim(S)>0 \}$.
\end{proposition}
\begin{remark}\label{remark:comp_regions}
    If we instead consider the open sets $C_{\alpha}$ for $\alpha \in A$ a vertex of the subdivision given by $C_\alpha\coloneqq\{x \in M_\R \colon L_h(x)=l_{\alpha}(x), L_h(x)>l_{\alpha'}(x) \textnormal{ for }\alpha'\neq\alpha\}$, we can see that $L_h$ is smooth on $C_{\alpha}$. This provides us with a labelling of complementary regions of $\mA_{\trop}$ that is compatible with the poset structure of Proposition \ref{proposition:amoeba_strata} (in the sense that the resulting polyhedral decomposition of $M_{\R}$ is anti-equivalent to the face lattice of $\mS$). 
\end{remark}

Tropical hypersurfaces in $M_{\R}$ can be related to algebraic hypersurfaces in $M_{\Cs}$ through the technique called \emph{tropicalisation}, which we shall briefly review here. Recall that we have a logarithmic map $\Log \colon  M_{\Cs} \rightarrow M_{\R}$ given by $\Log(z_1,\dots,z_n)=(\log(|z_1|),\dots,\log(|z_n|))$ in coordinates or intrinsically as $\Log= \textnormal{id} \otimes \log(|\cdot|)\colon M\otimes \Cs \rightarrow M\otimes \R$. For $\beta>0$, we also consider $\Log_{\beta}\coloneqq\frac{1}{\beta}\Log$. 

\begin{definition}
    The \emph{tropical amoeba} of a hypersurface $X \subset M_{\Cs}$ is defined by $\mA_X=\Log(X) \subset M_{\R}$.
\end{definition}

Given the choice of a function $h$ as above and a choice of constants $c_{\alpha} \in \Cs$ for $\alpha \in A$, we can also consider a one-parameter family of hypersurfaces $H_{\beta}$ for $\beta>0$ in $M_{\Cs}$ given by 
\begin{equation*}
    H_{\beta}\coloneqq\left\{z \in M_{\Cs}\colon \sum_{\alpha \in A}{c_{\alpha}e^{-\beta h(\alpha)}z^{\alpha}}=0 \right\}.
\end{equation*}
Note that we suppress the choice of $c$ and $h$ from our notation in order to make it lighter, since we will always be concerned with the behaviour of such a family with fixed values of $c$, $h$ and a varying parameter $\beta$. To simplify notation, we shall denote the corresponding monomial terms as $f_{\alpha,\beta}(z)\coloneqq c_{\alpha}e^{-\beta h(\alpha)}z^{\alpha}$ and the defining polynomial as $f_{\beta}(z)\coloneqq \sum_{\alpha \in A}{f_{\alpha,\beta}(z)}$. We can then describe the limit of this family as $\beta$ becomes large:

\begin{theorem}[{\cite[Corollary 6.4]{Mikh}}]\label{theorem:mikh_trop}
    Denote the tropical amoebae of $H_{\beta}$ by $\mA_{\beta}=\Log_{\beta}(H_{\beta})$. Then we have $\mA_{\beta} \rightarrow \mA_{\trop}$ as $\beta \rightarrow \infty$ in the Gromov-Hausdorff metric.
\end{theorem}

\begin{remark}\label{remark:dom_terms}
    In the light of this result, the partition of $C_{S}$ acquires a new interpretation: for $\beta \gg 0$, if $z_0 \in H_{\beta}$ and its projection $u_0=\Log_{\beta}(z_0) \in \mA_{h,\beta}$ is close to the stratum $C_{S}$, then the dominant terms in the defining equation for $H_{\beta}$ around $z$ are precisely $f_{\alpha,\beta}(z)$ for $\alpha \in S$, since we have 
    \begin{equation*}
        |f_{\alpha,\beta}(z)|=|c_{\alpha}|e^{\beta l_{\alpha}(u)}.
    \end{equation*}
    Therefore, we can locally approximate $H_{\beta}$ by $\left\{ z \in M_{\Cs}\colon  \sum_{\alpha \in S}{f_{\alpha,\beta}(z)}=0\right\}$ near $z_0$. 
\end{remark}

To simplify further calculations, we now fix the notation for coordinates on $M_{\Cs}$ and $M_{\R}$: suppose that $z=(z_1,\dots,z_n)$ is a choice of global coordinates on $M_{\Cs}$ identifying it with $(\Cs)^n$, then we write $\rho\coloneqq\Log(z)$ for the induced coordinates on $M_{\R}$ (so that the standard polar coordinates on $M_{\Cs}$ are $z_j=e^{\rho_j+i\theta_j}$) and $u=\Log_{\beta}(z)=\frac{1}{\beta}\rho$ for their rescaling that is natural in the context of Theorem \ref{theorem:mikh_trop}.

The process of taking $\beta \rightarrow \infty$ has another nice feature under the extra assumption that $\mS$ is a refined triangulation -- it guarantees smoothness of the respective hypersurface:

\begin{lemma}\label{lemma:sm_nontailored}
    Suppose that $P=\conv(A)$ is a lattice polytope in $N$, $h\in \R^A$ a refined triangulating function on $P$, $c_{\alpha} \in \Cs$ are non-zero constants for $\alpha \in A$, and let $H_{\beta}=\V(f_{\beta})$ be the associated family of hypersurfaces. Then there exists a constant $\beta_0>0$ such that for all $\beta>\beta_0$, $H_{\beta}$ is a smooth complex hypersurface in $M_{\Cs}$. 
\end{lemma}
\begin{proof}
    The result is analogous to \cite[Proposition 4.5.1]{MS} here, but we are working in a slightly different setting, so we explicitly spell out a proof following a similar strategy for the sake of clarity (see also \cite[Proposition B.1]{integrality}). 

    Suppose that $z_0 \in H_{\beta}$ is a singular point, then we must have $\frac{\partial f_{\beta}}{\partial z_j}(z_0)=0$ for all $j=1$, \dots, $n$. In particular, for any $v \in M$, it must be true that 
    \begin{equation}
    \label{equation:sing_sing}
        \sum_{\alpha \in A}{f_{\alpha,\beta}(z_0)\langle \alpha, v \rangle} =\sum_{j=1}^n{v_jz_j\frac{\partial f_{\beta}}{\partial z_j}(z_0)}= 0.
    \end{equation}
    Denote $u_0=\Log_{\beta}(z_0)$ and let $S$ be the set of lattice points $\alpha \in A$ such that\footnote{The awkward choice of asymptotics of $O(\beta^{-\frac12})$ instead of the obvious choice $O(\beta^{-1})$, which would also clearly yield stronger bounds, is for the sake of consistency with the tailored case in Section \ref{section:tailoring}.} $l_{\alpha}(u_0) \geq L_h(u_0)-\beta^{-\frac12}$ if and only if $\alpha \in S$. By Lemma \ref{lemma:aff-nbhd-structure} and Proposition \ref{proposition:amoeba_strata}, for sufficiently large $\beta$, $\conv(S)$ is a polytope in the regular subdivision $\mS$ of $P$ induced by $h$. Moreover, by the assumption that $h$ is a refined triangulating function, $S$ must be equal to the vertex set of a simplex $\conv(S)$. Now, pick some vertex $\alpha_0 \in S$ such that $l_{\alpha_0}(u_0)=L_h(u_0)$. By multiplying the defining equation for $H_{\beta}$ by $z^{-\alpha}$ for some $\alpha \in S \backslash \{ \alpha_0 \}$, we can assume that $0$ is a vertex of $S$ different from $\alpha_0$ without loss of generality. 
    
    Also, observe that for $\alpha \in A \backslash S$, we must have $l_{\alpha}(u_0)<l_{\alpha_0}(u_0)-\beta^{-\frac12}$, which means that such $\alpha$ satisfy
    \begin{equation}
    \label{equation:estsing}
        |e^{-\beta h(\alpha)}z_0^{\alpha}| < e^{-\sqrt{\beta}}|e^{-\beta h(\alpha_0)}z_0^{\alpha_0}|.
    \end{equation}
    Now, note that since $S$ is a simplex, we can find a vector $v \in M_{\R}$ such that $\langle \alpha_0, v \rangle = 1$ and $\langle \alpha , v \rangle = 0$ for all $\alpha \neq \alpha_0$ in $S$ (if the other vertices of $S$ are $\alpha_1=0$, \dots, $\alpha_k$, we can pick any vector perpendicular to $\alpha_2$, \dots, $\alpha_k$, but not to $\alpha_0$). 

    Plugging this vector into Equation \ref{equation:sing_sing}, it follows that if $z_0$ is a singular point, we must have:  \begin{equation*}
    \begin{split}
        |c_{\alpha_0}|\left|e^{-\beta h (\alpha_0)}z_0^{\alpha_0}\right|&=\left| \sum_{\alpha \in A \backslash S} c_{\alpha}e^{-\beta h(\alpha)}\langle \alpha, v \rangle z_0^{\alpha} \right| \\ &\leq
        \sum_{\alpha \in A \backslash S} |c_{\alpha}|\left| e^{-\beta h(\alpha)}\langle \alpha, v \rangle z_0^{\alpha} \right| \\
        &< \left(\sum_{\alpha \in A \backslash S}{|c_{\alpha}||\langle \alpha, v \rangle|} \right) |e^{-\beta h(\alpha_0)} z_0^{\alpha_0}| e^{-\sqrt{\beta}}, \\
    \end{split}
    \end{equation*}
    where we have applied the triangle inequality and the estimate from Equation \ref{equation:estsing}. By simplifying, we see that in order for $z_0$ to be a singular point, the other quantities have to satisfy  
    \begin{equation*}
        e^{\sqrt{\beta}} < \frac{1}{|c_{\alpha_0}|}\sum_{\alpha \in A \backslash S}{|c_{\alpha}||\langle \alpha, v \rangle|},
    \end{equation*}
    which can only happen for $\beta<K(S)$, where $K(S)$ is a constant that depends only on the simplex $S$. Therefore, if $\beta\geq\sup_{S} K(S)$, $H_{\beta}$ must be a smooth complex hypersurface. 
\end{proof}

\subsection{Complete intersections} \label{section:tropical-ci}

In this section, we shall generalise some of the theory from Section \ref{section:tropical-hypersurfaces} to complete intersections inside $M_{\Cs}$. 

Suppose that we have $r$ lattice polytopes $P_1=\conv(A_1)$, \dots, $P_r=\conv(A_r)$ in $N$ equipped with functions $h_j \in \R^{A_j}$ that induce regular subdivisions $\mS_j$ of $P_j$ and constants $c_{\alpha,j} \in \Cs$ with $\alpha \in A_j$ for all $j=1$, \dots, $r$. Then we can consider a one-parameter family of $r$ polynomials given by
\begin{equation*}
    f_{\beta,j}(z)\coloneqq\sum_{\alpha \in A_j}{c_{\alpha,j}e^{-\beta h_j(\alpha)}z^{\alpha}},
\end{equation*}
which consist of monomials $f_{\alpha,\beta,j}(z)\coloneqq c_{\alpha,j}e^{-\beta h_j(\alpha)}z^{\alpha}$ for $\alpha \in A_j$ and give rise to the following one-parameter family of complete intersections inside $M_{\Cs}$:
\begin{equation*}
    Z_{\beta}\coloneqq\V(f_{\beta,1},\dots,f_{\beta,r}).
\end{equation*}
Denote the amoeba of $Z_{\beta}$ as $\mA_{\beta}=\Log(Z_{\beta})$ and the amoebae of the individual hypersurfaces $H_{\beta,j}=\V(f_{\beta,j})$ as $\mA_{\beta,j}$ for $j=1, \dots, r$. By Theorem \ref{theorem:mikh_trop}, we know that $\mA_{\beta,j} \rightarrow \mA_{\trop,j}$ as $\beta \rightarrow \infty$. Therefore, it is a natural guess that $\mA_{\beta}$ should tend to the tropical complete intersection $\mA_{\trop}=\bigcap_{j=1}^r{\mA_{\trop,j}}$ as $\beta\rightarrow \infty$. 

Before proving the result (Theorem \ref{theorem:trop_ci}), we first give a brief overview of the combinatorial model for tropical complete intersections following \cite[Chapter 4.6]{MS}. In order to simplify our discussion, we shall restrict our attention to complete intersection that are transverse in the following sense:

\begin{definition}\label{definition:tcci}
    We say that $\mA_{\trop}=\bigcap_{j=1}^r{\mA_{\trop,j}}$ is a \emph{codimension $r$ transverse tropical complete intersection (TTCI)} if for any $r$-tuple of cells $C_1 \subset \mA_{\trop,1}$, \dots, $C_r \subset \mA_{\trop,r}$, the intersection $\bigcap_{j=1}^r \relint(C_j)$ is transverse.
\end{definition}

Note that since there are finitely many cells, transversality is a generic condition that can be achieved by perturbing the triangulating functions $h_j$. 

To study the combinatorics of TTCI's, we shall employ the \enquote{Cayley trick} and add $r$ auxiliary dimensions to the entire set-up. This leads to the following key definition:

\begin{definition}\label{definition:cayley-polytope}
    The \emph{Cayley polytope} $C(P_1,\dots,P_r)$ associated to $\mA_{\trop}$ is defined as 
    \begin{equation*}
        C(P_1,\dots,P_r)\coloneqq\conv(P_1\times e_1, \dots, P_r \times e_r) \subset N_{\R} \times \R^r. 
    \end{equation*}
\end{definition}

The triangulating functions $h_j$ then induce a function $h_\tot \in \R^{A_\tot}$ for $A_\tot\coloneqq \bigsqcup_{j=1}^r A_j \times \{e_j\}$ by setting $h_\tot(\alpha\times e_j)=h_j(\alpha)$ for any lattice point $\alpha \in A_j$, and since $C(P_1,\dots,P_r)=\conv(A_\tot)$, we get a subdivision $\mS$ of $C(P_1,\dots,P_r)$. The Minkowski sum $P\coloneqq P_1+\dots+P_r$ can be recovered from the Cayley polytope as the intersection of its rescaling with a certain affine subspace:
\begin{equation*}
    P \cong r \cdot C(P_1,\dots,P_r) \cap \{ (x,\lambda_1,\dots,\lambda_r) \in N_\R \times \R^r \colon \lambda_j=1 \textnormal{ for all } j=1,\dots,r\},
\end{equation*}
so we can view $P$ as a subset of $\cone(C(P_1,\dots,P_r))$. One can naturally extend the function $h_\tot$ from the Cayley polytope to a piecewise liner function $\psi'_{h_\tot}$ on $\cone(C(P_1,\dots,P_r))$ by setting $\psi'_{h_\tot}(t\cdot x)\coloneqq t\cdot \psi_{h_\tot}(x)$ for all $x \in C(P_1,\dots,P_r)$ and $t>0$ (recall that $\psi_{h_\tot}$ can be defined as the largest convex piecewise affine function on the affine subspace $\{(x,\lambda_1,\dots,\lambda_r) \in N_\R \times \R^r \colon \sum_{j=1}^r\lambda_j=1\}$ satisfying $\psi_{h_\tot}(\alpha) \leq h_\tot(\alpha)$ for all $\alpha \in A_{\tot}$), which leads to the following definition:
\begin{definition}
    The \emph{regular mixed subdivision} $\mS'$ of $P$ is the subdivision induced by $h' \in \R^{A'}$ obtained by restricting $\psi'_{h_\tot}$ to the set $A' \coloneqq A_1+\dots+A_r$.
\end{definition}

One can associate a cell $S' \in \mS'$ of the regular mixed subdivision to a cell $S \in \mS$ as follows: denote by $S_j$ the face of $S$ consisting of the vertices whose $e_j$-coordinate is $1$, then the cell corresponding to $S$ is $S'=S_1+\dots+S_r$ (note that this is only non-empty if $S_j \neq \emptyset$ for all $j$). Conversely, we can recover $S$ from $S'$ as the minimal cell of $\mS$ that contains $\frac1r S'$. We call the cell $S' \in \mS'$ \emph{mixed} if $\dim(S_j) \geq 1$ for all $j=1,\dots,r$. 

Note that since every point $p \in \mA_{\trop}$ in a TTCI belongs to the relative interior of a unique cell $C_{S_j}$ of $\mA_{\trop,j}$ for all $j=1$, \dots, $r$, there is a natural labelling of the cells of $\mA_{\trop}$ by $r$-tuples $(S_1,\dots,S_r)$ and the corresponding cell $C_{(S_1,\dots,S_r)}$ will be the relative interior of a polyhedron of the expected dimension (by transversality). Therefore, the only remaining question is what $r$-tuples actually correspond to non-empty strata. This is answered by the following result, which is a reformulation of \cite[Theorem 4.6.9]{MS} that is better suited for our purposes:

\begin{theorem}\label{theorem:amoeba_strata_ci}
    The poset of closed cells $\overline{C}_{(S_1,\dots,S_r)}$ of a transverse tropical complete intersection $\mA_{\trop}=\bigcap_{j=1}^r{\mA_{\trop,j}}$ associated to polytopes $P_j$ with functions $h_j$ is dual to the complex of mixed cells $S'=S_1+\dots+S_r$ in the induced subdivision $\mS'$ of $P=P_1+\dots+P_r$. 
\end{theorem}

\begin{proof}
    First, suppose that the cells $C_{S_1} \subset \mA_{\trop,1}$, \dots, $C_{S_r} \subset \mA_{\trop,r}$ have a non-empty intersection, where $S_j$ is a cell of the decomposition $\mS_j$ of $P_j$ induced by $h_j$. Then it is clear that $S\coloneqq \conv(S_j\times e_j \colon  1 \leq j \leq r)$ will have a non-empty intersection with the affine plane $\{(x,\lambda)\colon  \lambda=\frac1r(1,\dots,1) \}$, which means that there is a cell $S'=S_1+\dots+S_r$ of the mixed subdivision corresponding to $S$, and $\dim(S_j) \geq 1$ is forced by $S_j$'s actually labelling cells of a tropical hypersurface, thus $S'$ is mixed. Observe that transversality of the intersection allows us to compute the dimension of the cell as $\dim(S')=\dim(\sum_j S_j)=\sum_j \dim(S_j)$, since the affine span of $S_j$ can be identified with the conormal bundle of $C_{S_j}$ at an interior point.
     
    Conversely, suppose that we have a mixed cell $S'=S_1+\dots+S_r \subset P$ coming from a cell $S \subset C(P_1,\dots,P_r)$, then we want to find some $x \in M_{\R}$ such that $x \in C_{S_j} \subset \mA_{\trop,j}$ for $j=1,\dots,r$ (since the cell is mixed, $\dim(S_j)\geq1$, so they indeed label cells of the hypersurfaces). Recall that the subdivision $\mS$ of $C(P_1,\dots,P_r) \subset N_{\R} \times \R^r$ is the regular subdivision induced by $h_\tot$, so it has an associated tropical hypersurface $\mA_{\trop,\tot} \subset M_{\R} \times (\R^r)^\vee$ that contains a cell $C_S$ dual to $S$. Splitting the vertices of $C(P_1,\dots,P_r)$ based on which $P_j \times e_j$ they come from tells us that the the Legendre transform of $h_\tot$ can be expressed as
    \begin{equation*}
    \begin{split}
        L_{h_\tot}(x,\lambda)&=\sup_{\alpha \in A_{\tot}}\left\{\langle \alpha, (x,\lambda) \rangle - h(\alpha)\right\}=\max_{1 \leq j \leq r} \left\{\lambda_j+\sup_{\alpha_j \in A_j}{\langle \alpha_j, x\rangle - h_j(\alpha_j)}\right\} \\
        &= \max_{1 \leq j \leq r} \left\{ \lambda_j + L_{h_j}(x) \right\}.
    \end{split}
    \end{equation*}
    Pick a point $(x,\lambda) \in C_S$. By definition, we then have $L_{h_\tot}(x,\lambda)=\langle \alpha,(x,\lambda) \rangle-h(\alpha)$ for all $\alpha \in S$ and $L_{h_\tot}(x,\lambda)>\langle \alpha,(x,\lambda) \rangle-h(\alpha)$ otherwise. Through the above expression, this implies that $L_{h_\tot}(x,\lambda)-\lambda_j=L_{h_j}(x)=l_{\alpha_j}(x)$ if $\alpha_j \in S_j$ and $L_{h_j}(x)>l_{\alpha_j}(x)$ otherwise for all $j$. Hence, by definition, $x$ lies in the cell $C_{S_j}$ of $\mA_{\trop,j}$, which is the tropical hypersurface associated to $h_j$, so $C_S=\bigcap_j C_{S_j}$ is non-empty, as desired. 
\end{proof}

Now that we have some understanding of how tropical complete intersections behave, we can state and prove an extension of Theorem \ref{theorem:mikh_trop}.

\begin{theorem}\label{theorem:trop_ci}
    Let $Z_{\beta}$ be the one-parameter family of complete intersections associated to lattice polytopes $P_1=\conv(A_1)$, \dots, $P_r=\conv(A_r)$ and functions $h_j\in \R^{A_j}$. Suppose that the associated tropical complete intersection $\mA_{\trop}$ is transverse. Then the rescaled amoebae $\mA_{\beta}=\frac{1}{\beta}\Log(Z_{\beta})$ of the family of complete intersections inside $M_{\Cs}$ converge to a subset of the tropical complete intersection $\mA_{\trop}$ as $\beta \rightarrow \infty$ in the Gromov-Hausdorff metric. 
\end{theorem}
\begin{proof}
    By the proof of \cite[Corollary 6.4]{Mikh}, the amoeba $\Log_{\beta}(H_{\beta,j})$ is contained in a neighbourhood $U(\mA_{\trop,j},\delta)$ of the hypersurface $\mA_{\trop,j}$ of radius $\delta=O(\beta^{-1})$ (with respect to some arbitrarily chosen norm on $M_{\R}$). We also clearly have $\Log_{\beta}(Z_{\beta}) \subset \bigcap_j\Log_{\beta}(H_{\beta,j})$, therefore it follows that $\Log_{\beta}(Z_{\beta}) \subset \bigcap_jU(\mA_{\trop,j},\delta)$, hence by Corollary \ref{corollary:nbhd-intersections}, there exists a constant $K>0$ such that for all $\beta>0$ large enough, $\Log_{\beta}(Z_{\beta}) \subset U(\mA_{\trop},K\delta)$ with $\delta=O(\beta^{-1})$, so we are done.
\end{proof}

Analogously, we can provide a generalisation of Lemma \ref{lemma:sm_nontailored} to the case of complete intersections:

\begin{lemma}\label{lemma:sm_nontailored_ci}
    Let $Z_{\beta}$ be the one-parameter family of complete intersections associated to lattice polytopes $P_1=\conv(A_1)$, \dots, $P_r=\conv(A_r)$ and refined triangulating functions $h_j\in \R^{A_j}$ inducing triangulations $\mT_j$. Suppose that associated tropical complete intersection $\mA_{\trop}$ is transverse. Then there exists a constant $\beta_0>0$ such that all $Z_{\beta}$ with $\beta>\beta_0$ are smooth complete intersections of complex dimension $n-r$ inside $M_{\Cs}$. 
\end{lemma}
\begin{proof}
    We proceed analogously to Lemma \ref{lemma:sm_nontailored}: suppose that $z_0 \in Z_{\beta}$ is a singular point. Then the collection of $r$ vectors $df_{\beta,1}(z_0)$, \dots, $df_{\beta,r}(z_0)$ are not linearly independent. 

    Let $T_j$ be the set of vertices in $\mT_j$ such that $l_{\alpha}(u_0) \geq L_{h_j}(u_0)-\beta^{-\frac12}$ holds. For $\beta$ large enough, it is guaranteed that $T_j$ will be the vertex set of a simplex in $\mT_j$ (similarly to Lemma \ref{lemma:sm_nontailored}). Moreover, Lemma \ref{lemma:aff-nbhd-structure} tells us that $(T_1,\dots,T_r)$ produced this way will label a non-empty cell of $\mA_{\trop}$, which means that it corresponds to a mixed cell by Theorem \ref{theorem:trop_ci} and, in particular, $\dim(T_j) \geq 1$ for all $j$. 

    Pick vertices $\alpha_{j} \in T_j$ such that $L_{h_j}(u_0)=l_{\alpha_{j}}(u_0)$. Then, analogously to Equation \ref{equation:estsing}, we get
    \begin{equation}
    \label{equation:estsing2}
        |e^{-\beta h_j(\alpha)}z_0^{\alpha}| < e^{-\sqrt{\beta}}|e^{-\beta h_j(\alpha_{j})}z_0^{\alpha_{j}}|,
    \end{equation}
    for all $\alpha \in P_j \backslash T_j$ and $j=1$, \dots, $r$. 
    
    By multiplying each defining equation $f_{\beta,j}$ by $z^{-\alpha'_j}$ for some $\alpha'_j \in T_j\backslash \{ \alpha_{j} \}$, we can assume that $0$ is a vertex of $T_j$ different from $\alpha_{j}$ without loss of generality.
    By transversality of the intersection, we know that $\dim(\sum_jT_j)=\sum_j\dim(T_j)$, hence $\bigcup_j T_j$ is the vertex set of a simplex of the expected dimension (because the set is affinely independent, thanks to each of the sets  $T_j$ being affinely independent as a vertex set of a simplex in $\mT_j$ and the sum $\sum_j \R T_j$ of the associated vector subspaces being direct by transversality), hence we can pick a vector $v_j \in M_{\R}$ for each $j=1$, \dots, $r$ such that $\langle \alpha_{j}, v_k \rangle = \delta_{jk}$ and $\langle \alpha, v_j \rangle = 0$ for all other vertices $\alpha$ in $\bigcup_jT_j$. In coordinates, we have an expression for the differentials:
    \begin{equation*}
        df_{\beta,j}(z_0)=\sum_{\alpha \in A_j}{f_{\alpha,\beta,j}(z_0)\langle\alpha_j,d\rho+id\theta \rangle },
    \end{equation*}
    hence pairing this with a vector $v_k=\sum_l v_{k,l}\partial_{\rho_l}$ gives:
    \begin{equation*}
        df_{\beta,j}(z_0)(v_k)=\delta_{jk}f_{\alpha_j,\beta,j}(z_0)+\sum_{\alpha \in A_j\backslash T_j}{f_{\alpha,\beta,j}(z_0)\langle \alpha, v_k  \rangle}.
    \end{equation*}
    Therefore, the $(r \times r)$ matrix $A_{jk}=df_{\beta,j}(z_0)(v_k)$ is of the form 
    \begin{equation*}
        A=\begin{pmatrix}
         f_{\alpha_1,\beta,1}(z_0)& 0 &  \hdots & 0\\
        0 & f_{\alpha_2,\beta,2}(z_0) & \hdots & 0 \\
        \vdots & \ddots & & \vdots \\
        0 & 0 & \vdots & f_{\alpha_r,\beta,r}(z_0) \\
        \end{pmatrix} \cdot \left(I_r+D \right),
    \end{equation*}
    where the matrix $D$ has all entries bounded by $Ke^{-\sqrt{\beta}}$ for some constant $K>0$ that only depends on the tuple $(T_1,\dots,T_r)$. But if $df_{\beta,1}(z_0)$, \dots, $df_{\beta,r}(z_0)$ are linearly dependent, the matrix $A$ can not be full rank, hence $I_r+D$ must fail to be invertible, which can only happen if $Ke^{-\sqrt{\beta}}\geq 1$, so for $\beta > \log(K)^2$, $z_0$ must actually be a smooth point of the complete intersection. Varying the tuple $(T_1,\dots,T_r)$ over all the possibilities and taking the maximal value of $K$ in this lower bound on $\beta$, we obtain the desired result. 
\end{proof}

In fact, the same proof gives us a slightly stronger result, saying that the map $(f_{\beta,1},\dots,f_{\beta,r})$ will be a submersion not only along $Z_{\beta}$, but also at all the points projected by $\Log_\beta$ onto a certain neighbourhood $\mA_\trop \subset \mA_{\trop}(\beta)$. More explicitly, $\mA_{\trop}(\beta)$ is the \emph{closed affine neighbourhood} (see Appendix \ref{section:appendix-polyhedra}) of $\mA_{\trop}$ of size $\beta^{-\frac12}$ defined as
\begin{equation*}
    \mA_{\trop}(\beta)\coloneqq \{ u \in M_{\R}\colon L_{h_j}(u) \leq l_{\alpha}(u)+\beta^{-\frac12} \textnormal{ for at least two } \alpha \in P_j \textnormal{ for }j=1,\dots,r\}.
\end{equation*}
This is true since the proof only relies on the facts that simplices $T_j$ are at least one-dimensional (which is true for points projecting to $\mA_\trop(\beta)$ by definition) and that $(T_1,\dots,T_r)$ labels a cell of $\mA_\trop$ (which is true by the same argument invoking Lemma \ref{lemma:aff-nbhd-structure}). We record this fact for future use:

\begin{corollary}\label{corollary:weird_generalisation}
    In the setting of Lemma \ref{lemma:sm_nontailored_ci}, any point $z_0 \in \Log^{-1}_\beta(\mA_{\trop}(\beta))$ will be a regular point of the map $(f_{\beta,1},\dots,f_{\beta,r})\colon M_{\Cs} \rightarrow \C^r$ for all $\beta>\beta_0$.
\end{corollary}

\subsection{Batyrev--Borisov complete intersections} \label{section:tropical-bbci}

In this section, we specialise the setup of Section \ref{section:tropical-ci} to the case that will mostly be of interest to us. After reviewing some basic facts about translating nef partitions into tropical geometry, we provide an alternative description for the combinatorics of $\mA_{\trop}$ that will be useful for us in the rest of the paper. 

Suppose that we start from the data of a nef partition $\nabla=\nabla_1+\dots+\nabla_r$ of a reflexive polytope $\nabla$ inside $N_{\R}$ and a centred refined triangulating function $h$ on $\Delta^\vee$ (recall that $\Delta^\vee=\conv(\nabla_j)$, so this also induces centred refined triangulating functions $h_j\coloneqq h|_{\nabla_j}$). Denote the triangulation of $\partial\Delta^{\vee}$ coming from $h$ as $\mT$ and, analogously, denote the corresponding triangulations of the faces of $\nabla_j$ not containing the origin by $\mT_j$. We can obtain a star-shaped triangulation $\mT * 0$ of $\Delta^\vee$ by taking its maximal simplices to be $\overline{T}\coloneqq \conv(0,T)$ for maximal $T \in \mT$, and star-shaped triangulations $\mT_j*0$ of $\nabla_j$. We shall denote the tropical hypersurfaces associated to the functions $h_j$ as $\mA_{\trop,j}$ and the tropical hypersurface associated to $h$ itself as $\mA_{\trop,\tot}$ (note that unlike in the general case, these live inside $M_{\R}$, not inside $M_{\R} \times (\R^r)^\vee$).

\begin{definition}\label{definition:batbor_trop_ci}
    We call $\mA_\trop=\bigcap_{j=1}^r \mA_{\trop,j}$ a \emph{tropical Batyrev--Borisov complete intersection} associated to the nef partition $\nabla=\nabla_1+\dots+\nabla_r$ and a centred refined triangulating function $h$. 
\end{definition}

Observe that $h$ being refined on $\Delta^\vee$ automatically implies that any tropical BBCI is transverse, so the results from the previous section apply. Following the general procedure for TTCI's, we can assign a family of open Batyrev--Borisov complete intersections $Z_{\beta}$ to this set-up (allowing for the extra freedom of the choice of constants $c_{\alpha,j}\in\Cs$ as before). Note that Lemma \ref{lemma:sm_nontailored_ci} guarantees that for $\beta$ sufficiently large, this will indeed define a smooth complete intersection inside $M_{\Cs}$, which proves Theorem \ref{theorem:main-thm}(\ref{main-thm-1}). 

\begin{example}\label{example:running-trop}
    In the setting of Example \ref{example:running-nef}, consider the function $h \colon \Delta^\vee \cap N \rightarrow \R$ given by $h(0)=0$, $h(\pm \eta_i)=1$. This will be a centred refined triangulating function, since the faces of $\Delta^\vee$ are already unimodular simplices. An example of an open BBCI associated to these data is $Z_\beta=\{f_1=f_2=0\}$ given by equations $f_1(z)=e^{-\beta}(z_1+z_1^{-1}+z_3)-1$, $f_2(z)=e^{-\beta}(z_2+z_2^{-1}+z_3^{-1})-1$. The associated Legendre transforms $L_1(u)=\max\{u_1-1,-u_1-1,u_3-1,0\}$ and $L_2(u)=\max\{u_2-1,-u_2-1,-u_3-1,0 \}$ define tropical hypersurfaces $\mA_{\trop,1}$ and $\mA_{\trop,2}$ respectively. Their intersection $\mA_\trop$ is a tropical BBCI, which is the tropicalisation of $Z_\beta$ (see Figure \ref{fig:trop-ci}). The hypersurface $\mA_{\trop,\tot}$ is associated to the function $L_\tot(u)=\max\{L_1(u),L_2(u)\}$. 
\end{example}
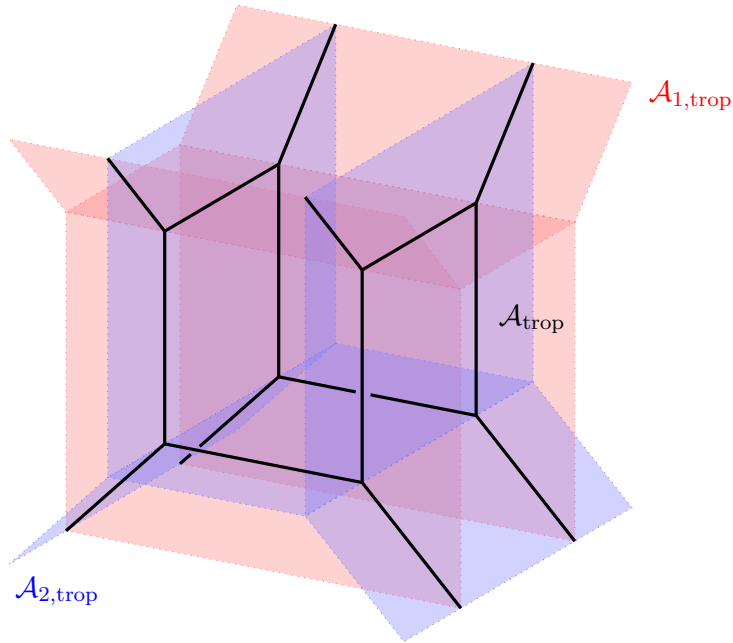
\begin{figure}[ht]
\centering
\tdplotsetmaincoords{70}{120} 
\begin{tikzpicture}[tdplot_main_coords, scale=1.5]

    \filldraw[red!50, fill opacity=0.35, dotted] (1,-2,1) -- (2,-2,2) -- (2,2,2) -- (1,2,1) -- cycle;
    \filldraw[red!50, fill opacity=0.35, dotted] (-1,-2,1) -- (-2,-2,2) -- (-2,2,2) -- (-1,2,1) -- cycle;
    
    \filldraw[red!50, fill opacity=0.35, dotted] (1,-2,1) -- (1,-2,-2) -- (1,2,-2) -- (1,2,1) -- cycle;
    
    \filldraw[red!50, fill opacity=0.35, dotted] (-1,-2,1) -- (-1,-2,-2) -- (-1,2,-2) -- (-1,2,1) -- cycle;
    
    \filldraw[red!50, fill opacity=0.35, dotted] (-1,-2,1) -- (1,-2,1) -- (1,2,1) -- (-1,2,1) -- cycle;

    \filldraw[blue!50, fill opacity=0.35, dotted] (-2,1,-1) -- (-2,2,-2) -- (2,2,-2) -- (2,1,-1) -- cycle;
    \filldraw[blue!50, fill opacity=0.35, dotted] (-2,-1,-1) -- (-2,-2,-2) -- (2,-2,-2) -- (2,-1,-1) -- cycle;
    
    \filldraw[blue!50, fill opacity=0.35, dotted] (-2,1,-1) -- (-2,1,2) -- (2,1,2) -- (2,1,-1) -- cycle;
    
    \filldraw[blue!50, fill opacity=0.35, dotted] (-2,-1,-1) -- (-2,-1,2) -- (2,-1,2) -- (2,-1,-1) -- cycle;
    
    \filldraw[blue!50, fill opacity=0.35, dotted] (-2,-1,-1) -- (-2,1,-1) -- (2,1,-1) -- (2,-1,-1) -- cycle;

    \draw[very thick, black] (1,1,1) -- (-1,1,1);
    \draw[very thick, black] (-1,1,-1) -- (-1,1,1);
    \draw[very thick, black] (1,1,1) -- (1,1,-1);
    \draw[very thick, black] (-1,-1,1) -- (1,-1,1);
    \draw[very thick, black] (-1,-1,1) -- (-1,-1,-1);
    \draw[very thick, black] (1,-1,1) -- (1,-1,-1);
    \draw[very thick, black] (-1,-0.22,-1) -- (-1,-1,-1);
    \draw[very thick, black] (-1,1,-1) -- (-1,-0.07,-1);
    \draw[very thick, black] (1,1,-1) -- (1,-1,-1);
    
    \draw[very thick,black] (1,1,1) -- (2,1,2);
    \draw[very thick,black] (1,-1,1) -- (2,-1,2);
    \draw[very thick,black] (-1,1,1) -- (-2,1,2);
    \draw[very thick,black] (-1,-1,1) -- (-2,-1,2);
    \draw[very thick,black] (1,1,-1) -- (1,2,-2);
    \draw[very thick,black] (1,-1,-1) -- (1,-2,-2);
    \draw[very thick,black] (-1,1,-1) -- (-1,2,-2);
    \draw[very thick,black] (-1,-1,-1) -- (-1,-1.8,-1.8);
    \draw[very thick,black] (-1,-1.89,-1.89) -- (-1,-2,-2);
    \node[label={$\color{red}\mathcal{A}_{\textnormal{trop},1}$}] at (0,3.75,2.5) {};
    \node[label={$\color{blue}\mathcal{A}_{\textnormal{trop},2}$}] at (3.75,-0.5,-1.75) {};
    \node[label={$\mathcal{A}_{\textnormal{trop}}$}] at (-0.25,2,0) {};
\end{tikzpicture}
    \caption{Tropical Batyrev--Borisov complete intersection from Example \ref{example:running-trop}}
    \label{fig:trop-ci}
\end{figure}

Following \cite{BB96} and \cite{HZ}, we now describe a combinatorial model for tropical Batyrev--Borisov complete intersections that is an alternative to Theorem \ref{theorem:amoeba_strata_ci} and works directly inside $M_\R$ without relying on the Cayley trick. 

Let $\nabla=\nabla_1+\dots+\nabla_r$ be an irreducible nef partition and $h$ a centred refined triangulating function inducing a triangulation $\mT$ of $\partial \Delta^\vee$. For a simplex $T \in \mT$, we shall denote $T_j\coloneqq T \cap \nabla_j$ and say that a simplex $T \in \mT$ is \emph{transversal} if $T_j \in \mT_j$ is non-empty for all $j=1, \dots, r$. Note that the collection of all the transversal simplices in $\mT$ forms an \emph{order ideal} $\mT^\trans$ of the face lattice $\mT$ (i.e. if $T$ is transversal, then any $T'$ containing $T$ is also transversal). For $T \in \mT^\trans$, we can consider the following polytope inside $\nabla$:
\begin{equation*}
    T_{\nabla}\coloneqq T_1+\dots+T_r,
\end{equation*}
which gives us a simplicial complex $\mT_{\nabla}\coloneqq\{T_{\nabla} \colon  T \in \mT^\trans \}$ and its geometric realisation $|\mT_\nabla| \coloneqq \bigcup_{T \in \mT^\trans} T_\nabla$. 

\begin{proposition}[{\cite[Proposition 2.2]{HZ}}]
    The geometric realisation of $\mT_{\nabla}$ is given by $|\mT_{\nabla}|=\partial\left(r\Delta^{\vee}\right) \cap \nabla$, which is a polyhedral subcomplex of $\partial\nabla$, and its faces can be written as $T_{\nabla}=rT \cap \nabla$. In particular, the face lattice of $T_{\nabla}$ is isomorphic to the poset $\mT^\trans$. 
\end{proposition}
\begin{remark}\label{remark:subcpx-intrinsic}
    This proposition demonstrates that although the definition of $\mT_{\nabla}$ depends on the choice of a triangulating function $h$, its geometric realisation is purely an invariant of the nef partition. 
\end{remark}
Note that we can also repeat the same process with the triangulation $\check{\mT}$ of $\partial\nabla^{\vee}$ to get a polytopal subcomplex $\check{\mT}_{\Delta}$ inside $\partial\Delta$. 

Recall that the cells of the tropical Batyrev--Borisov complete intersection $\mA_{\trop}$ are naturally labelled by the $r$-tuples of simplices $(\overline{T}_1,\dots,\overline{T}_r) \in (\mT_1*0)\times\dots\times(\mT_r*0)$. The following two results tell us which tuples actually correspond to a non-empty cell. 

\begin{proposition}[{\cite[Proposition 2.4]{HZ}}]\label{proposition:bdd_cells}
    The complex of bounded cells of $\mA_{\trop}=\bigcap_{j=1}^{r}{\mA_{\trop,j}}$ is equal to the complex of cells $C_{\overline{T}}$ of $\mA_{\trop,\tot}$ indexed by $\overline{T} \in \mT*0$ for $\overline{T}\coloneqq T*0$ with $T \in \mT^\trans$. 
\end{proposition}

\begin{remark}\label{remark:bdd_cells1}
    In the language of our discussion for general TTCI, this says that a tuple of simplices $(\overline{T}_1,\dots,\overline{T}_r)$ labels a non-empty bounded cell of $\mA_{\trop}$ if and only if there exists a simplex $T \in \mT^\trans$ such that $\overline{T}_j=(T*0)\cap\nabla_j$ for $j=1,\dots,r$. This also allows us to explicitly describe the cells of $\mA_{\trop}$ as intersections of cells of $\mA_{\trop,j}$'s through the equality $C_{\overline{T}}=C_{(\overline{T}_1,\dots,\overline{T}_r)}=\bigcap_{j=1}^{r}{C_{\overline{T}_j}}$. 
\end{remark}

\begin{corollary}\label{corollary:unbdd_cells}
    The unbounded cells of $\mA_{\trop}$ are indexed by tuples $(\overline{T_1},\dots,\overline{T_r}) \in (\mT_1*0)\times \dots \times (\mT_r*0)$ such that all the simplices $\overline{T}_j$ are positive-dimensional, $\overline{T}=\conv(\overline{T}_1,\dots,\overline{T}_r)$ is a simplex in $\mT^\trans*0$ and $0 \notin \overline{T}_j$ for some $j$. 
\end{corollary}
\begin{proof}
    First, we show that in order to label an unbounded non-empty cell of $\mA_{\trop}$, the tuple $(\overline{T}_1,\dots,\overline{T}_r)$ must be of the given form. Suppose that $x \in \mA_{\trop}$ lies inside an unbounded cell $C_{(\overline{T}_1,\dots,\overline{T}_r)}=\bigcap_{j=1}^{r}{C_{\overline{T}_j}}$. Since the cell is unbounded and convex, there must exist some direction $v \in M_{\R}\backslash \{0\}$ such that the ray $x+tv$, $t \in \R_{\geq 0}$, lies fully inside the cell $C_{(\overline{T}_1,\dots,\overline{T}_r)}$. Observe that under such circumstances, we must have $L_{h_j}(x+tv)=l_{\alpha}(x)+t\langle\alpha,v\rangle$ for $\alpha \in \overline{T}_j$ and $L_{h_j}(x+tv)>l_{\alpha}(x)+t\langle\alpha,v\rangle$ for $\alpha \in \nabla_j \backslash \overline{T}_j$. This means that there exist some non-negative constants $k_1$, \dots, $k_r$ such that $\langle \alpha, v \rangle =k_j$ for $\alpha \in \overline{T}_j$ and $\langle \alpha, v \rangle < k_j$ for $\alpha \in \nabla_j \backslash \overline{T}_j$ for all $j=1, \dots, r$. Observe that if $k_j=0$ for all $j$, then all the functions $L_{h_j}$ are constant along the ray, hence $L_h(x)=\sup_j L_{h_j}(x)$ is constant along the ray, which is impossible, since the Legendre transform of a function defined on a full-dimensional convex polytope is necessarily proper. Hence at least one $k_j$ must be positive (this also shows that $0 \notin \overline{T}_j$ for some $j$ is a necessary condition). 

    Consider the opposite ray $x-tv$ and note that at some point, it must leave the cell $C_{(\overline{T}_1,\dots,\overline{T}_r)}$, since for some index $j$ with $k_j>0$, $\alpha \in \overline{T}_j$ and $t$ large enough $l_{\alpha}(x-tv)=l_{\alpha}(x)-t k_j<0 \leq L_{h_j}(x-tv)$, so the ray leaves $C_{\overline{T}_j}$ eventually. Let $t_0$ be the smallest value such that $x-t_0v \notin C_{(\overline{T}_1,\dots,\overline{T}_r)}$. Note that the only way this could happen is that for some $j=1$, \dots, $r$ (possibly for multiple such $j$), we get a new maximiser $l_{\alpha}(x-t_0v)$ with $\alpha \in \nabla_j \backslash \overline{T}_j$ for $L_{h_j}(x-t_0v)$. Therefore, $x-t_0v \in C_{(\overline{T}'_1,\dots,\overline{T}'_r)}$ with $\overline{T}_j \subset \overline{T}'_j$. Observe that if the cell $C_{(\overline{T}'_1,\dots,\overline{T}'_r)}$ is bounded, then we are done, since by Proposition \ref{proposition:bdd_cells}, $\conv(\overline{T}'_1,\dots,\overline{T}'_r)$ is a simplex in $\mT^\trans*0$ and $\conv(\overline{T}_1,\dots,\overline{T}_r) \subset \conv(\overline{T}'_1,\dots,\overline{T}'_r)$ is its subsimplex. If the cell that we land in is still unbounded, we can repeat the process outlined above until we get to a bounded cell (because we make one of the $\overline{T}_j$'s bigger every time, it must eventually terminate). 

    Conversely, suppose that we are given a tuple $(\overline{T}_1,\dots,\overline{T}_r)$ that satisfies the hypotheses of this Proposition, then we need to show that the cell $C_{(\overline{T}_1,\dots,\overline{T}_r)} \subset \mA_{\trop}$ is non-empty. Let $\overline{S}=\conv(0,\overline{T}_1,\dots,\overline{T}_r)$, then $\overline{S}=S * 0 $ for some $S \in \mT^\trans$, so by Proposition \ref{proposition:bdd_cells}, $C_{\overline{S}}$ is a non-empty bounded cell, which can also be indexed by a tuple $(\overline{S}_1,\dots,\overline{S}_r)$ with $\overline{S}_j=\conv(0,\overline{T}_j)$. By considering an appropriate outwards pointing normal of $\partial \Delta^{\vee}$ along $S$, we can find a vector $v \in M_{\R}$ and constants $k_j \geq 0$ such that $k_j=0$ if $\overline{T}_j=\overline{S}_j$, $k_j>0$ if $\overline{T}_j \subsetneq \overline{S}_j$, $\langle \alpha, v \rangle = k_j$ if $\alpha \in \overline{S}_j$ and $\langle \alpha, v \rangle < k_j$ if $\alpha \in \nabla_j \backslash \overline{S}_j$ ($k_j\geq 0$ and the final two inequalities are true thanks to $v$ being an outwards pointing normal to $\nabla_j$ along $S_j$, while vanishing of appropriate $k_j$ can be guaranteed by Lemma \ref{lemma:nef-summands-disjoint}). By the analysis from the first part of the proof, if we pick a point $x \in C_{(\overline{S}_1,\dots,\overline{S}_r)}$ and some $t>0$ large enough, we will have $x+tv \in C_{(\overline{T}_1,\dots,\overline{T}_r)}$. 
\end{proof}
\begin{remark}\label{remark:live_in_amoeba}
    Note that unlike the bounded cells of $\mA_{\trop}$, the unbounded cells do not necessarily coincide with the cells of $\mA_{\trop,\tot}$, but they still lie inside the total tropical hypersurface (for example, the rays in $\mA_{\trop}$ in Figure \ref{fig:trop-ci} are different from the rays of $\mA_{\trop,\tot}$ and lie inside $2$-dimensional strata). 
\end{remark}

Finally, we shall also extend the terminology of this section to the fan $\Sigma$ and call $\sigma \in \Sigma$ \emph{transversal} if the set $\sigma_j(1)\coloneqq \sigma(1) \cap \Sigma_j(1)$ is non-empty for all $j=1,\dots,r$. In particular, for each such $\sigma$, we get cones $\sigma_j \subset \partial \sigma$ that satisfy $\sigma=\conv(\sigma_1,\dots,\sigma_r)=\sigma_1+\dots+\sigma_r$. Denote the poset of all transversal cones as $\Sigma_\trans$ and its geometric realisation inside $N_\R$ obtained by taking the union corresponding cones of $\Sigma$ as $|\Sigma_\trans|$ (note that this is not a fan, since a transversal cone has non-transversal subcones). Similarly to the setting of faces, this will be an upper order ideal in $\Sigma$ (in other words, intersection of two transversal cones is transversal and a cone containing a transversal subcone is itself transversal), and its realisation will be an open subset of $N_\R$, since each transversal cone is contained in an $n$-dimensional transversal cone.   

\subsection{Toric compactifications}\label{section:toric-cpct1}
To wrap up this section, we discuss the toric compactifications of tropical complete intersections and complete intersections inside $(\Cs)^n$. In addition to showing that some particular complete intersections one can consider in Conjecture \ref{conjecture:bb_hms} near the `tropical limit' are smooth, this allows us to show that $Z_\beta$ equipped with a homogeneous toric potential defines a Liouville domain for $\beta \gg 0$, and that the Liouville domain is independent of this auxiliary choice of the potential. In particular, we aim to extend Lemma \ref{lemma:sm_nontailored_ci} to this setting, along the lines of the discussion in \cite[Section 2.3]{Mikh} for the case of hypersurfaces. We will follow the treatment of \cite{Khovanskii77}, along with the more modern approach of \cite{MS} via tropical geometry. We adopt the following notion of a smooth compactification from \cite{ghhps}, along with a generalisation to the orbifold setting:

\begin{definition}\label{definition:sm-tor-cpct}
    An \emph{orbifold smooth toric compactification} of a submanifold $Y \subset M_{\Cs}$ is a smooth orbifold $\overline{Y} \subset \mX_{\check{\Sigma}}$, where $\check{\Sigma}$ is a complete simplicial fan in $M_\R$ and $\mX_{\check{\Sigma}}$ the associated toric stack, such that $Y=\overline{Y} \cap M_{\Cs}$ and the intersection of $\overline{Y}$ with every stratum of the toric boundary is transverse. We say that $Y$ admits a \emph{smooth toric compactification} if, in addition to the above, $\overline{Y}$ avoids all the orbifold strata of $\mX_{\check{\Sigma}}$ (and hence is a smooth variety). 
\end{definition}

\begin{example}\label{example:running-cpct}
    In Example \ref{example:running-trop}, we can read off that $Z_\beta$ is an elliptic curve with 8 punctures from the combinatorics of its tropicalisation $\mA_\trop$. It will admit a smooth toric compactification inside the smooth toric variety $X_{\check{\Sigma}}$ that is a crepant resolution of $X_\nabla$ from Example \ref{example:running-nef}. 
\end{example}

Our first task is to describe what the intersections of a toric compactification $\overline{Z} \subset \mX_{\check{\Sigma}}$ with boundary strata look like for a general smooth toric stack $\mX_{\check{\Sigma}}$ arising from a fan $\check{\Sigma}$ (as before, we equip every ray $\rho \in \check{\Sigma}(1)$ with its primitive generator $b_\rho$ to give it the structure of a stacky fan). In order to do that, we define the \emph{trivial tropicalisation} $\trop(H_\triv)$ of a very affine hypersurface $H \subset M_{\Cs}$ given as the vanishing locus of a Laurent polynomial $f(z)=\sum_{\alpha \in A} c_\alpha z^\alpha$ as the tropical hypersurface inside $M_\R$ associated to the function $L_\triv(x)\coloneqq \max_{\alpha \in A} \{ \langle \alpha,x\rangle \}$. For a complete intersection $Z=\bigcap_{j=1}^r H_j$, the trivial tropicalisation is given as $\trop(Z_\triv) \coloneqq  \trop(H_{1,\triv}) \cap_{st} \dots \cap_{st} \trop(H_{r,\triv})$ (where $\cap_{st}$ denotes the \emph{stable intersection} of tropical varieties \cite[Definition 3.6.5]{MS} that one needs to consider when dealing with non-transverse intersections in tropical geometry). The triviality of such a tropicalisation means that we are taking $h \equiv 0$ in the language of the previous sections, which corresponds to thinking of $Z$ as a \enquote{trivial} one-parameter family of complete intersections that does not depend on $\beta$. 

Suppose that the Laurent polynomial $f_j(z)=\sum_{\alpha \in A_j} c_{\alpha,j} z^\alpha$ has Newton polytope $P_j=\conv(A_j) \subset N_\R$ for $j=1,\dots,r$. Moreover, suppose that for a generic choice of coefficients $c_{\alpha,j}$, the equations $\{f_1(z)=0,\dots,f_r(z)=0\}$ define a non-empty codimension $r$ very affine complete intersection (by \cite[p.41]{Khovanskii1978}, this is equivalent to requiring $\dim(\sum_{i \in I} P_i) \geq |I|$ for all $I \subset \{1,\dots,r\}$, which can also be shown to be equivalent to $\trop(Z_\triv) \neq \emptyset$). Then, for any toric variety $X_{\check{\Sigma}}$, we can use the trivial tropicalisation to detect which boundary strata intersect $\overline{Z} \subset X_{\check{\Sigma}}$:

\begin{proposition}[{\cite[Theorem 6.3.4]{MS}}]\label{proposition:int-strata1}
     The compactification $\overline{Z}$ intersects the torus orbit $O(\sigma)$ for $\sigma \in \check{\Sigma}$ if and only if $\trop(Z_\triv)$ intersects $\relint(\sigma)$.
\end{proposition}

Similarly to the case of hypersurfaces, we can describe a particularly nice compactification of any given complete intersection: for a CI $Z$ whose defining equations have Newton polytopes $P_1$, \dots $P_r$, consider the polytope $P\coloneq P_1+\dots+P_r$ and denote its normal fan as $\check{\Sigma}_0$ (by \cite[Proposition 6.2.13]{CLS}, this is the same as the coarsest common refinement of the normal fans $N(P_1)$, \dots, $N(P_r)$). Then $X_0 \coloneqq X_{\check{\Sigma}_0}$ will be a projective toric variety with an associated polytope $P$, so there is a dimension-preserving poset isomorphism between the boundary strata of $X_0$ and the faces of $P$. We call the compactification of $Z$ inside $X_0$ (denoted $\overline{Z}_0$) its \emph{canonical compactification}. As a simple application of Proposition \ref{proposition:int-strata1}, we obtain a generalisation of \cite[Proposition 2.18]{Mikh}:

\begin{corollary}\label{corollary:can-cpctification}
    The canonical compactification $\overline{Z}_0 \subset X_0$ intersects the orbit $O(F)$ for $F \subset \partial P$ if and only if $F=F_1+\dots+F_r$ for some faces $F_j \subset \partial P_j$ of positive dimension. In particular, $\overline{Z}_0$ avoids all the orbits of dimension less than $r$.
\end{corollary}

To describe the intersections $Z(\sigma) \coloneqq Z \cap O(\sigma)$ in terms of equations, we will need a further assumption on the fan $\check{\Sigma}$ (introduced in \cite{Khovanskii77}).
\begin{definition}\label{definition:suff-fine}
    We say that $\check{\Sigma}$ is \emph{sufficiently fine with respect to $P=P_1+\dots+P_r$} if it is a refinement of the fan $\check{\Sigma}_0=N(P)$
\end{definition}
Note that this condition is equivalent to requiring that the restriction of the piecewise linear function $L_{j,\triv}$ to each cone $\sigma \in \check{\Sigma}$ is linear for $j=1,\dots,r$. 

The significance of this condition for us will be the following: when $\check{\Sigma}$ is sufficiently fine, then for any cone $\sigma \in \check{\Sigma}$, there exists a unique maximal face $P_j^\sigma \subset P_j$ such that $\sigma$ lies inside the normal cone to $P_j$ at $P_j^\sigma$. We pick a lattice point $\alpha^\sigma_j \in A_j \cap P_j^\sigma$ for every such face and then define the \emph{truncation of $f_j$ over $\sigma$} as
    \begin{equation*}
        f^\sigma_j(z) \coloneqq \sum_{\alpha \in A_j \cap P^\sigma_j} c_{\alpha,j} z^{\alpha-\alpha_j^\sigma},
    \end{equation*}
for $j=1,\dots,r$. The character lattice of the torus orbit $O(\sigma)$ is equal $\sigma^\perp \subset M$ and $P^\sigma_j$ is the intersection of $P_j$ with the translate $\alpha^\sigma_j+\sigma^\perp$, so all the characters appearing in $f^\sigma_j$ can be viewed as characters on the orbit $O(\sigma)$. In fact, the truncation is picking out the terms minimising $\langle \alpha, v_\rho \rangle$ for all ray generators $v_\rho$ of $\rho \in \sigma(1)$, which is their order of vanishing along a divisor $D_\rho$ containing $O(\sigma)$, so it is keeping track of the monomials that are dominant in the expression for $f_j(z) \cdot z^{-\alpha^\sigma_j}$ as we approach $O(\sigma)$. This idea motivates the following result.

\begin{proposition}[{\cite[Theorem 2]{Khovanskii77}}]\label{proposition:int-strata2}
    If $\sigma \in \check{\Sigma} \backslash \{0\}$ is a cone such that $Z(\sigma) \neq \emptyset$, then $Z(\sigma)=\mathbb{V}(f^\sigma_1,\dots,f^\sigma_r)$. Moreover, if the cone $\sigma$ is smooth, then the intersection $\overline{Z} \cap O(\sigma)$ will be transverse for a generic choice of coefficients $c_{\alpha,j}$, and so $Z(\sigma)$ is generically smooth. 
\end{proposition}

In particular, note that the definition of stable intersection ensures that $\dim(\sum_{i \in I} P^\sigma_i) \geq |I|$ for all $I \subset \{1,\dots,r\}$, so $Z(\sigma)$ is indeed non-empty. 

\begin{remark}
    Note that Propositions \ref{proposition:int-strata1} and \ref{proposition:int-strata2} extend to the setting of toric stacks $\mX_{\check{\Sigma}}$ when $\check{\Sigma}$ is simplicial, with the understanding that the closure $\overline{Z}$ has to be viewed as a stack rather than a variety (by working in affine charts and passing to quotients by finite groups, as we do in Corollary \ref{corollary:smooth-cpctification}). The condition of $\sigma$ being smooth in Proposition \ref{proposition:int-strata2} can also be relaxed to $\sigma$ being simplicial in that setting.  
\end{remark}

We now turn our attention back to canonical compactifications: note that the polytope $P$ can fail to be simple, which means that $X_P$ might not be an orbifold, so we can not directly use the canonical compactification to demonstrate the existence of orbifold smooth toric compactifications. However, we can always find a refinement of the normal fan $\check{\Sigma}\rightarrow \check{\Sigma}_0$ to a simplicial fan with $\check{\Sigma}(1)=\check{\Sigma}_0(1)$ by \cite[Proposition 11.1.7]{CLS}, which gives us a partial resolution $X_{\check{\Sigma}} \rightarrow X_{\check{\Sigma}_0}$ with $X_{\check{\Sigma}}$ an orbifold, and hence $\mX_{\check{\Sigma}}$ is a smooth Deligne--Mumford stack. One can immediately obtain a generalisation of Lemma \ref{lemma:sm_nontailored_ci} to compactifications inside such stacks $\mX_{\check{\Sigma}}$:

\begin{corollary}\label{corollary:smooth-cpctification}
    Let $Z_\beta$ be a family of complete intersections satisfying the hypotheses of Lemma \ref{lemma:sm_nontailored_ci}. Then $Z_\beta$ admits an orbifold smooth toric compactification for all $\beta>0$ large enough. 
\end{corollary}

\begin{proof}
    Consider any desingularisation $\mX_{\check{\Sigma}}\rightarrow X_{\check{\Sigma}_0}$ as above, and the compactifications $\overline{Z}_\beta$ inside it. As $\check{\Sigma}$ is a refinement of $\check{\Sigma}_0$ by definition, Proposition \ref{proposition:int-strata2} provides equations for the intersections $Z_\beta(\sigma)\coloneqq\overline{Z}_\beta \cap O(\sigma)$. In fact, $Z_\beta(\sigma)$ is the proper transform of the intersection with some orbit $O(F) \subset X_{\check{\Sigma}_0}$ for $F=F_1+\dots+F_r$ inside the canonical compactification from Corollary \ref{corollary:can-cpctification}, and the faces for the truncation are given simply as $P^\sigma_j=F_j$. 
    
    Let $\sigma\in \check{\Sigma}$ be a cone with $Z_\beta(\sigma)\neq \emptyset$, then since $\sigma$ is simplicial, the associated orbit has an affine Zariski neighbourhood of the form $U(\sigma) \cong O(\sigma) \times [\C^{\dim(\sigma}/G]$ for a finite group $G \cong (M \cap \bigoplus_{\rho \in \sigma(1)} \R b_\rho)/ \bigoplus_{\rho \in \sigma(1)} \Z b_\rho$ acting on $\C^{\dim(\sigma)}$ as a subgroup of the dense torus. Hence, to prove the smoothness of the affine stack $\overline{Z}_\beta \cap U(\sigma)$, we can lift it to a subvariety of $(\Cs)^{n-\dim(\sigma)}\times \C^{\dim(\sigma)} \rightarrow (\Cs)^{n-\dim(\sigma)}\times [\C^{\dim(\sigma)}/G]$ given by the same equations, and check its smoothness as a variety. 

    Therefore, without loss of generality, we can assume that $\sigma$ is a smooth cone, and so the orbit $O(\sigma)$ has an affine neighbourhood of the form $U(\sigma) \cong \C^{\dim(\sigma)} \times (\Cs)^{n-\dim(\sigma)}$. Since the restriction of a refined triangulating function $h_j$ to any face $F_j \subset P_j$ is still refined, Lemma \ref{lemma:sm_nontailored_ci} applies to the very affine complete intersection $Z_\beta(\sigma)$, so they are smooth for $\beta$ large enough, while the proof of the Lemma shows a stronger statement that $0$ is a regular value $\overline{f}_\beta\coloneqq(f_{1,\beta}\cdot z^{-\alpha^\sigma_1},\dots,f_{r,\beta}\cdot z^{-\alpha^\sigma_r}) \colon U(\sigma) \rightarrow \C^r$ extending a defining function $f_\beta \colon M_{\Cs} \rightarrow \C^r$ for $Z_\beta$, so we also get that $\overline{Z}_\beta \cap U(\sigma)$ is a smooth variety that intersects $O(\sigma)$ transversely. Since there are finitely many charts covering $\mX_{\check{\Sigma}}$, this shows that $\overline{Z}_\beta$ indeed defines an orbifold smooth toric compactification of $Z_\beta$ for all $\beta \gg0$. 
\end{proof}
Finally, we explain what happens in the setting of Batyrev--Borisov mirror construction that is of particular interest to us: in that case, the polytope $\nabla=\nabla_1+\dots+\nabla_r$ is reflexive and the variety $X_\nabla$ is a Gorenstein Fano variety. As per Remark \ref{remark:resolutions}, the choice of a centred refined triangulating function $\check{h}$ gives us a star-shaped regular triangulation $\check{\mT}$ of $\nabla^\vee$ and hence a desingularisation $X_{\check{\Sigma}} \rightarrow X_\nabla$. By \cite[Theorem 2.2.24]{Batyrev94}, this birational morphism is crepant and projective, and when the vertex set of $\check{\mT}$ is equal to $\nabla^\vee \cap M$, it gives a \emph{maximal projective crepant partial (MPCP) desingularisation} of $X_\nabla$. Therefore, by Corollary \ref{corollary:smooth-cpctification}, the choices of centred refined triangulating functions $h$ and $\check{h}$ indeed give us a one-parameter family $\overline{Z}_\beta \subset \mX_{\check{\Sigma}}$ of smooth orbifolds for sufficiently large $\beta$.  
\begin{remark}\label{remark: mpcs condition}
    If one wishes to avoid working with orbifolds (and hence get smooth compactifications rather than orbifold smooth ones), it is possible to generalise Batyrev's definition for hypersurfaces to Batyrev--Borisov complete intersections $Z \subset M_{\Cs}$ and ask that the \emph{MPCS condition} holds for $Z$, i.e. there exists an MPCP desingularisation $\check{\Sigma} \rightarrow \check{\Sigma}_0$ such that the compactification $\overline{Z}$ of $Z$ inside $\mX_{\check{\Sigma}}$ avoids all the orbifold strata in $\partial \mX_{\check{\Sigma}}$. In the light of Proposition \ref{proposition:int-strata1}, this is equivalent to asking that every cone $\sigma \in \check{\Sigma}$ contained in the relative interior of the normal cone of a `mixed' face $F=F_1+\dots+F_r \subset \partial \nabla$ with $F_j \subset \partial \nabla_j \backslash\{0\}$ and $\dim(F_j) \geq 1$ is smooth. Even more concretely, by Corollary \ref{corollary:can-cpctification}, it is sufficient for the variety $X_{\check{\Sigma}}$ to be smooth along torus orbits of dimension at least $r$ in order for the MPCS condition to hold.
\end{remark}
\section{Adapted K{\"a}hler potentials} \label{section:potentials}
\enlargethispage{3\baselineskip}
Our end goal is to relate the skeleton of a complete intersection with the \emph{FLTZ skeleton} $\eL_\Sigma=\bigcup_{\sigma \in \Sigma} \sigma^\perp \times \sigma \subset M_{S^1} \times N_\R \cong T^*M_{S^1}$. In order to make the desired comparison, we need to construct an identification $\Phi\colon M_{\Cs}\rightarrow T^*M_{S^1}$. Since $M_{\Cs}=M_{S^1}\times M_{\R}$, the easiest way to proceed is to pick an identification $M_\R \rightarrow N_\R$ between the pair of dual vector spaces. In the original work \cite{GS22}, the authors used a conveniently chosen inner product to accomplish that goal. However, existence of such an inner product is a restrictive extra condition (see Remark \ref{remark:pc_case_easy}), which was removed in \cite{Zhou} by working with a wider class of \emph{adapted} homogeneous degree $2$ potentials $\varphi \colon M_\R \rightarrow \R_{\geq 0}$. This section is devoted to introducing the notion of \emph{strongly adapted potentials} and explaining their importance for the main results of the paper. Our approach is also similar to the constructions that appear in \cite[Section 2]{han}. 

We construct a Liouville structure on $M_{\Cs}$ from an adapted potential as follows: given a smooth strictly convex function $\varphi \colon M_\R \rightarrow \R_{\geq 0}$, we can pull it back along $\Log\colon M_{\Cs}\rightarrow M_{\R}$ to a function on $M_{\Cs}$ (we shall abuse notation and also denote this pullback by $\varphi$) and use it as a Kähler potential there: let
\begin{equation*}
    \lambda_{\varphi}\coloneqq-d^{c}\varphi,\ \omega_{\varphi}\coloneqq-dd^{c}\varphi \textnormal{ and } g_{\varphi}(X,Y)\coloneqq\omega_{\varphi}(X,JY),
\end{equation*}
where $J$ is the standard complex structure on $M_{\Cs}$ given by $J\partial_{\rho_j}=\partial_{\theta_j}$, $J\partial_{\theta_j}=-\partial_{\rho_j}$. In polar coordinates, we can provide explicit expressions for all these:
\begin{equation*}
    \begin{split}
        \lambda_{\varphi}(\rho,\theta)&=\sum_{j=1}^{n}{\partial_j \varphi(\rho)d\theta_j}, \\
        \omega_{\varphi}(\rho,\theta)&=\sum_{1 \leq j < k \leq n}{\partial_{jk} \varphi(\rho)d\theta_j\wedge d\rho_k}, \\
        g_{\varphi}(\rho,\theta)&=\sum_{1 \leq j , k \leq n}{\partial_{jk}\varphi(\rho)(d\theta_j\otimes d\theta_k+d\rho_j\otimes d\rho_k)}.\\
    \end{split}
\end{equation*}
In such a set-up, we have the following simple description of the desired identification: observe that the differential $d\varphi$ can be viewed as a map from $M_\R$ to $N_\R$ via $u \in M_\R \mapsto d\varphi(u) \in T_{\varphi(u)}N_\R \cong N_\R$. We shall denote this map as $\Phi=d\varphi$, and also use the same notation for the complexification 
\begin{equation*}
    \Phi=\textnormal{id}\times d\varphi \colon M_{\Cs}=M_{S^1}\times M_{\R} \rightarrow M_{S^1} \times N_{\R}=T^*M_{S^1}.
\end{equation*}
\begin{proposition}[\cite{Zhou}]\label{proposition:capital-phi-iso}
For any smooth strictly convex potential $\varphi \colon M_\R \rightarrow \R_{\geq 0}$ homogeneous of degree $2$, the map $\Phi \colon (M_{\Cs},\lambda_\varphi) \rightarrow (T^*M_{S^1},\lambda_0)$ is a strict isomorphism of Liouville manifolds, where $\lambda_0$ is the standard Liouville form on $T^*M_{S^1}$. 
\end{proposition}

\begin{remark}\label{remark:scaling-weird}
    Due to the rescaling present in the process of tropicalisation, we usually work in coordinates $(u,\theta)$ rather than $(\rho,\theta)$. Therefore, it could seem more natural to consider a one-parameter family of potentials $\varphi \circ \Log_\beta$ instead of a single potential $\varphi \circ \Log$, which comes at an extra cost of having to deal with a family of Kähler forms. We shall mostly work with potentials $\varphi$ that are homogeneous of degree $2$, which means that these two differ by a factor $\beta^{-2}$ and it is straightforward to convert between the expressions for $\lambda$, $\omega$ and $g$ above and the expressions using $du_j=\beta^{-1} d\rho_j$. 
\end{remark}

Recall that a polytope $P$ is said to be \emph{simple} if each of its vertices is an intersection of precisely $n$ facets (this is clearly equivalent to the normal fan $N(P)$ being simplicial). Unlike the aforementioned papers, we restrict our attention to simple polytopes for technical reasons. This degree of generality suffices for our purposes, since the polytopes that we obtain via tropical geometry in the main part of the paper turn out to be simple anyway, because if $\mA_\trop$ is a tropical hypersurface dual to a triangulation, all the complementary regions and cells will be simple polytopes. This gives us an advantage of $P$ and all of its faces being manifolds with corners, which allows us to operate within the convenient framework developed by Joyce in \cite{joyce}. 

Unless specified otherwise, $P$ will be a full-dimensional simple polytope in $M_{\R}$ containing the origin in its interior. Denote its normal fan $\Sigma$ and for a cone $\sigma \in \Sigma$, let $F_{\sigma}$ be its dual face (with the convention that $P=F_0$). Pick a generator $v_{\rho}$ for every ray $\rho \in \Sigma(1)$ so that the polytope is given as $P=\{ x \in M_{\R} \colon \langle v_{\rho},x\rangle \leq 1 \textnormal{ for all } \rho \in \Sigma(1)\}$.
\begin{definition}\label{definition:bar_subdiv}
    The \emph{barycentric subdivision of a cone $\sigma$} is the Minkowski sum of segments from $0$ to $v_{\rho}$ for all $\rho \in \sigma(1)$; this will be a cube $[0,1]^{\dim(\sigma)}$ and we shall denote it $\sigma^{\ba}$. The \emph{barycentric subdivision of $\Sigma$} is defined as the union of barycentric subdivisions of all its cones, denoted $\Sigma^{\ba}$. 
\end{definition}
\begin{figure}[ht]
    \centering
    \begin{tikzpicture}[scale=2, every node/.style={font=\small}]

  \draw[->] (0,0) -- (-1,0) node[left] {$v_{\rho_4}$};
  \draw[->] (0,0) -- (0,-1) node[below] {$v_{\rho_3}$};
  \draw[->] (0,0) -- (1,1) node[right] {$v_{\rho_2}$};
  \draw[->] (0,0) -- (0,1) node[above left] {$v_{\rho_1}$};

  \fill[lightgray, opacity=0.3] (0,0) -- (-1,0) -- (-1,-1)  -- (0,-1) -- cycle;
  \draw[dashed] (0,0) -- (-1,0) -- (-1,-1) -- (0,-1) -- cycle;

  \fill[lightgray, opacity=0.3] (0,0) -- (0,-1) -- (1,0) -- (1,1) -- cycle;
  \draw[dashed] (0,0) -- (0,-1) -- (1,0) -- (1,1) -- cycle;

  \fill[lightgray, opacity=0.3] (0,0) -- (1,1) -- (1,2) -- (0,1) -- cycle;
  \draw[dashed] (0,0) -- (1,1) -- (1,2) -- (0,1) -- cycle;

  \fill[lightgray, opacity=0.3] (0,0) -- (0,1) -- (-1,1) -- (-1,0) -- cycle;
  \draw[dashed] (0,0) -- (0,1) -- (-1,1) -- (-1,0) -- cycle;
  \node at (-0.5,-0.5){$\sigma_3^{\ba}$};
  \node at (0.5,0){$\sigma_2^{\ba}$};
  \node at (0.5,1){$\sigma_1^{\ba}$};
  \node at (-0.5,0.5){$\sigma_4^{\ba}$};
  \node at (1.25,0) {$\color{black}\Sigma^\ba$};
\end{tikzpicture}

    \caption{Barycentric subdivision of the fan of the first Hirzerbruch surface}
    \label{fig:bar-subdiv}
\end{figure}
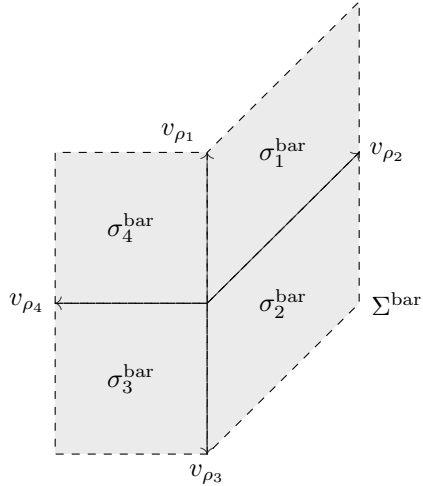
Note that we have $\sigma^{\ba} \cap \tau^{\ba} = (\sigma \cap \tau)^{\ba}$ for any pair of cones in $\Sigma$, so the barycentric subdivision of $\Sigma$ is a complex of cubes glued along cubes. In fact, $P$ admits a similar decomposition: let $x_0=0$ be the origin and for positive-dimensional cones $\sigma \in \Sigma$, pick an arbitrary point $x_{\sigma} \in \relint(F_{\sigma})$. Then we can observe that $P(\sigma)\coloneqq\conv(\{x_{\tau} \colon \tau \subset \sigma\})$ will be a cube of the same dimension as $\sigma^\ba$. Clearly, this gives a decomposition of $P$ into a complex of cubes that is combinatorially isomorphic to the decomposition of $\Sigma^{\ba}$. By inductively gluing homeomorphisms between the cubes in a compatible way, one obtains:

\begin{lemma}\label{lemma:bar_subdiv_homeo}
    There exists a homeomorphism $\psi \colon\Sigma^\ba \rightarrow P$ that induces homeomorphisms $\psi_\sigma \colon \sigma^\ba \rightarrow P(\sigma)$ for all $\sigma \in \Sigma$. 
\end{lemma}

For a cone $\sigma$, we shall denote the part of $\Sigma^\ba$ identified with the face $F_\sigma \subset \partial P$ under this homeomorphism as $\Sigma^\ba_\sigma$. More explicitly, we can describe $\Sigma^\ba_\sigma$ as the set of points in the boundary $\partial \Sigma^\ba$ that contain the term $\sum_{\rho \in \sigma}v_\rho$ in their Minkowski sum representation from Definition \ref{definition:bar_subdiv}. 

\begin{remark}\label{remark:pc_case_easy}
    In the paper \cite{GS22}, the authors exclusively consider potentials of the form $\varphi(x)=\lVert x \rVert ^2$ coming from an inner product. However, this is only possible under an extra assumption that the polytope $P$ is \emph{perfectly centred} (see their Definitions 6.2.1 and 6.22). However, there are polytopes that are known to not be perfectly centred, for example, \cite{non-pc-polytopes} provides an example of a $4$-dimensional polytope that is not even combinatorially isomorphic to a perfectly centred polytope, so can not arise from a perfectly centred fan\footnote{We thank Andrew Hanlon for pointing out this reference.}. Among other things, we sketch a proof of the fact that if we start from a perfectly centred polytope, then the quadratic potential satisfies all the conditions that we impose. Similarly, one can observe that a lot of the arguments we give in Section \ref{section:combinatorics} simplify in the perfectly centred case (which already encompasses a wide variety of interesting examples, such as complete intersections in products of weighted projective spaces). Therefore, one of the main goals of this section is essentially to make our results independent of assumptions on perfectly centredness. 
\end{remark}

\subsection{Strongly adapted potentials}\label{section:strongly-adapted-potentials}
In this subsection, we give a definition of what it means for a potential to be strongly adapted. We start out by discussing the situation for general embeddings $P \hookrightarrow N_\R$ and later specialise to the case when the embedding is given by the differential $\Phi=D\varphi$ of a convex function $\varphi$.

\begin{definition}\label{definition:w_standard_embedding}
    We say that an embedding $G \colon P \hookrightarrow N_{\R}$ is \emph{standard with respect to $\Sigma$} if the following conditions hold:
    \begin{enumerate}
        \item For every $\sigma \in \Sigma$, we have $G(F_{\sigma}) \subset \st(\sigma)$.
        \item For every pair of cones $\sigma \subset \tau$, the intersection $G(\relint(F_{\sigma})) \cap \tau$ is non-empty and transverse. In addition, we require that $G(\relint(F_\sigma))\cap \sigma$ is a single point. 
    \end{enumerate}
\end{definition}

Note that in the language of the second condition, the first condition can be equivalently rephrased as $G(\relint(F_\sigma)) \cap \tau=\emptyset$ for $\sigma \not\subset \tau$.To get better control of what the smooth topology of intersections $G(F_{\sigma}) \cap \tau$ looks like, we will need a more restrictive notion of being standard. In particular, we will fix a particularly convenient standard embedding of $P$ into $N_{\R}$ and work within its isotopy class. 

Recall that the barycentric subdivision of $\Sigma$ defined above provides us with a space $\Sigma^\ba \subset N_\R$ that is homeomorphic to $P$ (Lemma \ref{lemma:bar_subdiv_homeo}). There is a slight issue that $\Sigma^\ba$ might fail to be a manifold with corners, due to the presence of bends in directions normal to $\sigma$ along $\sigma$ (as it is the case on Diagram \ref{fig:bar-subdiv}), so we can not ask $G(P)$ to be isotopic to it. However, in Appendix \ref{section:appendix-smoothing}, we explain that $\Sigma^\ba$ admits an essentially unique smoothing relative to $\Sigma$, denoted as $\widetilde{\Sigma}^\ba$, that is a manifold with boundary (obtained by smoothing its boundary $\partial \Sigma^\ba$). We also show that manifolds with corners intersecting the cones of $\Sigma$ transversely, such as $G(P)$ for a standard embedding $G$, admit unique smoothings relative to $\Sigma$. Therefore, the following definition makes sense:

\begin{definition}\label{definition:standard_embedding}
    We say that an embedding $G \colon P \hookrightarrow N_{\R}$ is \emph{strongly standard with respect to $\Sigma$} if it is standard with respect to $\Sigma$ and there exists an ambient isotopy $\{H_t \}_{t \in [0,1]}$ that preserves $\Sigma$ and deforms the smoothing $\widetilde{G(P)}$ to a smoothing $\widetilde{\Sigma}^\ba$.
\end{definition}

\begin{definition}\label{definition:adapted_potential}
    We call a function $\varphi \colon M_{\R} \rightarrow \R_{\geq 0}$ a \emph{potential (strongly) adapted to $P$} if it is strictly convex, homogeneous of degree $2$, smooth away from the origin and its differential $\Phi\coloneqq D\varphi \colon M_{\R} \rightarrow N_{\R}$ gives us an embedding of $P$ that is (strongly) standard with respect to $\Sigma$. 
\end{definition}

\begin{remark}\label{remark:off_centred_potentials}
    Note that the origin does not really play a special role among the points in the interior of $P$, other than that we choose to scale from it. Therefore, for any $c \in \inte(P)$, we will say that $\varphi$ is a \emph{potential (strongly) adapted to $P$ centred at $c$} if $\varphi(x-c)$ is (strongly) adapted to $P-c$. Unless otherwise specified, we will still assume that all the potentials are centred at the origin. 
\end{remark}

The notion of strong adaptedness is necessary for our purposes to get stronger control of the smooth topology of the singular Lagrangians that naturally arise as skeleta of complete intersections. This is becomes very relevant, e.g., when we study the topology of skeleta of open BBCI's in Section \ref{section:combi-ci}.

\begin{remark}\label{remark:non_adapted_guys}
    It is also very natural to ask whether each adapted potential is strongly adapted; at the moment of writing, the answer to this question is unknown. One could hope that an argument along the lines of \cite[Proposition 2.10]{Zhou} could be used to show that the space of adapted potentials is at least connected, which would then imply a positive answer to this question through similar techniques to the ones we use in Section \ref{section:existence-adapted}. However, our condition (1) does not interact well with taking convex combinations, which makes things significantly more complicated. 
\end{remark}

The following lemma relates our definition to the existing notions of an adapted potential by establishing that being adapted implies adaptedness from \cite{Zhou} (the definition given there does not involve our first condition, even though the construction provided there satisfies it) and is essentially equivalent to the notion considered in \cite{han} in the context of subdivisions adapted to a fan (his Proposition 2.37 proves the existence of what we call adapted potentials).  

\begin{lemma}\label{lemma:adapted_criterion}
    A smooth, strictly convex function $\varphi$ homogeneous of degree $2$ is adapted to $P$ if and only if it satisfies $\Phi(F_{\sigma}) \subset \st(\sigma)$ and the function $\varphi$ attains its minimum over $F_{\sigma}$ at some point in $\relint(F_{\sigma})$ for all $\sigma \in \Sigma$.
\end{lemma}

\begin{proof}
    Let $\sigma \in \Sigma$ be a non-zero cone and $F_{\sigma}$ the associated face in $\partial P$, denote the linear subspace of $M_{\R}$ identified with the tangent space of $F_\sigma$ as $L_{\sigma}$. Then for every $u \in F_{\sigma}$, we know that $D_u\Phi(L_\sigma) \cap \spann_{\R}(\sigma) = \{0\}$: suppose that $D_u\Phi(v) \in \spann_{\R}(\sigma)$ for $v \in L_\sigma$, since $\spann_{\R}(\sigma)=\ann(L_{\sigma})$, this implies that $\langle D_u\Phi(v),v\rangle=0$, so by strict convexity of $\varphi$, we indeed must have $v=0$. By counting dimensions, it follows that $D_u\Phi(L_{\sigma})\oplus \spann_{\R}(\sigma)=N_{\R}$ and, therefore, whenever $\Phi(u) \in \tau$ for some cone $\tau$ containing $\sigma$, the intersection $\Phi(\relint(F_{\sigma})) \cap \tau$ is automatically transverse at $\Phi(u)$. 
    
    Similarly to \cite[Proposition 2.13]{Zhou}, we observe that $x_{\sigma}$ is a critical point of $\varphi$ restricted to the affine space containing $F_{\sigma}$ if and only if $\Phi(x_{\sigma}) \in \sigma$. Therefore, the second part of the condition for adaptedness is satisfied if and only if all the intersections $\Phi(\relint(F_{\sigma})) \cap \sigma$ are non-empty. To conclude the proof, we observe that the restriction of $\varphi$ to $\relint(F_{\sigma})$ is still strictly convex, hence $\Phi|_{F_\sigma}\colon F_{\sigma} \rightarrow N_{\R}/\sigma$ is a diffeomorphism from some neighbourhood of $x_{\sigma} \in \relint(F_{\sigma})$ to a neighbourhood of $0 \in N_{\R}/\sigma$, implying that the intersections $\Phi(\relint(F_{\sigma})) \cap \tau$ are non-empty for all $\sigma \subset \tau$ if $\Phi(\relint(F_{\sigma})) \cap \sigma$ is non-empty. 
\end{proof}

The main result of this section is the existence of strongly adapted potentials with respect to any simple polytope:

\begin{proposition}\label{proposition:adapted_existence}
    Let $P$ be a simple polytope inside $M_{\R}$ containing the origin in its interior. Then there exists a potential $\varphi$ strongly adapted to $P$. 
\end{proposition}

However, it might be difficult to directly produce smooth isotopies of $N_\R$ preserving $\Sigma$. To circumvent this issue, we introduce a particular class of piecewise smooth isotopies that can be constructed one cone at a time:

\begin{definition}\label{definition:gluable_data}
    A \emph{stratified pre-diffeomorphism} is a collection of smooth maps $\{H_{\sigma}\}_{\sigma \in \Sigma}$, where $H_{\sigma} \colon F_{\sigma} \times \sigma \rightarrow \sigma$, the map $H_\sigma(p,\cdot)$ is a diffeomorphism $\sigma \rightarrow \sigma$ mapping each boundary cone $\tau \subset \partial \sigma$ to itself for all $p \in F_\sigma$, and whenever $\tau \subset \sigma$ and $(p,x) \in F_{\sigma} \times \tau$, then we have $H_{\sigma}(p,x)=H_{\tau}(p,x)$\footnote{Note that both sides make sense due to canonical inclusions $\tau \hookrightarrow\sigma$ and $F_{\sigma} \hookrightarrow F_{\tau}$.}. A \emph{stratified pre-isotopy} is a one-parameter family of stratified pre-diffeomorphisms $H^t=\{H^t_{\sigma}\}_{\sigma \in \Sigma}$, which can get viewed as the same data, except with faces $F_{\sigma}$ replaced by $F_{\sigma} \times [0,1]$, so that we now have isotopies $H_{\sigma}^t(p,\cdot) \colon  \sigma \times [0,1]\rightarrow \sigma$ for $p \in F_\sigma$. 
\end{definition}

\subsection{Tented potentials}\label{section:tented-potentials}

In this subsection, we review a construction of a class of adapted potentials from \cite{Zhou}, which we will eventually show are strongly adapted. As before, throughout this section, we let $P$ be a simple polytope in $M_{\R}$ with a normal fan $\Sigma$ containing the origin in its interior. 

Recall that a convex function $\varphi$ that is homogeneous of degree $2$ is uniquely specified by a convex set $V=\{x \in M_\R \colon \varphi(x) \leq 1\}$ with smooth boundary that contains the origin in its interior. In fact, we can translate the property of being adapted to a polytope to properties of the sub-level set:

\begin{lemma}\label{lemma:adapt_via_sublevel}
    Suppose that we are given a compact convex set $V \subset M_{\R}$ containing the origin in its interior with a smooth boundary $S\coloneqq\partial V$ divided into regions $S_\sigma\coloneqq S \cap \cone(F_\sigma)$ for $F_\sigma \subset \partial P$ a face dual to $\sigma \in \Sigma \backslash 0$. Then the function given by $\varphi(x)\coloneqq \inf \{ \lambda>0 \colon \frac{1}{\sqrt{\lambda}}x \in V \}$ is a potential adapted to $P$ if and only if $V$ satisfies the following conditions:
    \begin{enumerate}
        \item[(1)] the boundary $S$ has everywhere positive curvature, i.e. the set $V$ is strictly convex;
        \item[(2)] for every non-zero $\sigma \in \Sigma$ and every $x \in S_\sigma$, the outwards pointing conormal to $S$ at $x$ lies in $\st(\sigma)$;
        \item[(3)] for every non-zero $\sigma \in \Sigma$, there is some $x_\sigma \in \relint(S_\sigma)$ such that the outwards pointing conormal to $S$ lies in $\sigma$. 
    \end{enumerate}
\end{lemma}

\begin{proof}
    The proof is essentially identical to the argument given in \cite[Proposition 2.13]{Zhou}.
\end{proof}

Now, we shall describe the construction of a level set for an adapted potential in steps, as outlined in \cite[Proposition 2.10]{Zhou}. We start from the piecewise linear defining function for $P$, denoted as $\phi_0$ (if $P$ is given by inequalities $\langle \nu, x \rangle \leq 1$ for some set of ray generators $\nu \in \Sigma(1)$, then $\phi_0(x)=\max_{\nu \in \Sigma(1)} \langle \nu, x \rangle$). For convenience, denote $B\coloneqq\prod_{F \subset \partial P} \relint(F)$, so that a point $u \in B$ determines a point $u_F \in \relint(F)$ for each face of $P$. Observe that a point in $u \in B$ induces a barycentric subdivision $\mP=\mP_u$ of $\partial P$ with simplices $S_C\coloneqq\conv(u_F \colon F \in C)$ indexed by chains $C$ in the face poset $\mS$ of $P$. We also denote by $\mS^m$ the set of maximal chains, corresponding to maximal simplices in the subdivision. 

As the first input of our construction, we fix a point $u \in B$, which gives us a barycentric subdivision $\mP$. Then, we pick a parameter $\varepsilon_1>0$ and define $\phi_1$ to be a piecewise linear function that is linear on $\cone(S_C)$ for any simplex $S_C \in \mP$ and satisfies $\phi_1(u_F)=1-\dim(F)\varepsilon_1$ for all faces $F$. Note that for $\varepsilon_1$ small enough, this function will be convex, $\varepsilon_1$-close to $\phi_0$ and its sub-level set will be $V_1=\conv(u'_F \colon F \subset \partial P)$ for $u'_F\coloneqq\frac{1}{1-\dim(F)\varepsilon_1}u_F$. Moreover, if we denote by $\nu_C$ the linear function obtained by restricting $\phi_1$ to a maximal simplex $S_C$ indexed by $C \in \mS^m$, we get an explicit description $\phi_1(x)=\max_{C \in \mS^m} \langle \nu_C, x \rangle$. These functions are specifically designed to attain their minima over each $F$ at $u_F \in \relint(F)$.

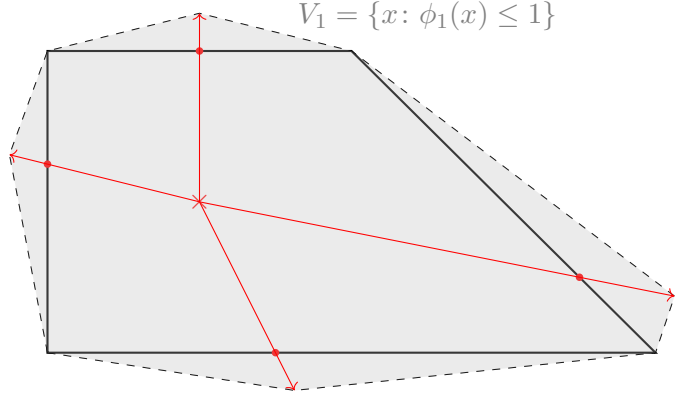
\begin{figure}[ht]
    \centering
\begin{tikzpicture}[scale=2]
    \coordinate (A) at (-1,-1);
    \coordinate (B) at (3,-1);
    \coordinate (C) at (1,1);
    \coordinate (D) at (-1,1);
    \coordinate (O) at (0,0);
    
    \draw[thick] (A) -- (B) -- (C) -- (D) -- cycle;
    
    \coordinate (MAB) at (0.5,-1); 
    \coordinate (MBC) at (2.5,-0.5); 
    \coordinate (MCD) at (0,1); 
    \coordinate (MDA) at (-1,0.25); 

    \node[fill,circle,inner sep=1pt, color=red] at (MAB) {};
    \node[fill,circle,inner sep=1pt, color=red] at (MBC) {};
    \node[fill,circle,inner sep=1pt, color=red] at (MCD) {};
    \node[fill,circle,inner sep=1pt, color=red] at (MDA) {};
    \node[color=red] at (O) {$\times$};

    \coordinate (VAB) at ($(MAB)+1/4*(MAB)-1/4*(O)$);
    \coordinate (VBC) at ($(MBC)+1/4*(MBC)-1/4*(O)$);
    \coordinate (VCD) at ($(MCD)+1/4*(MCD)-1/4*(O)$);
    \coordinate (VDA) at ($(MDA)+1/4*(MDA)-1/4*(O)$);
    
    \draw[dashed] (A) -- (VAB) -- (B) -- (VBC) -- (C) -- (VCD) -- (D) -- (VDA) -- cycle;
    \fill[lightgray, opacity=0.3] (A) -- (VAB) -- (B) -- (VBC) -- (C) -- (VCD) -- (D) -- (VDA) -- cycle;

    \draw[<-, color=red] (VAB) -- (O);
    \draw[<-, color=red] (VBC) -- (O);
    \draw[<-, color=red] (VCD) -- (O);
    \draw[<-, color=red] (VDA) -- (O);

    \node[label={$\color{black!50}V_{1}=\{x \colon \phi_1(x) \leq 1\}$}] at (1.5,1) {};
\end{tikzpicture}
    \caption{Level set of a tented potential $\phi_1$ before smoothing and convexifying for $P$ dual to the fan from Figure \ref{fig:bar-subdiv} (tent points are in red and the \enquote{tents} have height $\varepsilon_1$)}
    \label{fig:level-set}
\end{figure}

Second, we want to consider a smoothing of $\phi_1$. It is possible to just smooth by convolution, but to get better control on the directions in which the function is strictly convex, we describe an alternative construction, inspired by \cite[Section 4.3]{borshevar}. Pick a smooth convex cut-off function $q\colon \R \rightarrow \R$ that satisfies $q(x)=0$ if $x \leq 0$, $q(1)=1$ and $q$ is increasing and strictly convex on $[0, +\infty)$. Then we can define a family of functions $q_{\varepsilon}(x)\coloneqq\varepsilon^{-1}q(x-1+\varepsilon)$ for all $\varepsilon>0$ (one can think of $q_\varepsilon$ as a strictly convex smoothing of $\varepsilon^{-1}\cdot\max\{0,x-1+\varepsilon\}$). Then, given a sufficiently small parameter $\varepsilon_2>0$, we set $\phi_2(x)\coloneqq\sum_{C \in \mS^m}q_{\varepsilon_2}\left( \langle \nu_C, x \rangle\right)$. By construction, this function is smooth and convex, so the set $V_2\coloneqq\{ x \in M_{\R} \colon \phi_2(x) \leq 1\}$ is convex with smooth boundary. We can also note that if $\phi_2(x)=1$, then at least one of the summands $q_\varepsilon(\langle \nu_C,x\rangle)$ must be non-negative, which means that $\phi_1(x)=\max_C \langle \nu_C,x\rangle \geq 1-\varepsilon_2$, but all the summands must be smaller than or equal to $1$, so $\phi_1(x) \leq 1$. Therefore, we have $(1-\varepsilon_2)V_1 \subset V_2 \subset V_1$, so even though the function $\phi_2$ might not globally be close to $\phi_1$, it still gives us a smoothing of the sub-level set that is $\varepsilon_2$-close to $V_2$. 

Finally, we need to ensure strict convexity. In order to attain that, we pick a strictly convex function $\psi \colon M_{\R} \rightarrow \R$ whose $C^{\infty}$-norm on some large compact set containing $P$ is bounded (in practice, we shall be content with $\psi(x)=\lVert x-x_0 \rVert^2$ for some $x_0 \in M_{\R}$ and some inner product). After that, we pick $\varepsilon_3>0$ small and consider $\phi_3(x)\coloneqq\phi_2(x)+\varepsilon_3\psi(x)$. This is the function we use to define a sub-level set $V\coloneqq\{ x \in M_{\R} \colon \phi_3(x) \leq 1\}$ (which will, by construction, be $O(\varepsilon_3)$-close to $V_2$) and hence also a homogeneous function $\varphi\coloneqq\varphi(u,\psi,\varepsilon)$ for $u=\{u_F\}_{F \subset \partial P} \in B$, $\psi$ a strictly convex function and $\varepsilon=(\varepsilon_1, \varepsilon_2, \varepsilon_3) \in \R_{>0}^3$ small. We call the functions of such form \emph{potentials tented at $u$}; dependence on the other parameters will mostly be suppressed from our notation, since their choices usually play a secondary role. By combining some results from the literature, we get a weaker version of Proposition \ref{proposition:adapted_existence}:

\begin{lemma}\label{lemma:almost_ad_existence}
    For any $u \in B$ and any strictly convex $\psi$ that is $C^{\infty}$-bounded on a large compact set containing $P$, the function $\varphi(u,\psi,\varepsilon_1, \varepsilon_2, \varepsilon_3)$ will be a potential adapted to $P$, provided all the values of $\varepsilon_j>0$ are chosen to be sufficiently small. Moreover, the minimum of $\varphi$ over $F$ is attained at some point $\widetilde{u_F} \in \relint(F)$ with $\lVert u_F-\widetilde{u_F}\rVert\leq C(u,\psi)(\varepsilon_2+\varepsilon_3)$ for a positive constant $C(u,\psi)$ depending on the choice of $u$ and $\psi$. 
\end{lemma}

\begin{proof}
    We would like to show that Lemma \ref{lemma:adapt_via_sublevel} applies in the present case. This essentially follows by combining the methods from the proofs of \cite[Proposition 2.10]{Zhou} (our construction is essentially the same, aside from the smoothing procedure) and \cite[Proposition 2.37]{han} (to verify that condition (2) of the Lemma is satisfied).

    The given estimate on minimisers of $\varphi$ follows from the argument in \cite{Zhou} and the above discussion of how closely the sets $V_2$ and $V$ approximate the polytope $V_1$. 
\end{proof}

There is a also slight generalisation of tented potentials along the lines of Remark \ref{remark:off_centred_potentials}, which will be useful later: for $c \in \inte(P)$ say that $\varphi=\varphi(c,u,\psi,\varepsilon)$ is a \emph{$c$-centred potential tented at $u$} if $\varphi(x-c)$ is a potential tented at $u$. 

\subsection{Families of adapted potentials}\label{section:adapted-families}

In the proof of Proposition \ref{proposition:adapted_existence}, it will be crucial to consider families of tented potentials to construct various isotopies. The purpose this subsection sets up a basic framework for such arguments. Throughout the section, $K$ will be a compact manifold with corners parameterising a family of standard embeddings of $P$, which is the same as a smooth map $G\colon K \times P \rightarrow N_\R$ with each $G(p,\cdot)$ an embedding of $P$, so we denote them by $\{ G(p,\cdot)\}_{p \in K}$. In particular, from a family of adapted potentials $\{ \varphi(p,\cdot)\}_{p \in K}$, we can get a family of embeddings $\{\Phi(p,\cdot)\}_{p \in K}$.

We begin by introducing a notion of an embedding $G_\sigma \colon P \hookrightarrow N_\R$ that is only standard over some cone $\sigma \in \Sigma$:

\begin{definition}\label{definition:standard_over_cone}
    We say that an embedding $G_\sigma \colon P \hookrightarrow N_\R$ is \emph{standard over $\sigma$} if
    \begin{enumerate}
        \item The set $G_\sigma^{-1}(\sigma) \cap \partial P$ is covered by $\bigcup_{\tau \subset \sigma}{F_\tau}$.
        \item For every pair of cones $\tau \subset \tau' \subset \sigma$, the intersection $G(\relint(F_\tau)) \cap \tau'$ is non-empty and transverse. Moreover, $G(\relint(F_\tau)) \cap \tau$ is a point for all $\tau \subset \sigma$.
    \end{enumerate}
    We call a smooth, strictly convex, degree $2$ homogeneous function $\varphi_\sigma$ a \emph{potential adapted over $\sigma$} if $\Phi_\sigma\coloneqq D\varphi_\sigma$ is standard over $\sigma$.
\end{definition}

It follows from the definition that $G$ is standard if and only if it is standard over each $\sigma \in \Sigma$, and an embedding that is standard over $\sigma$ will also be standard over any $\tau \subset \sigma$. As above, we shall denote families of embeddings adapted over $\sigma$ as $\{G_\sigma(p,\cdot) \}_{p \in K}$.  For such a family, let $Y_\sigma\coloneqq\{(p,x) \in K \times P \colon G_\sigma(p,x) \in \sigma\} \subset K \times P$ and, more generally, for any pair of cones $\tau \subset \tau' \subset \sigma$, let $Y_{\tau,\tau'}\coloneqq\{(p,x) \in K \times F_\tau \colon G_\sigma(p,x) \in \tau'\} \subset K \times F_\tau$.

\begin{lemma}\label{lemma:cornery_stuff1}
    For any $\sigma \in \Sigma$, $Y_\sigma$ is a manifold with corners of dimension $\dim(K)+\dim(\sigma)$ and, for any pair of cones $\tau \subset \tau' \subset \partial \sigma$, $Y_{\tau,\tau'}$ is a manifold with corners of dimension $\dim(K)+\dim(\tau')-\dim(\tau)$ inside $\partial Y_\sigma$.
\end{lemma}

\begin{proof}
    Since its codomain is an ordinary manifold, the map $(x,p) \in K \times P \mapsto G_\sigma(p,x) \in N_{\R}$ is a smooth map between manifolds with corners $K \times P \rightarrow N_{\R}$ in the sense of \cite{joyce}. The inclusion $\sigma \rightarrow N_{\R}$ is also a smooth map. In fact, the two maps are transverse: suppose that $\tau \subset \sigma$ is a stratum of $\sigma$, $y \in \relint(\tau')$ and $(p,x) \in K \times P$ a point satisfying $G_\sigma(p,x)=y$. Since $G_\sigma(p,\cdot)$ is standard over $\sigma$, this means that $x \in \relint(F_{\tau})$ for some $\tau \subset \tau'$ and that $G_\sigma(p,\relint(F_{\tau}))$ intersects $\tau'$ transversely (and the intersection is non-empty), which is sufficient for the two maps to be transverse, since their codomain does not have a boundary. Therefore, by \cite[Theorem 6.4]{joyce}, the fibre product associated to these maps exists in the category of manifolds with corners and is, by definition, equal to $Y_\sigma$. The statement that $Y_{\tau,\tau'}$ is a manifold with corners follows from an analogous argument. Moreover, \cite[Proposition 6.7]{joyce} tells us that the boundary of $Y_\sigma$ decomposes into $Y_\tau$ for $\tau \subset \partial \sigma$ and the preimage of the boundary of $K\times P$, which means that $Y_{\tau,\tau'} \subset \partial Y_\sigma$. 
\end{proof}

We also have a projection onto the first factor $K \times P \xrightarrow{\pi_1} K$, which can be restricted to a continuous map $Y_\sigma \xrightarrow{\pi} K$. For the following discussion, \emph{submersions of manifolds with corners} should be understood in the sense of \cite[Definition 3.2(iv)]{joyce}.

\begin{lemma}\label{lemma:cornery_stuff2}
    The map $\pi \colon Y_\sigma \rightarrow K$ is a surjective submersion of manifolds with corners.
\end{lemma}

\begin{proof}
    Since $\pi$ is a composition of the projection map and the smooth inclusion, by \cite[Theorem 3.4]{joyce}, it is smooth. It is also clearly surjective. Therefore, it remains to prove that it is a submersion. To begin with, note that the projection $\pi_1 \colon P \times K \rightarrow K$ is a submersion by Example 3.5(b) of loc. cit.. In the rest of the proof, we follow Joyce's notation. 

    Let $p \in S^l(K)$ be a point in a stratum of depth $l$ and $y=(p,x) \in \pi^{-1}(p) \subset Y_\sigma$, suppose that $x \in \relint(F_\tau)$ and $G_\sigma(p,x) \in \relint(\tau')$ for some cones $\tau \subset \tau'$. In the above notation, we then have $y \in Y_{\tau,\tau'}$, with the depth of $y$ inside $Y_{\tau,\tau'}$ also being equal to $l$. As $Y_{\tau,\tau'}$ is a codimension $\dim(\sigma)-\dim(\tau')+\dim(\tau)$ submanifold inside $\partial Y_\sigma$, the depth of $y$ in $Y_{\sigma}$ is equal $k=l+\dim(\sigma)-\dim(\tau')+\dim(\tau)$. Being a submersion at $y$ is thus equivalent to surjectivity of the maps $D_y\pi\colon T_yY_\sigma \rightarrow T_pK$ and $D_y\pi \colon T_yS^k(Y_\sigma) \rightarrow T_pS^l(K)$. 

    Since $\pi_1$ is a submersion, the map $D_y\pi_1\colon T_y(K \times P) \rightarrow T_pK$ is surjective and its kernel is an $n$-dimensional space $V\cong T_xP \subset T_y(K \times P)$. The space $T_yY_{\sigma}$ has dimension $\dim(K)+\dim(\sigma)$ and $Y_{\sigma} \cap \pi_1^{-1}(p)=G_\sigma(p,\cdot)^{-1}(\sigma)$, which means that $\dim(T_yY_\sigma \cap V)=\dim(T_y(Y_{\sigma}\cap \pi_1^{-1}(p)))=\dim(\sigma)$, hence
    \begin{equation*}
        \dim(T_yY_\sigma+V)=\dim(T_yY_\sigma)+\dim(V)-\dim(T_yY_\sigma \cap V)=n+\dim(K).
    \end{equation*}

    So, the spaces $V=\ker(D_y\pi_1)$ and $T_yY_\sigma$ intersect transversely and the restriction of $D_y\pi_1$ to $T_yY_\sigma$ is still surjective, as desired. 

    We can repeat a similar argument for the second point: first, observe that the depth of $y$ in $K\times P$ is $k'=l+\dim(\tau)$, so the map $D_y\pi_1\colon T_yS^{k'}(K \times P) \rightarrow T_pS^l(K)$ is a surjective from $\pi_1$ being a submersion. This means that if we denote the kernel of this map by $W$, one has $\dim(W)=n-\dim(\tau)$ by rank-nullity and $W=T_yS^{\dim(\tau)}(\pi_1^{-1}(p))$. By definition of depth, $\dim(T_yS^k(Y_\sigma))=\dim(K)-l+\dim(\tau')-\dim(\tau)$ and $S^k(Y_\sigma)=S^l(Y_{\tau,\tau'})$ near $y$, thus $S^k(Y_\sigma)\cap S^{\dim(\tau)}(\pi^{-1}_1(p))=\{p\}\times (F_\tau \cap G_\sigma(p,\cdot)^{-1}(\tau'))$ near $y$ and taking tangent spaces gives $\dim(T_yS^k(Y_\sigma) \cap W)=\dim(F_\tau \cap G_\sigma(p,\cdot)^{-1}(\tau'))=\dim(\tau')-\dim(\tau)$ (by standardness of $G_\sigma(p,\cdot)$ over $\sigma$, the intersection is non-empty and transverse). Repeating the same rank-nullity calculation as above, we see that $T_yS^k(Y_\sigma)$ and $W$ intersect transversely inside $T_yS^{k'}(K \times P)$, so the restriction of $D_y\pi_1$ to $T_yS^k(Y_\sigma)$ is still surjective, hence we have shown that $\pi$ is a submersion. 
\end{proof}

We record the following version of Ehresmann's theorem for manifolds with corners, which is essentially folklore (using the same notion of submersion as above): 

\begin{theorem}\label{theorem:ehresmann_with_corners}
    Let $\pi \colon M \rightarrow N$ be a proper submersion of manifolds with corners, then it is a locally trivial fibre bundle of manifolds with corners. 
\end{theorem}

\begin{proof}
    The standard proof of Ehresmann's theorem via lifting vector fields in the base (see, e.g., \cite[Theorem 9.3]{Voisin}) works for this case if we use  \cite[Proposition 5.1]{joyce} instead of the usual statement about the local canonical form of submersions (such as \cite[Theorem 4.29]{Lee2012}) to construct the lift. 
\end{proof}

\begin{corollary}\label{corollary:cornery_stuff1}
    $Y_\sigma\xrightarrow{\pi}K$ is a locally trivial bundle of manifolds with corners.
\end{corollary}

Finally, we state an extension of Lemma \ref{lemma:almost_ad_existence} to a family setting.

\begin{lemma}\label{lemma:family_almost_ad}
    Given a continuously varying family of data $\{(c_p,u_p,\psi_p)\}_{p \in K}$ with $K$ a compact manifold with corners, one can pick $\varepsilon>0$ sufficiently small so that $\varphi(p,\cdot)\coloneqq\varphi(c_p,u_p,\psi_p,\varepsilon)$ is a family of adapted potentials parameterised by $K$. Moreover, the bound from Lemma \ref{lemma:almost_ad_existence} holds with some constant $C=C(\{(c_p,u_p,\psi_p)\}_{p \in K})$. 
\end{lemma}

\begin{proof}
    The proof of Lemma \ref{lemma:almost_ad_existence} goes through verbatim, since all the relevant bounds on $\varepsilon $ vary continuously with $(c,u,\psi)$, so we can use compactness of $K$ to get a uniform bound on the constants $\varepsilon$ over the entire family.
\end{proof}

\subsection{Existence of strongly adapted potentials}\label{section:existence-adapted}

In this section, we strengthen Lemma \ref{lemma:almost_ad_existence} and prove that tented potentials are strongly adapted. The main idea behind the proof is to construct well-behaved families of tented potentials by varying the tenting point (and the convexifying function) appropriately to interpolate between an arbitrary tented potential and one that looks \enquote{standard} near $F_\sigma$, which exploits the framework of Section \ref{section:adapted-families} to construct stratified pre-isotopies. 

\begin{lemma}\label{lemma:potential_deformation_isotopy}
    Let $K$ be a compact manifold with corners, $\sigma \in \Sigma$ a cone and $\{\varphi_\sigma(p,t,\cdot)\}_{(p,t) \in K\times[0,1]}$ a smooth family of potentials (whose centres are also allowed to vary) adapted over $\sigma$ with $\varphi_\sigma(p,0,\cdot)=:\widetilde{\varphi}_\sigma$ independent of $p\in K$. Then this induces a family $\{H^t_\sigma(p,\cdot)\}_{p \in K}$ of ambient isotopies of $\sigma$ between $\widetilde{\Phi}_\sigma(P) \cap \sigma$ and $\Phi_\sigma(p,1,P)\cap \sigma$. 
\end{lemma}

\begin{proof}
    Recall that $Y_\sigma=\{(p,t,x) \in K \times [0,1] \times P \colon \Phi_\sigma(p,t,x) \in \sigma\}$ is a manifold with corners by Lemma \ref{lemma:cornery_stuff1} and also has a structure of a fibre bundle over $K \times [0,1]$ by Corollary \ref{corollary:cornery_stuff1}. 
    
    Moreover, since $\widetilde{\varphi}_\sigma=\varphi_\sigma(p,0,\cdot)$ does not depend on $p$ by assumption, the sub-bundle $Y^0_\sigma$ lying over $K\times\{0\}$ is trivial, so there is a diffeomorphism $G^0_\sigma\colon K \times (\widetilde{\Phi}_\sigma(P) \cap \sigma) \xrightarrow{\sim} Y^0_\sigma$. 

    Now, we proceed analogously to Theorem \ref{theorem:ehresmann_with_corners}: consider a vector field $V=\partial_t$ on the base $K\times[0,1]$. Since $\pi$ is a submersion, there exists a lift $\Hat{V}$ to $Y_\sigma$ such that $\pi_*(\Hat{V})=V$. The time-$t$ flow of this vector field is a diffeomorphism $\Psi_t$ between $Y^0_\sigma$ and $Y^t_\sigma=Y_\sigma \cap \pi^{-1}(K\times \{t\})$, which yields a global trivialisation $G_\sigma\colon K\times [0,1] \times (\widetilde{\Phi}_\sigma(P) \cap \sigma) \rightarrow Y_\sigma$ given by $G_\sigma(p,t,x)=\Psi_t(G_\sigma^0(p,x))$. Note that $Y_\sigma$ can be also viewed as a sub-bundle of the trivial bundle $K\times [0,1] \times \sigma$ (explicitly, the embedding is given by $(p,t,x)\mapsto (p,t,\Phi_\sigma(p,t,x))$), so by extending the vector field $\Hat{V}$ to this bigger bundle, we also obtain a map $G'_\sigma \colon K \times [0,1] \times \sigma \rightarrow K \times [0,1] \times \sigma$ that extends $G_\sigma$. 

    Finally, the isotopy $H^t_\sigma(p,\cdot)\colon \sigma \times [0,1] \rightarrow \sigma$ is given by projecting $G'_\sigma(p,t,\cdot)$ onto the $\sigma$ component. By construction, this is an isotopy of $\sigma$ that takes $\widetilde{\Phi}_\sigma(P) \cap \sigma$ to $\Phi_\sigma(p,1,P) \cap \sigma$. It will map each boundary cone $\tau$ to itself, since the restrictions of $\pi$ to the strata $Y_\tau$ are still submersions and the construction of the lift $\Hat{V}$ sketched out in the proof of Theorem \ref{theorem:ehresmann_with_corners} ensures that $\Hat{V}$ is tangent to these boundary strata. It is also clear that we can extend the vector field to the entirety of the trivial bundle $\sigma \times K \times [0,1]$ while preserving this property, which yields the desired isotopy after integrating. 
\end{proof}

\begin{corollary}\label{corollary:potential_deformation_isotopy}
    Let $\varphi$ be an adapted potential and suppose that we have a family of functions $\{ \varphi(p,t,\cdot)\}_{(p,t) \in P \times [0,1]}$, with the restriction of the family to $F_\sigma$ being denoted as $\{\varphi_\sigma(p,t,\cdot)\}_{(p,t) \in F_\sigma \times [0,1]}$. Suppose that
    \begin{enumerate}
        \item for every cone $\sigma \in \Sigma$, $p \in F_\sigma$ and $t \in [0,1]$, the function $\varphi_\sigma(p,t,\cdot)$ is a potential adapted over $\sigma$ centred at some $c(p,t) \in \inte(P)$;
        \item for all $p \in P$, we have $\varphi(p,0,\cdot)=\varphi$;
        \item for any $\sigma \in \Sigma$ and $p\in F_\sigma$, the image of $\varphi_\sigma(p,1,\cdot)$ intersects $\sigma$ at $\sigma^{\ba}$.
    \end{enumerate}
    Then $\varphi$ is a strongly adapted. 
\end{corollary}

\begin{proof}
    By Lemma \ref{lemma:potential_deformation_isotopy}, the condition i) tells us how to use the data of $\varphi_\sigma(p,t,\cdot)$ to produce a family of isotopies $H^t_\sigma \colon F_\sigma \times [0,1] \times \sigma \rightarrow \sigma$. From conditions ii) and iii), this isotopy will be between $\Phi(P) \cap \sigma$ and $\sigma^\ba$. Moreover, since the data comes from such a bigger family, we can ensure that the lifts of vector fields from Lemma \ref{lemma:potential_deformation_isotopy} agree on overlaps by constructing them inductively, so the collection $\{H^t_\sigma\}_{\sigma \in \Sigma}$ satisfies the required compatibility property to define a stratified pre-isotopy that exhibits $\varphi$ as a strongly adapted potential (by Corollary \ref{corollary:smoothing-isotopies2}, we can use such a pre-isotopy to construct a genuine isotopy). 
\end{proof}

Therefore, the only ingredient that we are missing from the proof is a construction of the deformation data in a form we have described above for the case of tented potentials. The idea of our construction is roughly as follows: suppose that $v$ is a vertex of $P$ and let $g_v$ be the inner product making the unit normals of facets intersecting at $v$ into an orthonormal basis. Then a tented potential centred at a point $v'$ very close to $v$ whose tent points on faces adjacent to $v$ are given by orthogonal projections of $v'$ (and with distance from $v'$ measured using $g_v$ as a convexifying function) will make the image of $P$ inside the maximal cone $\sigma$ dual to $v$ look standard. Such a local model at $v$ may also be smoothly deformed to an arbitrary tented potential $\varphi$ by varying the tenting data and performing small modifications in the tenting construction. Therefore, this idea should give us the inputs for Corollary \ref{corollary:potential_deformation_isotopy} over the full-dimensional cones of $\Sigma$, so we just need to check that it has the desired properties and construct the gluing data over the lower-dimensional cones.
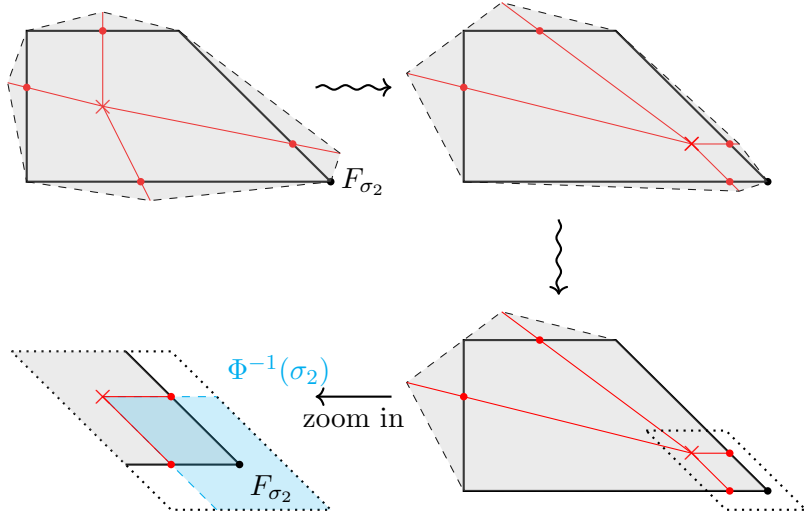
\begin{figure}[htbp]
    \centering
    \begin{tikzpicture}[scale=1]

    \coordinate (A) at (-1,-1);
    \coordinate (B) at (3,-1);
    \coordinate (C) at (1,1);
    \coordinate (D) at (-1,1);
    \coordinate (O) at (0,0);
    \node [fill,circle,inner sep=1pt] at (B){};
    \node [right] at (B) {$F_{\sigma_2}$};
    
    \draw[thick] (A) -- (B) -- (C) -- (D) -- cycle;
    
    \coordinate (MAB) at (0.5,-1); 
    \coordinate (MBC) at (2.5,-0.5); 
    \coordinate (MCD) at (0,1); 
    \coordinate (MDA) at (-1,0.25); 

    \node[fill,circle,inner sep=1pt, color=red] at (MAB) {};
    \node[fill,circle,inner sep=1pt, color=red] at (MBC) {};
    \node[fill,circle,inner sep=1pt, color=red] at (MCD) {};
    \node[fill,circle,inner sep=1pt, color=red] at (MDA) {};
    \node[color=red] at (O) {$\times$};

    \coordinate (VAB) at ($(MAB)+1/4*(MAB)-1/4*(O)$);
    \coordinate (VBC) at ($(MBC)+1/4*(MBC)-1/4*(O)$);
    \coordinate (VCD) at ($(MCD)+1/4*(MCD)-1/4*(O)$);
    \coordinate (VDA) at ($(MDA)+1/4*(MDA)-1/4*(O)$);    
    
    \draw[red] (O) -- (VAB);
    \draw[red] (O) -- (VBC);
    \draw[red] (O) -- (VCD);
    \draw[red] (O) -- (VDA);
    
    \draw[dashed] (A) -- (VAB) -- (B) -- (VBC) -- (C) -- (VCD) -- (D) -- (VDA) -- cycle;
    \fill[lightgray, opacity=0.3] (A) -- (VAB) -- (B) -- (VBC) -- (C) -- (VCD) -- (D) -- (VDA) -- cycle;

    \draw[->, thick, decorate, decoration={snake, amplitude=0.3mm, segment length=3mm}] (2.8,0.25) -- (3.8,0.25);

    \coordinate (S1) at (5.75,0);
    
    \draw[thick] ($(A)+(S1)$) -- ($(B)+(S1)$) -- ($(C)+(S1)$) -- ($(D)+(S1)$) -- cycle;
    \coordinate (O1) at (2,-0.5);
    \coordinate (MAB1) at (2.5,-1);
    \coordinate (MBC1) at (2.5,-0.5);
    \coordinate (MCD1) at (0,1); 
    \coordinate (MDA1) at (-1,0.25);

    \node[fill,circle,inner sep=1pt, color=red] at ($(MAB1)+(S1)$) {};
    \node[fill,circle,inner sep=1pt, color=red] at ($(MBC1)+(S1)$) {};
    \node[fill,circle,inner sep=1pt, color=red] at ($(MCD1)+(S1)$) {};
    \node[fill,circle,inner sep=1pt, color=red] at ($(MDA1)+(S1)$) {};
    \node[color=red] at ($(O1)+(S1)$) {$\times$};
    \node [fill,circle,inner sep=1pt] at ($(B)+(S1)$){};

    \coordinate (VAB1) at ($(MAB1)+1/4*(MAB1)-1/4*(O1)$);
    \coordinate (VBC1) at ($(MBC1)+1/4*(MBC1)-1/4*(O1)$);
    \coordinate (VCD1) at ($(MCD1)+1/4*(MCD1)-1/4*(O1)$);
    \coordinate (VDA1) at ($(MDA1)+1/4*(MDA1)-1/4*(O1)$);

    \draw[red] ($(O1)+(S1)$) -- ($(VAB1)+(S1)$);
    \draw[red] ($(O1)+(S1)$) -- ($(VBC1)+(S1)$);
    \draw[red] ($(O1)+(S1)$) -- ($(VCD1)+(S1)$);
    \draw[red] ($(O1)+(S1)$) -- ($(VDA1)+(S1)$);

    \draw[dashed] ($(A)+(S1)$) -- ($(VAB1)+(S1)$) -- ($(B)+(S1)$) -- ($(VBC1)+(S1)$) -- ($(C)+(S1)$) -- ($(VCD1)+(S1)$) -- ($(D)+(S1)$) -- ($(VDA1)+(S1)$) -- cycle;
    \fill[lightgray, opacity=0.3] ($(A)+(S1)$) -- ($(VAB1)+(S1)$) -- ($(B)+(S1)$) -- ($(VBC1)+(S1)$) -- ($(C)+(S1)$) -- ($(VCD1)+(S1)$) -- ($(D)+(S1)$) -- ($(VDA1)+(S1)$) -- cycle;

    \draw[->, thick, decorate, decoration={snake, amplitude=0.3mm, segment length=3mm}] (6,-1.5) -- (6,-2.5);

    \coordinate (S2) at (5.75,-4.1);
    \draw[thick] ($(A)+(S2)$) -- ($(B)+(S2)$) -- ($(C)+(S2)$) -- ($(D)+(S2)$) -- cycle;
    \draw[dashed] ($(A)+(S2)$) -- ($(B)+(S2)$) -- ($(C)+(S2)$) -- ($(VCD1)+(S2)$) -- ($(D)+(S2)$) -- ($(VDA1)+(S2)$) -- cycle;
    \fill[lightgray, opacity=0.3] ($(A)+(S2)$) -- ($(B)+(S2)$) -- ($(C)+(S2)$) -- ($(VCD1)+(S2)$) -- ($(D)+(S2)$) -- ($(VDA1)+(S2)$) -- cycle;

    \node[fill,circle,inner sep=1pt, color=red] at ($(MAB1)+(S2)$) {};
    \node[fill,circle,inner sep=1pt, color=red] at ($(MBC1)+(S2)$) {};
    \node[fill,circle,inner sep=1pt, color=red] at ($(MCD1)+(S2)$) {};
    \node[fill,circle,inner sep=1pt, color=red] at ($(MDA1)+(S2)$) {};
    \node[color=red] at ($(O1)+(S1)$) {$\times$};\node[color=red] at ($(O1)+(S2)$) {$\times$};
    \node [fill,circle,inner sep=1pt] at ($(B)+(S2)$){};

    \draw[red] ($(O1)+(S2)$) -- ($(MAB1)+(S2)$);
    \draw[red] ($(O1)+(S2)$) -- ($(MBC1)+(S2)$);
    \draw[red] ($(O1)+(S2)$) -- ($(VCD1)+(S2)$);
    \draw[red] ($(O1)+(S2)$) -- ($(VDA1)+(S2)$);

    \draw[thick, dotted] ($(O1)+(S2)+(-0.6,0.3)$) -- ($(O1)+(S2)+(0.45,0.3)$) -- ($(O1)+(S2)+(1.5,-0.75)$) -- ($(O1)+(S2)+(0.45,-0.75)$) -- cycle;

    \draw[->, thick] (3.8,-3.85) -- (2.8,-3.85) node[midway,below] {zoom in};

    \coordinate (O3) at (0,-3.85);
    \draw[thick, dotted] ($(O3)+(-1.2,0.6)$) -- ($(O3)+(0.9,0.6)$) -- ($(O3)+(3,-1.5)$) -- ($(O3)+(0.9,-1.5)$) -- cycle;

    \draw[thick] ($(O3)+(0.9,-0.9)+(-0.6,0)$) -- ($(O3)+(1.8,-0.9)$) -- ($(O3)+(0.3,0.6)$);

    \fill[lightgray, opacity=0.3] ($(O3)+(-1.2,0.6)$) -- ($(O3)+(0.3,0.6)$) -- ($(O3)+(1.8,-0.9)$) -- ($(O3)+(0.3,-0.9)$) --cycle;

    \fill[cyan, opacity=0.2] (O3) -- ($(O3)+(1.5,0)$) node[opacity=1, above right] {$\color{cyan}\Phi^{-1}(\sigma_2)$} -- ($(O3)+(3,-1.5)$) -- ($(O3)+(1.5,-1.5)$) --cycle;
    \draw[dashed, cyan] ($(O3)+(1.5,-1.5)$) -- (O3) -- ($(O3)+(1.5,0)$);

    \node[color=red] at (O3) {$\times$};
    \node[fill,circle,inner sep=1pt, color=red] at ($(O3)+(0.9,0)$) {};
    \node[fill,circle,inner sep=1pt, color=red] at ($(O3)+(0.9,-0.9)$) {};
    \node [fill,circle,inner sep=1pt] at ($(O3)+(1.8,-0.9)$){};
    \node [below right] at ($(O3)+(1.8,-0.9)$) {$F_{\sigma_2}$};

    \draw[red] (O3) -- ($(O3)+(0.9,0)$);
    \draw[red] (O3) -- ($(O3)+(0.9,-0.9)$);
\end{tikzpicture}
    \caption{Deformations of the piecewise linear functions (from the setting associated to Fig. \ref{fig:bar-subdiv} and \ref{fig:level-set}) to make $\varphi$ look standard near the vertex $v=F_{\sigma_2}$ in Lemma \ref{lemma:pre_isotopy_exists}}
    \label{fig:tent-deformation}
\end{figure}
\begin{lemma}\label{lemma:pre_isotopy_exists}
    Let $\varphi$ be a tented potential that is adapted to $P$. Then, provided the values of parameters $\varepsilon=(\varepsilon_1,\varepsilon_2,\varepsilon_3)$ are small enough, there exists a pre-isotopy $\{H^t_{\sigma}\}_{\sigma \in \Sigma}$ taking $\Phi(P) \cap \sigma$ to $\sigma^{\ba}$. 
\end{lemma}

\begin{proof}   
     Denote the starting adapted potential by $\varphi=\varphi(u,\psi,\varepsilon)$ for $u=\{u_F\}_{F \subset \partial P} \in B$, $\psi$ some strictly convex function and $\varepsilon \in \R_{>0}^3$. 

    Pick some small constant $\varepsilon_0>0$ and for each $\sigma \in \Sigma$, denote by $R_{\sigma}$ the set of points $x \in P$ such that $\langle \nu, x \rangle \geq 1-\varepsilon_0$ holds if and only if $\nu \in \sigma(1) \subset \Sigma(1)$. This gives us a decomposition of $P$ into polytopes, the most notable one being $P'\coloneqq R_0=(1-\varepsilon_0)P$. We denote its faces $F'_{\sigma}$ and for a point in $u \in B$, we can get a point $u'\coloneqq(1-\varepsilon_0)u \in \prod_{\sigma}\relint (F'_{\sigma})$. Also, since $\varepsilon_0$ was chosen to be small enough, we can assume that $R_{\sigma}\cong \relint(F_{\sigma}) \times [0,1]^{\dim(\sigma)}$.

    For every vertex $v$ of $P$, let $g_v$ be the unique inner product on $N_{\R}$ that makes the defining functionals of the facets meeting at $v$ an orthonormal basis. We can also pick an open cover $\mU=\{U_v\}$ of $P$ indexed by its vertices such that every face only intersects the opens associated to its vertices. This open cover admits some smooth partition of unity $\{\lambda_v\}$, which allows us to define a family of inner products $g_x\coloneqq\sum_v \lambda_v(x)g_v$ on $N_{\R}$ indexed by the points $x \in P$ (and a family of dual inner products on $M_{\R}$, which we shall also denote by $g_x$). By construction, if the point $x$ is close enough to some face $F$, then the defining functionals of facets meeting at $F$ are orthonormal with respect to $g_x$. In particular, by making $\varepsilon_0$ small enough, we can ensure that this holds for all $x \in R_{\sigma}$. 

    For a point $x \in \inte(P)$, denote its $g_x$-orthogonal projection onto the affine space containing $F_\sigma$ as $\pi^g_\sigma(x)$. Thanks to smallness of $\varepsilon_0$, we can guarantee that for all $x \in R_\sigma$, one has $\pi^g_\sigma(x) \in \relint(F_\sigma)$. So, we can pick some smooth cut-off functions $\lambda'_\sigma \colon M_\R \rightarrow [0,1]$ that are equal to $1$ on $R_\sigma$ and to zero whenever $\pi^g_\sigma(p) \notin F_\sigma$, and then define $u_\sigma(p)\coloneqq \lambda'_\sigma(p)\pi^g_\sigma(p)+(1-\lambda'_\sigma(p))u_\sigma$, which will always be a point in $\relint(F_\sigma)$. Note that since $F'_\sigma$ lies inside the closure $\overline{R_\tau}$ for any $\tau \subset \sigma$, the functions $\lambda'_\tau$ will be equal to $1$ over it, hence $u_\tau(p)=\pi^g_\tau(p)$ for all $p \in F'_\sigma$. Putting all the points $u_\sigma(p)$ together gives us a point $u(p) \in \prod_\sigma \relint(F_\sigma)$.
    
    In our family indexed by $P'\times [0,1]$, we will now want to interpolate between $\varphi$ and a tented potential $\varphi(p,1,\cdot)$ centred at $p$ arising from the data $(p,u(p),\lVert x-p \rVert^2_{g_p},\varepsilon)$. More formally, for $t\in[0,1]$, denote $p_t\coloneqq tp$, $u(p,t)\coloneqq tu(p)+(1-t)u$ and $\psi(p,t)(x)\coloneqq t\lVert x-p_t \rVert^2_{g_p}+(1-t)\psi(x)$. Note that all the parameter spaces for tented potentials are convex, so these are valid inputs for our tenting construction. Therefore, we obtain a smooth family of tented potentials $\varphi(p,t,\cdot)$ with parameters $(p_t,u(p,t),\psi(p,t),\varepsilon)$. By Lemma \ref{lemma:family_almost_ad}, we can pick the parameters $\varepsilon$ small enough to ensure that all these potentials are adapted to $P$. 

    Now, it is straightforward to check that the first two conditions of Corollary \ref{corollary:potential_deformation_isotopy} are satisfied for this family. However, the last one might not quite be satisfied: while $\varphi(p,1,\cdot)$ being centred at a point $p \in F'_\sigma$ close to $F_\sigma$ and tented at its orthogonal projections over all the neighbouring faces guarantees that the images look almost standard, the contributions from the neighbouring cones not contained in $\sigma$ might cause some small discrepancies. 
    
    We can change that as follows: recall that in the construction of tented potentials in Section \ref{section:tented-potentials}, we started by defining a characteristic function $\phi_1$ for a polytope that is obtained by putting \enquote{tents} over all faces of $P$ by scaling the points $u_F$ outwards from the centre by a factor $(1-\dim(F)\varepsilon_1)^{-1}$. The main purpose of this procedure was to gain control over the minima of the resulting potential $\varphi$ over each face (and, in particular, guarantee that they lie in the faces' relative interiors). 

    We can do things slightly differently by scaling each point $u_\sigma(p) \in \relint(F_\sigma)$ by a factor of $(1-(1-s\lambda'_\sigma(p))\dim(F)\varepsilon_1)^{-1}$ from the central point $p$ for some $s \in [0,1]$. For $s=0$, this recovers the original construction, while for $s=1$, it tells us that for $p \in F'_\sigma$, the points $u_\tau(p)$ for $\tau \subset \sigma$ are not scaled at all, since $\lambda'_\tau(p)=1$. Varying our construction leads to a different family of PL functions $\phi_{1,s}(p,\cdot)$ and, carrying through with the other steps with convexifier $\lVert x-p \rVert^2_{g_p}$, this produces a new family $\varphi_s(p,\cdot)$ parametrised by $p \in P$, $s \in [0,1]$ with $\varphi_0(p,\cdot)=\varphi(p,1,\cdot)$. We can observe that our proof of adaptedness of tented potentials still goes through in this setting: the only thing that requires extra care is that when $s\lambda'_\tau(p)$ becomes close to $1$, we no longer get a unique local minimum of $\phi_{2,s}$ over the face $F_\tau$. However, this can only happen if $\lambda'_\tau(p) \neq 0$, meaning that $\pi^g_\tau(p) \in \relint(F_\tau)$, so our convexifier has a minimum in $\relint(F_\tau)$ and the third auxiliary function $\phi_{3,s}$ that we produce will still have the right properties.

    Therefore, by extending the deformation $\varphi(p,t,\cdot)$ to $t \in [0,2]$ via $\varphi(p,t,\cdot)\coloneqq\varphi_{\chi(t-1)}(p,\cdot)$ for $t \geq 1$, where $\chi \colon \R \rightarrow [0,1]$ is some cut-off function with $\chi(t)=0$ for $t \leq 0$ and $\chi(t)=1$ for $t\geq 1$, we get a family of adapted potentials indexed by $P' \times [0,2]$. Clearly, this still satisfies the conditions i) and ii) of Corollary \ref{corollary:potential_deformation_isotopy}, and we shall now describe how to perform one further modification to fulfil condition iii).
    
    The level set of a potential $\varphi(p,2,\cdot)$ near $F_\sigma$ now looks like the level set of a function $\phi(p,x)=\sum_{\nu \in \sigma(1)}{q_{\varepsilon_2}(\langle \nu, x-p \rangle)}+\varepsilon_3 \lVert x-p \rVert ^2_{g_p}$ (since we have ensured that there is no tenting going on over all facets $F_\nu$ for $\nu \in \sigma(1)$ adjacent to $F_\sigma$, hence they no longer contribute by many different linear terms coming from the subdivision, but just by a single linear form coming from their defining equations). In particular, this means that if we pick a $g_p$-orthonormal basis $w_1$, \dots, $w_n$ for $M_\R$ such that $w_1$, \dots, $w_k$ are dual to the unit ray generators of $\sigma$, one has $\Phi(p,2,\cdot)^{-1}(\sigma) = p+ \cone(w_1,\dots,w_k)$\footnote{The cone $\cone(w_1,\dots,w_k)$ can also be intrinsically described as $\sigma^\vee \cap \ann(\sigma)^\perp$, with $\perp$ taken with respect to $g$.}. In particular, this is exactly the same as for $\lVert x-p \rVert^2_{g_p}$, so by extending the deformation to $t \in [0,3]$ via $\varphi(p,t,x)\coloneqq(1-\chi(t-2))\varphi(p,2,x)+\chi(t-2)\lVert x-p \rVert^2_{g_p}$ for $t \geq 2$, we get a larger family indexed by $P' \times [0,3]$. The potentials defined this way might no longer be adapted, however, when $p \in F_\sigma$, the above observation guarantees that $\Phi(p,t,\cdot)^{-1}(\sigma)=\Phi(p,2,\cdot)^{-1}(\sigma)$ for all $t \in [2,3]$. In particular, since the definition of adaptedness over $\sigma$ only depends on the preimage of $\sigma$, this means that $\varphi(p,t,\cdot)$ remains adapted over $\sigma$, hence this larger family still satisfies condition i) of Corollary \ref{corollary:potential_deformation_isotopy}. The conditions ii) and iii) are satisfied by design, so after introducing some obvious linear re-scalings (so that our family is parametrised by $P \times [0,1]$ rather than $P' \times [0,3]$ and the intersections are isotopic to $\sigma^{\ba}$ rather than $\varepsilon_0\sigma^{\ba}$), we are done by the aforementioned Corollary. 
\end{proof}

\begin{remark}\label{remark:norm_strongly_adapted}
    Note that a simplified version of this construction, where we interpolate by quadratic potentials instead, can be applied in the case of perfectly centred polytopes to show that if $\varphi(x)=\lVert x \rVert^2$ has a minimum over each face $F$ at some unique point in $\relint(F)$, then it is strongly adapted in our sense, which means that the original construction from \cite{GS22} fits within our framework. 
\end{remark}

\begin{proof}[Proof of Proposition \ref{proposition:adapted_existence}] 
    By Lemma \ref{lemma:almost_ad_existence}, tented potentials are adapted to $P$. By Lemma \ref{lemma:pre_isotopy_exists}, we conclude that such potentials are also strongly adapted. 
\end{proof}

\subsection{Potentials adapted to amoeba complements}\label{section:adapted-to-complements}

Suppose that $\mA_{\trop}$ is a smooth tropical hypersurface in $M_{\R}$ that is associated to a centred refined triangulating function $h\colon P \rightarrow \R$, then there is a unique compact complementary region $C_0$ that contains the origin in its interior, so we can consider potentials $\varphi$ adapted to $C_0$. In this subsection, we shall prove a number of estimates that hold for such $\varphi$, focusing on their behaviour on the non-bounded cells. 

For convenience, let $\langle \cdot,\cdot\rangle$ be some inner product on $M_{\R}$. Recall that a function $f\colon M_{\R} \rightarrow \R$ is called \emph{$m$-strongly convex} for $m>0$ if we have 
\begin{equation*}
    f(y) \geq f(x)+\langle \nabla f(x),y-x \rangle + \frac{m}{2}\lVert y-x \rVert^2,
\end{equation*}
for all $x, y \in M_{\R}$ and the gradient taken with respect to the chosen inner product. 

\begin{lemma}\label{lemma:stronk_cvx}
    For any potential $\varphi$ adapted to $C_0$, there exists a constant $m>0$ such that $\varphi$ is $m$-strongly convex.
\end{lemma}

\begin{proof}
    By \cite[Lemma 2.5]{Zhou}, the Hessian of $\varphi$ is $0$-homogeneous, i.e. we have $H_{\varphi}(\lambda x)=H_{\varphi}(x)$ for $\lambda > 0$. Consider the function $G\colon S^{n-1} \rightarrow \R_{> 0}$ that sends a point $x \in M_{\R}$ with $\lVert x \rVert=1$ to the smallest eigenvalue of $H_{\varphi}(x)$. Then $G$ is continuous and positive, hence there exists some constant $m_0>0$ such that $G(x) \geq m_0$ for all $x \in S^{n-1}$. By $0$-homogeneity, we therefore know that the eigenvalues of the Hessian $H_{\varphi}(x)$ are bounded below by $m_0$ for any $x \in M_{\R} \backslash \{0\}$, which implies that $\varphi$ is $m_0$-strongly convex on any convex subset $V \subset M_{\R} \backslash \{0\}$. 
    
    So we just need to get strong $m_1$-convexity over any line containing the origin for some $m_1>0$ (since $\varphi$ will generally fail to be twice differentiable at $0$, the Hessian criterion does not apply), then the function will be $m$-strongly convex for $m=\min\{m_0,m_1\}$. Let $R(t)=tv$ be a line through the origin with $\lVert v \rVert=1$, then we know that $\varphi(tv)=t^2\varphi(v)$ for $t\geq 0$ and $\varphi(tv)=t^2\varphi(-v)$ for $t<0$, so $\varphi|_R$ is $m_1$-strongly convex for $m_1=\min \{ \varphi(v), \varphi(-v) \}$. Therefore, if we take $m_1=\inf \{\varphi(v) \colon \lVert v \rVert =1 \}$, all the restrictions will be $m_1$-strongly convex and we are done. 
\end{proof}

Suppose that $C_{T}$ is an unbounded cell of $\mA_{\trop}$, then it must correspond to a simplex $T \in \mT$ not containing the origin. Consider $\overline{T}=\conv\{0,T\}$, then $C_{\overline{T}}$ is a face of $C_0$, so there exists a unique minimum of $\varphi$ over the face that we denote $u_{\overline{T}} \in \relint(C_{\overline{T}})$. By definition, for a point $u \in C_{\overline{T}}$, we have $\langle \alpha, u \rangle=h(\alpha)$ for all $\alpha \in T$. This simple observation allows us to make the following estimate on distance of points that are close to $C_T$ from $C_{\overline{T}}$:

\begin{lemma}\label{lemma:adapted_bound}
    Suppose that $m>0$ is a constant such that $\varphi$ is $m$-strongly convex. Then there exists a positive constant $c$ such that the inequality
    \begin{equation*}
        \langle d\varphi(u),u-u_{\overline{T}} \rangle \geq m\cdot\lVert u-u_{\overline{T}} \rVert^2  + c\cdot d_{\aff}(u,C_{\overline{T}}),
    \end{equation*}
    holds for all $u \in M_{\R}$ satisfying $l_{\alpha}(u) \geq 0$ for all $\alpha \in T$.                                   
\end{lemma}

\begin{proof}
    From the discussion in Appendix \ref{section:appendix-polyhedra}, the tropical hypersurface $\mA_{\trop}$ has a canonical parameterisation, so the notion of affine distance makes sense. Since $\varphi$ is $m$-convex, we have:
    \begin{equation*}
        \langle d\varphi(u),u-u_{\overline{T}} \rangle \geq \langle d\varphi(u_{\overline{T}}),u-u_{\overline{T}} \rangle + m\lVert u-u_{\overline{T}} \rVert^2,
    \end{equation*}
    so it suffices to relate  $\langle d\varphi(u_{\overline{T}}),u-u_{\overline{T}} \rangle$ to the affine distance of $u$ from $C_{\overline{T}}$. Denote the lattice points in $T$ as $\{\alpha_1,\dots,\alpha_k\}$, then $u_{\overline{T}}$ being a minimum over $C_{\overline{T}}$ is equivalent to $d\varphi(u_{\overline{T}}) \in \relint(\nc(C_{\overline{T}}))$ by \cite[Proposition 2.13]{Zhou}, so there exist positive constants $t_1, \dots, t_k$ such that $d\varphi(u_{\overline{T}})=\sum_{i=1}^k{t_i \alpha_i}$. Since $u_{\overline{T}} \in C_{\overline{T}}$, we have $\langle \alpha, u_{\overline{T}}\rangle = h(\alpha)$ for all $\alpha \in T$. Therefore, by assumption that $l_{\alpha_i}(u) \geq 0$ and the definition of affine distance, we get 
    \begin{equation*}
        \langle d\varphi(u_{\overline{T}}),u-u_{\overline{T}} \rangle= \sum_{i=1}^k{t_il_{\alpha_i}(u)} \geq \min_{1 \leq i \leq k}\{t_i\} \cdot \max_{1 \leq i \leq k}\{l_{\alpha_i}(u) \}=\min_{1 \leq i \leq k}\{t_i\} \cdot d_{\aff}(u,C_{\overline{T}}).
    \end{equation*}
    Therefore, we can take $c=\min_{1 \leq i \leq k}\{t_i\}$ (which only depends on $d\varphi(u_{\overline{T}})$) to get the desired inequality. 
\end{proof} 
\section{Tailoring} \label{section:tailoring}
In this section, we extend the tailoring methods for hypersurfaces introduced in \cite{Abou} (originally, the process of tailoring was introduced in \cite{Mikh}, but we follow the more quantitative approach through specifying explicit cut-off functions). Our construction of bump functions is also heavily inspired by the ones appearing in \cite{Zhou}. 

Recall that \emph{tailoring} (also called \emph{localisation}) is a method introduced by Mikhalkin in the context of studying very affine hypersurfaces. It allows to us to make the heuristic outlined in Remark \ref{remark:dom_terms} precise, in the sense that all the non-dominant terms of the defining polynomial are set to be zero by introducing some cut-off functions $\chi_{\alpha,\beta}$  and modifying the defining equation for a hypersurface $H_{\beta}$ to get a tailored hypersurface $\widetilde{H}_{\beta}$. 

The expected behaviour for the functions is that $\chi_{\alpha,\beta}=1$ on the region $C_{\alpha}$ and $\chi_{\alpha,\beta}=0$ outside its small neighbourhood (which converges to $C_\alpha$ as $\beta \rightarrow \infty$). This manifestly forces the cut-off functions to be non-holomorphic, so $\widetilde{H}_{\beta}$ is no longer a complex submanifold of $M_{\Cs}$. This means that some care needs to be taken in establishing that the result of tailoring is a symplectic submanifold. 

We first illustrate our tailoring method on the case of a hypersurface, emphasising the key differences from previously employed tailoring methods, and then proceed with the generalisation for complete intersections. 

Let $\varphi\colon M_{\R} \rightarrow \R$ be a strictly convex, proper, bounded-below function. Then the pullback of $\varphi$ to $M_{\Cs}$ under $\Log$ (that will also be denoted by $\varphi$) is the potential for a Kähler form $\omega=d\lambda$ with a Liouville one-form $\lambda=-d^c\varphi$. We call Kähler potentials of that form \emph{toric}. In this section, unless otherwise specified, we shall state our results for Kähler forms coming from a potential $\varphi$ that is toric and homogeneous of degree $2$. A simple compactness argument shows that the Riemannian metrics $g_1$ and $g_2$ induced by two such potentials $\varphi_1$ and $\varphi_2$ are Lipschitz-equivalent, so it does not matter which particular $\varphi$ we choose to measure lengths. Note that, in particular, all the potentials considered in Section \ref{section:potentials} fall into this category. 

\subsection{Hypersurfaces}\label{section:tailor-hyp}

Suppose that $P=\conv(A) \subset N$ is a lattice polytope, let $h\in \R^A$ be a refined triangulating function inducing a regular triangulation $\mT$ of $P$ and $c_{\alpha}$ are non-zero constants for $\alpha \in A$. Recall that to such a set-up, we can associate a one-parameter family of hypersurfaces $H_{\beta}$ inside $M_{\Cs}$. 

\begin{definition}\label{definition:compatible_cutoffs}
    We say that $\chi=\{ \chi_{\alpha,\beta} \colon M_{\Cs} \rightarrow [0,1]\}_{\alpha \in A, \beta>0}$ is a \emph{collection of $h$-compatible cut-off functions} if the following conditions hold:
    \begin{enumerate}
        \item the functions are toric, i.e. $\chi_{\alpha,\beta}$ is a function of the radial coordinate $u$;
        \item there exists a constant $C>0$ such that $\lVert d \chi_{\alpha,\beta}(u) \rVert_g \leq C$ for all $\alpha$, $\beta$;
        \item we have $\chi_{\alpha,\beta}(u)=1$ whenever $L_h(u)\leq l_{\alpha}(u)+\beta^{-\frac12}$.
    \end{enumerate}

    We call a collection of $h$-compatible functions \emph{localising} if there also exists a constant $C'>0$ such that $\chi_{\alpha,\beta}(u)=0$ whenever $L_h(u) \geq l_\alpha(u)+\beta^{-\frac12}+C'\beta^{-1}$.
\end{definition}

The condition of being localising serves to rule out the cases such as $\chi_{\alpha,\beta} \equiv 1$ in the later calculations where it becomes important that we only have a small number of non-zero terms. 

\begin{remark}\label{remark:weird-defn}
    As remarked earlier, the condition (2) is independent of the auxiliary choice of metric $g$ induced by a homogeneous potential. The conditions (2) and (3) might seem contradictory for localising cut-off functions, since we require that the function $\chi_{\alpha,\beta}$ changes its value from $0$ to $1$ over a region of size $O(\beta^{-1})$, but has bounded derivative. This is illusory, since the metric $g$ is defined using coordinates $\rho=\Log(z)$ rather than $u=\beta^{-1}\rho$, so the size of the region is actually $O(1)$ from the perspective of the metric. One could, equivalently, work with the potentials pulled back under $\Log_\beta$ (as in Remark \ref{remark:scaling-weird}), which would require us to replace the constant bound in (2) by one linear in $\beta$.
\end{remark}

\begin{lemma}\label{lemma:cutoff-fns-space}
    The spaces of compatible cut-off functions and localising compatible cut-off functions are convex. Moreover, for any refined triangulating function $h$, there exists a localising collection of cut-off functions compatible with $h$. 
\end{lemma}
\begin{proof}
    The first part is straightforward, since if $\chi_{\alpha,\beta}$ and $\chi_{\alpha,\beta}'$ are compatible with $h$, so is $t\cdot\chi_{\alpha,\beta}+(1-t)\cdot \chi'_{\alpha,\beta}$ for all $t \in [0,1]$. If the two functions are localising with constants $C$, $C'$, their convex combinations are going to be localising with a constant $\max\{C,C'\}$. 
    
    Our construction for the second part is a slight generalisation of the one appearing in \cite{Zhou}, which we review for completeness (the extra degrees of freedom give us access to results like Corollary \ref{corollary:bdary_conv}). We start by picking a function $\chi\colon\R\rightarrow \left[0,1\right]$ with the following properties:
        \begin{itemize}
            \item $\chi(x)=1$ for $x\geq 0$ and $\chi(x)=0$ for $x \leq -2$;
            \item $\chi(x)\exp(x)$ is convex.
        \end{itemize}
    From this, we can define auxiliary functions $\chi_{\alpha,\alpha',\beta}$ for any vertices $\alpha$, $\alpha'$ adjacent in $\mT$ and $\beta\geq0$: pick a constant $K_{\alpha,\alpha'} \geq 0$, then set    
    \begin{equation*}
        \chi_{\alpha,\alpha',\beta}(u)\coloneqq\chi \left( \beta (l_{\alpha}(u)-l_{\alpha'}(u))+\sqrt{\beta}+K_{\alpha,\alpha'}\right).
    \end{equation*}
    
    By design, this function satisfies $\chi_{\alpha,\alpha',\beta}(u)=1$ whenever $l_{\alpha}(u) \geq l_{\alpha'}(u) - \beta^{-\frac12}-K_{\alpha,\alpha'}\beta^{-1}$ and $\chi_{\alpha,\alpha',\beta}(u)=0$ is $l_{\alpha}(u) \leq l_{\alpha'}(u)-\beta^{-\frac12}-(2+K_{\alpha,\alpha'})\beta^{-1}$.
    
    With this notation, we can take our tailoring functions to be of the form
    
    \begin{equation*}
        \chi_{\alpha,\beta,K}(u)\coloneqq\prod_{\alpha' \in \Nbhd(\alpha)}{\chi_{\alpha,\alpha',\beta}(u)},
    \end{equation*}
    where the product is taken over $\Nbhd(\alpha)$, the set of all vertices $\alpha'$ connected to $\alpha$ in the triangulation $\mT$. Finally, by precomposing with $\Log_{\beta}$, these maps lift to the desired cut-off functions on $M_{\Cs}$.

    If we pick $C_0=\max_{\alpha,\alpha'}K_{\alpha,\alpha'}$, then the constructed functions $\chi_{\alpha,\beta,K}$ clearly satisfy the desired properties requiring them to be $0$ or $1$ in certain regions. The bound on derivatives follows from the discussion in \cite[Proposition 1.7]{Zhou}.
\end{proof}

Given a collection of compatible cut-off functions $\chi$, the first part of the Lemma tells us that $1-s+s\chi$ is compatible for $s \in [0,1]$, so one can deform the individual terms $f_{\alpha,\beta}(z)$ as follows:

\begin{equation*}
    f_{\alpha,s,\beta}(z)\coloneqq\left(1-s+s\chi_{\alpha,\beta}(z)\right)f_{\alpha,\beta}(z).
\end{equation*}

This corresponds to a deformation of the entire Laurent polynomial to its tailored version:
\begin{equation*}
    f_{s,\beta}(z)\coloneqq\sum_{\alpha \in A}{f_{\alpha,s,\beta}(z)}.
\end{equation*}

These can be used to define a two-parameter family of hypersurfaces $H_{s,\beta}$ with $s \in [0,1]$, $\beta>0$ as 

\begin{equation*}
    H_{s,\beta}\coloneqq\left\{z \in N_{\Cs}:f_{s,\beta}(z)=0 \right\}.
\end{equation*}

We also introduce the following notation for the fully tailored hypersurfaces and its amoeba:

\begin{equation*}
    \begin{split}
    \widetilde{H}_{\beta}&\coloneqq H_{1,\beta}, \widetilde{f}_{\beta}(z)\coloneqq f_{1,\beta}(z), \\
    \mA_{s,\beta}&\coloneqq\Log\left(H_{s,\beta}\right), \widetilde{\mA}_{\beta}\coloneqq\mA_{1,\beta}. \\
    \end{split}
\end{equation*}

With this setup in mind, we are ready to state the main result of this section:

\begin{proposition}\label{proposition:isotopy}
    Let $\omega$ be a Kähler form on $M_{\Cs} \cong (\Cs)^n$ coming from a degree $2$ homogeneous toric potential. Then for any choice of compatible cut-off functions, there exists a $\beta_0$ such that the hypersurface $H_{s,\beta}$ is a symplectic submanifold of $\left(M_{\Cs},\omega\right)$ for all $\beta>\beta_0$ and $s\in[0,1]$. 
\end{proposition}

This is analogous to \cite[Proposition 1.7]{Zhou} and \cite[Proposition 4.2]{Abou}, but with a different class of cut-off functions, so we review the proof that is essentially identical to the one in \cite{Zhou} in order to state a few key lemmas explicitly. 

For future reference, we shall also introduce notation for the modulus of the monomials as $F_{\alpha,s,\beta}(z)\coloneqq|f_{\alpha,s,\beta}(z)|$, $F_{s,\beta}(z)\coloneqq\sum_{\alpha\in A}{F_{\alpha,s,\beta}(z)}$. Moreover, we let $h=h_\omega$ be the Hermitian metric induced by $\omega$ (so that $\omega=\IM(h)$ and $g=\RE(\omega)$ is the associated Riemannian metric). Then the following estimates on holomorphic and antiholomorphic parts of the derivative of $f_{s,\beta}(z)$ are instrumental in proving \ref{proposition:isotopy}:

\begin{lemma}\label{lemma:tailoring_estimate}
    For any choice of $h$-compatible cut-off functions, there exist constants $C_1$, $C_2$, $C_3>0$ such that the following inequalities hold for all $\beta>0$, $s \in [0,1]$ and $z \in M_{\Cs}$:
        \begin{equation*}
            \begin{split}
                \lVert\overline{\partial}f_{s,\beta}(z)\rVert_h &\leq F_{s,\beta}(z) C_1 e^{-\sqrt{\beta}},\\
                \lVert\partial f_{s,\beta}(z)\rVert_h &\geq F_{s,\beta}(z)(C_2-C_3 e^{-\sqrt{\beta}}).
            \end{split}
        \end{equation*}
\end{lemma}
\begin{proof}
    These bounds are precisely the intermediate steps in the proof of \cite[Proposition 1.7]{Zhou}, since the estimates used there rely only on the axiomatic properties of compatible cut-off functions rather than the explicit construction of $\chi_{\alpha,\beta,K}$ provided earlier. 
\end{proof}

Another ingredient in the proof is a classical result due to Donaldson that relates the above inequalities to properties of symplectic subspaces inside a vector space equipped with a compatible triple:

\begin{lemma}[\cite{Donaldson}]\label{lemma:donaldson}
    Let $V$ be an $n$-dimensional vector space equipped with a compatible triple $(J,\omega,h)$. Then if $a\colon V \rightarrow \C$ and $b\colon V \rightarrow \C$ are, respectively, complex linear and complex antilinear maps satisfying $\lVert a \rVert_h > \lVert b \rVert_h$, the real vector subspace $\ker(a+b) \subset V$ is symplectic with respect to $\omega$. 
\end{lemma}

\begin{proof}[Proof of Proposition \ref{proposition:isotopy}](following \cite{Abou})
    To show that $H_{s,\beta}=f_{s,\beta}^{-1}(0)$ is a smooth symplectic submanifold in $M_{\Cs}$ of codimension $2$, we need to prove that $\ker(df_{s,\beta}(z)) \cong T_zH_{s,\beta}$ is a symplectic codimension $2$ subspace of $T_zM_{\Cs}$. By Lemma \ref{lemma:donaldson} and the decomposition $df_{s,\beta}=\partial f_{s,\beta}+\overline{\partial}f_{s,\beta}$, it suffices to prove that 
    \begin{equation*}
        \lVert\partial f_{s,\beta}(z)\rVert_h > \lVert\overline{\partial}f_{s,\beta}(z)\rVert_h
    \end{equation*}
    holds for all $z \in H_{s,\beta}$ for $\beta$ large enough. However, this clearly follow from Lemma \ref{lemma:tailoring_estimate}. 
\end{proof}

As advertised earlier, the tailoring process makes the heuristic from Remark \ref{remark:dom_terms} exact, so locally, the hypersurface $\widetilde{H}_{\beta}$ will be given as the vanishing locus of a function $\sum_{\alpha \in T}{\widetilde{f}_{\alpha,\beta}(z)}$ for some simplex $T \in \mT$. This shall greatly simplify the calculations that we have to do while we compute the skeleton, since it allows us to do them locally and neglect all the non-dominant defining terms. Recall that by Remark \ref{remark:comp_regions}, the tropical amoeba $\mA_\trop$ has complementary regions $C_\alpha$ labelled by $\alpha \in A$, so we denote the corresponding regions for logarithmic tailored amoebae as $\widetilde{C}_{\alpha,\beta}$. In fact, one can provide explicit equations for boundaries of these regions:

\begin{proposition}[{\cite[Proposition 1.9]{Zhou}}]\label{proposition:bdary_eqn}
    The boundary of the connected component $\widetilde{C}_{\alpha,\beta}$ of $M_{\R}\backslash\widetilde{\mA}_{\beta}$ is given by the equation 
    \begin{equation*}
        \partial\widetilde{C}_{\alpha,\beta}=\left\{\sum_{\alpha' \in \Nbhd(\alpha)}{\widetilde{F}_{\alpha',\beta}(u)}=|c_{\alpha}|\right\}.
    \end{equation*}
\end{proposition}

In the most relevant case, where $P$ contains the origin and $h$ induces a star-shaped triangulation centred at the origin, there is a particularly convenient choice of the constants $K_{\alpha,\alpha'}$ in our tailorings. Namely, we can guarantee that the single compact complementary region $\widetilde{C}_{0,\beta}$ is convex, which plays an important role in determining the skeleton. Then we call a collection of $h$-compatible cut-off functions \emph{centred} if the complementary region $\widetilde{C}_{0,\beta}\subset M_{\R}\backslash\widetilde{\mA}_{\beta}$ containing the origin is convex. 

\begin{corollary}\label{corollary:bdary_conv}
    For any centred refined triangulating function $h$ on the polytope $P$, there exists a collection of localising centred $h$-compatible cut-off functions.
\end{corollary}
\begin{proof}
    We will show that for an appropriate choice of constants $K_{\alpha,\alpha'}$ in the construction from Lemma \ref{lemma:cutoff-fns-space}, the resulting collection will be centred. By Proposition \ref{proposition:bdary_eqn}, we know that the boundary of the region $\widetilde{C}_{0,\beta}$ is a level set of the function $G(u)=\sum_{\alpha \in A \backslash \{ 0\}}{\widetilde{F}_{\alpha,\beta}(u)}$. Moreover, we can see that 
    \begin{equation*}
        \widetilde{C}_{0,\beta}=\{u \in M_{\R} \colon G(u) \leq |c_0| \}. 
    \end{equation*}
    Hence it suffices to show that the function $G$ is convex on some region containing $\widetilde{C}_{0,\beta}$. 

    Denote $K\coloneqq2+\log|c_0|-\min_{\alpha \in A}{\log|c_{\alpha}|}$ and consider the set 
    \begin{equation*}
        V=\{ u \in M_{\R} \colon L_h(u) \leq K\beta^{-1}\}.
    \end{equation*} 

    We shall make the following choice for constants $K_{\alpha,\alpha'}$:

    \begin{equation*}
        K_{\alpha,\alpha'}=     \begin{cases}
      0, & \textnormal{if}\ \alpha'= 0\\
      K, & \textnormal{otherwise}.
    \end{cases}
    \end{equation*}
    
    Observe that we also have $\widetilde{C}_{0,\beta} \subset V$: if $u \notin V$, then any $\alpha$ satisfying $L_h(u)=l_\alpha(u)$ must be non-zero, since we have $L_{h}(u)>K\beta^{-1}>0$. Any such $\alpha$ will will also satisfy $\chi_{\alpha,\beta}(u)=1$, since $u \in C_\alpha$, and $|c_{\alpha}|e^{\beta l_{\alpha}(u)} > |c_0|$ by the choice of $K$, thus $G(u)\geq\widetilde{F}_{\alpha,\beta}(u)>|c_0|$ and $u \notin \widetilde{C}_{0,\beta}$.  
    
    We shall prove that for $u \in V$ and any $\alpha \neq 0$, we have $\chi_{\alpha,\beta}(u)=\chi_{\alpha,0,\beta}(u)$: clearly, if $\chi_{\alpha,0,\beta}(u)=0$, then $\chi_{\alpha,\beta}(u)=0$ as well. When $\chi_{\alpha,0,\beta}(u)>0$, we get the inequality $l_{\alpha}(u)>-(\beta^{-\frac12}+2\beta^{-1})$ from the support estimates from Lemma \ref{lemma:cutoff-fns-space}. By the assumption $u \in V$, we also know know that $l_{\alpha'}(u) \leq K\beta^{-1}$ holds for all $\alpha' \in \Nbhd(\alpha) \backslash \{0\}$. Therefore, summing these two inequalities yields $l_{\alpha'}(u)-l_{\alpha}(u)<\beta^{-\frac12}+(K+2)\beta^{-1}$, so the aforementioned support estimates tell us that $\chi_{\alpha,\alpha',\beta}(u)=1$ holds for all $\alpha' \in \Nbhd(\alpha) \backslash \{0\}$, which means that $\chi_{\alpha,\beta}(u)=\prod_{\alpha' \in \Nbhd(\alpha)}{\chi_{\alpha,\alpha',\beta}(u)}=\chi_{\alpha,0,\beta}(u)$, as required. 

    Finally, recall that we have chosen the function $\chi$ so that $\chi(x)e^x$ is convex on $\R$, so by composing with a linear function and multiplying by a positive constant, we get that $\widetilde{F}_{\alpha,\beta}(u)=e^{\beta l_{\alpha}(u)}\chi_{\alpha,0,\beta}(u)$ is convex. But since a sum of convex functions is convex, $G$ is indeed convex on $V$, so we are done. 
\end{proof}

\subsection{Complete intersections}\label{section:tailor-ci}

Throughout this section, fix lattice polytopes $P_1=\conv(A_1)$, \dots, $P_r=\conv(A_r)$ in $N_{\R}$ with $P=P_1+\dots+P_r$, refined triangulating functions $h_j$ on $P_j$ inducing triangulations $\mT_j$ of $P_j$ and non-zero constants $c_{\alpha,j}$ for all $\alpha \in A_j$ and all $j=1$,\dots,$r$. Suppose that the induced tropical complete intersection $\mA_\trop=\bigcap_{j=1}^r \mA_{\trop,j}$ is transverse, so that we are in the setting of Section \ref{section:tropical-ci}. Let $Z_{\beta}$ be the associated one-parameter family of complete intersections in $M_{\Cs}$ with defining polynomials $f_{\beta,1}(z)$, \dots, $f_{\beta,r}(z)$. 

\begin{definition}\label{definition:compatible-cutoff-ci}
    We call $\chi=(\chi_1,\dots,\chi_r)$ a \emph{collection of $(h_j)$-compatible cut-off functions} if $\chi_j=\{ \chi_{\alpha,\beta}\}_{\alpha \in A_j,\beta>0}$ is $h_j$-compatible for all $j=1,\dots,r$. We say that $\chi$ is \emph{localising} if every $\chi_j$ is localising. 
\end{definition}

We can now deform all the polynomials by picking a $(h_j)$-compatible collection $\chi$ and then introducing a factor of $\chi_{\alpha,\beta,j}$ in front of the corresponding monomial terms in $f_{\beta,j}(z)$ as if we were only tailoring the individual hypersurface $H_{\beta,j}=\V(f_{\beta,j})$, which leads us to the following definition of the tailoring of $Z_{\beta}$:

\begin{equation*}
    \begin{split}
        f_{\alpha,s,\beta,j}(z)&\coloneqq\left(1-s+s\chi_{\alpha,\beta,j}(z)\right)e^{-\beta h(\alpha)}z^{\alpha}, \\
        f_{s,\beta,j}(z)&\coloneqq\sum_{\alpha \in A_j}{f_{\alpha,s,\beta,j}(z)}, \\
        Z_{s,\beta}&\coloneqq\V(f_{s,\beta,1},\dots,f_{s,\beta,r}),\\
        \widetilde{f}_{\beta,j}(z)&\coloneqq f_{1,\beta,j}(z) \textnormal{ and } \widetilde{Z}_{\beta}\coloneqq Z_{1,\beta}.\\
    \end{split}
\end{equation*}

The following proposition confirms that this is indeed a sensible way of tailoring:

\begin{proposition}\label{proposition:isotopy_ci}
    Let $\omega$ be a toric Kähler form on $M_{\Cs}$ coming from a degree $2$ homogeneous potential. Then for any choice of a collection of $(h_j)$-compatible cut-off functions $\chi$, there exists a $\beta_0>0$ such that the complete intersection $Z_{s,\beta}$ is a smooth $2(n-r)$-dimensional symplectic submanifold of $\left(M_{\Cs},\omega\right)$ for all $\beta>\beta_0$ and $s \in [0,1]$. 
\end{proposition}

The proof relies heavily on the estimates from Lemma \ref{lemma:tailoring_estimate} and a generalisation of Lemma \ref{lemma:donaldson}. In order to provide a more systematic treatment of the linear algebra of compatible triples, we first recall some definitions from \cite{Donaldson} and \cite{CielebakMohnke}. 

Given a finite dimensional real inner product space $(V,\langle \cdot, \cdot \rangle)$, define the \emph{angle between two non-zero vectors} $x,y \in V$ as

\begin{equation*}
    \angle (x,y) \coloneqq \cos^{-1}\left(\frac{\langle x, y \rangle}{\lVert x \rVert \lVert y \rVert} \right) \in [0,\pi].
\end{equation*}

The \emph{angle between a linear subspace $X \subset V$ and a vector $y \in V$} is 

\begin{equation*}
    \angle(X, y) \coloneqq \inf_{x \in X \backslash \{ 0\}}{\angle(x,y)},
\end{equation*}

which is always going to be equal to $\angle (\pi^{\perp}(y),y )$. Finally, the \emph{angle between two linear subspaces} $X,Y \subset V$ is 

\begin{equation*}
    \angle ( X, Y) \coloneqq \sup_{y \in Y\backslash \{0\}}{\angle(x,Y)}=\sup_{y \in Y\backslash\{0\}}{\inf_{x \in X \backslash \{0\}}}{\angle(x,y)} \in [0,\frac12\pi].
\end{equation*}

Note that the definition is manifestly not symmetric, in the sense that one might have $\angle(X,Y)\neq \angle(Y,X)$. We are mostly interested in the case when the inner product comes from a compatible triple $(J,\omega,\langle\cdot,\cdot\rangle)$ (that also automatically has a Riemannian metric $g$ and a Hermitian metric $h$). There, one defines the \emph{angle of a real subspace} $X \subset V$ as $\angle(X)\coloneqq\angle(X,JX) \in [0,\frac12\pi]$. For example, it is immediate that $X$ is a complex subspace if and only if $\angle(X)=0$, or that we have the following relationship with symplectic forms:

\begin{lemma}\label{lemma:approx_hol_symp}
    A subspace $X \subset V$ satisfies $\angle(X)<\frac12\pi$ if and only if it is a symplectic subspace.
\end{lemma}
\begin{proof}
    The angle condition is equivalent to $X \cap (JX)^\perp=\{0\}$, which is the same as $X \cap X^\omega=\{0\}$.
\end{proof}

Therefore, the following definition provides a notion that is stronger than being a symplectic submanifold, but weaker than being a complex submanifold:

\begin{definition}\label{definition:ahol-submfd}
    For $0\leq \delta\leq\frac12\pi$, we say that a submanifold $X$ of a Kähler manifold $M$ is \emph{$\delta$-approximately holomorphic} if $\angle(TM)\leq\delta$.
\end{definition}

\begin{lemma}\label{lemma:angles_estimate}
    Let $V$ be an $n$-dimensional complex vector space equipped with a complex structure $J$ and Kähler metric $h$. Then for every $\varepsilon,\varepsilon'>0$, there exists a $\delta>0$ with the following property: suppose that $f=(f_1,\dots,f_r)\colon V \rightarrow \C^r$ is a real linear map whose components can be decomposed into complex linear and complex anti-linear parts as $f_j(z)=\langle z, a_j \rangle_h + \langle b_j, z \rangle_h$, and let $A$ be an $(r\times r)$-matrix with entries $A_{jk}=\frac{\langle a_j, a_k \rangle_h}{\Vert a_j\rVert_h\lVert a_k\rVert_h}$. If $f$ satisfies $\left|\det(A)\right| \geq \varepsilon$ and $\frac{\lVert b_j\rVert_h}{\lVert a_j\rVert_h}<\delta$ for all $j=1,\dots,r$, then $\ker(f) \subset V$ is a $2(n-r)$-dimensional subspace with angle at most $\varepsilon'$. 
\end{lemma}
\begin{proof}
    Denote the $J$-compatible Riemannian metric associated to $h$ as $g=\RE(h)$. By \cite[Lemma 8.3(c)]{CielebakMohnke}, we know that $\angle(X)=\angle(X^{\perp})$, thus it is sufficient to prove the result for $V_0\coloneqq\ker(f)^{\perp}$ instead. Also, note that if we rescale a map $f_j$ by some positive real constant, it does not change the ratio $\frac{\lVert b_j \rVert_h}{\lVert a_j \rVert_h}$ nor the matrix $A$, so, without loss of generality, assume that $\lVert a_j \rVert = 1$ and $\lVert b_j \rVert < \delta$. For $f_j(z)=\langle z, a_j \rangle_h + \langle b_j, z \rangle_h$, we can rewrite it as $f_j(z)=\langle z,a_j+b_j \rangle_g-i\langle z,J(a_j-b_j)\rangle_g$, hence the orthogonal complement of the kernel has to be $\ker(f_j)^{\perp}=\spann_{\R}(a_j+b_j,J(a_j-b_j))$. Hence, if we denote $v_{2j-1}=a_j+b_j$, $v_{2j}=J(a_j-b_j)$ for $j=1$, \dots, $r$, these vectors give us a spanning set for $V_0$. 

    Denote the image of $A \in GL(r,\C)$ under the standard embedding $GL(r,\C) \hookrightarrow GL(2r,\R)$ as $\Hat{A}$, then $\det(\Hat{A})=\det(A)^2 \geq \varepsilon^2$ ($A$ is clearly Hermitian, so its determinant is real). Observe that in coordinates, we can write $\Hat{A}_{jk}=g(\Hat{v}_j,\Hat{v}_k)$ for $\Hat{v}_{2l-1}=a_l$, $\Hat{v}_{2l}=Ja_l$ with $l=1$, \dots, $r$, since $\langle a_j,a_k\rangle_h=g(a_j,a_k)+ig(a_j,Ja_k)$ maps to a $(2\times2)$-block 
    \begin{equation*}
        \begin{pmatrix}
        g(a_j,a_k) & -g(a_j,Ja_k) \\
        g(a_j,Ja_k) & g(a_j,a_k) \\
    \end{pmatrix}=\begin{pmatrix}
        g(a_j,a_k) & g(Ja_j,a_k) \\
        g(a_j,Ja_k) & g(Ja_j,Ja_k) \\
    \end{pmatrix},
    \end{equation*}
    under the  embedding. 
    
    In particular, if we consider the $(2r\times2r)$ real matrix $B$ with entries $B_{lm}=g(v_l,v_m)$, the assumption on norms $\lVert b_j \rVert_h$ being $\delta$-small means that the matrix $B$ is $3\delta$-close to $\Hat{A}$ for $\delta<1$ (by Cauchy--Schwarz, one can get estimates like $|g(\Hat{v}_{2j},\Hat{v}_{2k})-g(v_{2j},v_{2k})|=|g(a_j,b_k)+g(b_j,a_k)+g(b_j,b_k)| \leq 2\delta+\delta^2 \leq 3\delta$). 
    
    Since the set of possible values of the matrix $\hat{A}$ is a compact subset of $GL(2r,\R)$ (we have a lower bound on the determinant by assumption and an upper bound on entries of the matrix by Cauchy--Schwarz), its $3\delta$-neighbourhood will still lie inside the open set $GL(2r,\R) \subset \textnormal{Mat}_{2r\times 2r}(\R)$, which forces $\det(B)>0$, so $V_0$ is a real $2r$-dimensional subspace of $V$ for all sufficiently small $\delta$. Finally, observe that the angle function is continuous on the non-compact Stiefel manifold $St(2r,V)$ (i.e. the set of linearly independent $2r$-tuples of vectors in $V$), while the above discussion tells us that $V_0$ lies $3\delta$-close to a tuple $(a_1,Ja_1,\dots,Ja_r)$ that spans a complex subspace $\Hat{V_0}\coloneqq\spann_\R(a_1,Ja_1,\dots,Ja_r)$ and is constrained to lie in a compact subset $K \subset St(2r,V)$. By the above discussion, we know that our basis lies in a $3\delta$-neighbourhood of $K$ in $St(2r,V)$, so it has to lie in a slightly bigger compact set $K \subset K'$ for sufficiently small $\delta$. Finally, recall that being complex implies $\angle(\hat{V}_0)=0$, hence the conclusion $\angle(V_0) \leq \varepsilon'$ follows from uniform continuity of $\angle$ on $K'$ for all sufficiently small $\delta$. 
\end{proof}
The desired generalisation of Donaldson's lemma we get is the following:

\begin{corollary}\label{corollary:donaldson_on_steroids}
    In the setup of Lemma \ref{lemma:angles_estimate}, for any $\varepsilon>0$, there exists a $\delta>0$ such that $\ker(f)$ is a symplectic subspace of $(V,\omega)$.
\end{corollary}

\begin{proof}
    Follows from Lemma \ref{lemma:angles_estimate} and Lemma \ref{lemma:approx_hol_symp}. 
\end{proof}

\begin{remark}
    In Donaldson's original result, Lemma \ref{lemma:donaldson}, normalisation forces that the $(1\times 1)$-matrix $A$ is always equal to $1$, so the determinant condition becomes vacuous. Unlike our Corollary \ref{corollary:donaldson_on_steroids}, it also provides a sharp bound on $\delta$, since it allows $\delta=1$, which is clearly optimal. 
\end{remark}

With this in mind, the following couple of computational lemmas demonstrate that the linear algebra result can be applied in our setting.

\begin{lemma}\label{lemma:partials_lin_indep}
    For any choice of compatible cut-off functions, there exists a $\beta_0$ such that for all $\beta>\beta_0$, $s\in[0,1]$ and $z \in Z_{s,\beta}$, the vectors $\partial f_{s,\beta,1}(z)$, \dots, $\partial f_{s,\beta,r}(z)$ are linearly independent $\C$-linear functionals . 
\end{lemma}

\begin{proof}
    We can use the same strategy as in the proof of Lemma \ref{lemma:sm_nontailored_ci} with $Z_{\beta}$ replaced by $Z_{s,\beta}$. 

    Like there, we begin by defining the simplices $T_j \in \mT_j$, pick $\alpha_j \in T_j$ such that $L_{h_j}(u)=l_{\alpha_j}(u)$ and without loss of generality shift the defining equations so that $0 \in T_j \backslash \{\alpha_j\}$ for $j=1$, \dots, $r$, which means that we can find vectors $v_1$,\dots,$v_r$ satisfying $\langle \alpha_j,v_k \rangle=\delta_{jk}$ and $\langle \alpha,v_k \rangle = 0$ for all other $\alpha \in \bigcup_j T_j$. Note that the estimate from Equation \ref{equation:estsing2} still holds in this setting.
    
    The coordinate expression for the differential in this case becomes somewhat more complicated: writing $f_{\alpha,s,\beta,j}(z)=(1-s+s\chi_{\alpha,\beta,j}(z))f_{\alpha,\beta,j}(z)$ and applying the product rule gives
    \begin{equation*}
        \partial f_{s,\beta,j}(z)=\sum_{\alpha \in A_j}{f_{\alpha,s,\beta,j}(z)\langle \alpha_j,d\rho+id\theta\rangle+s\partial\chi_{\alpha,\beta,j}(z)f_{\alpha,\beta,j}(z)}.
    \end{equation*}
    Similarly to Lemma \ref{lemma:sm_nontailored_ci}, the terms of the sum with $\alpha \in A_j \backslash T_j$ can be bounded through Equation \ref{equation:estsing2}. The extra term involving $\partial \chi$ is bounded by a constant multiple of $F_{s,\beta,j}(z)\cdot e^{-\sqrt{\beta}}$: by definition of compatible cut-off functions, $\lVert d\chi_{\alpha_j,\beta}(z) \rVert_g$ is bounded and is non-zero only when the corresponding $f_{\alpha,0,\beta,j}(z)$ term is smaller or equal to $F_{s,\beta,j}(z)\cdot e^{-\sqrt{\beta}}$ (cf. Lemma \ref{lemma:tailoring_estimate}). Therefore, if we consider the matrix $A$ with entries $A_{jk}=\partial f_{s,\beta,j}(x)(v_k)$, it will be of the form 
        \begin{equation*}
        A=\begin{pmatrix}
         f_{\alpha_1,s,\beta,1}(z)& 0 &  \hdots & 0\\
        0 & f_{\alpha_2,s,\beta,2}(z) & \hdots & 0 \\
        \vdots & \ddots & & \vdots \\
        0 & 0 & \vdots & f_{\alpha_r,s,\beta,r}(z) \\
        \end{pmatrix} \cdot \left(I_r+D \right),
    \end{equation*}
    with the remainder $D$ having the modulus of all entries bounded by $Ke^{-\sqrt{\beta}}$ for some constant $K$ (we have chosen $\alpha_j$ so that $l_{\alpha_j}(u)=L_{h_j}(u)$ and $F_{\alpha_j,s,\beta,j}$ is the largest term in the sum defining $F_{s,\beta,j}$, so this indeed follows from the above bounds). Since $\chi_{\alpha_j,\beta,j}(u)=1$ by $L_{h_j}(u)=l_{\alpha_j}(u)$, we get $f_{\alpha_j,s,\beta,j}(x)=f_{\alpha_j,\beta,j}(x)$ and hence this is the exact same expression as in Lemma \ref{lemma:sm_nontailored_ci}, so we can deduce the desired statement analogously. 
\end{proof}

\begin{remark}\label{remark:weird_generalisation}
    Note that Corollary \ref{corollary:weird_generalisation} also generalises to this setting, since the cut-off functions only depend on $u$. 
\end{remark}

\begin{lemma}\label{lemma:cursed_determinant}
     Define the function $A\colon M_{\Cs} \rightarrow \textnormal{Mat}_{r\times r}(\C)$ by setting 
    \begin{equation*}
        A_{jk}(z)\coloneqq \frac{\langle\partial f_{s,\beta,j}(z) ,\partial f_{s,\beta,k}(z)\rangle_{h}}{\lVert \partial f_{s,\beta,j}(z) \rVert_{h}\lVert \partial f_{s,\beta,k}(z) \rVert_{h}},
    \end{equation*}
    for $1 \leq j,k \leq r$. Then for any choice of compatible cut-off functions, there exist constants $\beta_0$ and $\varepsilon>0$ such that the inequality
    \begin{equation*}
        \det(A(z)) \geq \varepsilon,
    \end{equation*}
    holds for all $\beta>\beta_0$, $s\in[0,1]$ and $z \in Z_{s,\beta}$. 
\end{lemma}

\begin{proof}  
    We prove that the bound holds over the set $\Log_\beta^{-1}(\mA_\trop(\beta))$ from Corollary \ref{corollary:weird_generalisation}, which contains $Z_{s,\beta}$ for all $s \in [0,1]$ and all sufficiently large $\beta>0$ by our proof of Lemma \ref{lemma:partials_lin_indep}. By Remark \ref{remark:weird_generalisation}, we know that the vectors $\partial f_{s,\beta,j}$ are linearly independent and the determinant is non-zero for all $z \in \Log^{-1}_\beta(\mA_\trop(\beta))$, so we just need to produce a uniform lower bound ($A(z)$ is positive semidefinite by design, so the determinant is positive). By the characterisation of $\mA_\trop(\beta)$ from Lemma \ref{lemma:aff-nbhd-structure}, we know that it admits a decomposition into closed full-dimensional polyhedral cells indexed by the cells of $\mA_\trop$. By Theorem \ref{theorem:combi-ci}, we know that these are indexed by the mixed cells $T=T_1+\dots+T_r$ in the subdivision of $P$. Denote the full-dimensional cell of $\mA_\trop(\beta)$ associated to such a tuple as $C_{(T_1,\dots,T_r)}(\beta)$.

    We shall also introduce a generalised version of the notation from the statement of this Lemma: for any Hermitian metric $h'$ on $\C^n$, we can identify all the tangent spaces $T_z M_{\Cs}$ with $\C^n$ by picking a trivialisation of the tangent bundle and then consider $A' \in \textnormal{Mat}_{r\times r}(\C)$ with entries
    \begin{equation*}
        A'_{jk}(z,h')=\frac{\langle\partial f_{s,\beta,j}(z) ,\partial f_{s,\beta,k}(z)\rangle_{h'}}{\lVert \partial f_{s,\beta,j}(z) \rVert_{h'}\lVert \partial f_{s,\beta,k}(z) \rVert_{h'}}.
    \end{equation*}
    In particular, the matrix $A_{jk}$ can be recovered as $A_{jk}(z)=A'_{jk}(z,h(z))$.  
    
    For a cell $C_{(T_1,\dots,T_r)}$ of $\mA_\trop$, denote the set of possible values of the expression $\sum_{\alpha \in T_j} {w_{\alpha,j} \cdot \alpha}$ for complex weights satisfying $\sum_\alpha |w_{\alpha,j}|=1$ as $K_{T_j} \subset M_\C$ and let $K_{(T_1,\dots,T_r)} \coloneqq K_{T_1} \times \dots \times K_{T_r}$. Suppose that $u \in C_{(T_1,\dots,T_r)}(\beta) \subset \mA_{\trop}(\beta)$, then by definition of $(h_j)$-compatible cut-off functions, we have $\chi_{\alpha,\beta,j}(u)=1$ for all $\alpha \in T_j$, and hence also $f_{s,\alpha,\beta,j}(z)=f_{\alpha,\beta,j}(z)$. Our estimates from Lemma \ref{lemma:partials_lin_indep} then give:
    \begin{equation}\label{equation:cursed1}
        F_{s,\beta,j}(z)^{-1}\cdot \partial f_{s,\beta,j}(z)=F_{s,\beta,j}(z)^{-1}\cdot\left(\sum_{\alpha \in T_j}{f_{\alpha,\beta,j}(z)\langle \alpha,d\rho+id\theta\rangle}\right) + O(e^{-\sqrt{\beta}}),
    \end{equation}
    where the magnitude of the remainder is measured using $h(z)$. In other words, an appropriate rescaling $w_j= F_{s,\beta,j}(z)^{-1} \cdot \partial f_{s,\beta,j}(z)$ of $\partial f_{s,\beta,j}(z)$ lies $O(e^{-\sqrt{\beta}})$-close to $K_{T_j}$, since the coefficients $w_{\alpha,j}=F_{s,\beta,j}(z)^{-1} \cdot f_{\alpha,\beta,j}(z)$ have moduli summing up to $1+O(e^{-\sqrt{\beta}})$. Note that the second part of the proof of that Lemma shows that for sufficiently large $\beta$, any tuple $w=(w_1,\dots,w_r)$ that is $O(e^{-\sqrt{\beta}})$-close to $K_{(T_1,\dots,T_r)}$ has to be linearly independent. Combined with compactness of $K_{(T_1,\dots,T_r)}$, this observation yields a lower bound on the determinant of the matrix $B(w)$ with entries $B_{jk}(w)\coloneqq\frac{\langle w_j,w_k \rangle_h}{\lVert w_j \rVert_h \rVert w_k \rVert_h}$ for $w$ varying over a neighbourhood of $K_{(T_1,\dots,T_r)}$ and any fixed Hermitian metric $h'$. Therefore, by taking $w_j=\partial f_{s,\beta,j}(z)$ for $z\in \Log_\beta^{-1}(C_{(T_1,\dots,T_r)}(\beta))$, we can cancel out the scaling factors $F_{s,\beta,j}(z)$ and get a lower bound on the determinant of $A(z,h')$. By repeating this argument for each of the finitely many cells of $\mA_\trop(\beta)$, it follows that for any fixed Hermitian metric $h'$, there exists a constant $\varepsilon>0$ such that the bound $|\det(A'(z,h))| \geq \varepsilon$ holds for all sufficiently large $\beta >0$ and all $z \in \Log^{-1}_\beta(\mA_\trop(\beta))$. 
    
    Finally, we observe that this approach still works if we replace $h'$ by a continuously varying family of Hermitian metrics $\{h'_p\}_{p \in V}$ parameterised by a compact set $V$, and the bound does not depend on the precise value of $\beta>0$, as long as it is sufficiently large. Since the potential $\varphi$ is toric and homogeneous of degree $2$, \cite[Proposition 2.4]{Zhou} tells us that $h$ is homogeneous of degree $0$, so we have $h(z)=h(\lVert u\rVert_{g_0}^{-1} \cdot u)$. Therefore, the bound $\det(A'(z,h(u)))>\varepsilon$ for all $z \in \mA_\trop(\beta)$, sufficiently large $\beta$ and $u \in M_\R$ satisfying $\lVert u \rVert_{g_0}=1$ follows from the discussion above and implies the desired statement. 
\end{proof}

\begin{corollary}\label{corollary:cursed-determinant}
    In the set-up of Lemma \ref{lemma:cursed_determinant}, the smallest eigenvalue of $A(z)$ is also bounded below by some constant $\varepsilon'>0$ for all $s \in [0,1]$, $z \in Z_{s,\beta}$ and sufficiently large $\beta>0$.
\end{corollary}

\begin{proof}
    By the Cauchy-Schwarz inequality, the entries of $A(z)$ are bounded above by $1$. Therefore, we have $|\sum_{k=1}^rA_{jk}v_k| \leq \sum_{k=1}^r|v_k|$ for all $v=(v_1,\dots,v_r) \in \R^r$ and $j=1,\dots,r$, so if $v$ is an eigenvector with eigenvalue $\lambda$, summing over all $j$ gives $\lambda \leq r$ (we know that $\lambda \geq 0$, since $A$ is positive definite). Therefore, the eigenvalues of $A(z)$ are bounded above by $r$. This means that the smallest eigenvalue must be bounded by $\varepsilon'=\varepsilon\cdot r^{-r+1}$, as desired. 
\end{proof}

\begin{remark}\label{remark:cursed-determinant}
    The same calculations show that the $(2r \times 2r)$-matrix $\hat{A}$ that is the image of $A$ under the embedding $GL(r,\C) \xhookrightarrow{} GL(2r,\R)$ (see Lemma \ref{lemma:angles_estimate}) also admits an analogous global lower bound on the determinant and on the smallest eigenvalue, since $\det(\hat{A})=\det(A)^2$. 
\end{remark}

\begin{corollary}\label{corollary:approx_hol_ci}
     For any choice of $(h_j)$-compatible cut-off functions $\chi=(\chi_1,\dots,\chi_r)$ and any $\varepsilon'>0$, there exists a $\beta_0>0$ such that the complete intersection $Z_{s,\beta}$ is a smooth $2(n-r)$-dimensional $\varepsilon'$-approximately holomorphic submanifold of $M_{\Cs}$ for all $\beta>\beta_0$ and $s \in [0,1]$. 
\end{corollary}

\begin{proof}
    It suffices to prove that for all $z \in Z_{s,\beta}$, $\beta_0>\beta$ and $s\in[0,1]$, the intersection
    \begin{equation*}
        \bigcap_{j=1}^{r}{\ker(df_{s,\beta,j}(z))} \subset T_zM_{\Cs},
    \end{equation*}
    is a subspace of dimension $2(n-r)$ with angle at most $\varepsilon'$. 

    We can decompose the maps into holomorphic and antiholomorphic parts as $df_{s,\beta,j}=\partial f_{s,\beta,j} + \overline{\partial} f_{s,\beta,j}$. By Lemma \ref{lemma:partials_lin_indep}, we know that the for $\beta$ large enough, the map $(\partial f_{s,\beta,1},\dots,\partial f_{s,\beta,r})\colon T_zM_{\Cs} \rightarrow \C^r$ is full rank for $\beta>\beta_0$. By Lemma \ref{lemma:cursed_determinant}, we know that if we increase $\beta_0$ enough, there exists a constant $\varepsilon>0$ such that $|\det(A(z))| \geq \varepsilon$ holds for all $\beta>\beta_0$, $z$ and $s$. Therefore, by Lemma \ref{lemma:angles_estimate}, there exists a constant $\delta$ such that if 
    \begin{equation*}
        \frac{\lVert \overline{\partial} f_{s,\beta,j} \rVert_h} {\lVert\partial f_{s,\beta,j}\rVert_h} \leq \delta,
    \end{equation*}
    holds for all $z$, $s$, $j$ and $\beta>\beta_0$, the desired conclusion follows. However, this can be accomplished by increasing $\beta_0$, since Lemma \ref{lemma:tailoring_estimate} tells us that all the ratios are bounded by $Ce^{-\sqrt{\beta}}$ for some constant $C>0$, so we are done.

\end{proof}

\begin{proof}[Proof of Proposition \ref{proposition:isotopy_ci}]
    This follows from Corollary \ref{corollary:approx_hol_ci} and Lemma \ref{lemma:approx_hol_symp}. 
\end{proof}

\subsection{Batyrev--Borisov complete intersections}\label{section:tailor-bbci}

Suppose now that we specialise to the case of open Batyrev--Borisov complete intersections: we start from the data of a nef partition $\nabla=\nabla_1+\dots+\nabla_r$ and a centred refined triangulating function $h$ that induces a regular triangulation $\mT$ of $\Delta^{\vee}$ and hence also of all $\nabla_j$'s. Suppose that we also normalise the constants $c_{\alpha,j}$ so that $c_{0,j}<0$ for all $j$. Recall that from a choice of $(h_j)$-compatible cut-off functions $\chi=(\chi_1,\dots,\chi_r)$, we have the hypersurfaces $H_{s,\beta,j}\coloneqq\V(f_{s,\beta,j})$ for all $j=1$,\dots,$r$. We also consider the following hypersurface:

\begin{equation*}
\begin{split}
    f_{s,\beta,\tot}(z)&\coloneqq \sum_{j=1}^{r}{f_{s,\beta,j}(z)}=\sum_{j=1}^{r}\sum_{\alpha \in \nabla_j}{f_{\alpha,s,\beta,j}(z)},\\
    H_{s,\beta,\tot}&\coloneqq \V(f_{s,\beta,\tot}) \subset M_{\Cs}.
\end{split}
\end{equation*}

By Lemma \ref{lemma:nef-summands-disjoint}, we know that $\nabla_j \cap \nabla_k = \{0 \}$ for all $j\neq k$, therefore the above equation implies that the Newton polytope of the untailored polynomial $f_{\beta,\tot}(z)$ is $\Delta^{\vee}=\conv\{\nabla_1,\dots,\nabla_r\}$ and its monomials are $f_{0,\beta,\tot}=c_{0,\tot}$ (where $c_{0,\tot}=\sum_{j=1}^r{c_{0,j}}$) and $f_{\alpha,\beta,\tot}(z)=f_{\alpha,\beta,j}(z)$ for $\alpha \in \nabla_j \backslash \{0\}$. Note that this gives us a tailoring of the complex hypersurface $H_{\beta,\tot}$ in a way that is compatible with the tailoring $Z_{s,\beta}$. We record the properties of this tailoring $H_{s,\beta,\tot}$ as a proposition:

\begin{proposition}\label{proposition:total_tailoring}
    For any choice of compatible cut-off functions, there exists a $\beta_0>0$ such that $H_{s,\beta,\tot}$ is a symplectic $(2n-2)$-dimensional submanifold of $(M_{\Cs},\omega)$ for all $\beta>\beta_0$ and $s \in [0,1]$.
\end{proposition}

\begin{proof}
    By assumption, the function $h$ defines a refined triangulation of $\Delta^{\vee}$ (for a general complete intersection, we would not have a good way of obtaining a regular subdivision of $\conv\{P_1,\dots,P_r\}$). Therefore, to apply Proposition \ref{proposition:isotopy}, it suffices to check that the collection of cut-off functions given by

    \begin{equation*}
        \begin{split}
        \chi_{\alpha,\beta,\tot}(u)&\coloneqq\chi_{\alpha,\beta,j}(u) \textnormal{ whenever } \alpha \in \nabla_j \backslash \{0\},\\
        \chi_{0,\beta,\tot}(u)&\coloneqq \frac{1}{|c_{0,\tot}|}\sum_{j=1}^{r}|c_{0,j}|\cdot\chi_{0,\beta,j}(u),
        \end{split}
    \end{equation*}
    are $h$-compatible on $\Delta^\vee$. The first two conditions are clear and for the final one, we can observe that $L_h(u)=\max\{L_{h_1}(u),\dots,L_{h_r}(u)\}$ and the collections $\{ \chi_{\alpha,\beta,j}\}$ are $(h_j)$-compatible for $j=1$,\dots,$r$, hence when $L_h(u) \leq l_{\alpha}(u)+\beta^{-\frac12}$ for some $\alpha \in \Delta^\vee$, then also $L_{h_j}(u) \leq l_{\alpha}(u)+\beta^{-\frac12}$ for all $j$, which implies that $\chi_{\alpha,\beta,j}(u)=1$ and thus $\chi_{\alpha,\beta,\tot}(u)=1$, as desired.
\end{proof}

By analogy with the hypersurface setting, we call the collection of $(h_j)$-compatible cut-off functions $\chi=(\chi_1,\dots,\chi_r)$ \emph{centred} if the complementary regions $\widetilde{C}_{0,\beta,j}\subset M_{\R}\backslash\widetilde{\mA}_{\beta,j}$ for $j=1,\dots,r$ and $\widetilde{C}_{0,\beta,\tot}\subset M_{\R}\backslash\widetilde{\mA}_{\beta,\tot}$ are all convex.

\begin{proposition} \label{proposition:all_shall_be_convex}
    For any nef partition $\nabla=\nabla_1+\dots+\nabla_r$ and any centred refined triangulating function $h\colon \Delta^\vee \rightarrow \R$, there exists a collection of localising centred $(h_j)$-compatible functions $\chi=(\chi_1,\dots,\chi_r)$. 
\end{proposition}

\begin{proof}
    We proceed analogously to Corollary \ref{corollary:bdary_conv} to show that the cut-off functions can be made centred, except we pick a different constant
    \begin{equation*}
        K=2+\log\left(|c_{0,\tot}|\right)-\min_{\alpha \in \Delta^{\vee}}\log|c_{\alpha,\tot}|.
    \end{equation*}
    By assumption, we have $c_{0,j}<0$, which means that $|c_{0,\tot}|=\sum_{j=1}^r|c_{0,j}|>|c_{0,j}|$ holds for all $j=1,\dots,r$, so the proof of the aforementioned Corollary goes through without any modifications for each $\widetilde{C}_{0,\beta,j}$: we can consider the defining function $G_j(u)$ and show that it is convex on a region $V_j=\{ u \in M_\R \colon l_\alpha(u) \leq K\beta^{-1} \textnormal{ for } \alpha \in \nabla_j\}$ containing $\widetilde{C}_{0,\beta,j}$ by analysing the tailored monomials individually. The region $\widetilde{C}_{0,\beta,\tot}$ will then have a defining function $G_\tot(u)=\sum_{j=1}^r G_j(u)$, which is a sum of convex functions on the region $V= \bigcap_{j=1}^r V_j$. Our choice of $K$ guarantees that $V$ must contain $\widetilde{C}_{0,\beta,\tot}$, hence it is a sub-level set of a convex function and we are done. 
\end{proof}

\subsection{Toric compactifications}\label{section:toric-cpct2}

To conclude this section, we shall study the consequences of our tailoring process for the toric compactifications of $Z_{s,\beta}$ inside toric stacks (along the lines of Section \ref{section:toric-cpct1}). In the rest of the section, we consider desingularisations $\mX_{\check{\Sigma}} \rightarrow X_P$ that correspond to simplicial refinements $\check{\Sigma} \rightarrow \check{\Sigma}_0=N(P)$ and give us an orbifold smooth toric compactifications of $Z_\beta$ for all $\beta \gg 0$ by Corollary \ref{corollary:smooth-cpctification}. Recall that Proposition \ref{proposition:int-strata2} also gives us equations for intersections $Z_\beta(\sigma)\coloneqq Z_\beta \cap O(\sigma)$ for $\sigma \in \check{\Sigma}$ in terms of truncations of $f_{\beta,j}$ over appropriate faces $P^\sigma_j \subset P_j$.

\begin{definition}\label{definition:cutoffs_cpct}
    We say that a collection of $(h_j)$-compatible cut-off functions $\chi=(\chi_1,\dots,\chi_r)$ is \emph{$\check{\Sigma}$-proper} for a simplicial refinement $\check{\Sigma} \rightarrow \check{\Sigma}_0$ if every $\chi_{\alpha,\beta,j}\colon M_{\Cs} \rightarrow \R$ extends to a smooth function $\overline{\chi}_{\alpha,\beta,j}\colon \mX_{\check{\Sigma}} \rightarrow \R$ and, in addition, for every cone $\sigma \in \check{\Sigma}$, the cut-off functions $\chi^\sigma=(\chi^\sigma_1,\dots,\chi^\sigma_r)$ on $O(\sigma)\cong (M/\sigma)_{\Cs}$ defined as $\chi^\sigma_j\coloneqq\{ \overline{\chi}_{\alpha,\beta,j}|_{O(\sigma)} \}_{\alpha\in P^\sigma_j,\beta>0}$ are $(h^\sigma_j)$-compatible with $h^\sigma_j\coloneqq h_j|_{P_j^\sigma}$ as cut-off functions on $(M/\sigma)_\R$. If there exists a simplicial refinement $\check{\Sigma} \rightarrow \check{\Sigma}_0$ for which $\chi$ is $\check{\Sigma}$-proper, we call $\chi$ \emph{proper}. 
\end{definition}

\begin{lemma}\label{lemma:cutoffs_cpct}
    For any refined triangulating functions $(h_j)$ and any choice of a simplicial refinement $\check{\Sigma} \rightarrow \check{\Sigma}_0$ as in Corollary \ref{corollary:smooth-cpctification}, there exists a collection of $\check{\Sigma}$-proper localising $(h_j)$-compatible cut-off functions.
\end{lemma}

\begin{proof}
    Let $\sigma \in \check{\Sigma}$ be a simplicial cone with an associated torus orbit $O(\sigma)$ and an affine chart $U(\sigma)$, it clearly suffices to check that the functions $\chi$ extend to smooth functions satisfying the requirements over $U(\sigma)$ for all such cones. Analogously to our proof of Corollary \ref{corollary:smooth-cpctification}, passing to a suitable cover shows that it is sufficient to treat the case when $\sigma$ is smooth and $U(\sigma) \cong \C^{\dim(\sigma)} \times (\Cs)^{n-\dim(\sigma)}$.
    
    Consider the construction from Lemma \ref{lemma:cutoff-fns-space} and pick standard coordinates $(w_\sigma,z_\sigma) \in \C^{\dim(\sigma)} \times (\Cs)^{n-\dim(\sigma)}$, so that $z_\sigma$ is a coordinate on $O(\sigma)=\{w_\sigma=0\}$. Analogously to \cite[Lemma 2.3]{ghhps}, we observe that our cut-off functions are independent of $w_\sigma$ near $w_\sigma=0$, so they admit smooth extensions to $\mX_\Sigma$. Moreover, these extensions can be explicitly computed as 
    \begin{equation*}
        \chi_{\alpha,\beta,j}^\sigma(u_\sigma)= 
        \begin{cases}
            \prod_{\alpha' \in \Nbhd_{P_j^\sigma}(\alpha)}\chi^\sigma_{\alpha,\alpha',\beta,j}(u_\sigma), \textnormal{ if } \alpha \in P^\sigma_j \\
            0, \textnormal{ otherwise},
        \end{cases} 
    \end{equation*}
    where we are using our chart to explicitly identify $O(\sigma)$ with $(M/\sigma)_{\Cs}$ to get real coordinates $u_\sigma=\Log_\beta(z_\sigma)$ on $(M/\sigma)_\R$, $\Nbhd_{P^\sigma_j}(\alpha)\coloneqq \Nbhd(\alpha)\cap P^\sigma_j$ and 
    \begin{equation*}
        \chi^\sigma_{\alpha,\alpha',\beta,j}(u)\coloneqq\chi(\beta(l_{\alpha}(u)-l_{\alpha'}(u))+\sqrt{\beta}+K_{\alpha,\alpha',j}),
    \end{equation*}
    with the constants $K_{\alpha,\alpha',j}$ and the functions $h^\sigma_j$ being induced from $A_j$. These expressions make it clear that these cut-off functions are indeed $(h^\sigma_j)$-compatible whenever $Z_\beta(\sigma)\neq \emptyset$. 
\end{proof}

With this in mind, we can define a tailored version of the truncations $f^\sigma_{\beta,j}(z)$ from Section \ref{section:toric-cpct1} as
\begin{equation*}
    f_{s,\beta,j}^\sigma(z)\coloneqq\sum_{\alpha \in P^\sigma_j}{(1-s+s\chi^\sigma_{\alpha,\beta,j}(u))f_{\alpha,\beta,j}(z)},
\end{equation*}
for $j=1,\dots,r$, with the corresponding two-parameter family of complete intersections inside $(M/\sigma)_{\Cs}$. Analogously to Proposition \ref{proposition:int-strata2}, one checks that $f^\sigma_{s,\beta,j}(z) \cdot z^{-\alpha^\sigma_j}$ are the defining equations for $Z_{s,\beta}(\sigma)\coloneqq \overline{Z}_{s,\beta} \cap O(\sigma)$.

\begin{theorem}\label{theorem:cpctification_ci_tailored}
     Consider a one-parameter family $Z_{s,\beta}$ of complete intersections associated to the data of $P=P_1+\dots+P_r$, refined triangulating functions $(h_j)$ and $\check{\Sigma}$-proper $(h_j)$-compatible cut-off functions $\chi$ for a desingularisation given by a simplicial refinement $\check{\Sigma}$. Then there exists a $\beta_0>0$ such that $Z_{s,\beta}$ admits an orbifold smooth toric compactification inside $\mX_{\check{\Sigma}}$ for all $\beta>\beta_0$ and $s \in [0,1]$. 
\end{theorem}

\begin{proof}
    The proof of Corollary \ref{corollary:smooth-cpctification} goes through with the changes outlined above if we replace Lemma \ref{lemma:sm_nontailored_ci} by Proposition \ref{proposition:isotopy_ci}.
\end{proof}

Combining this theorem with Corollary \ref{corollary:approx_hol_ci} and some general results from \cite[Appendix B]{ghhps} gives us a way to construct a Liouville domain from the family $Z_{s,\beta}$ and also shows that, up to isomorphism, it does not depend on the ad hoc choices we have made along the way, such as the tailoring functions or the adapted potential $\varphi$.

\begin{remark}\label{remark:orbifold extensions}
    Although the results from \cite[Appendix B]{ghhps} are phrased in the setting of smooth compactifications, their proofs extend to orbifold smooth compactifications: in fact, smoothness is only used in loc. cit. to perform an explicit coordinate computation in Lemma B.18 to show that manifolds admitting smooth toric compactifications are `asymptotically cylindrical', which generalises to the more general case as follows. For a non-smooth simplicial cone $\sigma$ with an affine chart $U(\sigma) \cong [\C^{\dim(\sigma)}/G] \times (\Cs)^{n-\dim(\sigma)}$, we have a quotient map $\C^{\dim(\sigma)} \times (\Cs)^{n-\dim(\sigma)} \rightarrow U(\sigma)$. On the dense tori, this corresponds to a degree $|G|$ cover $(\Cs)^n \rightarrow (\Cs)^n$, with the property that it intertwines the Euler vector fields (the cover is induced by an integral linear map in $GL(M_\R)$ with determinant $|G|$, so the pushforward of the radial vector field is the radial vector field). Therefore, one can use the coordinate computation on a smooth chart to show that a suitable lift of the orbifold smooth compactification is asymptotically cylindrical, and asymptotic cylindricity of the quotient follows. 
\end{remark}

\begin{corollary}\label{corollary:tailored_ci_LD}
    Let $\varphi$ be the degree $2$ homogeneous Kähler potential for the symplectic form $\omega$ on $M_{\Cs}$. Then for $\beta>0$ and $C>0$ sufficiently large, the compact manifold with boundary $A_{s,\beta,C}=Z_{s,\beta} \cap \{ \varphi \leq C \}$ equipped with the restriction Liouville form $\lambda|_{A_{s,\beta,C}}$ is a Liouville domain for all $s \in [0,1]$. In particular, its Liouville completion $\Hat{A}_{s,\beta,C}$ is independent of $C$, so we shall omit it from the notation. 
\end{corollary}

\begin{proof}
    By Corollary \ref{corollary:approx_hol_ci}, we get that for $\beta$ large enough, $Z_{s,\beta}$ is going to be $\frac17\pi$-almost holomorphic. Hence, by Theorem \ref{theorem:cpctification_ci_tailored}, we can apply \cite[Lemma B.3]{ghhps} to get the desired conclusion. 
\end{proof}

\begin{corollary}\label{corollary:isotopy_ci}
    The isotopy $Z_{s,\beta}$ between $Z_{\beta}=Z_{0,\beta}$ and $\widetilde{Z}_{\beta}=Z_{1,\beta}$ provides an exact symplectomorphism between Liouville manifolds $\Hat{A}_{0,\beta}$ and $\Hat{A}_{1,\beta}$ for all $\beta>0$ large enough. 
\end{corollary}

\begin{proof}
    As above, note that for large $\beta$, all the submanifolds $\frac17\pi$-approximately holomorphic, so we can apply \cite[Lemma B.4]{ghhps}. 
\end{proof}

Finally, we can also relate the potentials that we consider to the natural \emph{compactification potential} $\varphi_c$, which is a term from \cite{ghhps} synonymous with Seidel's notion of a potential associated to a divisor complement from \cite{based-seidel}. 

\begin{corollary}[Theorem \ref{theorem:main-thm}(\ref{main-thm-2})]\label{corollary:tailored-LD-cpctification-LD}
    The Liouville manifold $\widehat{(A_{1,\beta},\lambda)}$ associated to $\widetilde{Z}_\beta$ with the degree two homogeneous potential $\varphi$ is isomorphic to $\widehat{(A_{c},\lambda_c)}$, the Liouville manifold associated to $Z_\beta$ with the compactification potential $\varphi_c$. 
\end{corollary}

\begin{proof}
    Combine Corollary \ref{corollary:isotopy_ci} with \cite[Lemma B.2]{ghhps} (which applies because of Corollary \ref{corollary:smooth-cpctification}). 
\end{proof}  
\section{Computations of skeleta} \label{section:skeleta}
In this section, we perform to the main calculation and describe the skeleton of open Batyrev--Borisov complete intersections. Our strategy is similar to the approach taken in \cite{Zhou} for the case of hypersurfaces, but some of the details differ. In particular, by setting $r=1$ throughout the section, we obtain an alternative proof of Zhou's main theorem in the case when $\Delta^\vee$ is reflexive (see Remark \ref{remark:generality} for an explanation of why we chose this level of generality). Also, note that we do not require the potential $\varphi$ to be \emph{strongly} adapted in this section, the extra control becomes useful only when we study the combinatorics of the skeleton in Section \ref{section:combinatorics}. 

\subsection{Set-up and notation}\label{section:skeleta-setup}
Throughout this section, suppose that we have fixed an irreducible nef partition of a reflexive polytope $\nabla=\nabla_1+\dots+\nabla_r$, a centred refined triangulating function $h$ on $\Delta^\vee$ inducing triangulations $\mT$ of $\partial\Delta^\vee$ and $\mT_j$ of $\partial \nabla_j$, constants $c_{\alpha,j}$ such that $c_{0,j}<0$ and $c_{\alpha,j}>0$ for all $\alpha \in \nabla_j\backslash\{0\}$ and all $j$, and a collection of $(h_j)$-compatible cut-off functions $\chi=(\chi_1,\dots,\chi_r)$ that is proper, localising and centred.

These data determine a family of open Batyrev--Borisov complete intersections $Z_{s,\beta}$ and, in particular, a family of tailored complete intersections $\widetilde{Z}_{\beta}$. We also have associated hypersurfaces $\widetilde{H}_{\beta,j}$ for $j=1,\dots,r$ that are used to define the complete intersections and the hypersurface $\widetilde{H}_{\beta,\tot}$ defined in Section \ref{section:tailor-bbci}. Denote the associated hypersurface amoebae as $\widetilde{\mA}_{\beta,j}$, $\widetilde{\mA}_{\beta,\tot}$ respectively, and the amoeba of $\widetilde{Z}_{\beta}$ as $\widetilde{\mA}_{\beta}$. Consider also the tropical hypersurfaces $\mA_{\trop,j}$ for $j=1,\dots,r$ (coming from functions $h|_{\nabla_j}$), $\mA_{\trop,\tot}$ (coming from $h$) and the tropical Batyrev--Borisov complete intersection $\mA_{\trop}=\bigcap_{j=1}^r \mA_{\trop,j}$. We fix a Kähler potential $\varphi$ adapted with respect to the amoeba complement $C_{0,\trop,\tot}$, which gives us a strict isomorphism of Liouville domains $\Phi \colon M_{\Cs} \rightarrow T^*M_{S^1}$ by Proposition \ref{proposition:capital-phi-iso}. 

By Proposition \ref{proposition:isotopy_ci}, we know that $\widetilde{Z}_{\beta}$ will be a smooth $(2n-2r)$-dimensional exact submanifold of $(M_{\Cs},\lambda_{\varphi})$ for all $\beta$ large enough, hence it makes sense to speak of its skeleton. 

Recall that the \emph{FLTZ skeleton} associated to this set-up is given by 
\begin{equation*}
   \eL_{\Sigma}=\bigcup_{\sigma \in \Sigma}{\sigma^\perp\times\sigma} \subset T^*M_{S^1},
\end{equation*}
where $\Sigma$ is a simplicial fan over $\mT$ with cones $\sigma=\cone(S)$ for $S \in \mT$. Note that we shall sometimes write the Lagrangian $\sigma^\perp \times \sigma$ as $S^\perp \times \cone(S)$ instead, to emphasise that we also have a canonical choice of ray generators as the vertices $\alpha \in S$. For convenience, we will also denote 
\begin{equation*}
\begin{split}
    B_{\beta,j}&\coloneqq \Log_{\beta}^{-1}(\partial \widetilde{C}_{0,\beta,j}) \textnormal{ for } j=1,\dots, r,\\ B_{\beta,\tot}&\coloneqq\Log^{-1}_{\beta}(\partial \widetilde{C}_{0,\beta,\tot}) \textnormal { and } B_{\beta}\coloneqq \bigcap_{j=1}^r B_{\beta,j}.
\end{split}
\end{equation*}
With this in mind, we are ready to state the main result:
\begin{theorem}[Theorem \ref{theorem:main-thm}(\ref{main-thm-3})]\label{theorem:skeleton_ci}
    For all $\beta>0$ large enough, the Liouville skeleton of $\widetilde{Z}_\beta$ is given as:
    \begin{equation*}
        \skel(\widetilde{Z}_{\beta})=\Phi^{-1}(\eL_{\Sigma})\cap B_{\beta}.
    \end{equation*}
\end{theorem}

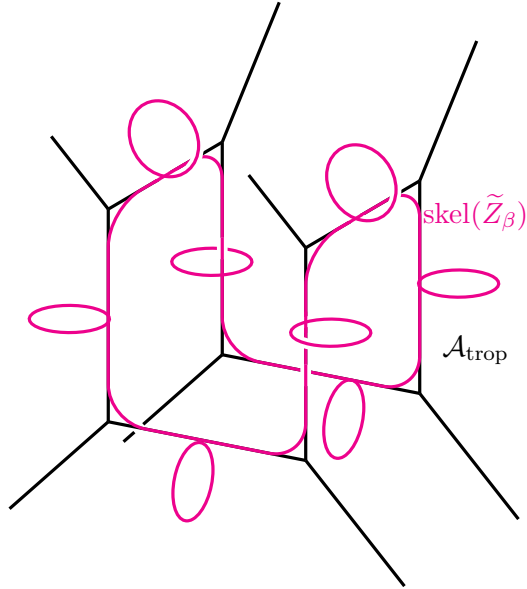
\begin{figure}[ht]
    \centering
    \tdplotsetmaincoords{70}{120} 
    \begin{tikzpicture}[tdplot_main_coords, scale=1.5]
        \draw[very thick, black] (1,1,1) -- (-1,1,1);
        \draw[very thick, black] (-1,1,-1) -- (-1,1,1);
        \draw[very thick, black] (1,1,1) -- (1,1,-1);
        \draw[very thick, black] (-1,-1,1) -- (1,-1,1);
        \draw[very thick, black] (-1,-1,1) -- (-1,-1,-1);
        \draw[very thick, black] (1,-1,1) -- (1,-1,-1);
        \draw[very thick, black] (-1,1,-1) -- (-1,-1,-1);
        \draw[very thick, black] (1,1,-1) -- (1,-1,-1);
        
        \draw[very thick,black] (1,1,1) -- (2,1,2);
        \draw[very thick,black] (1,-1,1) -- (2,-1,2);
        \draw[very thick,black] (-1,1,1) -- (-2,1,2);
        \draw[very thick,black] (-1,-1,1) -- (-2,-1,2);
        \draw[very thick,black] (1,1,-1) -- (1,2,-2);
        \draw[very thick,black] (1,-1,-1) -- (1,-2,-2);
        \draw[very thick,black] (-1,1,-1) -- (-1,2,-2);
        \draw[very thick,black] (-1,-1,-1) -- (-1,-1.8,-1.8);
        \draw[very thick,black] (-1,-1.89,-1.89) -- (-1,-2,-2);
    
        \draw[very thick, magenta, rounded corners = 5mm] (1,1,1) -- (-1,1,1) -- (-1,1,-1) -- (-1,-1,-1) -- (-1,-1,1) -- (1,-1,1) -- (1,-1,-1)-- (1,1,-1) -- cycle;
        \draw[line width=4pt, white] (1,1,-0.1)--(1,1,-0.3);
        \draw[very thick, magenta] (1,1,-0.1)--(1,1,-0.3);
        \begin{scope}[canvas is xy plane at z=0.5]
            \draw[line width=4pt, color=white] (1.6,1.6) circle (10pt);
            \draw[very thick, color=magenta] (1.6,1.6) circle (10pt);
        \end{scope}
        \draw[very thick,magenta] (1,1,0.4)--(1,1,0.2);
        
        \begin{scope}[canvas is xy plane at z=0]
            \draw[very thick, color=magenta] (-1.25,1.25) circle (10pt);
        \end{scope}
        \begin{scope}[canvas is xy plane at z=-0.25]
            \draw[very thick, color=magenta] (-1.25,-1.25) circle (10pt);
        \end{scope}
        \draw[line width=4pt, white] (-1,-1,-0.1)--(-1,-1,0.1);
        \draw[very thick, magenta] (-1,-1,-0.1)--(-1,-1,0.1);
        \begin{scope}[canvas is xy plane at z=0]
            \draw[very thick, color=magenta] (1.25,-1.25) circle (10pt);
        \end{scope}
        \begin{scope}[canvas is yz plane at x=0.37]
            \draw[line width=4pt, white] (1.2,1.45) circle (10pt);
            \draw[very thick, color=magenta] (1.2,1.45) circle (10pt);
        \end{scope}
        \draw[very thick, magenta] (0,1,1)--(0.5,1,1);
        \begin{scope}[canvas is yz plane at x=0.37]
            \draw[line width=4pt, white] (-0.8,1.5) circle (10pt);
            \draw[very thick, color=magenta] (-0.8,1.5) circle (10pt);
        \end{scope}
        \draw[very thick, magenta] (0.2,-1,1)--(-0.35,-1,1);
        \begin{scope}[canvas is xz plane at y=-0.2]
            \draw[very thick, color=magenta] (0.9,-1.45) circle (10pt);
        \end{scope}
        \begin{scope}[canvas is xz plane at y=0]
            \draw[very thick, color=magenta] (-1.4,-1.55) circle (10pt);
        \end{scope}
    
        \node[label={$\mathcal{A}_{\textnormal{trop}}$}] at (-0.25,2,-0.5) {};
        \node[label={$\color{magenta}\textnormal{skel}(\widetilde{Z}_\beta)$}] at (-0.25,2,0.75) {};
    \end{tikzpicture}
    \caption{The skeleton associated to the open BBCI $Z_\beta=\{z \in (\Cs)^3 \colon e^{-\beta}(z_1+z_1^{-1}+z_3)-1=0, e^{-\beta}(z_2+z_2^{-1}+z_3^{-1})-1=0\}$ from Example \ref{example:running-trop}}
    \label{fig:skeleton}
\end{figure}

\begin{remark}\label{remark:main-thmiii-follows}
    To deduce the exact statement of Theorem \ref{theorem:main-thm}(\ref{main-thm-3}) from this, it suffices to pin down the topology of $\Phi(B_\beta)$, which is a simple application of our results on smoothing from Appendix \ref{section:appendix-smoothing}: by definition, we have $\Phi(B_\beta)\cong M_{S^1} \times \Phi(\partial\widetilde{C}_{0,\beta})$ for $\partial\widetilde{C}_{0,\beta}\coloneqq \bigcap_{j=1}^r \partial\widetilde{C}_{0,\beta,j}$. Combining our Proposition \ref{proposition:smoothing-bbci} with \cite[Theorem 2.5]{HZ} (saying that $\partial C_{0,\trop} \coloneqq \bigcap_{j=1}^r \partial C_{0,\trop,j}$ is a topological $(n-r)$-sphere) proves that $\partial \widetilde{C}_{0,\beta}$ is a smooth $(n-r)$-dimensional sphere inside $M_\R$, which is precisely what we need. See Figure \ref{fig:skeleton} for an illustration of what such skeleta look like. 
\end{remark}

\subsection{Outline of the proof}\label{section:skeleta-outline}

We give an outline of the proof of Theorem \ref{theorem:skeleton_ci}, highlighting where some of the assumptions from the set-up are used. Note that it is also instructive to see how our arguments work in the case $r=1$ where the entire setting is more standard, even though some of them become unnecessarily complicated. To simplify notation in our calculations, denote $\Lambda_\beta\coloneqq\Phi^{-1}(\eL_{\Sigma})\cap B_{\beta}$.

The proof is based on a few observations: first, it is particularly important that the cut-off functions $\chi$ are \emph{localising} to determine which terms in the defining equations $\widetilde{f}_{\beta,j}$ are non-zero (this allows us to perform calculations locally, as the name suggests). More precisely, given a point $z \in \widetilde{Z}_\beta$, denote 
\begin{equation*}
    \begin{split}
        T_j(z) &\coloneqq \{ \alpha \in \nabla_j \colon \chi_{\alpha,\beta,j}(z)>0\},\\
        T'_j(z)&\coloneqq \{\alpha \in \nabla_j \colon \chi_{\alpha,\beta,j}(z)=1 \},
    \end{split}
\end{equation*}
for $j=1,\dots,r$. Then $\alpha \in T_j$ are the only terms that actually contribute to the defining equation $\widetilde{f}_{\beta,j}$ and we also have the following:

\begin{lemma}\label{lemma:localising-important}
     Then there exists a constant $\beta_0>0$ such that for all $\beta>\beta_0$ and all $z \in \widetilde{Z}_\beta$, the tuples $(T_1(z),\dots,T_r(z))$ and $(T_1'(z),\dots,T'_r(z))$ label non-empty cells of $\mA_\trop$. In particular, $T\coloneqq \conv(T_1,\dots,T_r)$ and $T'\coloneqq\conv(T_1',\dots,T_r')$ are simplices in $\mT*0$.
\end{lemma}

\begin{proof}
    Since the functions $\chi$ are localising, there exists a constant $K>0$ such that the inequality $L_{h_j}(u)-l_\alpha(u) \leq \beta^{-\frac12}+K\beta^{-1}$ holds for all $\alpha \in T_j$ and $j=1,\dots,r$. If we denote the set of $\alpha \in \nabla_j$ satisfying this inequality as $\widetilde{T}_j$, then one can show that $(\widetilde{T}_1,\dots,\widetilde{T}_r)$ labels a non-empty cell of $\mA_\trop$ analogously to the beginning of the proof of Lemma \ref{lemma:sm_nontailored_ci}. Since $\mT_j*0$ is a triangulation, $T'_j \subset T_j \subset \widetilde{T}_j$ by definition and $|T'_j(z)| \geq 2$ by $\widetilde{f}_{\beta,j}(z)=0$, it follows that both $(T_1,\dots,T_r)$ and $(T'_1,\dots,T_r')$ label non-empty cells. Finally, using the characterisation of cells of $\mA_\trop$ from Proposition \ref{proposition:bdd_cells} and Corollary \ref{corollary:unbdd_cells} gives the last part of the statement. 
\end{proof}

Second, to prove that $\Lambda_\beta$ is preserved by the Liouville flow, we show that it is stratified by strongly exact Lagrangians (recall that a half-dimensional submanifold $L \subset (Y,\lambda)$ inside an exact symplectic manifold is called a \emph{strongly exact Lagrangian} if $\lambda|_L \equiv 0$). The key observation is the following standard fact:

\begin{lemma}\label{lemma:strongly-exact-in-skeleton}
    Suppose that $(Y,\lambda)$ is a Liouville manifold and $L \subset Y$ is a closed, strongly exact Lagrangian. Then $L \subset \skel(Y,\lambda)$.
\end{lemma}
\begin{proof}
    Denote the Liouville vector field of $Y$ as $X_\lambda$. Pick a point $p \in L$ and a tangent vector $v \in T_pL$, then strong exactness implies $0=\lambda_p(v)=\omega_p(X_\lambda(p),v)$. Varying $v$ tells us that $X_\lambda(p) \in (T_pL)^\omega$, but since $L$ is a Lagrangian, this is equivalent to $X_\lambda(p)$ being tangent to $L$. Therefore, the flow of $X_\lambda$ preserves $L$, which implies that it does not escape to infinity and, in particular, has to be a subset of the skeleton.
\end{proof}

One can check that the FLTZ Lagrangian $\eL_\Sigma$ is strongly exact and the property $\lambda|_L\equiv 0$ is clearly inherited by any submanifold of $L$, which makes this result relevant to the problem at hand.

\begin{proof}[Proof of Theorem \ref{theorem:skeleton_ci}]
    We start by checking that $B_\beta \cong \partial\widetilde{C}_{0,\beta} \times M_{S^1}$ is indeed a manifold of the expected dimension (Lemma \ref{lemma:ctct-hyp}). 
    
    Since the potential $\varphi$ is \emph{adapted} to $C_{0,\trop,\tot}$ and $\widetilde{C}_{0,\beta}$ is close to the tropical complete intersection $\partial C_{0,\trop} \subset \partial C_{0,\trop,\tot}$ by construction, we can use the condition $\Phi(C_{(T_1,\dots,T_r)}) \subset \st(\cone(T_1,\dots,T_r))$ from adaptedness to prove that the summands $\widetilde{f}_{\alpha,\beta,j}(z)$ are all non-negative for $\alpha \in \nabla_j \backslash \{0\}$ and $z \in \Lambda_\beta$, which means that $\Lambda_\beta \subset \widetilde{Z}_\beta$ (Lemma \ref{lemma:contains}). 
    
    Closeness to the tropical complete intersection also implies that $\Phi(\partial\widetilde{C}_{0,\beta})$ intersects the cones $\sigma \in \Sigma$ transversely and avoids non-transversal cones. Therefore, $\Lambda^0_\beta \subset M_{\Cs}$ can be stratified by $\Lambda_\beta(\sigma)\coloneqq \Phi^{-1}(\sigma^\perp \times \sigma) \cap B_\beta \cong \sigma^\perp \times (\Phi^{-1}(\sigma) \cap \partial \widetilde{C}_{0,\beta})$ for $\sigma \in \Sigma^\trans$. Since $\eL_\Sigma$ is a strongly exact Lagrangian, each $\Lambda_\beta(\sigma)$ will be a strongly exact Lagrangian with corners, so by applying Lemma \ref{lemma:strongly-exact-in-skeleton}, we see that $\Lambda_\beta \subset \skel(\widetilde{Z}_\beta)$ (see Corollary \ref{corollary:contains}).

    Finally, we need to check that the Liouville vector field of $\widetilde{Z}_\beta$ has no zeroes away from $\Lambda_\beta$ (which is an analogue of \cite[Proposition 5.9]{Zhou}). This is the most computationally involved part of the proof, presented as Lemma \ref{lemma:no-zeroes}. The idea behind it is as follows: to demonstrate that $X_\lambda^{\widetilde{Z}}(z) \neq 0$ for $z \notin \Lambda_\beta$, we produce a vector $v \in T_z\widetilde{Z}_\beta$ such that $\omega(X_\lambda^{\widetilde{Z}}(z),v)=\langle\lambda(z),v\rangle \neq 0$. Consider a \enquote{projection} $\pi_{\aff}(u)$ of $u\in \mA_\trop(\beta)$ onto an appropriate cell $C \subset \mA_\trop$ that is chosen depending on the exact position of $u$ (following Corollary \ref{corollary:aff-nbhd-structure}). Then, we connect $\pi_{\aff}(u)$ to a minimum $u_{\min}$ over some other cell $C'$, so that the vector $v\coloneqq \frac{\pi_\aff(u)-u_{\min}}{\lVert\pi_\aff(u)-u_{\min} \rVert}$ is tangent to a cell of $\mA_\trop$. We can then apply the estimate from Lemma \ref{lemma:adapted_bound} to get a lower bound on the quantity $\langle d\varphi(u),v\rangle$. The cells are picked so that the lift $Jv$ is almost tangent to $\widetilde{Z}_\beta$ at $z$, and we use uniform estimates similar to Lemma \ref{lemma:cursed_determinant} to find a nearby tangent vector. This strategy is applied in Case \ref{case:easy-noncpt}, \ref{case:hard-noncpct} and \ref{case:easy-cpct} with minor variations. Finally, in the remaining Case \ref{case:hard-cpct}, we can observe that the defining functions for $\widetilde{Z}_\beta$ will be holomorphic near $z$, which allows us to perform a direct calculation showing that $\lambda(z)|_{\widetilde{Z}} \neq 0$.
\end{proof}

\subsection{Calculations}\label{section:skeleta-calculations} In this section, we repay the technical debt from Section \ref{section:skeleta-outline} and fill in the details from the proof of Theorem \ref{theorem:skeleton_ci}.

\begin{lemma} \label{lemma:ctct-hyp}
    $B_{\beta}$ is a smooth $(2n-r)$-dimensional submanifold of $M_{\Cs}$ nowhere tangent to the Liouville vector field $X_{\lambda}$.
\end{lemma}

\begin{proof}
    We first show that $\partial \widetilde{C}_{0,\beta}$ is a smooth codimension $r$ submanifold of $M_{\R}$: by Proposition \ref{proposition:bdary_eqn}, we can write $\widetilde{C}_{0,\beta,j}=\{u \in M_{\R} \colon G_j(u) \leq|c_{0,j}|\}$ using the functions $G_j(u)=\widetilde{F}_{\beta,j}(u)-\widetilde{F}_{0,\beta,j}(u)\coloneqq\sum_{\alpha \in \nabla_j\backslash\{ 0\}}{|\widetilde{f}_{\alpha,\beta,j}(u)|}$. 

    For a point $u \in \partial \widetilde{C}_{0,\beta}$, we can use the real equations $\widetilde{F}_{\beta,j}$ in place of the complex $\widetilde{f}_{\beta,j}$ to define $T_j(u)$ and $T'_j(u)$ for $j=1,\dots,r$ as in Lemma \ref{lemma:localising-important}, show that they are both simplices in $\mT_j*0$ and that $T=\conv(T_1,\dots,T_r)$, $T'=\conv(T_1',\dots,T'_r)$ are simplices in $\mT$ (note that $G_j(u)=|c_{0,j}|$ automatically implies that $\chi_{0,\beta,j}(u)=1$ for all $j$, so we always have $0 \in T'_j(u)$). Also, by definition of $T_j$, we can write $G_j(u)$ as a sum over $\alpha \in T_j(u) \backslash \{0\}$ rather than over the entire polytope $\nabla_j \backslash \{0\}$.

    We can use estimates from Lemma \ref{lemma:partials_lin_indep} (with a slight variation that we are computing differentials of real functions rather than holomorphic differentials of complex functions, so the situation simplifies) to prove that $D_uG_j$ is $O(e^{-\sqrt{\beta}})$-close to a point in $\relint(\cone(T'_j))$ by bounding the contributions from $T_j \backslash T'_j$. This allows us to conclude the the vectors $D_u
    G_1$, \dots, $D_uG_r$ are linearly independent analogously to the aforementioned Lemma. Therefore, the point $(|c_{0,1}|,\dots,|c_{0,r}|)$ is a regular value of $(G_1,\dots,G_r)$, which means that $\partial\widetilde{C}_{0,\beta}$ is a manifold of codimension $r$. Analogously, we observe that $|c_{0,\tot}|$ is a regular value of the function $G=\sum_j G_j$ (which is proper and convex by the assumption that $\chi$ is centred), so $\widetilde{C}_{0,\beta,\tot}$ is a full-dimensional convex set containing the origin in its interior with smooth boundary $\partial \widetilde{C}_{0,\beta,\tot}$. 
    
    By \cite[Proposition 2.7]{Zhou}, we know that $(\Log_\beta)_*(X_\lambda)$ is a positive rescaling of the radial vector field on $M_\R$. Since $\widetilde{C}_{0,\beta,\tot}$ is convex, any ray starting from the origin intersects its boundary at most once, so $(\Log_\beta)_*(X_\lambda)$ must be positively transverse to $\partial \widetilde{C}_{0,\beta,\tot}$. In particular, since $\partial \widetilde{C}_{0,\beta} \subset \partial\widetilde{C}_{0,\beta,\tot}$, the vector field is also nowhere tangent to $\partial \widetilde{C}_{0,\beta}$. The desired statements now follow by taking preimages of the set-up inside $M_\R$ under the submersion $\Log_\beta \colon M_{\Cs} \rightarrow M_\R$. 
\end{proof}

\begin{lemma} \label{lemma:contains}
    For sufficiently large $\beta>0$, $\Lambda_\beta$ is contained in $\widetilde{Z}_{\beta}$.
\end{lemma}

\begin{proof}
    Suppose that we have a point $z = (u, \theta) \in M_{\Cs}$ such that $u \in \partial\widetilde{C}_{0,\beta}$ and $\Phi(z) \in \eL_{\Sigma}$. 

    Define $T_j(u)$, $T(u)$ analogously as in the proof of Lemma \ref{lemma:ctct-hyp} and let $\tau\coloneqq\cone(T)$. Recall that $\widetilde{F}_{\alpha,\beta,j}(z)$ for $\alpha \in T_j(u)$ are the only non-zero monomials in $\widetilde{F}_{\beta,j}(z)$, hence we have:
    \begin{equation*}
        \sum_{\alpha \in T_j(u)\backslash \{ 0\}}{\widetilde{F
        }_{\alpha,\beta,j}(z)}=|c_{0,j}|,
    \end{equation*}
    for all $j=1$, \dots, $r$. 
    
    By Proposition \ref{proposition:bdd_cells}, we know that the cell $C_{(T_1,\dots,T_r)}$ of $\mA_{\trop}$ is just a face $C_{T,\tot}$ of the polytope $C_{0,\trop,\tot}$. Since the potential $\varphi$ is adapted to $C_{0,\trop,\tot}$, we see that $d\varphi(C_{T,\tot}) \subset \st(\tau)$. The inequalities $L_{h_j}(u)-l_\alpha(u) \leq \beta^{-\frac12}+K\beta^{-1}$ for $\alpha \in T_j(u) \backslash \{0\}$, all $j$ and some constant $K>0$ coming from $\chi_{\alpha,\beta,j}(u)>0$ tell us that the affine distance (cf. Appendix \ref{section:appendix-polyhedra}) of $u$ from $C_{(T_1,\dots,T_r)}$ is bounded by a constant multiple of $\beta^{-\frac12}$. Therefore, by Lemma \ref{lemma:aff-vs-norm}, we have $d(u,C_{T,\tot})=O(\beta^{-\frac12})$, so a compactness argument shows that $d\varphi(u) \in \st(\tau)$ for $\beta \gg0$. From $\Phi(z) \in \eL_{\Sigma}$, this means that $\Phi(z) \in S^{\perp} \times \sigma$ for some $\sigma=\cone(S)$ and $T(u) \subset S$, which gives $\theta \in S^{\perp} \subset T(u)^\perp$. Therefore, $z^{\alpha}=
    |z^{\alpha}|$ for all $\alpha \in \nabla_j \backslash \{ 0 \}$, so $\widetilde{f}_{\alpha,\beta,j}(z)=|\widetilde{f}_{\alpha,\beta,j}(z)|$ and
    \begin{equation*}
        \widetilde{f}_{\beta,j}(z) = \left(\sum_{\alpha \in T_j(u)\backslash\{0\}}{\widetilde{f}_{\alpha,\beta,j}(z)}\right)+c_{0,j} =\left(\sum_{\alpha \in T_j(u)\backslash \{0\}}{\widetilde{F}_{\alpha,\beta,j}(z)}\right)-|c_{0,j}|= 0,
    \end{equation*}
    which means that $z \in \widetilde{Z}_{\beta}$. 
\end{proof}

\begin{remark}\label{remark:int_fltz}
    Our proof shows that only the strata $\sigma^{\perp} \times \sigma$ where $\sigma \in \Sigma$ is transversal contribute to the skeleton of an open BBCI, so we could use a smaller Lagrangian
    \begin{equation*}
        \eL_{\Sigma}^{\trans}\coloneqq\bigcup_{\sigma \in \Sigma^\trans}{\sigma^{\perp}\times\sigma} \subset \eL_{\Sigma},
    \end{equation*}
    that only accounts for transversal cones instead of the entire FLTZ skeleton. 
\end{remark}

\begin{remark}\label{remark:positive_locus_ci}
    The proof also shows that $\Lambda_\beta$ lies inside the \textit{positive real locus of $\widetilde{Z}_{\beta}$} $\widetilde{Z}^+_\beta\coloneqq\widetilde{Z}_\beta \cap B_{\beta}$, i.e. the set of all $z \in \widetilde{Z}_{\beta}$ such that $\widetilde{f}_{0,\beta,j}(z)=c_{0,j}$ and $\widetilde{f}_{\alpha,\beta,j}(z) \geq 0$ for all $j=1,\dots,r$ and $\alpha \in \nabla_j \backslash \{0\}$.
\end{remark}

We shall now study the dynamics of the Liouville vector field $X^{\widetilde{Z}}_{\lambda}$: it is a standard fact that for $z \in \widetilde{Z}_{\beta}$, one can write it in terms of the ambient Liouville vector field $X_\lambda$ as
\begin{equation*}
    X_{\lambda}(z)=X^{\widetilde{Z}}_{\lambda}(z)+X^{\perp}_{\lambda}(z),
\end{equation*}
with $X^{\widetilde{Z}}_{\lambda}(z)$ equal to the symplectic orthogonal projection of $X_{\lambda}(z)$ onto $T_z\widetilde{Z}_{\beta}$ and $X^{\perp}_{\lambda}(z)$ symplectically orthogonal to $T_z\widetilde{Z}_{\beta}$. We record the following result that generalises \cite[Proposition 5.3]{Zhou} by expressing the Liouville vector field over the positive locus of $\widetilde{Z}_\beta$ in terms of of Hamiltonian vector fields of imaginary parts of $\widetilde{f}_{\beta,1}$, \dots, $\widetilde{f}_{\beta,r}$:

\begin{lemma}\label{lemma:liou-dynamics}
    For all $z\in \widetilde{Z}^+_{\beta}$ in the positive locus, there exist real numbers $c_j$ for $j=1,\dots,r$ such that at least one of them is non-zero and they satisfy
    \begin{equation*}
        X^{\perp}_{\lambda}(z)=\sum_{j=1}^r c_j X_{\IM(\widetilde{f}_{\beta,j})}(z).
    \end{equation*}
\end{lemma}

\begin{proof}
    Observe that the map
    \begin{equation*}
        d\widetilde{f}_{\beta}=(d\widetilde{f}_{\beta,1}, \dots, d\widetilde{f}_{\beta,r})\colon (T_z\widetilde{Z}_{\beta})^{\omega} \rightarrow \C^r
    \end{equation*}
    is an isomorphism by Proposition \ref{proposition:isotopy_ci} and both $X^{\perp}_{\lambda}$ and $X_{\IM(\widetilde{f}_{\beta,j})}$ for $j=1,\dots,r$ are symplectically orthogonal to $T_z\widetilde{Z}_{\beta}=\ker(d\widetilde{f}_{\beta}\colon T_zM_{\Cs}\rightarrow \C^r)$, so it suffices to prove the desired statement for pushforwards of the two vector fields.  
    
    Note that at such $z$, we have $\widetilde{f}_{\alpha,\beta,j}(z)=\widetilde{F}_{\alpha,\beta,j}(z)$ by the reasoning from Lemma \ref{lemma:contains}, so 
    \begin{equation*}
    \begin{split}
        d\widetilde{f}_{\beta,j}(z)&=\sum_{\alpha \in \nabla_j}{\widetilde{f}_{\alpha,\beta,j}(z)\langle \alpha,d\rho+id\theta\rangle+f_{\alpha,\beta,j}(z)d\chi_{\alpha,\beta,j}(z)}\\
        &=d\widetilde{F}_{\beta,j}(z)+i\sum_{\alpha \in \nabla_j}{\widetilde{F}_{\alpha,\beta,j}(z)\langle \alpha,d\theta\rangle}.
    \end{split}
    \end{equation*}
    Therefore, we can use the Riemannian metric $g$ induced by $\varphi$ to conveniently express the Hamiltonian vector fields as
    \begin{equation*}
        X_{\IM(\widetilde{f}_{\beta,j})}=g^{-1}\left(\sum_{\alpha \in \nabla_j} \widetilde{F}_{\alpha,\beta,j}(z)\langle \alpha,d\rho\rangle\right).
    \end{equation*}
    By bounds similar to Lemma \ref{lemma:partials_lin_indep}, there exists a constant $K_1>0$ such that the inequality 
    \begin{equation*}
        \left\lVert\nabla_g \widetilde{F}_{\beta,j}(z)-X_{\IM(\widetilde{f}_{\beta,j})}\right\rVert_g=\left\lVert\sum_{\alpha \in \nabla_j}{F_{\alpha,\beta,j}(z)d\chi_{\alpha,\beta,j}(z)}\right\rVert_{g} \leq K_1 e^{-\sqrt{\beta}},
    \end{equation*}
    holds for all $z \in \widetilde{Z}_\beta^+$ (the bound given there involves an extra term $\widetilde{F}_{\beta,j}(z)$, but this expression is equal to $2|c_{0,j}|$ along $\widetilde{Z}^+_{\beta}$, hence we can absorb it into the constant). Consider two real $(r \times r)$-matrices $M_{jk} \coloneqq \langle \nabla_g \widetilde{F}_{\beta,j}(z),\nabla_g\widetilde{F}_{\beta,k}(z) \rangle_g $ and $M'_{jk} \coloneqq \langle d\widetilde{f}_{\beta,j}(z),X_{\IM(\widetilde{f}_{\beta,k})} \rangle$. We can combine the above expressions to show that
    \begin{equation*}
    \begin{split}
        |M_{jk}(z)-M'_{jk}(z)|&= |\langle d\widetilde{F}_{\beta,j}(z),\nabla_g\widetilde{F}_{\beta,k}(z)-X_{\IM(\widetilde{f}_{\beta,k})}\rangle|\\
        &\leq \lVert d\widetilde{F}_{\beta,j}(u) \rVert_g \cdot K_1 e^{-\sqrt{\beta}} \leq K_2 e^{-\sqrt{\beta}},
    \end{split}
    \end{equation*}
    where we are using $z \in \widetilde{Z}^+_\beta$ to get an upper bound on $\lVert d\widetilde{F}_{\beta,j} \rVert_g$ by a constant (all the monomials in the expression are bounded along the positive locus). Since the entries of $M(z)$ are uniformly bounded over $\widetilde{Z}_\beta^+$, an argument analogous to Corollary \ref{corollary:cursed-determinant} tells us that the eigenvalues of $M(z)$ are bounded below by some $K_3>0$ for $z \in \widetilde{Z}^+_\beta$ (the lack of normalisation does not matter, because we are working over a compact set). Therefore, for all sufficiently large $\beta \gg0$, the remainder term $K_2 e^{-\sqrt{\beta}}$ will be smaller than $K_3$, so the matrix $M'$ will be invertible. This means that if we let $v \in \R^n$ be a vector with entries $v_j=\langle d\widetilde{f}_{\beta,j}(z), X^{\perp}_{\lambda} \rangle$, there will exist a unique vector $c=(c_1,\dots,c_r) \in \R^r$ such that $M'c=v$. By definition, it will then also satisfy
    \begin{equation*}
        \left\langle d\widetilde{f}_{\beta}, \sum_{k=1}^r{c_kX_{\IM(\widetilde{f}_{\beta,k})}} \right\rangle = \langle d\widetilde{f}_{\beta}, X^{\perp}_{\lambda} \rangle,
    \end{equation*}
    as desired. To see that $c \neq 0$, we observe that $v \neq 0$, because we have
    \begin{equation*}
    \begin{split}
        v_j&=\langle d\widetilde{f}_{\beta,j}(z),X^\perp_\lambda \rangle=\langle d\widetilde{f}_{\beta,j}(z),X_\lambda \rangle = \langle d\widetilde{f}_{\beta,j}(z),\nabla_g\varphi \rangle \\
        &= \langle d\widetilde{F}_{\beta,j}(z),\nabla_g \varphi \rangle=K_4\cdot \langle d\widetilde{F}_{\beta,j}(z),\sum_{a=1}^n \rho_a \partial_{\rho_a} \rangle >0,
    \end{split}
    \end{equation*}
    for some $K_4=K_4(z)>0$ by \cite[Proposition 2.7]{Zhou}, and the derivative of $\widetilde{F}_{\beta,j}$ along the radial vector field is positive by convexity of $\widetilde{C}_{0,\beta,j}$ observed in Lemma \ref{lemma:ctct-hyp}, so $v \in \R^n_{>0}$, which shows that $v\neq0$.
\end{proof}

As remarked in Section \ref{section:skeleta-outline}, we do not use this result in the main proof. Instead, we record the following Lemma, which is straightforward to prove by induction on $\dim(X)$.

\begin{lemma}\label{lemma:flows_corners}
    Let $M$ be a manifold with a complete vector field $v \in \Gamma(TM)$ and $X \subset M$ a closed submanifold-with-corners of $M$. Suppose that for all $p \in X$, $v(p)$ is tangent to each stratum of $X$ containing $p$. Then $X$ is preserved under the flow of $v$, i.e. for $p \in X$, we have $\Phi^t_v(p) \in X$ for all $t\geq0$. 
\end{lemma}

The main application of this is that it allows us to use Lemma \ref{lemma:strongly-exact-in-skeleton} for compact Lagrangians with corners, which is precisely what we do in the following Corollary:

\begin{corollary} \label{corollary:contains}
    We have $\Lambda_\beta \subseteq \skel(\widetilde{Z}_{\beta})$.
\end{corollary}

\begin{proof}
    We shall show that $\Lambda^0_\beta(\sigma)\coloneqq\Phi^{-1}(\sigma^\perp \times \sigma)\cap B_{\beta} \subset \skel(\widetilde{H}_\beta)$ for all $\sigma \in \Sigma^\trans$, where $\sigma=\cone(S)$ for $S \in \mT^\trans$ (we have already seen that the intersection is empty for non-transversal simplices in Lemma \ref{lemma:contains}). For simplicity, denote $Y_\sigma=\sigma^\perp \times \sigma$.

    First, we prove that the intersection $\Phi^{-1}(Y_\sigma)\cap B_{\beta}$ is transverse. The transversality statement is clearly equivalent to $\sigma$ intersecting $\Phi(\partial \widetilde{C}_{0,\beta})$ transversely inside $N_\R$, which is the same as
    $\dim((\Phi)_*(\ker(dG_1(u),\dots,dG_r(u))) + \R \sigma)=n$ for all $u$ in the intersection (where $G_j$ are the defining functions for $B_\beta$ from Lemma \ref{lemma:ctct-hyp}). We can use the Hessian metric $g$ induced by $\varphi$ to simplify this through $(\Phi)_*(\ker(dG_1(u),\dots,dG_r(u))=(dG_1(u),\dots,dG_r(u))^\perp_g$. By the proof of Lemma \ref{lemma:contains}, $d\varphi(u) \in \sigma$ implies that $\widetilde{f}_{\beta,j}(z)$ can only-have non-zero terms corresponding to $\alpha \in S_j$. This also implies that $dG_j(u) \in \sigma$ for all $j$, so $(\R \sigma)^\perp_g \subseteq (dG_1(u),\dots,dG_r(u))^\perp_g$ and the desired statement follows from $(\R \sigma)^\perp_g\oplus \R \sigma \cong \R^ n$. Therefore, by \cite[Theorem 6.4]{joyce}, $\Lambda^0_\beta(\sigma)$ is an $(n-r)$-dimensional manifold with corners. 
    
    Note that $\lambda_0$ is zero when restricted to $Y_\sigma$, hence it is also zero when restricted to $\Phi(\Lambda^0_\beta(\sigma))$, so it must be a strongly exact Lagrangian submanifold of $\Phi(\widetilde{Z}_\beta)$ with corners. By construction the map $\Phi$ satisfies $\Phi^*\lambda_0=\lambda_\varphi$, so $\Lambda^0_\beta(\sigma)$ is a strongly exact Lagrangian in $\widetilde{Z}_{\beta}$. In particular, this means that the Liouville vector field $X^{\widetilde{Z}}_\lambda$ is tangent to $\Lambda_\beta^0(\sigma)$ over its interior by our proof of Lemma \ref{lemma:strongly-exact-in-skeleton}. 

    By \cite[Proposition 6.7]{joyce}, we know that the boundary strata of $\Lambda^0_\beta(\sigma)$ are given by $\Phi^{-1}(Y_{\sigma,\tau}) \cap B_\beta$, where $Y_{\sigma,\tau}\coloneqq \sigma^\perp \times \relint(\tau) \subset \partial Y_\sigma$ are the boundary strata of $Y_\sigma$ indexed by $\tau \subsetneq \sigma$. But these can also alternatively be described as $\Lambda^0_\beta(\sigma) \cap \inte(\Lambda^0_\beta(\tau))$. The above considerations for a strongly exact Lagrangian $\Lambda^0_\beta(\tau)$ tell us that for all $z$ in such a boundary stratum, $X^{\widetilde{Z}}_\lambda \in T_z \Lambda^0_{\beta}(\tau)=T_z(\Phi^{-1})_*(\tau^\perp \times \relint(\tau)) \cap T_z B_\beta$. But we have already seen that the Liouville vector field also lies in $T_z\Lambda^0_\beta(\sigma)=T_z(\Phi^{-1})_*(\sigma^\perp \times \sigma)) \cap T_z B_\beta$, so the desired statement $X^{\widetilde{Z}}_\lambda \in T_z\Phi^{-1}(Y_{\sigma,\tau}) \cap T_zB_\beta$ follows by simple linear algebra.

    Therefore, the Liouville vector field is tangent to all the boundary strata of the compact manifold with boundary $\Lambda^0_\beta(\sigma)$, so we can deduce the desired statement from Lemma \ref{lemma:flows_corners}.
\end{proof}

Finally, we verify that there are no other zeroes of the Liouville vector field. Our approach is similar to \cite[Proposition 5.9]{Zhou}, with certain details spelled out more explicitly to ensure that everything indeed generalises to the complete intersection case.

\begin{lemma}\label{lemma:no-zeroes}
    The Liouville vector field $X_{\varphi}^{\widetilde{Z}}$ on $\widetilde{Z}_{\beta}$ has no zeroes away from $\Lambda_{\beta}$. 
\end{lemma}
\begin{proof}
    Suppose that $z$ is a zero of the Liouville vector field $X^{\widetilde{Z}}_{\lambda}$ (equivalently, $\lambda|_{\widetilde{Z}_{\beta}}=0$). As in Lemma \ref{lemma:localising-important}, we introduce the simplices $T'_j(z) \subseteq T_j(z) \subset\nabla_j$ that label non-empty cells of $\mA_\trop$ and $T(z), T'(z) \in \mT*0$. We can also consider the sets 
    \begin{equation*}
        T''_j(z)\coloneqq \{ \alpha \in \nabla_j \colon L_{h_j}(u)-l_\alpha(u) \leq \beta^{-\frac12}\},
    \end{equation*}
    and $T''(z)\coloneqq \conv(T_1'',\dots,T_r'')$. Since $\chi$ are $(h_j)$-compatible cut-off functions, we have $T''_j(z) \subseteq T'_j(z)$, so $T''_j$ are simplices in $\mT_j*0$ and $T'' \subseteq T'$ is a simplex of $\mT*0$. By considering the stratification of $\mA_\trop(\beta)$ introduced in Lemma \ref{lemma:aff-nbhd-structure} and the start of the proof of Lemma \ref{lemma:cursed_determinant}, we know that the simplices are defined so that $u \in C_{(T''_1,\dots,T''_r)}(\beta)$ and, in particular, $(T''_1,\dots,T''_r)$ also labels a non-empty cell of $\mA_\trop$. 
    
    Analogously, we know that the cut-off functions are localising, so $\chi_{\alpha,\beta,j}(u)>0$ forces $L_{h_j}(u)-l_\alpha(u) \leq \sqrt{2}\beta^{-\frac12}$ for all sufficiently large $\beta$, so we also know that $u$ has affine distance at most $\sqrt{2}\beta^{-\frac12}$ from $C_{(T_1,\dots,T_r)}$. Therefore, we also have a tuple of simplices $T_j \subseteq T^0_j$ given as
    \begin{equation*}
        T^0_j(z)\coloneqq \{ \alpha \in \nabla_j \colon L_{h_j}(u)-l_\alpha(u) \leq \sqrt{2}\beta^{-\frac12}\},
    \end{equation*}
    along with a compound simplex $T^0\coloneqq \conv(T^0_1,\dots,T^0_r)$. As above, one can see that $T^0_j$'s are chosen precisely so that we have $u \in C_{(T^0_1(z),\dots,T^0_r(z))}(\frac12\beta) \subset \mA_\trop(\frac12\beta)$. 
    
    We shall now distinguish several cases based on what the tuples $(T''_1,\dots,T''_r) \subseteq(T^0_1,\dots,T^0_r)$ are (in the language of \cite[Proposition 5.9]{Zhou}, this roughly corresponds to treating good regions and bad regions differently).
    \begin{case}\label{case:easy-noncpt}
        Assume that the cell $C_{(T^0_1,\dots,T^0_r)}$ is unbounded. Let $\overline{T}_j\coloneqq\conv\{0,T^0_j\}$, then by Corollary \ref{corollary:unbdd_cells}, we must have $T^0_j \subsetneq \overline{T}_{j}$ for some $j$ and $C_{(\overline{T}_1,\dots,\overline{T}_r)}$ is the largest bounded cell of $\mA_{\trop}$ that lies in the closure of $C_{(T^0_1,\dots,T^0_r)}$. By Proposition \ref{proposition:bdd_cells}, $C_{(\overline{T}_1,\dots,\overline{T}_r)}=C_{\overline{T}}$ is also a face of the polytope $C_{0,\tot}$ labelled by the simplex $\overline{T}\coloneqq\conv\{0,T\}$. Since $u \in C_{(T^0_1,\dots,T^0_r)}(\frac12\beta)$, Corollary \ref{corollary:aff-nbhd-structure} guarantees that we can find a point $u^0 \in C_{(T^0_1,\dots,T^0_r)} \subseteq \overline{C}_{(T_1,\dots,T_r)}$ that has affine distance at least $\sqrt{2}\beta^{-\frac12}$ from $\partial \overline{C}_{(T^0_1,\dots,T^0_r)}$ and satisfies $\rVert u-u^0 \rVert_{g_0} \leq K_1\cdot \beta^{-\frac12}$ for a constant $K_1$. Denote the unique minimum of $\varphi$ over $C_{\overline{T}}=C_{(\overline{T}_1,\dots,\overline{T}_r)}$ as $u_{\min}$. 
        
        Consider the vector $v\coloneqq\frac{u^0-u_{\min}}{\lVert u^0-u_{\min} \rVert}$, then we can apply Lemma \ref{lemma:adapted_bound} (since $u^0 \in \overline{C}_{T_j,j} \subset \mA_{\trop,j}$, we have $l_{\alpha}(u^0) \geq 0$ for all $\alpha \in T_j$ and $j=1,\dots,r$, so the Lemma applies) to get a bound
        \begin{equation*}
            \langle d\varphi(u^0),v \rangle = \frac{1}{\lVert u^0-u_{\min}\rVert} \langle d\varphi(u^0),u^0-u_{\min}\rangle \geq m\cdot\lVert u^0-u_{\min}\rVert + c\cdot\frac{d_{\aff}(u^0,C_{\overline{T}})}{\lVert u^0-u_{\min}\rVert},
        \end{equation*}
        for some positive constants $m$ and $c$ only depending on $\varphi$ and $\overline{T}$. By construction, we know that $d_\aff(u^0,C_{\overline{T}}) \geq d_{\aff}(u^0,\partial C_{(T^0_1,\dots,T^0_r)}) \geq \sqrt{2} \beta^{-\frac12}$. Therefore, by combining the above estimate with the AM-GM inequality, we get 
        \begin{equation*}
            \langle d\varphi(u^0),v \rangle \geq 2\sqrt{mcd_{\aff}(u^0,C_{\overline{T}})} \geq K_2 \cdot \beta^{-\frac14},
        \end{equation*}
        for some positive constant $K_2$. Since the map $d\varphi$ is Lipschitz (as it is $1$-homogeneous and continuously differentiable), $v$ is a unit vector by construction and we know that $u$ is $O(\beta^{-\frac12})$-close to $u^0$, it follows that there exists some constant $K_3>0$ such that 
        \begin{equation}\label{equation:no-zeroes}
            \langle d\varphi(u),v \rangle \geq K_2\beta^{-\frac14}-K_3\beta^{-\frac12}.
        \end{equation}
        Because $\lambda=d^c\varphi=d\varphi \circ J$, the inequality implies that $\langle \lambda(z), Jv \rangle \neq 0$ for $\beta$ large enough. Since $v$ is a vector connecting two points in $\overline{C}_{(T_1,\dots,T_r)}$, there must exist some constants $c_1, \dots, c_r$ such that $\langle \alpha, v \rangle = c_j$ for all $\alpha \in T_j$, which gives us 
        \begin{equation*}
            \langle d\widetilde{f}_{\beta,j}(z), Jv \rangle = i \sum_{\alpha \in T_j}{\widetilde{f}_{\alpha,\beta,j}(z)\langle \alpha, v \rangle} = ic_j\widetilde{f}_{\beta,j}(z)=0,
        \end{equation*}
        so $Jv \in \ker(d\widetilde{f}_{\beta,j}(z))$ for all $j=1,\dots,r$, which means that $Jv \in T_z\widetilde{Z}_{\beta}$, so $\lambda(z)|_{\widetilde{Z}_{\beta}} \neq 0$, a contradiction. 
    \end{case}
    \begin{case}\label{case:hard-noncpct}
        Assume that $C_{(T^0_1,\dots,T^0_r)}$ is bounded, but $C_{(T''_1,\dots,T''_r)}$ is an unbounded cell of $\mA_{\trop}$. Analogously to the previous case, we can take $\overline{T}_j\coloneqq\conv\{0,T''_j\}$, $\overline{T}=\conv\{0,T''\}$ and then project $u$ to a point $u'$ in $C_{(T''_1,\dots,T''_r)}$ to produce a unit vector $v$ that satisfies the same lower bound, since we know that the distance of $u$ from $C_{0,\trop,\tot}$ is at least $\beta^{-\frac12}$ (Equation \ref{equation:no-zeroes}). The key difference is that this time, we do not necessarily get $Jv \in T_z\widetilde{Z}_{\beta}$, since we only have $\langle \alpha, v \rangle =c_j$ for $\alpha \in T''_j(z)$ and there are some remainder terms corresponding to $\alpha \in T_j \backslash T''_j$:
        \begin{equation*}
        \begin{split}
            \langle d\widetilde{f}_{\beta,j}(z), Jv \rangle &= i \sum_{\alpha \in T_j}{\widetilde{f}_{\alpha,\beta,j}(z)\langle \alpha, v \rangle}\\ &= ic_j\left( \sum_{\alpha\in T_j''}{\widetilde{f}_{\alpha,\beta,j}(z)}\right)+i\sum_{\alpha\in T_j\backslash T''_j}{\widetilde{f}_{\alpha,\beta,j}(z)\langle \alpha, v \rangle} \\
            &= i\sum_{\alpha\in T_j\backslash T''_j}{\widetilde{f}_{\alpha,\beta,j}(z)(\langle \alpha, v \rangle-c_j)}.
        \end{split}
        \end{equation*}
        By definition of $T''_j(z)$, we have $L_{h_j}(u) \geq l_\alpha(u)+\beta^{-\frac12}$ for $\alpha \in T_j\backslash T''_j$, which gives us a bound $|\widetilde{f}_{\alpha,\beta,j}(z)| \leq e^{-\sqrt{\beta}}\widetilde{F}_{\beta,j}(z)$ (as in Lemma \ref{lemma:partials_lin_indep}). Therefore, there exists a constant $K_4>0$ independent of $z$ such that we have
        \begin{equation*}
            \left|\left\langle (\widetilde{F}_{\beta,j}(z))^{-1}\cdot d\widetilde{f}_{\beta,j}(z), Jv \right\rangle\right| \leq K_4 e^{-\sqrt{\beta}},
        \end{equation*}
        for all $z$ in this region and $j=1,\dots,r$. Consider the rescaling of $d\widetilde{f}_{\beta}$ given by $\ell(z)\coloneqq ((\widetilde{F}_{\beta,1}(z))^{-1}\cdot d\widetilde{f}_{\beta,1}(z),\dots,(\widetilde{F}_{\beta,r}(z))^{-1}\cdot d\widetilde{f}_{\beta,r}(z))$, which restricts to an isomorphism between $T_z\widetilde{Z}_\beta^\omega$ and $\C^r$ by Corollary \ref{corollary:approx_hol_ci}, so there exists a unique $Jv' \in T_z\widetilde{Z}_\beta^\omega$ such that $\ell(z)(Jv')=\ell(z)(Jv)$. Moreover, by repeating a discussion similar to Lemma \ref{lemma:cursed_determinant}, we can see that $\ell(z)$ can only take values in some fixed compact set $V \subset GL(\C,r)$. Therefore, the operator norm of $\ell(z)^{-1}$ with respect to $g_0$ is bounded by some constant independent of $z$, so the bound on $\| \ell(z)(Jv)\|$ implies that we must have $\lVert v' \rVert_{g_0} \leq K_5 e^{-\sqrt{\beta}}$ for some constant $K_5>0$. Finally, observe that we have $\langle d\widetilde{f}_{\beta,j}(z), J(v-v') \rangle=0$ by construction, so $J(v-v')\in T_z \widetilde{Z}_{\beta}$ and since $\lVert v' \rVert_{g_0}$ is bounded, we still get a lower bound analogous to Equation \ref{equation:no-zeroes} with an extra term of order $e^{-\sqrt{\beta}}$, so the desired contradiction follows from $\langle \lambda(z),J(v-v')\rangle=-\langle d\varphi(z),v-v'\rangle \neq 0$.
    \end{case}
    \begin{case}\label{case:easy-cpct}
        Suppose that $C_{(T^0_1,\dots,T^0_r)} \subsetneq C_{(T''_1,\dots,T''_r)}$ are both bounded cells. As in the previous cases, we can consider the \enquote{projection} of $u$ onto the cell $C_{(T^0_1,\dots,T^0_r)} \subset \partial C_{(T''_1,\dots,T''_r)}$ as $u^0$ and the minimum of $\varphi$ over the cell $C_{(T_1'',\dots,T''_r)}$ as $u_{\min}$, then  we know that $d(u,u^0) \leq K_1 \cdot \beta^{-\frac12}$ for a positive constant $K_1$. Moreover, since $u_{\min} \in \relint(C_{(T_1'',\dots,T''_r)})$, there exists a positive $\delta>0$ such that the inequality $d(u_{\min},\partial C_{(T_1'',\dots,T''_r)})>\delta$ holds for all such cells in $C_{0,\trop,\tot}$. Therefore, for $v\coloneqq \frac{u^0-u_{\min}}{\lVert u^0-u_{\min}\rVert}$, Lemma \ref{lemma:adapted_bound} gives us an estimate
        \begin{equation*}
            \langle d\varphi(u^0),v \rangle \geq m\cdot \lVert u^0-u_{\min} \rVert \geq m \cdot d(u_{\min},\partial C_{(T''_1,\dots,T''_r)}) \geq K_2, 
        \end{equation*}
        for some positive constant $K_2$. As in Case \ref{case:easy-noncpt}, this also implies
        \begin{equation*}
            \langle d\varphi(u),v \rangle \geq K_2-K_3\beta^{-\frac12}.
        \end{equation*}
        Now, we proceed similarly to Case \ref{case:hard-noncpct} to produce a tangent vector from $v$: since $u^0-u_{\min}$ is a vector between two points in $\overline{C}_{(T''_1,\dots,T''_r)}$, $v$ is tangent to that cell, hence $Jv$ is  \enquote{almost tangent} to $\widetilde{Z}_{\beta}$ with an error term of order $O(e^{-\sqrt{\beta}})$, so we can find an $O(e^{-\sqrt{\beta}})$-small vector $v'$ such that $J(v-v')$ is tangent to $\widetilde{Z}_\beta$ at $z$ and shows that $\lambda(z)|_{\widetilde{Z}_\beta} \neq 0$.
    \end{case}
    \begin{case}\label{case:hard-cpct}
        Finally, assume that $(T^0_1,\dots,T^0_r)=(T''_1,\dots, T''_r)$ and they label a bounded cell. The cut-off functions do not contribute to $\widetilde{f}_{\beta}$ near $z$, so we have $\widetilde{f}_{\alpha,\beta,j}(z)=f_{\alpha,\beta,j}(z)$ for $\alpha \in T_j$ and $\widetilde{f}_{\alpha,\beta,j}(z)=0$ otherwise. Since we are assuming that $d^c\varphi(z)$ vanishes on $T_z\widetilde{Z}_{\beta}$, there must exist real constants $a_j$, $b_j$ for $j=1,\dots,r$ such that 
        \begin{equation*}
        d^c\varphi(\rho)=\sum_{j=1}^r{a_j\RE(d\widetilde{f}_{\beta,j}(z))+b_j\IM(d\widetilde{f}_{\beta,j}(z))}.
        \end{equation*}
        By looking at the $d\rho_a$ components and using linear independence of the non-zero vectors in $T=\conv\{T_1,\dots,T_r\}$, we can see that phases of all the terms can only differ by a multiple of $\pi$, so they must be all real and at least one $f_{\alpha,\beta,j}(z)$ with $\alpha \in T_j \backslash \{0\}$ must be positive for each $j$. If all the terms are positive, then the point lies in the positive locus $\widetilde{Z}_{\beta}^+$ and looking at $d\theta$-components shows that $\Phi(z) \in \tau^\perp \times \tau$ for $\tau=\cone(T)$, which implies that $z \in \Lambda_\beta$, as desired. 

        Otherwise, suppose that there exists some $\alpha \in T_j$ such that $f_{\alpha,\beta,j}(z)<0$. Denote 
        \begin{equation*}
        \begin{split}
            T_j^+&\coloneqq\{ \alpha \in T_j \colon f_{\alpha,\beta,j}(z)>0 \},\\
            T_j^-&\coloneqq\{ \alpha \in T_j \colon f_{\alpha,\beta,j}(z)<0 \},
        \end{split}
        \end{equation*}
        then we know that $u$ lies in the complete intersection inside $M_\R$ given by $\{u \in M_{\R} \colon G_{j}^\pm(u)=|c_{0,j}| \textnormal{ for all } j=1,\dots,r\}$, where
        \begin{equation*}
            G_{j}^\pm(u)= \sum_{\alpha \in T_j^+}{c_{\alpha}e^{\beta l_{\alpha}(u)}}-\sum_{\alpha \in T_j^-}{c_{\alpha}e^{\beta l_{\alpha}(u)}}.
        \end{equation*}
        As we are over a holomorphic region, $d\varphi(z)$ must also vanish on $T_z\widetilde{Z}_{\beta}$, so there must exist constants $c_1 \dots, c_r$ such that 
        \begin{equation*}
            d\varphi(u)=\sum_{j=1}^r{c_jdG_{j}^\pm(u)}.
        \end{equation*}
        But by adaptedness of $\varphi$ and our proof of Lemma \ref{lemma:contains}, we know that $d\varphi(u) \in \relint(\sigma)$ for some $\tau \subset \sigma$. It is also straightforward to compute that $dG_{j}^\pm(u) \in \relint(\cone(T^+_j,-T^-_j))$, so $T^+_j$ being non-empty for all $j$ along with linear independence necessarily implies $c_j>0$ for all $j$, which contradicts the equation if $T^-_j$ is non-empty for some $j$, so this case can not occur and we are done. 
    \end{case}

\end{proof}

\begin{remark}\label{remark:phases}
    All our calculations also go through if we change the phases of $c_{\alpha,j} \in \nabla_j \backslash \{0\}$ from $0$ to $-\Theta(\alpha)$ for some function $\Theta\colon \Delta^{\vee}\rightarrow S^1$, which corresponds to rotating the FLTZ skeleton to 
    \begin{equation*}
        \eL_{\Sigma,\Theta}\coloneqq \bigcup_{\sigma \in \Sigma} \{(\theta,r) \in M_{S^1}\times N_{\R} \colon r \in \sigma \textnormal{ and } \langle \alpha,\theta \rangle = \Theta(\alpha) \textnormal{ for all } \alpha \in \sigma(1) \}.
    \end{equation*}
    For more details on the necessary changes, see \cite{Zhou}, where this more general case is treated. Deformation of Lagrangian skeleta of this form is also studied in more detail in \cite{zhou-thesis}.
\end{remark}
\section{Combinatorics of the skeleton} \label{section:combinatorics}
In this section, we shall outline how the combinatorics of $\skel(\widetilde{Z}_{\beta})$ on the A-side matches up with the combinatorics of the mirror. We also present a slightly different view of the case of hypersurfaces, revisiting the setting of \cite{GS22}. The results are motivated by generalising the following observation from our running example:

\begin{example}\label{example:running-ex-cover}
    The FLTZ skeleton $\eL_{\Pp^1} \subset T^*S^1$ mirror to $\Pp^1$ is the union of the zero section with the cotangent fibre over zero, and the FLTZ skeleton $\eL_{\bullet}$ mirror to a point is a point. Therefore, the skeleton depicted in Figure \ref{fig:skeleton} admits an open cover by eight copies of $\eL_{\Pp^1}$ overlapping over open intervals, which are the same as stabilisation $\eL_{\bullet} \times (0,1)$ (see Figure \ref{fig:mirror-covers} for an illustration of the cover). 
\end{example}

We operate within the same set-up as in Section \ref{section:skeleta}, with an extra requirement that the potentials that we consider are \emph{strongly} adapted. Throughout the entire section, we shall speak about diffeomorphisms between stratified spaces, which will always mean germs of diffeomorphisms of the ambient smooth manifolds that respect the stratifications. This is notion is a natural choice for our purposes, since all the spaces that we consider already come embedded inside $T^*M_{S^1}$ or $M_{\Cs}$, so we do not have to deal with the questions concerning a choice of the embedding.

\subsection{Covers of the FLTZ skeleton}\label{section:std-covers-skeleton} Before proving our main results, we shall review the material covered in \cite[Section 4.3]{GS22} and, in particular, explain what we mean by \emph{standard} inclusions of certain FLTZ skeleta in Theorem \ref{theorem:main-thm} and other places. 

As before, let $\Sigma \subset N_\R$ be a simplicial fan with an associated FLTZ skeleton $\eL_{\Sigma} \subset T^*M_{S^1}$. By picking ray generators for all rays in $\Sigma(1)$, we obtain diffeomorphisms of manifolds with corners $\sigma \cong \sigma/\tau \times \tau$ for all cones $\tau \subset \sigma$ that are also compatible with inclusions (i.e. for all $\tau \subset \sigma' \subset \sigma$, the identifications take the inclusions $\sigma' \hookrightarrow \sigma$ to $\sigma'/\tau \times \tau \hookrightarrow \sigma/\tau \times \tau$). In particular, if we restrict the identifications to $\sigma \cap \st(\tau) \cong \sigma/\tau \times \relint(\tau)$ and vary over the cones $\sigma$ containing $\tau$, these glue to a homeomorphism $\st(\tau) \cong (N/\tau)_{\R} \times \relint(\tau)$ that takes the stratification of the left hand side by $\Sigma$ to the stratification of the right hand side by $\Sigma/\tau$ (on the first factor) and whose restriction to $\relint(\sigma)$ is smooth for each $\tau \subset \sigma$. 

 However, the issue is that the piecewise linear identifications constructed this way are not smooth. However, observe that they can be fully reconstructed from the (piecewise linear) surjections $\pi_\sigma \colon \st(\sigma) \rightarrow \relint(\sigma)$ whose fibres are slices of the quotients $q_\sigma \colon \st(\sigma)\rightarrow (N/\sigma)_\R$ satisfying certain compatibilities. Therefore, we introduce the following class of maps:
 \begin{definition}\label{definition:coherent-projections}
    We say that a collection of surjective submersions $\pi_\sigma \colon \st(\sigma) \rightarrow \relint(\sigma)$ indexed by non-zero cones of $\Sigma$ is a \emph{collection of coherent projections} if
    \begin{enumerate}
        \item The fibres of $\pi_\sigma$ project diffeomorphically onto $(N/\sigma)_{\R}$ via the quotient map $q_\sigma\colon\st(\sigma)\rightarrow (N/\sigma)_\R$.
        \item For any pair of cones $\tau \subseteq \sigma$, we have $\pi_\tau \circ \pi_\sigma=\pi_\tau|_{\st(\sigma)}$. 
    \end{enumerate}
\end{definition}
Observe that, in particular, by taking $\tau=\sigma$ in the second condition, we get that the coherent projection must satisfy $\pi_\sigma|_{\relint(\sigma)}=\textnormal{id}_{\relint(\sigma)}$. The idea is that a choice of coherent projections $\{\pi_\sigma \}_{\sigma \in \Sigma \backslash\{0\}}$ gives us diffeomorphisms $\st(\sigma) \cong (N/\sigma)_{\R} \times \relint(\sigma)$ that take the stratification by $\Sigma$ to the stratification by $\Sigma/\sigma$ and such that an inclusion $\st(\sigma) \hookrightarrow \st(\tau)$  corresponds to an inclusion 
\begin{equation}\label{equation:comp-projections}
   (N/\sigma)_\R \times \relint(\sigma) \cong q_\tau(\st(\sigma)) \times \relint(\tau) \hookrightarrow (N/\tau)_\R \times \relint(\tau),
\end{equation}
for all pairs of cones $\tau \subset \sigma$ in $\Sigma$. 

In Appendix \ref{section:appendix-smoothing}, we explicitly construct smoothings of the piecewise linear maps considered as motivation earlier to show that there always exists a collection of coherent compatible projections (Corollary \ref{corollary:coherent-projections}) such that all the maps are also homogeneous of degree $1$. Given another choice $\{\pi'_\sigma \}_{\sigma \in \Sigma \backslash \{0\}}$, we obtain different identifications, but there will also exist canonical comparison maps $c^{\pi,\pi'}_{\sigma} \colon (N/\sigma)_{\R} \times \relint(\sigma) \rightarrow (N/\sigma)_{\R} \times \relint(\sigma)$ that lift the identity map of $(N/\sigma)_\R$, preserve the stratification by $\Sigma/\sigma$ and intertwine the two identifications given by $(q_\sigma,\pi_\sigma)$, $(q_\sigma,\pi'_\sigma)$. Therefore, we are free to pick a single collection of coherent projections, since any other one would yield essentially equivalent identifications. Note that this discussion also explains the relationship between these diffeomorphisms and the piecewise linear homeomorphisms that we started from: the maps obtained by gluing give piecewise smooth coherent projections, so we get piecewise smooth intertwining maps relating them with identifications from actual smooth coherent projections. 

We can pursue this strategy further to describe decompositions of FLTZ skeleta: consider the open subsets $\eL_{\Sigma}^\sigma\coloneqq\bigcup_{\sigma \subseteq \tau}{\tau^\perp\times\relint(\tau}) \subset \eL_\Sigma$, so that $\eL_\Sigma^0=\eL_\Sigma$. We then have $\eL_\Sigma^\sigma=\eL_\Sigma\cap(M_{S^1} \times \st(\sigma))=\eL_\Sigma \cap(\sigma^\perp\times\st(\sigma))$ for all $\sigma \in \Sigma$, where the first equality is by definition (and makes it clear that these are indeed open subsets), while the second one follows from the fact that $\tau^\perp \hookrightarrow \sigma^\perp$ for all $\sigma \hookrightarrow \tau$. A collection of coherent projections then gives us diffeomorphisms $\sigma^\perp \times \st(\sigma) \cong \sigma^\perp \times (N/\sigma)_\R \times \relint(\sigma)=T^*\sigma^\perp \times \relint(\sigma)$. Since these are compatible with stratifications by $\Sigma/\sigma$, they induce identifications $\eL_\Sigma^\sigma \cong \eL_{\Sigma/\sigma}\times\relint(\sigma)$ for all $\sigma \in \Sigma$ (as mentioned in the introduction to this Section, \enquote{identification} means a germ of an ambient diffeomorphism between a stratified subspace of $T^*M_{S^1}|_{\sigma^\perp}$ and a stratified subspace of $T^*\sigma^\perp \times \relint(\sigma)$ that respects the stratifications). The compatibility relations for coherent projections then tell us that an inclusion $\eL_\Sigma^\sigma \hookrightarrow \eL_\Sigma^\tau$ is identified with
\begin{equation}\label{equation:comp-projections-lag}
    \eL_{\Sigma/\sigma} \times \relint(\sigma) \cong \eL^{\sigma/\tau}_{\Sigma/\tau} \times \relint(\tau) \hookrightarrow \eL_{\Sigma/\tau} \times \relint(\tau),
\end{equation}
for all cones $\tau \subset \sigma$ in $\Sigma$, so we have provided a smooth analogue of \cite[Lemma 4.3.1]{GS22}. The maps $\eL_{\Sigma/\sigma} \times \relint(\sigma) \cong \eL^\sigma_\Sigma \hookrightarrow \eL_\Sigma$ constructed this way are precisely the \emph{standard inclusions} that we alluded to earlier. 

Finally, we show how this construction can be used to describe a cover of $\partial \eL_\Sigma$ (strengthening Lemma 4.3.2 in loc. cit.). The key observation is the following:

\begin{lemma}\label{lemma:lie-action-product}
    Let $G$ be a Lie group acting smoothly and properly on a manifold $X$ and also on its product with another manifold $Y$ in a way that lifts the action on $X$ (i.e. the projection $X\times Y \rightarrow X$ is $G$-equivariant). Moreover, suppose that the action on $X$ is free and admits a section $s\colon X/G \rightarrow X$. Then there exists a diffeomorphism $(X\times Y)/G \cong (X/G)\times Y$. 
\end{lemma}

In our context, the space $N_\R$ admits a scaling action $\R \curvearrowright N_\R$ given as $t\cdot v=e^tv$ that is free away from zero. We shall denote the \emph{sphere at infinity} as $\partial_\infty N_\R\coloneqq(N_\R\backslash\{0\})/\R \cong S^{n-1}$ and for a general conical subset $U \subset N_\R \backslash \{0\}$, let $\partial_\infty U \coloneqq C/\R \subset \partial_\infty N_\R$. In particular, non-zero cones of the simplicial fan $\Sigma$ define a triangulation $\partial_\infty\Sigma$ of the sphere at infinity by simplices $\partial_\infty\sigma$. Therefore, we have relatively open subsets $\relint(\partial_\infty\sigma)$ and their open neighbourhoods $\st(\partial_\infty\sigma)$ inside $\partial_\infty N_\R$.

Suppose that the coherent projections $\pi_\sigma$ are equivariant with respect to this action (this is the synonymous with being homogeneous of degree $1$, so such projections exist by Corollary \ref{corollary:coherent-projections}), then we can apply Lemma \ref{lemma:lie-action-product} with $G=\R$, $X=\relint(\sigma)$, $Y=(N/\sigma)_\R$ and the two different scaling actions in $\st(\sigma) \cong (N/\sigma)_\R \times \relint(\sigma)$. This yields diffeomorphisms $\st(\partial_\infty \sigma) \cong (N/\sigma)_\R \times \relint(\partial_\infty \sigma)$, which will satisfy a compatibility relation analogous to Equation \ref{equation:comp-projections}.

By applying Lemma \ref{lemma:lie-action-product} to the fibrewise scaling action on the cotangent bundle $T^*M_{S^1}$ in the discussion of the FLTZ skeleton, we obtain an open cover of $\partial \eL_\Sigma$ by the open sets $\partial \eL_\Sigma^\sigma \cong \eL_{\Sigma/\sigma} \times \relint(\partial_\infty \sigma)$, with compatibility relations precisely as in Equation \ref{equation:comp-projections-lag}, but with $\relint(\sigma)$ replaced by $\relint(\partial_\infty \sigma)$. 

\subsection{Hypersurfaces}\label{section:combi-hyp} In this section, we explain how one can construct an open cover of the skeleton of a very affine hypersurface from a slightly different perspective than \cite{GS22}, which we shall later generalise to the case of open Batyrev--Borisov complete intersections.

Suppose that we start from a polytope $\Delta^\vee$ with a centred refined triangulating function $h$ and a corresponding simplicial fan $\Sigma$. For brevity, denote $\Lambda\coloneqq\skel(\widetilde{H}_\beta)$. We shall also pick an inner product on $N_\R$, which gives us a unit sphere $S(N_\R) \subset N_\R$ that projects diffeomorphically onto $\partial_\infty N_\R$. 

\begin{lemma}\label{lemma:isotope-skeleton-hyp}
    There exists a smooth ambient isotopy that preserves $\Sigma$ and takes $\Phi(\partial \widetilde{C}_{0,\beta})$ to $S(N_\R)$. 
\end{lemma}
\begin{proof}
    Thanks to the requirement that the cut-off functions are centred, the region $\widetilde{C}_{0,\beta} \subset M_\R$ is convex and contains the origin, hence it is star-shaped, so its image under the homogeneous map $\Phi=d\varphi \colon M_\R \rightarrow N_\R$ is also star-convex. Therefore, the statement follows through an appropriate rescaling.
\end{proof}

If we apply an appropriate rescaling of this isotopy to the second factor in the decomposition $T^*M_{S^1}=M_{S^1} \times N_\R$, we can deform $\Phi(\Lambda)=\eL_\Sigma \cap \Phi(\Log_\beta^{-1}(\partial \widetilde{C}_{0,\beta}))$ to a standard boundary $\partial \eL_\Sigma$, so the discussion from Section \ref{section:std-covers-skeleton} immediately recovers the following results (in particular, note that we do not need strong adaptedness of $\varphi$, since radial scaling is sufficient to construct all the necessary isotopies):

\begin{proposition}[{\cite[Corollary 4.3.2]{GS22}}]\label{proposition:combi-hyp}
    There exists a cover of $\Lambda$ by open sets $\Lambda(\sigma)$ that are anti-indexed by the poset of non-zero cones $\sigma \in \Sigma$, each $\Lambda(\sigma)$ is diffeomorphic to $\eL_{\Sigma/\sigma}\times\relint(\partial_\infty \sigma)$ and the inclusions $\Lambda(\tau) \hookrightarrow \Lambda(\sigma)$ correspond to the standard inclusions $\eL_{\Sigma/\tau} \times \relint(\partial_\infty \tau) \hookrightarrow \eL_{\Sigma/\sigma} \times \relint(\partial_\infty \sigma)$ under these identifications for all $\sigma \subset \tau$.
\end{proposition}

\begin{proof}
    As explained above, if we take $\Lambda(\sigma)\coloneqq \Lambda \cap \Phi^{-1}(M_{S^1}\times\st(\sigma))$, the statement follows from Lemma \ref{lemma:isotope-skeleton-hyp} and Section \ref{section:std-covers-skeleton}.
\end{proof}

\subsection{Complete intersections}\label{section:combi-ci} Suppose that we are given the data of an irreducible nef-partition $\nabla=\nabla_1+\dots+\nabla_r$ of a reflexive lattice polytope $\nabla \subset M_\R$ with a refined centred triangulating function $h$ inducing regular triangulations $\mT$ of $\partial \Delta^\vee$, $\mT_j$ of $\partial\nabla_j\backslash\{0\}$ which also gives us a simplicial fan $\Sigma$. Recall that we have defined $\Sigma^\ba$, the \emph{barycentric subdivision} of a fan $\Sigma$, by gluing cubes $\sigma^\ba \subset \sigma$ (Definition \ref{definition:bar_subdiv}) and \emph{transversal cones} $\sigma \in \Sigma$ for which all the sets $\sigma_j(1)\coloneqq \sigma(1) \cap \Sigma_j(1)$ are non-empty (Section \ref{section:tropical-bbci}). Putting these two notions together, we define the \emph{transversal barycentric subdivision} as
\begin{equation*}
    \Sigma^\ba_\trans\coloneqq\bigcup_{\sigma \in \Sigma_\trans}{\Sigma^\ba_\sigma}.
\end{equation*}
Along with the combinatorial data, we fix some $(h_j)$-compatible cut-off functions $\chi$ satisfying the same properties as in Section \ref{section:skeleta} and a potential $\varphi$ that is strongly adapted to the polytope $C_{0,\trop,\tot}$. As before, we write $\Lambda\coloneqq\skel(\widetilde{Z}_\beta)$ for simplicity. 

The main difficulty of this case is that we no longer have an easily available identification between $\partial \widetilde{C}_{0,\beta}$ with a standard space, such as the sphere $S(N_\R)$ from Lemma \ref{lemma:isotope-skeleton-hyp}. Instead, we need to rely on a more involved result, Corollary \ref{corollary:smoothing-isotopies2}, to replace this with a lower-dimensional sphere obtained by smoothing $\Sigma^\ba_\trans$ relative to the fan $\Sigma$. 

We also need to generalise the ordinary scaling action that played a key role in the hypersurfaces case: recall that every cone $\sigma$ (not necessarily transversal) admits a unique decomposition as $\sigma=\sigma_1+\dots+\sigma_r$, where each $\sigma_j$ has ray generators in $\Sigma_j(1)$, since the ray generators form a simplex $T \in \mT$ that is partitioned into $T_1$,\dots,$T_r$ and we can take $\sigma_j\coloneqq \cone(T_j*0)$. Then $\sigma$ is transversal if and only if all $\sigma_j$ are non-zero. Since all the cones involved are simplicial, we also get a decomposition of the underlying vector spaces $\R\sigma=\bigoplus_{j=1}^r \R\sigma_j$ and hence a diffeomorphism $\sigma \cong \sigma_1 \times \dots \times \sigma_r$.

\begin{definition}
    For a cone $\sigma=\sigma_1+\dots+\sigma_r$, we define the \emph{nef partition scaling action} as the smooth $\R^r$-action given by $t \cdot v\coloneqq e^{t_1}v_1+\dots+e^{t_r}v_r$ for $v=v_1+\dots+v_r$ and $t=(t_1,\dots,t_r)$. 
\end{definition}

Observe this action is free on $\relint(\sigma)$ if and only if the cone $\sigma$ is transversal. In that case, the orbit space of $\relint(\sigma) \cong \relint(\sigma_1)\times\dots\relint(\sigma_r)$ can clearly be identified with the product of simplices $\relint(\partial_\infty \sigma_1)\times\dots\relint(\partial_\infty \sigma_r)$. We denote the corresponding quotient map as $q_{nef}\colon \relint(\sigma) \rightarrow \relint(\partial_\infty \sigma_1)\times\dots\relint(\partial_\infty \sigma_r)$. In fact, these actions extend to smooth actions on manifolds with corners $\overline{\relint(\sigma)}$, the closure of $\relint(\sigma)$ inside $|\Sigma_\trans|$ (see the end of Section \ref{section:tropical-bbci}), where the orbit space becomes diffeomorphic to the compact space $\partial_\infty \sigma_1 \times \dots \partial_\infty \sigma_r$. Since the extensions agree on overlaps, they glue to a continuous, piecewise smooth action $A_{nef}\colon\R^r \curvearrowright \lvert \Sigma_\trans\rvert$. We could directly glue the orbit spaces for individual cones following the same formulae to obtain an orbit space for this global action, but we instead make the following observation:

\begin{lemma}\label{lemma:orbit-space-barycentric}
    $\Sigma^\ba_\trans$ is a slice of the continuous action $A_{nef}$.
\end{lemma}

\begin{proof}
    Note that it suffices to prove that the quotient map $q_{nef}$ by the action restricted to $\overline{\sigma}^\ba_\trans$ (with the closure taken inside $\Sigma^\ba_\trans$) is a homeomorphism onto $\partial_\infty \sigma_1 \times \dots \partial_\infty \sigma_r$, since the required statement then follows by gluing. But this follows by definition of $\Sigma_\trans$: let $\sigma=\sigma_1+\dots+\sigma_r$ be the decomposition induced by the nef partition and for a ray $\rho \in \sigma(1)$, denote the linear functional defined on $\R \sigma$ and dual to the ray generator $v_\rho$ as $\eta^\rho_\sigma$. Then one can check that for $\eta_j\coloneqq\max\{\eta^\rho_\sigma \colon \rho \in \sigma_j(1)\}$, we have 
    \begin{equation*}
        \overline{\sigma}^\ba_\trans=\{v \in \sigma \colon \eta_j(v)=1 \textnormal{ for all } j=1,\dots,r \}.
    \end{equation*}
    Since $\eta_j(t\cdot v)=t_j\eta_j(v)$ for all $j$, the orbit of $A_{nef}$ through $v$ with $v \in \overline{\relint(\sigma)}$ (closure taken inside $|\Sigma^\trans|$) must intersect $\overline{\sigma}^\ba_\trans$ at exactly one point (that is explicitly given as $\frac{v_1}{\eta_1(v)}+\dots+\frac{v_r}{\eta_r(v)}$ for $v=v_1+\dots+v_r$). Therefore, restricting the quotient map to $\overline{\sigma}^\ba_\trans$ yields a continuous bijection with the orbit space $\partial_\infty \sigma_1 \times \dots \times \partial_\infty \sigma_r$, which is a homeomorphism by the topological inverse function theorem.
\end{proof}

We can also smooth the piecewise linear functions $\eta_1$, \dots, $\eta_r$ defining $\Sigma^\ba_\trans$ by convolution to obtain a family of smoothings $\widetilde{\Sigma}^\ba_\trans$ with a single parameter $\varepsilon>0$ (see Section \ref{section:ooga-booga-smoothing} and, in particular, the discussion around Corollary \ref{corollary:trop-smoothings} for more details on how the smoothing is defined). For such an object, we can do better and strengthen the previous Lemma to a statement about diffeomorphisms over each transversal cone:

\begin{corollary}\label{corollary:smoothed-ci-cone}
    For any $\varepsilon$-smoothing $\widetilde{\Sigma}^\ba_\trans$ with a sufficiently small parameter $\varepsilon>0$, the manifold with corners $\widetilde{\Sigma}^\ba_\trans \cap \overline{\sigma}$ will be a smooth section of $q_{nef}\colon \sigma \rightarrow \partial_\infty \sigma_1 \times \dots \times \partial_\infty \sigma_r$ for all transversal cones. 
\end{corollary}

\begin{proof}
    First, observe that by definition of $\Sigma$-smoothings, the space $\widetilde{\sigma}^\ba_\trans\coloneqq\widetilde{\Sigma}^\ba_\trans \cap \overline{\sigma}$ is indeed a manifold with corners (since $\widetilde{\Sigma}^\ba_\trans$ is a $\Sigma$-smoothing by Corollary \ref{corollary:trop-smoothings}). We shall prove that the orbits of the $\R^r$-action intersect $\widetilde{\sigma}^\ba_\trans$ transversely. From that, it follows that $q_{nef}|_{\widetilde{\sigma}^\ba_\trans}$ is a local diffeomorphism onto its image inside $\partial_\infty \sigma_1 \times \dots \times \partial_\infty \sigma_r$. Since the domain is compact, the image has to be closed, but it is also open, so by connectedness of the codomain, the map is surjective. Since $\widetilde{\tau}^\ba_\trans$ is a single point for a minimal transversal cone $\tau\subset\sigma$ (by properties of $\Sigma$-smoothings, it is ambient isotopic to $\tau^\ba_\trans$, which is a point), sizes of preimages being locally constant forces the map to be injective, hence a diffeomorphism $q_{nef}\colon \widetilde{\sigma}^\ba_\trans \rightarrow \partial_\infty \sigma_1 \times \dots \times \partial_\infty \sigma_r$. 
    
    By definition, $v \in \widetilde{\Sigma}^\ba_\trans$ if and only if $\widetilde{h}_{j,\varepsilon}(v)=1$ for all $j=1,\dots,r$, so it suffices to check that the $(r\times r)$ matrix $M$ with entries $M_{ij}=D_v\widetilde{h}_i(v_j)$ is invertible. By definition of the smoothing $\max_\varepsilon$, $D_v\widetilde{h}_i(v_j)$ is a convex combination of $D_v\widetilde{h}^\rho_\varepsilon(v_j)$ for $\rho \in \Sigma_i(1)$. By Corollary \ref{corollary:convolution-est-1}, we know that each of these terms satisfies $|D_v\widetilde{h}^\rho_\varepsilon(v_j)-\delta_{ij}|<C\varepsilon$ for some constant $C>0$, which tells us that $\lVert M-I_r\rVert<C\varepsilon$, so for $\varepsilon$ small enough, the matrix $M$ is indeed invertible. 
\end{proof}

Unlike in the hypersurface setting, we have one more issue to deal with -- the nef scaling action is smooth over the relative interior of every cone, but not globally smooth, so we can not directly repeat the calculation from Corollary \ref{corollary:smoothed-ci-cone} over $\st(\sigma)$ and apply Lemma \ref{lemma:lie-action-product}. The idea behind smoothing the action is as follows: observe that the infinitesimal nef scaling action at $x \in \relint(\sigma)$ is generated by the vectors $x_1$, \dots, $x_r$. Thanks to the choice of generators $v_\rho$ for all $\rho \in \sigma(1)$, we can express these in terms of the piecewise smooth functions $h^\rho$ (whose restriction to $\sigma$ is simply the dual functional $\eta^\rho_\sigma$) as $x_j=\sum_{\rho \in \sigma_j(1)} h^\rho(x)\cdot v_\rho$. Therefore, in order to smooth the action, it suffices to smooth the piecewise linear functions $h^\rho$ by convolution. This is done in detail in Section \ref{section:ooga-booga-smoothing} (see Definition \ref{definition:total-scaling-action-sm} and the discussion that follows after it). The construction also involves the freedom of choice of two extra smoothing parameters $(\varepsilon',\delta')$, in addition to the parameter $\varepsilon$ that we are using for smoothing $\Sigma^\ba_\trans$ to $\widetilde{\Sigma}^\ba_\trans$. It also forces us to work with smaller sets $\sigma^{fr} \subset \relint(\sigma)$ and $\st^{fr}(\sigma) \subset \st(\sigma)$ where the smoothed action is free, but it is clear that $\widetilde{\Sigma}^\ba_\trans$ is still covered by $\st^{fr}(\sigma)$ for transversal cones $\sigma$ for all sufficiently small $\varepsilon'$ and $\delta$. We denote such a smoothing $\widetilde{A}_{nef}$ and the associated quotient map $\widetilde{q}_{nef}$.

\begin{lemma}\label{lemma:transverse-action}
    For any sufficiently small $\varepsilon>0$, the $\varepsilon$-smoothing $\widetilde{\Sigma}^\ba_\trans$ will be transverse to the fibres of $\widetilde{q}_{nef}$ for all $\varepsilon'$ and $\delta'$ small enough.
\end{lemma}
\begin{proof}
    By Lemma \ref{lemma:vector-field-smoothings}, we know that the vector fields $\widetilde{V}_\rho$ generating the infinitesimal action associated to the total scaling action are $\varepsilon'$-close to the generators $V_\rho$. Therefore, the vector fields $\widetilde{V}_1$, \dots, $\widetilde{V}_r$ generating the nef scaling action are also $\varepsilon'$-close to the generators $V_1$, \dots, $V_r$ to the non-smooth action. By the calculation from Corollary \ref{corollary:smoothed-ci-cone}, we know that the matrix $M$ with entries $M_{ij}(x)=D_xh_{i,\delta}(V_j(x))$ is $\varepsilon$-close to the identity matrix for all $x\in \widetilde{\Sigma}^\ba_\trans$, which means that $\widetilde{M}_{ij}=D_xh_{i,\delta}(\widetilde{V}_j(x))$ will be $(\varepsilon+\varepsilon')$-close to $I_r$, so the matrix will be invertible for all $\varepsilon, \varepsilon'$ sufficiently small, which is equivalent to the desired transversality statement.  
\end{proof}

\begin{corollary}\label{corollary:smoothed-ci-cone-2}
    For every transversal cone, the orbit space $\widetilde{q}_{nef}(\sigma^{fr})$ is diffeomorphic to $\relint(\partial_\infty \sigma_1)\times \dots \times \relint(\partial_\infty \sigma_r)$.
\end{corollary}
\begin{proof}
    By Lemma \ref{lemma:transverse-action} and the discussion from the proof of Corollary \ref{corollary:smoothed-ci-cone}, we get that $\widetilde{\Sigma}^\ba_\trans \cap \sigma^{fr}$ is the orbit space of the action on $\sigma^{fr}$ (we can consider the closure of $\sigma^{fr}$ in the union of $\st^{fr}(\tau)$ over all $\tau \in \Sigma_\trans$ to reduce it to a statement about compact manifolds with corners as there and then pass to interiors). By definition of $\widetilde{\Sigma}^\ba_\trans$, it is clear that for appropriate values of smoothing parameters, one has $\widetilde{\Sigma}^\ba_\trans \cap \sigma^{fr}=\widetilde{\Sigma}^\ba_\trans \cap \relint(\sigma)$, so we are done by the aforementioned Corollary. 
\end{proof}

One can also define \emph{equivariant coherent projections} $\pi_\sigma \colon \st^{fr}(\sigma)\rightarrow\sigma^{fr}$ associated to the smoothing $\widetilde{A}_{nef}$ (see Definition \ref{definition:coherent-projections-equiv} and Corollary \ref{corollary:coherent-projections-equiv}), so we get identifications $\st^{fr}(\sigma) \cong (N/\sigma)_\R \times \sigma^{fr}$. Moreover, these identifications are also compatible with $\st^{fr}(\sigma)\hookrightarrow {\st^{fr}(\tau)}$ for all transversal $\tau \subset \sigma$ in the sense of Equation \ref{equation:comp-projections}. With that in mind, we are ready to prove the main result of this section: 

\begin{theorem}[Theorem \ref{theorem:main-thm}(\ref{main-thm-4})]\label{theorem:combi-ci}
    There exists a cover of the skeleton $\Lambda$ by open sets $\Lambda(\sigma)$ that are anti-indexed by the poset of transversal cones $\sigma \in \Sigma_\trans$, each $\Lambda(\sigma)$ diffeomorphic to $\eL_{\Sigma/\sigma} \times \relint(\partial_\infty \sigma_1) \times \dots \times \relint(\partial_\infty \sigma_r)$ and the inclusions $\Lambda(\tau) \hookrightarrow \Lambda(\sigma)$ correspond to the standard inclusions under these identifications for all $\sigma \subset \tau$. 
\end{theorem}

\begin{proof}
    Suppose that the smoothing parameters $\varepsilon, \varepsilon', \delta$ are picked so that all the previous results of this section apply. By Corollary \ref{corollary:smoothing-isotopies2}, we know that $\Phi(\partial \widetilde{C}_{0,\beta})$ is isotopic (relative to $\Sigma$) to a smoothing $\widetilde{\Sigma}^\ba_\trans$ constructed above via some isotopy $G_t$. By Theorem \ref{theorem:skeleton_ci}, this means that the isotopy $\overline{G}_t=\textnormal{id} \times G_t$ takes $\Phi(\Lambda)$ to $\Lambda_\trans\coloneqq(M_{S^1} \times \widetilde{\Sigma}^\ba_\trans) \cap \eL_\Sigma$. We shall define $\Lambda(\sigma)\coloneqq\Lambda \cap \Phi^{-1}(M_{S^1} \times \st^{fr}(\sigma))$ as an open cover, so that $\Phi(\Lambda(\sigma))$ gets taken to $\Lambda_\trans(\sigma)\coloneqq\Lambda_\trans \cap (M_{S^1} \times \st^{fr}(\sigma))$ by $\overline{G}_t$. 

    Note that by Lemma \ref{lemma:lie-action-product} and the definition of equivariant coherent projections, the orbit spaces $\widetilde{q}_{nef}(\st^{fr}(\sigma))$ are diffeomorphic to $\widetilde{q}_{nef}(\sigma^{fr}) \times (N/\sigma)_\R$, so by Lemma \ref{corollary:smoothed-ci-cone-2}, we get identifications $\widetilde{q}_{nef}(\st^{fr}(\sigma)) \cong (N/\sigma)_\R \times \relint(\partial_\infty \sigma_1) \times \dots \times \relint(\partial_\infty \sigma_r)$ that are compatible with stratifications by $\Sigma$ and $\Sigma/\sigma$. Therefore, on the cotangent bundles, they naturally induce diffeomorphisms $\Lambda_\trans(\sigma) \cong \eL_{\Sigma/\sigma} \times \relint(\partial_\infty \sigma_1) \times \dots \times \relint(\partial_\infty \sigma_r)$, and the compatibility with inclusions follows from properties of coherent projections. 
\end{proof}

We end this section by commenting on how to shrink the open sets considered in Theorem \ref{theorem:combi-ci} to a different cover that is somewhat more practical. One can replace the open sets $\st^{fr}(\sigma) \subset N_\R$ in the definition of $\Lambda(\sigma)$ by smaller neighbourhoods $\st^\varepsilon(\sigma)$ of $\sigma^{fr}$, which corresponds to restricting our attention to star-shaped open neighbourhoods $U^\varepsilon(\sigma)$ of $0\in (N/\sigma)_{\R}$ or radius at most $\varepsilon>0$, giving us open subsets $\Lambda^\varepsilon(\sigma)\coloneqq \Lambda \cap \Phi^{-1}(\sigma^\perp\times \st^\varepsilon(\sigma))$ instead. Note that the notation is chosen so that setting $\varepsilon=0$ can be seen as recovering the sets $\Lambda^0(\sigma)=\Lambda \cap \Phi^{-1}(\sigma^\perp \times \sigma)$ considered in Section \ref{section:skeleta} as the stratification of $\Lambda$ into compact smooth Lagrangians with corners. 

We still get inclusions $\st^\varepsilon(\sigma) \hookrightarrow \st^\varepsilon(\tau)$ for all pairs of cones $\tau \subset \sigma$ and since quotienting by a subspace maps a small star-shaped neighbourhood of the origin to a small star-shaped neighbourhood of the origin, we can arrange them so that $q_{\sigma/\tau}(U^\varepsilon(\tau))= U^\varepsilon(\sigma)$ and $U^\varepsilon(\tau)$ is a trivial fibre bundle over $U^\varepsilon(\sigma)$ through $q_{\sigma/\tau}$. The identifications $\sigma^\varepsilon \cong U^\varepsilon(\sigma) \times \sigma_{fr}$ from equivariant coherent projections then satisfy compatibilities
\begin{equation*}
    U^\varepsilon(\sigma) \times \sigma_{fr} \cong q_\tau(\sigma^\varepsilon) \times \tau_{fr} \hookrightarrow U^\varepsilon(\tau) \times \tau_{fr},
\end{equation*}
which can be seen as restrictions of Equation \ref{equation:comp-projections}. On the level of FLTZ skeleta, we know that $(\sigma^\perp \times U^\varepsilon(\sigma))\cap \eL_{\Sigma/\sigma} \cong \eL_{\Sigma/\sigma}$ for all $\sigma \in \Sigma_\trans$, since the inclusion $U^\varepsilon(\sigma) \hookrightarrow (N/\sigma)_\R$ is isotopic to identity via a scaling isotopy. Therefore, the conclusion of Theorem \ref{theorem:combi-ci} still holds for the open cover $\{\Lambda^\varepsilon(\sigma)\}_{\sigma \in \Sigma_\trans}$ for all sufficiently small $\varepsilon>0$. 
\section{Stabilising complete intersections} \label{section:stabilisation}
In this section, we provide a further refinement of our description of the skeleton from previous two sections: recall that in order to define the cosheaf $\msh_\Lambda$, it is necessary to embed a suitable stabilisation of the associated Liouville domain $A_\beta$ into the co-sphere bundle of some manifold. In general, such embeddings can be obtained from $h$-principles (see, e.g., \cite[Section 7]{gps3}). However, in the case at hand, we can exploit the fact that $A_\beta \subset M_{\Cs}$ is a subdomain of a very affine variety to construct an explicit embedding into $S^*M_{S^1}$. This will have some important consequences in the last part of Section \ref{section:hms}, where we show that our HMS equivalences are functorial under inclusions of open BBCIs and compatible with toric mirror symmetry. The main result of this section is the following statement, which is an analogue of \cite[Corollary 6.2.6]{GS22} (to avoid clutter, we are going to omit the ubiquitous index $\beta$ from our notation and just assume that $\beta \gg0$): 

\begin{theorem}\label{theorem:emb}
    For all sufficiently large $\beta>0$, there exists a Liouville subdomain $U \subset (T^*M_{S^1},\lambda_{0})$ completing to $T^*M_{S^1}$, a Liouville subdomain $A \subset \widetilde{Z}$ and a Liouville hypersurface embedding 
    \begin{equation*}
        i\colon A \times \C^{r-1} \hookrightarrow \partial U.
    \end{equation*}
    Moreover, $i$ can be chosen so that the image of $\skel(\widetilde{Z}\times(\C,\pm\infty)^{r-1})=\Lambda \times \R^{r-1}$ is a neighbourhood of $\Phi(\Lambda)$ inside $\partial\eL_\Sigma \subset \partial U$.
\end{theorem}

Before starting the proof, we remark that a similar set-up has previously been considered by Avdek in \cite[Section 5]{avdek}. Since we need more control over how skeleta interact with the entire construction, we can not use the results from there directly, so we explicitly adapt all the relevant arguments to the present setting for the sake of clarity. 

Recall that the Liouville domain associated to $\widetilde{Z}$ was previously defined as $A'=\widetilde{Z} \cap \{\varphi \leq C \}$, where $C>0$ is some sufficiently large constant. However, we shall shrink it to a small neighbourhood of $\Lambda$ by applying the Liouville flow of $X^{\widetilde{Z}}_{\lambda}$ to it and instead consider 
\begin{equation*}
    A\coloneqq \Phi^{-T}_{X}(A'),
\end{equation*}
for some sufficiently big $T>0$. Clearly $A' \hookrightarrow A$ is a trivial inclusion of Liouville domains for any $T$, so it does not matter which particular value we choose. The reason why we modify $A'$ is to have the following result at hand: 
\begin{lemma}\label{lemma:vf_not_tangent}
    For $T>0$ large enough, the Liouville vector field $X_{\lambda}$ of $M_{\Cs}$ is nowhere tangent to $\widetilde{H}_{\tot}$ along $A$.
\end{lemma}
\begin{proof}
    Analogously to Section \ref{section:skeleta}, consider the splitting $X_{\lambda}=X^{\perp}_{\lambda}+X^{\widetilde{H}}_{\lambda}$ obtained via symplectic orthogonal projection onto $T_z\widetilde{H}_{\tot}$. It suffices to check that $X^{\perp}_{\lambda} \neq 0$ along $A$.
    
    By a calculation similar to the one in Lemma \ref{lemma:liou-dynamics}, it follows that $X^{\perp}_{\lambda}(z) \neq 0$ along $\Lambda$, because we have:
    \begin{equation*}
    \begin{split}
        \langle d\widetilde{f}_{\tot}(z), X^{\perp}_{\lambda} \rangle &= \langle d\widetilde{f}_{\tot}(z), X_{\lambda} \rangle = \langle d\widetilde{f}_{\tot}(z), \nabla_g \varphi \rangle \\
        &= \langle d\widetilde{F}_{\tot}(z), \nabla_g \varphi \rangle=c_3\langle d\widetilde{F}_{\tot}(z), \sum_{a=1}^n{\rho_a\partial_{\rho_a}} \rangle > 0.
    \end{split}
    \end{equation*}
    Since being non-zero is an open condition, the same holds over a sufficiently small open neighbourhood of $\Lambda$ (by compactness of the skeleton, it will be true for an $\varepsilon$-neighbourhood for some $\varepsilon>0$). Because $A'$ is compact, we can guarantee that for any $T$ large enough, its image under time-$(-T)$ Liouville flow $A$ will lie inside the $\varepsilon$-neighbourhood, so the desired statement follows. 
\end{proof}

Therefore, we an use the same idea as in \cite[Corollary 6.2.6]{GS22} to construct a Liouville domain $U_0$ that contains a neighbourhood of $A$ inside $\widetilde{H}_{\tot}$ in its boundary and completes to $M_{\Cs}$: the Liouville flow gives us an identification $M_{\Cs}\backslash M_{S^1} \cong B_{\tot} \times \R$ with a projection $\pi \colon M_{\Cs}\backslash M_{S^1} \rightarrow B_{\beta,\tot}$. By Lemma \ref{lemma:vf_not_tangent}, $A$ is transverse to the fibres of $\pi$, so if we pick a sufficiently small neighbourhood $A_{\tot}$ of $A$, $\pi(A_{\tot})$ will be immersed inside $B_{\tot}$. Since it can also made to be arbitrarily close to the skeleton $\Lambda$ that is embedded, we can assume that $\pi(A_{\tot})$ is an embedded submanifold with boundary of codimension $1$. Therefore, $A_{\tot}$ can be written as a graph of some function defined over $\pi(A_{\tot}) \subset B_{\tot}$. We can arbitrarily extend this function to the rest of $B_{\tot}$ and declare the Liouville domain $U_0$ to be its graph. 

To get the desired result, it suffices to stabilise $\widetilde{Z}$ inside $\widetilde{H}_{\tot}$ appropriately. That requires extending the standard version of the symplectic neighbourhood form (see, e.g., \cite[Theorem 3.4.10]{mcbook}) from the closed setting to the setting of Liouville domains. Following the terminology from \cite{gps1}, recall that a \emph{Liouville inclusion} is a proper\footnote{More precisely, we mean that its natural extension to a map of completions via Liouville flows is a proper map of Liouville manifolds.}, codimension zero embedding $i\colon X \hookrightarrow Y$ such that there exists a function $f$ on $X$ supported away from the boundary such that $i^*\lambda_Y=\lambda_X+df$, and that a Liouville inclusion is \emph{trivial} if $i(X)$ can be smoothly deformed to $Y$ through Liouville inclusions. Since we will be interested in Liouville domains rather than Liouville manifolds, we introduce the following more restrictive notions that impose some extra conditions on what happens at the boundary.

\begin{definition}\label{definition:liouv-embeddings}
    We call a smooth map between Liouville domains $i\colon(X,\lambda_X) \rightarrow (Y,\lambda_Y)$ with $\dim(X)<\dim(Y)$ a \emph{positive codimension Liouville embedding} if it is a diffeomorphism onto its image, the boundaries satisfy $i(\partial_\infty X)= i(X) \pitchfork \partial_\infty Y$ and there exists a function $f$ supported away from the boundary such that $i^*\lambda_Y=\lambda_X+df$. We call a Liouville inclusion that is also a diffeomorphism a \emph{Liouville diffeomorphism}.
\end{definition}

Clearly, Liouville diffeomorphisms are the same as trivial Liouville inclusions $\hat{X} \hookrightarrow \hat{Y}$ of Liouville manifolds that map a boundary at infinity $\partial_\infty X$ to a specified boundary at infinity $\partial_\infty Y$. Also, observe that for any positive codimension Liouville embedding, we have $i^*\lambda_Y=\lambda_X$ near $\partial_\infty X$, which means that the Liouville vector field of $Z_X$ of $i(X)$ is positively transverse to $i(\partial_\infty X)$, so $(i(X),\lambda_Y|_{i(X)})$ is a Liouville domain and $i\colon X\rightarrow i(X)$ is a Liouville diffeomorphism. With that in mind, we shall sometimes drop the embedding $i$ from our notation (and write, e.g., $Z_X$ rather than $Z_{i(X)}$ for the Liouville vector field of $i(X)$ near the boundary). Before proving the desired neighbourhood theorem, we show that positive codimension Liouville embeddings can be \enquote{straightened} near the boundary. 

\begin{lemma}\label{lemma:liouv-inc}
    Let $i\colon X \hookrightarrow Y$ be a positive codimension Liouville embedding. Then there exist trivial Liouville inclusions $X \hookrightarrow X'$, $Y \hookrightarrow Y'$ and a positive codimension Liouville embedding $i' \colon X' \hookrightarrow Y'$ extending $i$ such that the Liouville vector field $Z_{Y'}$ is tangent to $i'(X')$ near $\partial_\infty X'$. 
\end{lemma}

\begin{proof}
    Inside the Liouville completion of $Y$ given by $\hat{Y}=Y \bigcup_{\partial_\infty Y} \partial_\infty Y \times [0,\infty)_r$, we can consider the piecewise-smooth submanifold $\hat{X}_0\coloneqq X \bigcup_{\partial_\infty X} \partial_\infty X \times [0,\infty)_r$ given by attaching a cylinder on $\partial_\infty i(X)$ to $i(X)$ at the boundary $\partial_\infty Y$ (as earlier, we shall drop $i$ from the notation). The pushforward of $Z_X$ under the Liouville embedding $i$ has to be positively transverse to $\partial_\infty Y$, so there exists a constant $c>0$ such that $\hat{X}_0(c)\coloneqq\hat{X}_0 \cap (\partial_\infty Y \times (-c,\infty))$ is a family of contact submanifolds parameterised by the cylindrical coordinate $r$. Denote the projection onto the $r$-coordinate as $\pi \colon \hat{X}_0(c) \rightarrow (-c,\infty)$ and the contact submanifold $\pi^{-1}(s)=\hat{X}_0 \cap (\partial_\infty Y \times \{s\})$ as $\Sigma_0^s$ for $s>-c$. Since all the fibres $\partial_\infty Y \times \{s\}$ are contactomorphic via the Liouville flow of $\hat{Y}$, we can equivalently consider $(\Sigma^s_0)_{s \in (-c,\infty)}$ as a family of contact submanifolds of $\partial_\infty Y$. 

    By employing appropriate cut-off functions, we can pick a small $\varepsilon>0$ and reparameterise the smooth family $(\Sigma^s_0)_{s \in (-c,0]} \subset \partial_\infty Y$ to a smooth family $(\Sigma^s)_{s \in (-c,0]} \subset \partial_\infty Y$ that satisfies $\Sigma^s=\Sigma^s_0$ whenever $s<-2\varepsilon$ and $\Sigma^s=\Sigma^0_0$ whenever $s>-\varepsilon$. Gluing this to the constant family over non-negative real numbers then yields a smooth family $(\Sigma^s)_{s \in (-c,\infty)}$ and, therefore, a smoothing $\hat{X}(c)\coloneqq\bigcup_{s \in (-c,\infty)} \Sigma^s \times \{s\}$ of $\hat{X}_0(c)$ inside $\hat{Y}$. By construction, this smoothing will satisfy $\hat{X}(c) \cap \pi^{-1}(s)=i(X) \cap \pi^{-1}(s)$ for $s<-2\varepsilon$ and $Z_Y=\partial_r$ will be tangent to it over the locus $s>-\varepsilon$. By replacing the non-smooth cylindrical end $\hat{X}_0(c)$ with $\hat{X}(c)$, we also obtain a smoothing $\hat{X}$ of $\hat{X}_0$.

    We claim that $\hat{X}$ is a symplectic submanifold of $\hat{Y}$: for a point $(x,s) \in \hat{X}$ with $s \in (-c,\infty)$ and $x \in \Sigma^s$, we have a decomposition $T_{(x,s)}\hat{X}=\xi_s(x) \oplus R_s(x) \oplus V_s(x)$, where $\xi_s$ is the contact distribution of $\Sigma^s$, $R_s$ is its Reeb vector field and $V_s(x)$ is any vector that is symplectically orthogonal to $\xi_s(x)$ and whose $\partial_r$-component is equal to $1$ (the choice of $V_s(x)$ is well-defined up to a multiple of $R_s(x)$). Since $\xi_s(x) \cap \xi_s(x)^\omega=\{0\}$ and $R_s(x) \in \xi_s(x)^\omega$, it suffices to check that $\omega(V_s,R_s) \neq 0$ holds to show that $T_{(x,s)}\hat{X} \cap T_{(x,s)}\hat{X} ^\omega=\{0\}$. However, since we are working over the cylindrisation of $\partial_\infty Y$, we know that $\omega_{(x,s)}=e^s(dr \wedge \lambda_x+d\lambda_x)$, where $\lambda$ is the contact form on $\partial_\infty Y$. Therefore, $\omega_{(x,s)}(R_s(x),\cdot)=-e^s\lambda(R_s(x)) dr=-e^sdr$ and so $\omega(V_s(x),R_s(x))>0$, as desired.

    The Liouville vector field $Z=Z_{\hat{X}}$ will then be a linear combination of $V$ and $R$, since it is symplectically orthogonal to $\xi$. By the relation $\omega(Z_s(x),R_s(x))=1$ and the above calculation, we see that the coefficient of $V$ in such an expression must be positive, therefore, $Z_s$ has a positive $\partial_r$-component along $\Sigma^s \times \{s\}$ for all $s>-c$. Therefore, $\hat{X}$ is isomorphic to the Liouville completion of the Liouville domain $X_{\varepsilon}= \hat{X} \backslash \bigcup _{s>-2\varepsilon} \Sigma^s\times \{s\}$. Since $X_{\varepsilon} \hookrightarrow X$ is a trivial Liouville inclusion (because $\hat{X}$ agrees with $X$ over the locus where $r<-2\varepsilon$ by construction), this also means that $\hat{X}$ is isomorphic to the Liouville completion of $X$. Therefore, it suffices to take a suitably large Liouville domain $Y' \subset \hat{Y}$ and $X'\coloneqq Y' \cap \hat{X}$ to get the desired extension $i'$, because we have arranged that $Z_Y$ is tangent to $\hat{X}$ outside of a compact set. 
\end{proof}

Suppose that we are given a Liouville embedding $i \colon X \hookrightarrow Y$ of codimension $2k$, then we can consider an open tubular neighbourhood $N(X) \subset Y$ of $X$ with $\partial_\infty N(X)\coloneqq N(X) \cap \partial_\infty Y$. Note that $N(X)$ naturally inherits the structure of an exact symplectic manifold, but its closure might fail to be a Liouville domain, since we can not guarantee that $Z_Y$ is outward pointing along the part of $\partial \overline{N(X)}$ that is not contained in $\overline{\partial_\infty N(X)}$. However, we still have the following neighbourhood theorem for positive codimension Liouville embeddings.

\begin{proposition}\label{proposition:liouville-nbhd-theorem}
    Let $i_0 \colon X_0 \hookrightarrow Y_0$ and $i_1 \colon X_1 \hookrightarrow Y_1$ be positive codimension Liouville embeddings. Suppose that $\varphi \colon X_0 \rightarrow X_1$ is a Liouville diffeomorphism and $\Phi$ an isomorphism of symplectic normal bundles covering $\varphi$. Then there exist trivial Liouville inclusions $X_j \hookrightarrow X'_j$, $Y_j \hookrightarrow Y'_j$ such that $\varphi$, $i_0$ and $i_1$ extend to make the following diagrams commute:
    \[\begin{tikzcd}
    {X_0} & {X_0'} \\
    {X_1} & {X_1'}
    \arrow[hook, from=1-1, to=1-2]
    \arrow["\varphi"', from=1-1, to=2-1]
    \arrow["\varphi'", from=1-2, to=2-2]
    \arrow[hook, from=2-1, to=2-2]
    \end{tikzcd}
    \quad\textnormal{and}\quad
    \begin{tikzcd}
    X_j \arrow[hook, r] \arrow[hook, d, "i_j"'] & X_j' \arrow[hook,d, "i_j'"] \\
    Y_j \arrow[hook,r] & Y'_j
    \end{tikzcd}
    \quad\textnormal{for } j=0,1,\]
    with $\Phi$ also extending to an isomorphism $\Phi'$ of symplectic normal bundles of $X_j'$ inside $Y'_j$. Moreover, there is a diffeomorphism $\psi \colon N(X_0') \rightarrow N(X'_1)$ between neighbourhoods of $X'_j$ inside $Y'_j$ that extends $\varphi'$ such that $d\psi$ induces $\Phi'$ and we also have $\psi^*\lambda_{Y'_1}=\lambda_{Y_0'}+df$ for some smooth function $f$ vanishing identically along $X'_0$ and near $\partial_\infty N(X_0')$. 
\end{proposition}

\begin{proof}
    By Lemma \ref{lemma:liouv-inc}, we can replace $X_j$, $Y_j$ by larger subdomains $X'_j$, $Y'_j$ to arrange that the Liouville vector field of $Y_j$ is tangent to the subdomain $i_j(X_j)$ near the boundary, so without loss of generality suppose that that was the case to begin with to avoid notational clutter (clearly, one can also perform these enlargements simultaneously in a way that is compatible with the extension of $\varphi$ to an isomorphism of Liouville domains via the Liouville flow, which guarantees that the first diagram commutes, and the extension $\Phi'$ is also constructed through the Liouville flow). As in the case of \cite[Theorem 3.4.10]{mcbook}, we can use the isomorphism $\Phi$ of symplectic normal bundles along with suitably chosen exponential maps to obtain a diffeomorphism $\varphi' \colon N(X_0) \rightarrow N(X_1)$ extending $\varphi$. For the sake of brevity, denote $\omega_0=d\lambda_0$ the forms on $N(X_0)$ obtained by restricting the ones from $Y_0$ and $\omega_1=d\lambda_1$ the ones obtained by pulling back the forms from $Y_1$ through $\varphi'$. 

    Analogously to the closed case, the forms $\omega_0$ and $\omega_1$ agree on $TY_0|_{X_0}$ by construction, and we also know that there exists a compactly supported function $f_0$ on $X_0$ such that $\lambda_1|_{X_0}=\lambda_0|_{X_0}+df_0$. In particular, this tells us that $\lambda_1|_{X_0}=\lambda_0|_{X_0}$ near infinity. Because the Liouville vector fields of $\lambda_{Y_j}$ are assumed to be tangent to $X_j$ near infinity, the Liouville vector fields $Z_j$ associated to $\lambda_j$ will be tangent to $X_0$ near infinity. Since Liouville diffeomorphisms preserve Liouville vector fields near infinity, this means that $Z_0|_{X_0}=Z_1|_{X_0}$ holds close to $\partial_\infty X_0$, so we get a stronger statement that the forms $\lambda_j=\omega(Z_j,\cdot)$ agree on $TY_0|_{X_0}$ near infinity. 

    Now, we proceed analogously to the proof of \cite[Proposition 11.8]{CE}, but working relative to a subdomain $X_0$: denote $\lambda_s=s\lambda_1+(1-s)\lambda_0$, then the above discussion implies that $\lambda_0$ and $\lambda_s$ agree on $T(\partial_\infty X_0)|_{\partial_\infty Y_0}$ along the closed contact submanifold $\partial_\infty X_0 \subset \partial_\infty Y_0$ for all $s \in [0,1]$. Therefore, by the neighbourhood theorem for contact submanifolds (\cite[Theorem 2.5.15]{geiges} along with a slight refinement analogous to the one suggested in \cite[Remark 2.5.12]{geiges}), there exists a family of strict contactomorphisms $\Phi_s \colon (\partial_\infty N(X_0),\lambda_0) \rightarrow (\partial_\infty N(X_0),\lambda_s)$ that map $\partial_\infty X_0$ to itself (possibly after shrinking the neighbourhood $N(X_0)$) and induce the identity along $T(\partial_\infty Y_0)|_{\partial_\infty X_0}$. Since the corresponding Liouville vector fields $Z_s$ are all transverse to $\partial_\infty N(X_0)$, these can be extended to diffeomorphisms $\Psi_s \colon Op(\partial_\infty N(X_0)) \rightarrow Op(\partial_\infty N(X_0))$ that satisfy $\Psi^* \lambda_s=\lambda_0$ and fix $X_0$, $TY_0|_{X_0}$ (the final two conditions are satisfied as $Z_s=sZ_1+(1-s)Z_0$ is independent of $s$ along $X_0$ and tangent to it near infinity). Outside of the neighbourhoods, we can extend them arbitrarily, to get a family of diffeomorphisms $\Psi_s \colon N(X_0) \rightarrow N(X_0)$ fixing $X_0$, inducing the identity on $TY_0|_{X_0}$ and satisfying $\Psi_s^*\lambda_s-\lambda_0=0$ near the boundary $\partial_\infty N(X_0)$. 

    The desired result then follows by Moser's stability theorem for exact symplectic forms \cite[Theorem 6.8]{CE}, with the issue of flows only existing locally resolved by shrinking the tubular neighbourhood even further (cf. \cite[Lemma 6.10]{CE}).
\end{proof}

We can also use the data of a positive codimension Liouville embedding to construct a model neighbourhood of $X$: consider the embedding $i_{std}$ of $X$ into the total space of the the symplectic normal bundle $TX^\omega$. By \cite[Theorem 5.3]{avdek}, there will then exist a small tubular neighbourhood of the zero section $N_{std}(X) \subset \Tot(TX^\omega)$ and a $1$-form $\lambda_{std}$ that makes $\overline{N_{std}(X)}$ a Liouville domain with corners and agrees with $\lambda_Y$ on $TY|_{X}$ (even though the form $\lambda_{std}$ constructed there is not unique, the space of forms satisfying these conditions is clearly contractible, so the resulting Liouville domains are deformation equivalent and it does not matter which one we choose). Applying Proposition \ref{proposition:liouville-nbhd-theorem} then tells us that this is indeed a model neighbourhood for any positive codimension Liouville embedding of $X$:

\begin{corollary}\label{corollary:liouville-nbhd-theorem}
    Let $i \colon X \hookrightarrow Y$ be a  positive codimension Liouville embedding, then there exist trivial Liouville inclusions $X \hookrightarrow X'$, $Y \hookrightarrow Y'$ with an extension $i' \colon X' \hookrightarrow Y'$ and a diffeomorphism $\psi \colon N_{std}(X') \rightarrow N(X')$ that satisfies $\psi^*\lambda_Y=\lambda_{std}+df$ for some smooth function $f$ vanishing identically along $X'$ and near $\partial_\infty N_{std}(X')$.
\end{corollary}

In order to apply these results and finish the proof, it now remains to understand the symplectic normal bundle of the complete intersection $\widetilde{Z}$ inside the hypersurface $\widetilde{H}_\tot$.

\begin{lemma}\label{lemma:normal-bundle-trivial}
    The symplectic normal bundle $T\widetilde{Z}^\omega$ of $\widetilde{Z}$ inside $T\widetilde{H}_{\tot}$ is trivial.
\end{lemma}

\begin{proof}
    For convenience, denote the normal bundle as $E \rightarrow\widetilde{Z}_\beta$, then by \cite[Theorem 2.6.3]{mcbook}, we know that it suffices to find an $\omega$-tame complex structure $J$ on $E$ such that $(E,J)$ is isomorphic to $\widetilde{Z} \times \C^{r-1}$ as a complex vector bundle. 

    Note that via the musical isomorphism induced by $\omega$, the symplectic normal bundle of a symplectic submanifold is isomorphic to the conormal bundle. Analogously to Lemma \ref{lemma:angles_estimate}, we can show that the images of differential 1-forms defined as $v_{2j-1}=\partial \widetilde{f}_{j}-\overline{\partial}\widetilde{f}_{j}$, $v_{2j}=J_0(\partial \widetilde{f}_{j}+\overline{\partial}\widetilde{f}_{j})$ for $j=1,\dots,r$ (where $J_0$ is the complex structure coming from $M_{\Cs}$) under the musical isomorphism define a trivialisation of the symplectic normal bundle $E'$ of $\widetilde{Z}$ inside of $M_{\Cs}$ as a real vector bundle.

    We can define a complex structure $J$ on $E'$ by declaring that it sends $v_{2j-1} \mapsto v_{2j}$ and $v_{2j} \mapsto -v_{2j-1}$ for all $j$. In particular, if all the antiholomorphic parts $\overline{\partial} \widetilde{f}_{j}$ vanish, this will agree with $J_0$. In general, we claim that $J$ will be $\omega$-tame for all $\beta>0$ large enough. Clearly, it suffices to show that the matrix $M(z)$ with entries $M_{jk}=\frac{\omega(v_j,Jv_k)}{\lVert v_j \rVert_g \lVert v_k \rVert_g}$ is positive definite at all $z \in \widetilde{Z}$, where we are using the Riemannian metric $g$ associated to $\omega$ and $J_0$ to normalise the entries. 

    From Lemma \ref{lemma:tailoring_estimate}, we know that there is a constant $C>0$ such that $\lVert \partial \widetilde{f}_{j}(z) \rVert >Ce^{-\sqrt{\beta}} \cdot \lVert \overline{\partial} \widetilde{f}_{j}(z) \rVert$ holds for all $z \in \widetilde{Z}$, $\beta>0$ and $j=1,\dots,r$. This means that $M$ is going to be $O(e^{-\sqrt{\beta}})$-close to the matrix $\hat{A}$ with entries $\hat{A}_{jk}=\frac{\langle \hat{v}_j,\hat{v}_k \rangle_g}{\lVert \hat{v}_j \rVert_g \lVert \hat{v}_k \rVert_g}$ for $\hat{v}_{2j-1}=\partial \widetilde{f}_{j}$, $\hat{v}_{2j}=J\partial \widetilde{f}_{j}$. By Corollary \ref{corollary:cursed-determinant} and Remark \ref{remark:cursed-determinant}, we know that the smallest eigenvalue of  $\hat{A}(z)$ is bounded below by some $\varepsilon'>0$ for all $z \in \widetilde{Z}_{\beta}$. This implies the inequality $w^T\hat{A}w \geq \varepsilon'\cdot w^Tw$ for all $w \in \R^{2r}$, so combining that with $M$ being close to $\hat{A}$ gives a bound $w^T M w \geq (\varepsilon'-C'e^{-\sqrt{\beta}})w^T w$ for some positive $\varepsilon'$, $C'$ and all $z \in \widetilde{Z}_\beta$, $w \in \R^{2r}$, hence $w^T M w>0$ for all $w\neq 0$ and all $\beta$ large enough, so $J$ is indeed $\omega$-tame on $E'$.
    
    To see that $E=E' \cap T\widetilde{H}_{\tot}$ is trivial, observe that $v_{2j-1}-v_{2r-1}$ and $v_{2j}-v_{2r}$ for $j=1,\dots,r-1$ are all tangent to $\widetilde{H}_{\tot}$ and linearly independent, hence they trivialise $E$ as a real vector bundle. By definition, $J_0$ preserves the span of these vectors, so it also defines an $\omega$-tame complex structure on $E$ that identifies it with a trivial rank $r-1$ complex vector bundle, as desired.
\end{proof}

\begin{corollary}\label{corollary:lag-subbundle-trivial}
    Let $F \subset E$ be the sub-bundle of the symplectic normal bundle of $\widetilde{Z}$ in $T\widetilde{H}_{\tot}$ spanned by Hamiltonian vector fields $w_j=X_{\Im(\widetilde{f}_{j})}-X_{\Im(\widetilde{f}_{r})}$ for $j=1,\dots,r-1$. Then there exists a symplectic trivialisation $\Phi \colon (E,\omega) \rightarrow \widetilde{Z} \times (\C^{r-1},\omega_{\C^{r-1}})$ that identifies the restriction of $F$ to the positive real locus $\widetilde{Z}^+$ with the trivial Lagrangian sub-bundle $\widetilde{Z}^+ \times \R^{r-1}$.
\end{corollary}

\begin{proof}
    From the local expression for $X_{\Im(\widetilde{f}_{j})}$ along $\widetilde{Z}^+$ from Lemma \ref{lemma:liou-dynamics}, it is clear that $\omega(w_j,w_k)=0$ along the locus, so the sub-bundle $F|_{\widetilde{Z}_\beta^+}$ is indeed Lagrangian. Since the complex structure $J$ defined in Lemma \ref{lemma:normal-bundle-trivial} is $\omega$-tame, this also means that the sub-bundle is also totally real with respect to it. We observe that the space of symplectic forms on $E$ that tame $J$ and make the sub-bundle $F|_{\widetilde{Z}^+}$ Lagrangian is contractible (by an argument analogous to \cite[Proposition 2.6.4]{mcbook} and \cite[Exercise 2.6.5]{mcbook}). Let $\Psi$ be the complex trivialisation constructed in Lemma \ref{lemma:normal-bundle-trivial}, then both $\omega_0=\omega$ and $\omega_1=\Psi^*\omega_{\C^{r-1}}$ lie in the space, hence there exists a homotopically unique path of forms $\omega_t$ between them. 

    Analogously to the proof of Theorem 2.6.3 in loc. cit., we can use Corollary 2.5.14 to obtain a symplectic trivialisation, if one can guarantee that the Lagrangian subspaces can be preserved throughout the deformation constructed there. This follows from an observation that for any smooth family of Lagrangian subspaces $(L_t)$ of a symplectic space $(V,\omega)$, there exists a smooth family of symplectomorphisms $\varphi_t$ satisfying $\varphi_t(L_0)=L_t$ (we can pick a compatible complex structure to define canonical Lagrangian complements $(K_t)$ varying smoothly with $t$ and we have, by definition, a smooth family of isomorphisms $\psi_t \colon L_0 \rightarrow L_t$, which then extend uniquely to smooth family symplectomorphisms $\varphi_t\colon L_0 \oplus K_0 \rightarrow L_t \oplus K_t$, as desired), and we can compose the smooth family of isomorphisms constructed there with this one. 
\end{proof}

\begin{proof}[Proof of Theorem \ref{theorem:emb}]
    Let $A \subset \widetilde{Z}$, $U_0 \subset M_{\Cs}$ be the Liouville domain constructed in Lemma \ref{lemma:vf_not_tangent} and the subsequent discussion. Note that there exists a Liouville subdomain $V_0 \subset H_{\tot}$ such that $A \hookrightarrow V_0$ is a positive codimension Liouville embedding: $A' \hookrightarrow V'$ for $V'\coloneqq \widetilde{H}_{\tot} \cap \{\varphi \le C\}$ will be a Liouville embedding, so it just suffices to shrink $V'$ along with $A'$ appropriately. Therefore, by Lemma \ref{lemma:normal-bundle-trivial} and Corollary \ref{corollary:liouville-nbhd-theorem}, we can slightly inflate the subdomains inside $\widetilde{H}_{\tot}$ to obtain an embedding $i_0 \colon A \times D^{2r-2}_\varepsilon(0) \hookrightarrow V_0$ that satisfies $i^*\lambda_{V}-\lambda_{\widetilde{Z}}-\lambda_{D}=df$ for some function $f$ that vanishes along $A$ and near $\partial_\infty A \times D^{2r-2}_\varepsilon(0)$. 

    We can also ensure that the embedding satisfies $i_0 (\Lambda \times (-\varepsilon, \varepsilon)^{r-1}) \subset \Phi^{-1}(\eL_\Sigma)$: first, we note that this can be guaranteed for a smooth embedding $\varphi'$ constructed via exponential maps Proposition \ref{proposition:liouville-nbhd-theorem}. That is because we can consider the trivialisation of $T\widetilde{Z}^\omega$ from Corollary \ref{corollary:lag-subbundle-trivial} and, in particular, use the flows of vector fields $w_1$, \dots, $w_{r-1}$ from there to trivialise the Lagrangian sub-bundle along $\widetilde{Z}^+$ (we have checked that their Poisson bracket is zero, so the flows also commute). Then, by the calculation from Lemma \ref{lemma:liou-dynamics}, we know that the Hamiltonian vector fields $X_{\Im(\widetilde{f}_j)}$ are tangent to $\Phi^{-1}(\eL_\Sigma)$ along $\Lambda$, so $w_1$, \dots, $w_{r-1}$ are also tangent to $\Phi^{-1}(\eL_\Sigma)$, which implies the desired statement. The deformation near the boundary of the Liouville domain does not affect the skeleta, so we can also guarantee that the maps $\Psi_s$ also fix $\Lambda \times (-\varepsilon,\varepsilon)^{r-1}$. Finally, it is straightforward to check that the smooth Lagrangian strata of $\Lambda \times (-\varepsilon,\varepsilon)^{r-1}$ are strongly exact with respect to both the Liouville forms $\lambda_0$ and $\lambda_1$, so the skeleton must also be preserved by the deformation obtained through Moser stability: the vector field $V_s$ whose flow is considered there satisfies $\omega_s(V_s,\cdot)=\lambda_1-\lambda_0$, so if $L$ is a strongly exact Lagrangian with respect to both $\lambda_j$'s, it is a strongly exact Lagrangian with respect to all $\lambda_s$ and $\omega_s(V_s,\cdot)|_L=0$, thus $V_s \in( TL)^{\omega_s}=TL$ and $V_s$ is tangent to $L$, hence its flow preserves $L$ (one can use Lemma \ref{lemma:flows_corners} to make this more precise for Lagrangians with corners, analogously to Corollary \ref{corollary:contains}). Moreover, we have shown that the primitive $f$ must vanish identically along $\Lambda$, since $df|_{\Lambda_i}=0$ by the observed strong exactness property for smooth strata $L_i \subset \Lambda \times (-\varepsilon,\varepsilon)^{r-1}$, $f|_{\widetilde{Z}}=0$ and connectedness of $\Lambda$.

    We can now finish analogously to \cite[Theorem 5.5]{avdek}: since $V_0 \subset \partial U_0$ is a Liouville hypersurface, it admits a neighbourhood of the form $([-\delta,\delta]_s\times V_0,ds+\lambda_{V_0}) \subset (\partial U_0,\lambda_{U_0})$ for a small $\delta>0$. By decreasing $\varepsilon>0$, we can assume that the primitive satisfies $|f|<\delta$. Therefore, if we compose the embedding $i$ with the diffeomorphism $V_0 \rightarrow \Gamma(f) \subset [-\delta,\delta]_s\times V_0$ given by flowing in the Reeb flow direction for time $f$, we get an embedding $i_f \colon A \times D^{2r-2}_\varepsilon(0) \hookrightarrow \Gamma(f)$ that satisfies $i_f^*\lambda_{\Gamma(f)}=\lambda_{\widetilde{Z}}+\lambda_{D}$. Since $f$ vanishes along $\widetilde{Z}$, the restrictions of $i_0$ and $i_f$ to it agree, so the map $i_f$ extends the tautological embedding $A\hookrightarrow \partial U_0$. Analogously, since $f$ vanishes along $\Lambda$, the image of $\Lambda$ will be a subset of $\Phi^{-1}(\eL_\Sigma) \cap \partial U_0$. Therefore, the desired statement follows for $i=\Phi \circ i_f$ and $U=\Phi(U_0)$.
\end{proof}

\section{HMS for open Batyrev--Borisov complete intersections}\label{section:hms}
In this section, we use the results of preceding sections to prove homological mirror symmetry for open Batyrev--Borisov complete intersections and explore some of its interesting features. Throughout the section, we assume that we have an irreducible nef partition $\nabla=\nabla_1+\dots+\nabla_r$, refined centred triangulating function $h$ on $\Delta^\vee$, potential $\varphi$, cut-off functions $\chi_{\alpha,\beta}$ and sufficiently large $\beta \gg0$ such that the results of Section \ref{section:skeleta} and \ref{section:combinatorics} apply for $Z\coloneqq \widetilde{Z}_\beta$ with a skeleton $\Lambda$. 

\subsection{B-side mirror}\label{section:b-side}

We can formalise the construction of the mirror stack $\check{Z}$ outlined in Section \ref{section:introduction} within this set-up: inside the ambient toric stack $\mX_\Sigma$, one can consider the closure of toric orbits $\overline{O(\sigma)}$ for all $\sigma \in \Sigma_\trans$ and take their union $\check{Z}\coloneqq\partial_\trans \mX_\Sigma \subset \partial \mX_\Sigma$. Note that when $r=1$, the transversal boundary becomes the entire toric boundary. By adapting the argument for that special case, we can also show that the gluing $\partial_\trans \mX_\Sigma$ is a colimit of stacks (where all the orbit closures are equipped with the reduced structure as closed substacks of $\mX_\Sigma$).

\begin{lemma}\label{lemma:mirror-space-colimit}
    In the category of algebraic stacks, we have $\partial_\trans \mX_\Sigma= \colim_{\sigma \in \Sigma_\trans} \overline{O(\sigma)}$.
\end{lemma}

\begin{proof}
    This is just a variation on \cite[Lemma 3.4.1]{GS22}: we start from a diagram $\{ \overline{O(\sigma)}\}$ indexed by $\Sigma^\trans$ with maps $\sigma \hookrightarrow \tau$ inducing closed embeddings $\overline{O(\tau)} \hookrightarrow \overline{O(\sigma)}$. Therefore, by \cite[Corollary 3.9]{schwede}, the colimit exists (we can compute it iteratively via pushouts, since the indexing set is a poset) and it admits a map into $\partial_\trans \mX_\Sigma$, so it suffices to check that it is an isomorphism.

    The question is étale local, we can check the statement over maximal affine toric charts $X_\sigma$ (that are varieties rather than stacks) indexed by full-dimensional simplicial cones $\sigma \in \Sigma$. In order to have non-empty intersection with $\partial _\trans \mX_\Sigma$, the cone has to be transversal, so we can write $\sigma=\sigma_1+\dots+\sigma_r$ for simplicial cones $\sigma_j$ of positive dimension. Recall that the ring of functions on $X_\sigma$ is just $\C[\sigma^\vee]$. For a subcone $\tau \subsetneq \sigma$, there is a dual face $\tau^\perp \subset \partial \sigma$ that corresponds to the characters that do not vanish along $\overline{O(\tau}) \subset X_\sigma$. In particular, this means that the vanishing ideal of $\overline{O(\tau)}$ is $I_\tau=\C[\sigma^\vee\backslash\tau^\perp]$, so its function ring is isomorphic to $\C[\sigma^\vee]/I_\tau \cong \C[\tau^\perp]$ (observe that if we have a semigroup $G$ and a semigroup ideal $I \leq G$, then $\C[I]$ is an ideal of the semigroup ring $\C[G]$). 

    Therefore, if we denote the vanishing ideal of $\partial_\trans X_\sigma$ as $I_\trans$, it suffices to check that $\C[\sigma^\vee]/I_\trans \rightarrow \lim_{\tau \in \sigma_\trans} \C[\sigma^\vee]/I_\tau$ is a ring isomorphism (where an inclusion $\tau_1 \rightarrow \tau_2$ induces a map $\sigma^\vee\backslash \tau_1^\perp \rightarrow \sigma^\vee \backslash \tau_2^\perp$ and hence a map between the quotient rings going in the same direction). By the general description of fibre products in the category of rings, an element of the RHS can be described as a tuple of elements $p_\tau \in \C[\sigma^\vee]/I_\tau$ for minimal transversal cones $\tau \in \sigma_\trans$ such that whenever $\tau_1$, $\tau_2$ are minimal transversal subcones of $\tau'$, then $p_{\tau_1}/I_{\tau'}=p_{\tau_2}/I_{\tau'} \in \C[\sigma^\vee]/I_{\tau'}$. Since the coordinate ring of $\overline{O(\tau)}$ is also isomorphic to $\C[\tau^\perp]$, each $p_\tau$ is just a finite linear combination $p_\tau=\sum_{\chi \in \tau^\perp} c_{\chi,\tau} \cdot \chi$. The coherence conditions guarantee that $c_{\chi,\tau_1}=c_{\chi,\tau_2}$ holds whenever $\chi \in \tau_1^\perp \cap \tau_2^\perp$ (by looking at the quotient by the maximal cone $\tau'$ such that $\chi \in (\tau')^\perp$), so every element on the RHS corresponds to a finite linear combination of characters $\chi \in \bigcup_{\tau \in \sigma_\trans} \tau^\perp$, which can be viewed as an element of $\C[\sigma^\vee]$, so the map from $\C[\sigma^\vee]$ into the limit is surjective. This description also tells us that the kernel of this map is given by $\bigcap_{\tau \in \sigma_\trans} I_\tau$, which is equal to the vanishing ideal of $\bigcup_{\tau \in \sigma_\trans} \overline{O(\sigma)}=\partial_\trans X_\sigma$, so the kernel is equal to $I_\trans$ by definition and we are done. 

\end{proof}

\begin{figure}[ht]
    \centering
    \tdplotsetmaincoords{70}{120} 
    \begin{tikzpicture}[tdplot_main_coords, scale=1.2]
    \draw[thick, rounded corners = 5mm] (1,1,1) -- (-1,1,1) -- (-1,1,-1) -- (-1,-1,-1) -- (-1,-1,1) -- (1,-1,1) -- (1,-1,-1)-- (1,1,-1) -- cycle;

    \draw[thick, black] (1,7,2) -- (-1,7,2) -- (-1,7,0) -- (-1,5.95,0);
    \draw[thick,black] (-1,5.75,0) -- (-1,5,0) -- (-1,5,2) -- (1,5,2) -- (1,5,0)-- (1,7,0) -- (1,7,2);
    \draw[line width=4pt, white] (1,1,-0.1)--(1,1,-0.3);
    \draw[thick] (1,1,-0.1)--(1,1,-0.3);
    \begin{scope}[canvas is xy plane at z=0.5]
        \draw[line width=4pt, color=white] (1.6,1.6) circle (10pt);
        \draw[thick] (1.6,1.6) circle (10pt);
    \end{scope}
    \draw[thick] (1,1,0.4)--(1,1,0.2);
    
    \begin{scope}[canvas is xy plane at z=0]
        \draw[thick] (-1.25,1.25) circle (10pt);
    \end{scope}
    \begin{scope}[canvas is xy plane at z=-0.25]
        \draw[thick] (-1.25,-1.25) circle (10pt);
    \end{scope}
    \draw[line width=4pt, white] (-1,-1,-0.1)--(-1,-1,0.1);
    \draw[thick] (-1,-1,-0.1)--(-1,-1,0.1);
    \begin{scope}[canvas is xy plane at z=0]
        \draw[thick, blue] (1.25,-1.25) circle (10pt);
    \end{scope}
    \begin{scope}[canvas is yz plane at x=0.37]
        \draw[line width=4pt, white] (1.2,1.45) circle (10pt);
        \draw[thick] (1.2,1.45) circle (10pt);
    \end{scope}
    \draw[thick] (-0.1,1,1)--(0.3,1,1);
    \begin{scope}[canvas is yz plane at x=0.37]
        \draw[line width=4pt, white] (-0.8,1.5) circle (10pt);
        \draw[thick] (-0.8,1.5) circle (10pt);
    \end{scope}
    \draw[thick, black] (0.2,-1,1)--(-0.38,-1,1);
    \begin{scope}[canvas is xz plane at y=-0.2]
        \draw[thick, red] (0.9,-1.42) circle (10pt);
    \end{scope}
    \begin{scope}[canvas is xz plane at y=0]
        \draw[thick] (-1.4,-1.55) circle (10pt);
    \end{scope}
    \draw[<->,dashed] (0,2.5,0.4) -- (0,4,0.65);
    \node[label={$\textnormal{skel}(Z)$}] at (-0.25,2,1) {};
    \node[label={$\check{Z}$}] at (-0.25,7.75,2) {};

    \draw[very thick,blue] (1,5,2) node[left] {$\color{blue}\mathbb{P}^1$} -- (1,5,0) ;
    \draw[very thick,red] (1,5,0) -- (1,7,0) node[below] {$\color{red}\mathbb{P}^1$};
    \node[fill,circle,inner sep=1.5pt,violet] at (1,5,0) {};
    \draw[|-|,thick, blue, rounded corners= 5mm] (0.2,-1,1) -- (1,-1,1) node[left] {$\color{blue}\mathbb{L}_{\mathbb{P}^1}$} -- (1,-1,-1) -- (1,-0.2,-1);
    \draw[|-|,thick, red, rounded corners=5mm] (1,-1,-0.2) -- (1,-1,-1) -- (1,1,-1) node[below] {$\color{red}\mathbb{L}_{\mathbb{P}^1}$} -- (1,1,0);
    \draw[line width=1.5pt, violet, rounded corners=5mm] (1,-1,-0.2) -- (1,-1,-1) node[below left] {$\color{violet}\mathbb{L}_{\bullet} \times I$} -- (1,-0.2,-1);
\end{tikzpicture}
    \caption{Idea behind our proof of HMS illustrated on the running example}
    \label{fig:mirror-covers}
\end{figure}
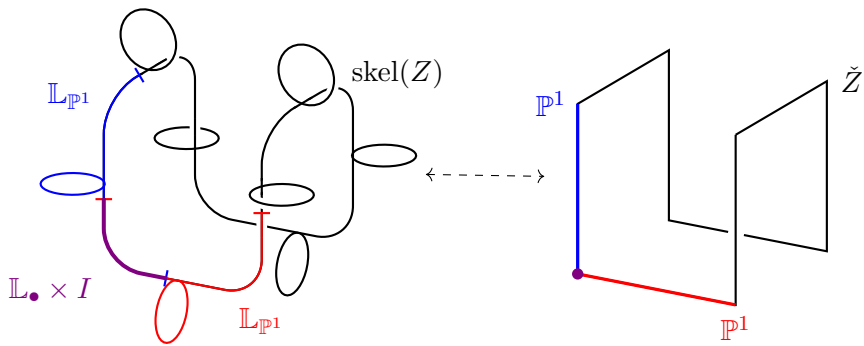

\begin{example}\label{example:running-b-side}
    In the setting of Example \ref{example:running-nef}, the B-side mirror will be given as the chain of 8 $\Pp^1$'s. Note that the corresponding subcomplex inside the moment polytope $\Delta$ of $X_\Delta$ is clearly isomorphic to the bounded part of $\mA_\trop$ from Example \ref{example:running-trop}. This matches the decomposition of the A-side skeleton (see Example \ref{example:running-ex-cover}, Fig. \ref{fig:skeleton}), as illustrated by Figure \ref{fig:mirror-covers}.
\end{example}

Since our definition of $\check{Z}$ is based on a somewhat ad hoc combinatorial gluing construction intended to match the description of the skeleton of $Z$, we finish the section by explaining how it fits within some well-known frameworks for producing mirror pairs. 

\begin{remark}\label{remark:gross-siebert-mirror}
    Since $\check{Z}$ is a singular stack obtained by gluing smooth toric stacks along their boundaries, it could also naturally arise as a central fibre of a toric degeneration. Such spaces play a central role in the \emph{Gross--Siebert programme} (see \cite{gs06}, \cite{gs11}), which is a framework for mirror symmetry that is more general than (resolutions of) Calabi--Yau complete intersections in Fano toric varieties. It is discussed in great detail in \cite{Gross} how to view Batyrev--Borisov mirror pairs as a special case of the Gross--Siebert programme. Through that lens, both \cite{GS22} and the present work can be seen as studying the A-side mirrors of the central fibres of particular Gross--Siebert toric degenerations. 
\end{remark}

\begin{remark}\label{remark:aak-mirror}
    Since our A-side space $Z$ is a complete intersection inside an affine toric variety, it falls within the scope of the Abouzaid--Auroux--Katzarkov mirror construction from \cite{aak}. The candidate B-side mirror that one would obtain by running their argument is a certain \emph{Landau-Ginzburg model} $(Y,W_Y)$, where $Y$ is a non-compact smooth toric stack and $W_Y$ a regular function on it. By using techniques similar to Sections \ref{section:tropical-ci} and \ref{section:tropical-bbci}, we can translate between the combinatorics of our set-up and the more general one involving the Cayley trick to prove that $Y=\Tot(\mL_1^\vee \oplus \dots \oplus \mL_r^\vee \rightarrow \mX_\Sigma)$, where the ample line bundles $\mL_j \rightarrow \mX_\Sigma$ satisfying $\mL_1 \otimes \dots \otimes\mL_r \cong -K_{\mX_{\Sigma}}$ come from the decomposition of rays $\Sigma(1)=\bigsqcup_{j=1}^k \Sigma_j(1)$ induced by the nef-partition, while the function $W_Y=\sum_{j=1}^r p_js_j$ is obtained by evaluating the defining sections $s_j \in \Gamma(\mX_\Sigma,\mL_j)$ to get fibrewise linear functions $p_js_j$ on $\Tot(\mL_j^\vee)$, see \cite{barbacovi-segal} for more details on the notation. The critical locus of $W_Y$ can then be naturally identified with $\check{Z} \subset Y$ inside the zero section, so by Orlov's Knörrer periodicity, we have an equivalence $\coh(\check{Z}) \cong MF(Y,W_Y)$, and Theorem \ref{theorem:hms-for-bbci} also implies a special case of homological mirror symmetry for the AAK mirrors (which is the set-up studied in \cite{msz}).
\end{remark}

\subsection{Categorical background}\label{section:categories}

Before presenting our proof of the main HMS equivalences, we give an overview of all the categories that appear in the statements and the proofs. We follow the standard conventions for $A_\infty$-categories from \cite{seidel-book} and \cite{sheel-thesis}, and all the categories that appear are assumed to be $\Z$-graded and $\C$-linear. In particular, we write $\perf \mC$ for the split-triangulated closure of an $A_\infty$-category $\mC$. Notions such as fullness, faithfulness or commutativity diagrams should always be understood at the level of the cohomology categories (we will always emphasise when one can make certain statements strict on the chain level).

For the algebraic categories involved, we will follow the conventions of \cite{grrrrr1} and work with derived pre-triangulated dg-categories (or, equivalently, $\C$-linear stable $\infty$-categories). In particular, \emph{derived category of an algebraic variety $X$} will be denoted as $\coh(X)$, and should understood to be the bounded derived category of abelian category $coh(X)$ of coherent sheaves on $X$, considered as a $\C$-linear dg-category. We will also briefly need the larger category of \emph{ind-coherent sheaves on $X$}, $\indcoh(X)$, defined by gluing the ind-completions of $\coh(U)$ over affine $U$ with an étale morphism $U \rightarrow X$ (see \cite{Gaitsgory2013} for more details; note that this is not equivalent to the ind-completion of $\coh(X)$ in general). In particular, one can recover $\coh(X) \subset \indcoh(X)$ as the full subcategory of compact objects. All the functors between derived categories should also be understood as derived, for example, the pushforward $f_* \colon \indcoh(X) \rightarrow \indcoh(Y)$ by a morphism $f\colon X \rightarrow Y$ will always denote the derived pushforward. We also make use of all these notions for Deligne--Mumford stacks, starting with the derived dg-category $\coh(\mX)$, which can be understood as equivariant versions of the scheme versions in suitable contexts. 

We will use the set-up considered in \cite{gps1} for partially wrapped Fukaya categories. The only difference is that we will also require the existence of extra data necessary to define $\Z$-gradings rather than $\Z/2$-gradings. We shall not discuss this in any detail, see \cite[Section 5.3]{gps3} for a review of gradings in this setting (we are interested in algebraic complete intersections in $(\Cs)^n$, which come with a canonical stable polarisation).

We finish this overview by recalling some generalities surrounding wrapped microlocal sheaves following the original treatment in \cite[Section 3]{nadler16}. We also mostly stick to the notation used there; in particular, we use the superscript $\diamondsuit$ whenever there are no restrictions on sizes of stalks. Note that we purposefully avoid this language in the statements of our main results for the sake of simplicity, but microlocal sheaves play a crucial role in the proofs, since they provides the necessary local-to-global principle on the A-side. 

For a Whitney stratification $\mS$ of a manifold $X$, we denote the category of $\mS$-constructible sheaves of complexes of $\C$-vector spaces on $X$ as $Sh^\diamondsuit_{\mS}(X)$ (recall that a sheaf is \emph{$\mS$-constructible} if its restriction to any stratum $S \in \mS$ is locally constant); this clearly has a natural structure of a pre-triangulated dg-category over $\C$. If $\mS'$ is a refinement of $\mS$, there is a natural restriction functor $Sh^\diamondsuit_{\mS}(X) \hookrightarrow Sh^\diamondsuit_{\mS'}(X)$, so we define the category of \emph{(large) constructible sheaves} on $X$ as $Sh^\diamondsuit(X)\coloneqq \bigcup_{\mS} Sh^\diamondsuit_{\mS}(X)$.

Every sheaf $\mF$ on $X$ also has a \emph{microsupport} (or \emph{singular support}; see \cite[Section 3.3]{nadler16}), which is a closed conical subset $ss(\mF) \subset T^*X$. For any subset $Z \subset T^*X$, we write $Sh_Z^\diamondsuit(X)$ for the constructible sheaves whose support is contained in $Z$. When $\mF$ is constructible, the microsupport will be a singular Lagrangian inside $T^*X$. Therefore, the categories $Sh^\diamondsuit_\Lambda(X)$ for a conical singular Lagrangian $\Lambda \subset T^*X$ will be of particular interest to us. For a stratification $\mS$, one can show that $Sh^\diamondsuit_{\mS}(X)=Sh^\diamondsuit_{N^*\mS}(X)$, so this is compatible with the notation introduced earlier. Moreover, for a singular Lagrangian $\Lambda \subset T^*X$, one can always find a triangulation $\mS$ such that $\Lambda \subset N^*\mS$, so there is a fully faithful embedding $Sh^\diamondsuit_\Lambda(X) \hookrightarrow Sh^\diamondsuit_{\mS}(X)$ and any sheaf in $Sh^\diamondsuit_\Lambda(X)$ is $\mS$-constructible, therefore, we have not lost anything by imposing the constructibility assumption from the start. When $U \subset X$ is an open set, we write $Sh^\diamondsuit_{\mS}(U)$ for sheaves on $U$ constructible with respect to $\mS \cap U$, and $Sh^\diamondsuit_\Lambda(U)$ as sheaves on $U$ with microsupport contained in $\Lambda \cap \pi^{-1}(U)$ (where $\pi \colon T^*X \rightarrow X$ is the projection onto the base). 

One can take the idea of considering $Sh^\diamondsuit(U)$ for different open subsets $U \subset X$ further and construct a sheaf of dg-categories on $T^*X$ whose global sections are $Sh^\diamondsuit(X)$ through a process called \emph{microlocalisation}: more precisely, the sheaf $\msh_\Lambda^\diamondsuit$ on $T^*X$ is obtained by sheafifying
\begin{equation*}
    \msh^{\diamondsuit,\,pre}_\Lambda(U)\coloneqq \Sh_\Lambda^\diamondsuit(X)/\Sh^\diamondsuit_{\Lambda\backslash (\Lambda \cap U)}(X).
\end{equation*}
This is, in fact, a pushforward of a sheaf of categories on $\Lambda$.  To be more precise, the sections of $\msh^\diamondsuit_\Lambda(U)$ for a conical singular Lagrangian $\Lambda$ are compactly generated dg-categories by \cite[Proposition 4.10]{coco}, and restriction functors admit both adjoints, so $\msh^\diamondsuit_\Lambda$ should be understood as a sheaf valued in the $\infty-$category $^*DG^*$ of compactly generated dg-categories with functors admitting both adjoints (sheafification of $\msh^{\diamondsuit,pre}_\Lambda$ is also taken inside $^*DG^*$, but is equivalent to sheafifying inside the $\infty-$category of dg-categories).

The \emph{large} constructible versions of constructible and microlocal sheaves introduced above are reminiscent of ind-coherent sheaves defined in the previous section. To obtain the natural analogues of $\coh$, we pass to compact objects and define the \emph{wrapped constructible sheaves} as $\Sh(X) \coloneqq \Sh^\diamondsuit(X)^c$ and \emph{wrapped microlocal sheaves} as $\msh_\Lambda(U) \coloneqq \msh_\Lambda^\diamondsuit(U)^c$ for all $U \subset T^*X$. This has the effect of reversing the arrows for the following formal reason: let $dg$ be the $\infty-$category of small pretriangulated dg-categories with exact functors, then taking ind-completion and passing to adjoints defines a fully faithful functor $dg \hookrightarrow {^{**}DG} \cong (^*DG^*)^{op}$, and its image are precisely the compact objects by \cite[Propositions 5.3.5.11, 5.4.2.4]{tapas}. Therefore, colimits in $dg$ are turned into limits in $^*DG^*$, and $\msh_\Lambda$ can be naturally viewed as a cosheaf valued in $dg$.  

Finally, in order to define the sheaf $\msh_\Lambda^\diamondsuit$ for $\Lambda$ a skeleton of a general Weinstein domain $Y$, one can follow the ideas of \cite{Shende21} (see also \cite[Section 7]{gps3} for a slightly different perspective) and use a Liouville hypersurface embedding of a suitable stabilisation $Y'$ of $Y$ into some cotangent bundle $T^*X$, which corresponds to a thickening of $\Lambda$ to a mostly Legendrian stop $\Lambda' \subset S^*X$, so we can take $\msh^\diamondsuit_\Lambda\coloneqq\msh^\diamondsuit_{\Lambda'}|_\Lambda$.

\subsection{Homological mirror symmetry}\label{section:proof-hms}

The descriptions of the cover of the skeleton provided by Theorem \ref{theorem:combi-ci} and the construction of $\check{Z}$ from Section \ref{section:b-side} are in line with the strategy for gluing mirror symmetry of toric stacks that first appeared in \cite{GS22} and was then developed more systematically in \cite{GS23}. We explain how to make this connection more precise and use it to prove homological mirror symmetry for open Batyrev--Borisov complete intersections. 

Recall that the skeleton $\Lambda$ of $Z$ is obtained by looking at a certain codimension $r$ sphere $S\coloneqq\partial\widetilde{C}_{0,\beta} \subset M_\R$ and taking $\Lambda\coloneqq\Phi^{-1}(\eL_\Sigma) \cap (M_{S^1} \times S) \subset M_{\Cs}$. We have also seen that for an adapted potential $\varphi$, $\Phi(S)$ intersects all the cones $\sigma \in \Sigma$ transversely (and that it has non-empty intersections only with transversal cones). Therefore, following \cite[Example 2.1]{GS23}, $S$ naturally comes with a structure of a fanifold $(S,\Sigma)$ from pulling back the stratification of $N_\R$ by $\Sigma$. 

The space $\check{Z}$ considered earlier can then be viewed as a special case of the B-side construction from Section 3 in loc. cit. associated to the fanifold $(S,\Sigma)$, as witnessed by the expression from Lemma \ref{lemma:mirror-space-colimit}. The situation is more subtle on the A-side -- even though one could follow the general recipe (Section  4 in loc. cit.) for a construction of a Weinstein domain $\mathbf{W}(S,\Sigma)$ from a fanifold for $(S,\Sigma)$, the space is of a fundamentally different nature than $Z$, so it is unclear how the two of them are related. One is led to the following conjecture:

\begin{conjecture}\label{conjecture:equivalence-with-fanifold}
    There exists an isomorphism of Weinstein domains $Z \cong \mathbf{W}(S,\Sigma)$.
\end{conjecture}

We make no attempt to prove Conjecture \ref{conjecture:equivalence-with-fanifold}, but merely comment on why one should expect a statement like that to hold and how it is relevant to the problem at hand. Note that our Theorem \ref{theorem:combi-ci} provides a diffeomorphism $\Lambda \cong \skel(\mathbf{W}(S,\Sigma))$ along with local gluing descriptions analogous to the ones in \cite[Theorem 4.1]{GS23}, but since the Lagrangians are singular, this is not sufficient to identify their neighbourhoods. The most natural way to proceed would therefore be to prove a Weinstein neighbourhood theorem for singular Lagrangians locally modelled on spaces like $\eL_\Sigma$ (i.e. positive conormal bundles of hypersurface arrangements inside of a cotangent bundle). The closest analogue of such a statement in the literature is perhaps \cite[Theorem 3.27]{agen}, which gives a uniqueness result for thickenings of Lagrangian skeleta with \emph{arboreal singularities}. We also note that our Theorem \ref{theorem:hms-for-bbci} would immediately follow from Conjecture \ref{conjecture:equivalence-with-fanifold} along with the above discussion identifying the B-side spaces and HMS for fanifolds \cite[Theorem 5.4]{GS23}.

Our strategy to prove Theorem \ref{theorem:hms-for-bbci} will rely on a slightly weaker statement: even though we can not identify $Z$ with a space constructed from a fanifold, the identifications between the open cover and standard charts from Theorem \ref{theorem:combi-ci} are sufficiently robust for all the gluing arguments used in \cite{GS23} to go through with minor changes. In particular, note that the identification of open sets $\Lambda(\sigma)$ with standard charts $\eL_{\Sigma/\sigma} \times V_\sigma$ for $V_\sigma\coloneqq\relint(\partial_\infty \sigma_1)\times \dots \times \relint(\partial_\infty \sigma_r)$ implies that $\inte(\Lambda^0(\sigma))=\Lambda \cap \Phi^{-1}(M_{S^1}\times\relint(\sigma)) \cong \sigma^\perp \times V_\sigma$. We can use this along with strong exactness of $\Lambda$ to produce standard symplectic neighbourhoods:

\begin{proposition}\label{proposition:symplectic-nbhd-charts}
    For each transversal cone $\sigma$, there exists a small neighbourhood $U_\sigma \subset Z$ of $\inte(\Lambda^0(\sigma))$ such that this identification lifts to an exact symplectomorphism $\eta_\sigma\colon U_\sigma \rightarrow W_\sigma$, where $W_\sigma$ is a neighbourhood of the zero section in $T^*\sigma^\perp \times T^*V_\sigma$, which also induces a diffeomorphism $\eta_\sigma(\Lambda \cap U_\sigma) \cong (\eL_{\Sigma/\sigma} \times V_\sigma) \cap W_\sigma$.
\end{proposition}

\begin{proof}
    First, observe that we can consider the closures $\overline{V}_\sigma=\partial_\infty \sigma_1 \times \dots \times \partial_\infty \sigma_r$ of $V_\sigma$ to avoid dealing with non-compact Lagrangians (and instead have Lagrangians with corners $\Lambda^0(\sigma)$). By restricting the diffeomorphisms constructed in Theorem \ref{theorem:combi-ci} to a small neighbourhood of $\sigma$, we obtain open sets $U'_\sigma \subset Z$, $W'_\sigma \subset T^*\sigma^\perp \times T^* \overline{V}_\sigma$ and a diffeomorphism $\eta_\sigma' \colon U'_\sigma \rightarrow W'_\sigma$ that satisfies $\eta'_\sigma(\Lambda \cap U'_\sigma)=(\eL_{\Sigma/\sigma} \times \overline{V}_\sigma) \cap W'_\sigma$. Therefore, it suffices to find an exact symplectomorphism $\psi$ from a smaller neighbourhood $U_\sigma \subset U'_\sigma$ with the restriction form $\lambda_0$ to another such neighbourhood equipped with the pullback form $\lambda_1=(\eta_\sigma')^*\lambda_{std}$ from the cotangent bundle such that $\psi(\Lambda) \subset \Lambda$. 
    
    Observe that $\Lambda \cap \sigma$ is a smooth Lagrangian with corners for both $\lambda_0$ and $\lambda_1$, so the standard Weinstein neighbourhood theorem produces such a symplectomorphism (we can deal with the fact that it has boundary the same way as we did in Corollary \ref{corollary:contains}). In order to check that such a deformation constructed via Moser stability sends $\Lambda$ to itself, we can exploit the strong exactness of $\Lambda$ with respect to both symplectic forms, analogously to what we did in the proof of Theorem \ref{theorem:emb}. 
\end{proof}

\begin{proof}[Proof of Theorem \ref{theorem:hms-for-bbci}]
    In line with the general discussion above, we use the embedding $Z\times \C^{r-1} \hookrightarrow S^*M_{S^1}$ from Theorem \ref{theorem:emb} to define the cosheaf $\msh_\Lambda$ by looking at sheaves whose microsupport lies inside the thickening $\Lambda'\cong \Lambda \times \R^{r-1} \subset S^*M_{S^1}$ of $\Phi(\Lambda)$. 

    By \cite[Theorem 1.4]{gps3}, we know that $\perf \mW(Z)^{op} \cong \msh_{\Lambda}(\Lambda)$, so we prove the equivalence $\msh_{\Lambda}(\Lambda) \cong \coh(\check{Z})$ instead. 

    We can consider two functors from the poset $\Sigma_\trans^{op}$: the functor $\indcoh$ defined as $\sigma \mapsto \indcoh(\overline{O(\sigma)})$ and the functor $\msh^\diamondsuit$ defined by sending $\sigma \mapsto \msh_{\Lambda}^\diamondsuit(\Lambda^\varepsilon(\sigma))$ for the open cover $\Lambda^\varepsilon(\sigma)= \Phi^{-1}(\sigma^\perp \times \st^\varepsilon(\sigma)) \cap \Lambda$ introduced at the end of Section \ref{section:combi-ci} for some sufficiently small $\varepsilon>0$. If we prove that these two functors are quasi-isomorphic, the desired assertion follows: since $\msh_{\Lambda}^\diamondsuit$ is a sheaf of categories and the open sets $\{ \Lambda^\varepsilon(\sigma) \}_{\sigma \in \Sigma_\trans}$ are a cover, we have $\msh_{\Lambda}^\diamondsuit(\Lambda) \cong \hocolim_{\sigma \in \Sigma_\trans} \msh(\sigma)$. By an inductive application of \cite[Corollary 2.5]{nadler16} (see also \cite[Theorem 8.A.1.2]{grrrrr2}) over our indexing poset, we also know that $\indcoh(\check{Z}) \cong \hocolim _{\sigma \in \Sigma_\trans} \indcoh(\sigma)$, so by \cite[Corollary 1.2.6]{Gaitsgory2013}, it suffices to pass to compact objects on both sides to recover the original statement. 

    The results from \cite[Section 7]{GS22} on functoriality of the coherent-constructible correspondence show that the functor $\indcoh$ is isomorphic to the functor $\Sh^\diamondsuit$ given by $\sigma \mapsto \Sh_{\eL_{\Sigma/\sigma}}^\diamondsuit(\sigma^\perp)$ for a sub-torus $\sigma^\perp \subset M_{S^1}$ and $\eL_{\Sigma/\sigma} \subset \sigma^\perp \times (N/\sigma)_\R \cong T^*\sigma^\perp$ (the precise statement is proved in their Theorem 7.4.1 for the larger indexing poset $\Sigma$ that contains $\Sigma_\trans$; note that it suffices to ensure that certain diagrams commute on homology without worrying about higher coherences, since the conditions of Lemma 7.1.1 from loc. cit. are satisfied), so it suffices to describe a quasi-equivalence between $\Sh^\diamondsuit$ and $\msh^\diamondsuit$. By \cite[Lemma 7.2.2]{GS22} (and \cite[Lemma 4.33]{GS23}), the arrow $\Sh^\diamondsuit_{\eL_{\Sigma/\sigma}}(\sigma^\perp) \rightarrow \Sh^\diamondsuit_{\eL_{\Sigma/\tau}}(\tau^\perp)$ in the diagram corresponding to an inclusion $\sigma \hookrightarrow \tau$ is the composition of Sato microlocalisation along $\tau^\perp$ and restriction to $p \times \tau^\perp$ for any point $p \in \relint(\tau/\sigma)$. 

    Finally, the gluing argument from \cite[Proposition 4.34]{GS23} yields an isomorphism between the two functors $\msh^\diamondsuit$ and $\Sh^\diamondsuit$, provided we replace their Theorem 4.1 by the fact that the open cover $\{\Lambda^\varepsilon(\sigma)\}_{\sigma \in \Sigma_\trans}$ satisfies Theorem \ref{theorem:combi-ci} and Proposition \ref{proposition:symplectic-nbhd-charts} for sufficiently small $\varepsilon>0$. 
\end{proof}

\subsection{Functoriality results}\label{section:functoriality}

In this section, we explain how to use our computations of skeleta to relate different very affine complete intersections. In particular, by constructing an embedding of $\perf \mW(Z)$ into the (split closure of a) wrapped Fukaya of a certain very affine hypersurface, we explain how to intertwine the equivalence proved in Section \ref{section:proof-hms} with toric mirror symmetry. 

Recall that a \emph{partition} $\mP=(P_1,\dots,P_k)$ of a finite set $S$ is simply a decomposition into disjoint subsets $S=\bigsqcup_{i=1}^k P_i$. For two partitions, we say that $\mP'$ is a \emph{refinement} of $\mP$ if for any $P' \in \mP'$, there exists some $P \in \mP$ such that $P' \subseteq P$. We can organise all the partitions of a finite set into a poset, with $\mP' \leq \mP$ if and only if $\mP'$ is a refinement of $\mP$. Therefore, one can define an associated category $\Part(S)$ of partitions of $S$, whose objects are partitions and there is a unique morphism $\mP' \rightarrow \mP$ whenever $\mP'$ is a refinement of $\mP$. We will be particularly interested in using $\Part_r\coloneqq \Part(\{1,\dots,r\})$ for $r \in \Z_{>0}$. 

Partition of sets induce the following natural operation on nef partitions of a reflexive polytope $\Delta \subset M_\R$:

\begin{definition}\label{definition:total-refinement}
    Given a nef partition $\Delta=\Delta_1+\dots+\Delta_r$ and a partition $\mP \in \Part_r$ into $r' \leq r$ sets $P_1,\dots,P_{r'}$, we define the \emph{$\mP$-grouping of $\{\Delta_j\}_{j=1}^r$} as a length $r'$ nef partition $\Delta=\Delta^\mP_1+\dots+\Delta_{r'}^\mP$ for $\Delta^\mP_i\coloneqq\sum_{j \in P_i} \Delta_j$. We call the result of this operation on the dual nef partition $\nabla^\mP=\nabla^{\mP}_1+\dots+\nabla^{\mP}_{r'}$ the \emph{$\mP$-cogrouping of $\{ \nabla_j\}_{j=1}^r$}.
\end{definition}

In fact, one can describe the effect of cogrouping explicitly:

\begin{lemma}\label{lemma:grouped-nef-partition}
    The summands of the $\mP$-cogrouping are given by $\nabla^\mP_i\coloneqq\conv\{ \nabla_j \colon j \in P_i\}$ for $1 \leq i \leq r'$. 
\end{lemma}

\begin{proof}
    Recall that to a polytope $Q \subset M_\R$ containing the origin, we can associate its \emph{characteristic function} $\psi \colon N_\R \rightarrow \R$ given by $\psi(y)\coloneqq\sup_{x \in Q} \langle y, x \rangle$. It is then straightforward to check that for polytopes $Q$, $Q'$ with characteristic functions $\psi$, $\psi'$, the characteristic function for $Q+Q'$ is $\psi+\psi'$ and the characteristic function for $\conv\{Q,Q'\}$ is $\max\{\psi,\psi'\}$. 

    Therefore, if we denote the characteristic functions for $\Delta_j$ as $\psi_j$, then the characteristic function of $\Delta^\mP_i$ is given as $\psi^\mP_i=\max_{j \in P_i} \{\psi_j\}$. By definition of the dual nef partition, $\nabla_i^\mP$ is the convex hull of $0$ and the points $y \in \partial \Delta^\vee$ satisfying $\psi^\mP_i(y)=1$, so the formula for $\psi^\mP_i$ implies the desired statement. 
\end{proof}

The operation of cogrouping might seem less natural, but it is the more natural operation to consider when one is interested in open Batyrev--Borisov complete intersections: suppose that $Z_\beta \subset M_{\Cs}$ is a family of open BBCI's corresponding to a nef partition $\nabla=\nabla_1+\dots+\nabla_r$ and a centred refined triangulating function $h$ on $\Delta^\vee$. Moreover, suppose that the coefficients of the defining polynomials $f_{\beta,j}(z)=\sum_{\alpha \in \nabla_j} c_{\alpha,j} e^{-\beta h(\alpha)} z^\alpha$ satisfy the sign conventions $c_{\alpha,j}>0$ for all $\alpha\ne0$ and $c_{0,j}<0$ for all $1 \leq j \leq r$. For any partition $\mP$, we can then look at the partial sums $f^\mP_{\beta,i}\coloneqq\sum_{j \in P_i} f_{\beta,j}$ for $1 \leq i \leq r'$, which yields an associated one-parameter family of open BBCI's $Z_\beta(\mP)$ of codimension $r'$. Due to the sign convention, it is guaranteed that the Newton polytope of $f^\mP_i$ is the cogrouping $\nabla^\mP_i$, so $Z_\beta(\mP)$ is the family associated to the co grouped nef partition. Clearly, taking $\mP=\{1\} \sqcup \dots \sqcup \{r\}$ recovers $Z_\beta$ itself, and if $\mP'$ is a refinement $\mP$, then we get an inclusion $Z_\beta(\mP') \hookrightarrow Z_\beta(\mP)$. We have previously considered another instance of this construction: the hypersurface $H_{\beta,\tot}$ is nothing else than $Z_\beta(\{1,\dots,r\})$.  The data of such a construction can be conveniently organised as follows:

\begin{lemma}\label{lemma:functor1}
    Given a nef partition $\{\nabla_j\}_{j=1}^r$, a centred refined triangulating function $h$, a choice of coefficients $c_{\alpha,j}$ and any sufficiently large $\beta \gg0$, there is a functor $Z$ from $\Part_r$ into the category of smooth affine schemes that takes $\mP \mapsto Z_\beta(\mP)$, where $Z_\beta(\mP)$ is an open BBCI associated to the $\mP$-cogrouping $\{ \nabla_i^\mP\}_{i=1}^{r'}$ with the same data.
\end{lemma}

\begin{remark}\label{remark:cpct-case-hard}
    The situation becomes more complicated when we try to compactify: since cogrouping using different partitions $\mP$ yields different polytopes $\nabla^\mP$, the natural toric compactifications of different varieties $Z(\mP)$ will lie inside different Fano toric varieties $X_{\nabla^\mP}$ (or resolutions thereof), so there is no obvious way of \enquote{compactifying} the functor from Lemma \ref{lemma:functor1}, which is perhaps the main reason why this operation has not been previously considered in the context of Batyrev--Borisov mirror symmetry. The situation is opposite on the mirror side: grouping does not change the polytope $\Delta$, so we are working inside a fixed toric variety $X_\Delta$. However, it is less clear how one should pick the defining Laurent polynomials associated to Minkowski sums of polytopes, our most natural guess being $\check{f}_{\beta,i}^\mP\coloneqq\prod_{j \in P_i}(\check{f}_{\beta,j}-1)+1$. 
\end{remark}

We can also set things up so that the computation of the skeleta goes through for all the complete intersections $Z_\beta(\mP)$ \emph{simultaneously}. First, observe that the tropical hypersurface $\mA_{\trop,\tot}$ is unaffected by the cogrouping process, since it does not depend on the nef partition. Therefore, it suffices to describe a construction of localising tailoring functions that agree on a neighbourhood of the skeleta. Note that if we try to simply sum up tailored polynomials (like we did when we were considering $\widetilde{H}_{\beta,\tot}$ earlier), the result might have too many non-zero coefficients $\chi_\alpha$ and hence will not be localising. 

\begin{lemma}\label{lemma:simultaneous-tailoring}
    There exists a collection of proper, localising and centred cut-off functions $\chi^\mP$  for $Z_{\beta}(\mP)$ indexed by $\mP \in \Part_r$ such that we have $\chi^\mP_{\alpha,\beta}(u)=\chi^{\mP'}_{\alpha,\beta}(u)$ for all $\mP$, $\mP' \in \Part_r$, $\beta\gg0$, $\alpha \in \Delta^\vee$ and all $u$ satisfying $L_h(u) \leq K\beta^{-1}$ for some constant $K>0$.
\end{lemma}

\begin{proof}
    The functions obtained via the construction from Proposition \ref{proposition:total_tailoring} (where $K$ is as specified there) will satisfy the requirement, since by the estimates from Corollary \ref{corollary:bdary_conv}, we know that regardless of $\mP$, the resulting cut-off functions $\chi^\mP_{\alpha,\beta}$ will be equal to $\chi_{\alpha,0,\beta}(u)$ for $\alpha \neq 0$ and to $1$ for $\alpha=0$ over the designated region. 
\end{proof}

Therefore, we have all the ingredients we need for the proof from Section \ref{section:skeleta} to go through for all $\widetilde{Z}_\beta(\mP)$ at once. Since the smoothings $\partial\widetilde{C}_{0,\beta}(\mP)=\bigcap_{i=1}^{r'} \partial \widetilde{C}_{0,\beta,i}(\mP)$ of associated tropical complete intersections lie inside the cell $C_{0,\trop,\tot}$, they also lie inside the region where the conclusion of Lemma \ref{lemma:simultaneous-tailoring} holds, so by inspecting the defining equations for $\partial \widetilde{C}_{0,\beta,i}(\mP)$ and Theorem \ref{theorem:skeleton_ci}, we get the following consequence of Lemma \ref{lemma:functor1}:

\begin{corollary}\label{corollary:functor-skel}
    Given a choice of data as above, there exists a functor $\skel$ from $\Part_r$ into closed subsets of $\partial \eL_\Sigma \subset S^*M_{S_1}$ given by $\skel \colon \mP \mapsto \Phi(\skel(\widetilde{Z}_\beta(\mP)))$. 
\end{corollary}

By performing the constructions inductively over $\Part_r$, we could use a process analogous to the proof of Theorem \ref{theorem:emb} to construct a functor $LD$ taking $\mP$ to a Liouville domain $A_\beta(\mP)$ associated to $\widetilde{Z}_\beta(\mP)$ (with the shrinking from sub-level set $\{ \varphi \leq C\}$ performed so that all the inclusion maps $LD(\mP') \rightarrow LD(\mP)$ are positive codimension Liouville embeddings). After that, we could use the thickening process described there to \enquote{fatten up} $LD$ to a functor $LD'$ that takes $\mP$ to the associated stabilised Liouville sector $LD'(\mP)\coloneqq LD(\mP) \times (D^2,\pm1)^{r'-1} \subset S^*M_{S^1}$ and the maps are embeddings of exact symplectic manifolds. This also gives us a way to compatibly thicken the functor $\skel$ to $\skel'$ such that $\skel'(\mP)$ is an open subset of $\partial \eL_\Sigma$ diffeomorphic to $\skel(\mP) \times D^{r'-1}$ (and the inclusions are inclusions of open subsets), which we record as a Corollary:

\begin{corollary}\label{corollary:functor-fat-skel}
    The functor $\skel$ can be `thickened up' to a functor $\skel'$ from $\Part_r$ into open subsets of $\partial \eL_\Sigma$ given by $\skel' \colon \mP \mapsto \skel(\Phi(\widetilde{Z}_\beta(\mP))\times(\C,\pm\infty)^{r'-1})$ for suitable thickened embeddings $LD'(\mP)$ of Liouville subdomains of $\widetilde{Z}_\beta(\mP)$ constructed via Theorem \ref{theorem:emb}. 
\end{corollary}

By noting that open inclusions of singular Legendrians $\skel'(\mP)$ induce functors on the associated microlocal sheaf categories, we get the following functoriality result about wrapped Fukaya categories of $LD(\mP)$:

\begin{corollary}\label{corrolary:functor-fuk}
    There is a (strict) functor $\mW$ from $\Part_r$ into the category of $A_\infty$-categories given by $\mW \colon \mP \mapsto \perf \mW(LD(\mP))^{op}$.
\end{corollary}

\begin{proof}
    We construct the functor as $\mW\colon \mP \mapsto \msh_{\skel'(\mP)}(\skel'(\mP))$, with the maps $\msh_{\skel'(\mP')} \rightarrow \msh_{\skel'(\mP)}$ given simply by restrictions of $\msh_{\partial \eL_{\Sigma}}$ to smaller open subsets, which is equivalent to a functor given in the statement through the identification from \cite{gps3} and the Künneth theorem of \cite{gps2} identifying the partially wrapped Fukaya category of the stabilisation $\perf \mW(LD'(\mP))$ with the wrapped Fukaya category $\perf \mW(LD(\mP))$ of the domain itself.
\end{proof}

\begin{remark}\label{remark:cooked-functoriality}
    Observe that the results from \cite{gps1} about functoriality of wrapped Fukaya categories and our construction of $LD'$ do not suffice for the Corollary, since we have not shown that the maps $LD'(\mP') \rightarrow LD'(\mP)$ are inclusions of Liouville sectors, so it is unclear how to give a direct geometric description of these functors without relying on the comparison results of \cite{gps3}. 
\end{remark}

It is even more straightforward to describe what happens to the B-side mirrors $\check{Z}(\mP)$ under this process: after cogrouping, the fan $\Sigma$ does not change, but $\Sigma_\trans(\mP)$ becomes larger (since we merge some of the sets $\Sigma_j(1)$ together, so being transversal becomes an easier condition to satisfy). Therefore, we also have a functor $\check{Z} \colon \mP \mapsto \partial_{\trans,\mP} \mX_{\Sigma}$ for $\partial_{\trans,\mP}\mX_{\Sigma}$ defined by gluing toric boundary strata associated to $\mP$-transversal cones, where $\check{Z}(\mP') \hookrightarrow \check{Z}(\mP)$ are closed embeddings. This also means that we get another functor $\coh \colon \mP \mapsto \coh(\check{Z}(\mP))$ and we can upgrade our mirror symmetry result to a family statement:

\begin{proposition}\label{proposition:mirror-symmetry-functors}
    There exists an equivalence of functors $\mW \cong \coh$.
\end{proposition}
\begin{proof}[Proof (sketch)]
    We use the description of $\mW(\mP)$ via microlocal sheaves on thickened skeleta $\skel'(\mP)$, then one can repeat the gluing procedure from our proof of Theorem \ref{theorem:hms-for-bbci} to obtain equivalences $\mW(\mP) \cong \coh(\mP)$ for all partitions $\mP$. The arguments from \cite{GS23} that we rely on also make the identifications functorial (cf. their Theorem 5.2), since we also have the corresponding diagram of fanifolds obtained from $\partial\widetilde{C}_{0,\beta}(\mP)$ (as constructed in Section \ref{section:proof-hms}; it does not matter that $LD'(\mP)$ are not obtained via the handle attachment construction from loc.cit., since that can be circumvented analogously to our earlier proof). 
\end{proof}

Perhaps the most important consequence of this result is Theorem \ref{theorem:comm-diag}, which can be obtained by specialising to the arrow $\{1\} \sqcup \dots \sqcup \{r\} \rightarrow \{1,\dots,r\}$ corresponding to the embedding of $Z$ inside the hypersurface $H_\tot$. 

\begin{proof}[Proof of Theorem \ref{theorem:comm-diag}]
    From Proposition \ref{proposition:mirror-symmetry-functors}, we obtain a commutative diagram 
    \[\begin{tikzcd}
    	{\perf \mW(Z)^{op}} & {\perf \mW(H_{\tot})^{op}} \\
        {\coh(\check{Z})} & {\coh(\partial \mX_\Sigma)}
    	\arrow[from=1-1, to=1-2]
        \arrow[from=2-1, to=2-2]
    	\arrow[sloped, "\sim", from=1-1, to=2-1]
    	\arrow[sloped, "\sim", from=1-2, to=2-2]
    \end{tikzcd}\]
    where we are using the fact that $H_\tot$ is a very affine hypersurface for which the results of \cite{Zhou} apply, so $\skel(H_\tot)=\partial \eL_\Sigma$. The desired statement follows by appending this to the commutative diagram
    \[\begin{tikzcd}
    	{\perf \mW(H_{\tot})^{op}} &  {\perf \mW(T^*M_{S^1},\partial \eL_\Sigma)^{op}}\\
        {\coh(\partial \mX_\Sigma)} & {\coh(\mX_\Sigma)}
    	\arrow[from=1-1, to=1-2]
        \arrow[from=2-1, to=2-2]
    	\arrow[sloped, "\sim", from=1-1, to=2-1]
    	\arrow[sloped, "\sim", from=1-2, to=2-2]
    \end{tikzcd}\]
    from Figure 6 of \cite{GS22}.
\end{proof}

\begin{remark}\label{remark:arrow-interpretation}
    Despite Remark \ref{remark:cooked-functoriality}, we can still provide some intuition behind the functor at the top row of Theorem \ref{theorem:comm-diag}: recall that $\mW(H_{\tot})\rightarrow \mW(T^*M_{S^1},\partial\eL_\Sigma)$ is known as the \emph{cup functor} associated to a Liouville hypersurface embedding $H_{\tot} \hookrightarrow S^*M_{S^1}$. Similar considerations for $Z$ along with Theorem \ref{theorem:emb} would give us a cup functor $\mW(Z) \rightarrow \mW(T^*M_{S^1},\Lambda\times \R^{r-1})$ for a suitable open embedding $\Lambda \times \R^{r-1} \hookrightarrow \partial \eL_{\Sigma}$. By \cite[Theorem 1.20]{gps2}, there is a stop removal functor $\mW(T^*M_{S^1},\partial \eL_{\Sigma}) \rightarrow \mW(T^*M_{S^1},\Lambda \times \R^{r-1})$ that is a localisation at a collection of linking discs, so our functor $\mW(Z) \rightarrow \mW(T^*M_{S^1},\partial \eL_{\Sigma})$ can be understood as a lift of the cup functor through this localisation. At the level of microlocal sheaves, the functor corresponds to a corestriction of $\msh$ along an open inclusion of $\Lambda \times \R^{r-1}$, followed by the left adjoint to microlocalisation $\msh_{\partial \eL_\Sigma}(\partial \eL_\Sigma) \rightarrow \Sh_{M_{S^1}}(\eL_\Sigma)$. 
\end{remark}
\appendix
\section{Neighbourhoods of polyhedra} \label{section:appendix-polyhedra}
In this appendix, we study properties of polyhedra and polyhedral complexes in $\R^n$ and their neighbourhoods. We shall work in some fixed $n$-dimensional real vector space $V$ equipped with an inner product $\langle \cdot,\cdot\rangle$. Unless specified otherwise, all the distances are computed with respect to the induced norm. 

Recall that a \emph{polyhedron} is a closed subset $P\subset V$ for which there exists a finite set of affine functions $\varphi_1$, \dots, $\varphi_k$ such that $P=\{x \in V \colon \varphi_i(x) \leq 0 \textnormal{ for all } i=1,\dots,k\}$. We call any such set $\{\varphi_i\}_{i=1}^k$ a \emph{defining collection for $P$}. Clearly, the choice of such a set is not unique, and defining collections form a poset under inclusion. Therefore, we call $P$ equipped with a particular choice of a minimal defining collection a \emph{parameterised polyhedron}.

Following \cite{Mikh}, we call a closed subset $\Sigma \subset V$ a \emph{polyhedral complex} if it can be presented as a finite union of polyhedra (called \emph{cells}) so that the faces in $\partial C$ for a $k$-dimensional cell are $(k-1)$-dimensional cells of $\Sigma$ and $\relint(C) \cap \relint(C')=\emptyset$ for $C\neq C'$. We call a polyhedral complex $\Sigma$ \emph{parameterised} if each cell $C$ is parameterised by $S_C$ and there are inclusions $i_{C,C'}\colon S_{C'} \hookrightarrow S_C$ for any pair of cells $C \hookrightarrow C'$ that satisfy $i_{C,C'}\circ i_{C',C''}=i_{C,C''}$. 

\begin{example}
    The tropical hypersurface $\mA_{\trop,h}$ associated to $h$ is a polyhedral complex. It also comes with a natural parameterisation, given by the affine functionals $l_{\alpha}$ of the Legendre transform. Analogously, since an intersection of polyhedral complexes is a polyhedral complex, tropical complete intersections naturally admit a structure of a parameterised polyhedral complex. 
\end{example}

The reason for introducing parameterisations is that they are necessary to make the following definition:

\begin{definition}
    For a parameteristed polyhedron $(P,S)$, the \emph{affine distance from $P$} is given as 
    \begin{equation*}
        d_\aff(x,P)\coloneqq\max_{\varphi \in S} \varphi(x).
    \end{equation*}
    For a parameterised polyhedral complex $\Sigma$, we let $d_\aff(x,\Sigma)\coloneqq\min_C d_\aff(x,C)$, where $C$ runs over polyhedral cells of $\Sigma$.
\end{definition}

We record a proof of the standard fact that affine distance is Lipschitz-equivalent to the distance function $d(x,P)\coloneqq\inf_{y \in P}\lVert x-y\rVert$ induced by the norm:

\begin{lemma}\label{lemma:aff-vs-norm}
    For any parameterised polyhedron $(P,S)$, there exists a constant $K>0$ such that the inequalities $K^{-1}\cdot d(x,P) \leq d_\aff(x,P) \leq K\cdot d(x,P)$ hold for all $x \in V$. 
\end{lemma}

\begin{proof}
    The upper bound for $d_\aff(x,P)$ is straightforward: pick a constant $K$ such that $\lVert \varphi \rVert <K$ holds for all $\varphi \in S$ (a norm of $V$ induces an operator norm on linear functionals on $V$, hence also on affine linear functionals), then the existence of a point $y \in P$ satisfying $\Vert x-y \rVert \leq \delta$ implies $\varphi(x)<K\delta$ for all $\varphi \in S$ and $\delta>0$, as desired. 

    For the lower bound, let us denote the orthogonal projection of $x$ onto $P$ as $\pi(x)$. Since $P$ is a polyhedron, $\pi(x)$ is equal to the orthogonal projection onto the smallest affine space containing some face $F$ of $P$ and we have $d(x,P)=\lVert v \rVert$ for $v=x-\pi(x)$ by definition. 
    
    Let $\varphi_1$, \dots, $\varphi_l$ be the functions in $S$ that are identically equal to $0$ on $F$, but do not vanish identically on $P$, and denote the functions vanishing on the entire polyhedron $P$ by $\varphi_{l+1}$, \dots, $\varphi_k$. By the Riesz representation theorem, the functions $\varphi_i(x)$ can then uniquely be written in the form $\varphi_i(x)=\langle x, w_i \rangle+c_i$ for $w_i \in V$, $c_i \in \R$. The (exterior) normal cone of $P$ at any $x \in \relint(F)$ is then given by 

    \begin{equation*}
        \nc_P(F)=\R_{\geq 0}w_1 \oplus \dots \oplus \R_{\geq 0}w_l\oplus \R w_{l+1} \oplus \dots \oplus \R w_k,
    \end{equation*}
    where directness of the sum (i.e. linear independence of $w_i$'s) follows by minimality of our defining set $S$. 

    Write $v=x-\pi(x)$, because $\pi(x)$ is the orthogonal projection of $x$ onto $F$, $v$ must lie in $\nc_P(F)$. By swapping the signs of $w_{l+1}$, \dots, $w_k$ if necessary (and splitting the problem into $2^{k-l}$ cases that way), we can assume that $v \in \bigoplus_{i=1}^k \R_{\geq 0} w_i$ holds without loss of generality. By definition of affine distance, we then have $\langle w_i, v \rangle \leq d_{\aff}(x,P)$ for all $i=1,\dots,k$. 

    Let $\Delta\coloneqq\conv\{w_1,\dots,w_k\}$ be the simplex spanned by $w_i$'s, since it does not contain the origin, the distance of $\Delta$ from the origin must be equal to $c_F>0$. Since $v \in \bigoplus_{i=1}^k \R_{\geq 0} w_i = \cone(\Delta)$, there exists some positive rescaling $\widetilde{v} \in \Delta$ of $v$. By linearity of inequalities, $\langle w_i, v \rangle \leq d_\aff(x,P)$ forces $\langle w, v \rangle \leq d_\aff(x,P)$ for all $w \in \Delta$. Therefore, since $\widetilde{v} \in \Delta$ and hence $\lVert \widetilde{v} \rVert \geq c_F$, we also get
    \begin{equation*}
        c_F \cdot d(x,P)=c_F\cdot\lVert v\rVert\leq\lVert \widetilde{v} \rVert \cdot \lVert v \rVert = \langle \widetilde{v},v \rangle \leq d_\aff(x,P), 
    \end{equation*} 
    which means that any constant satisfying $K>c_F^{-1}$ will work in the inequality. Since $P$ has finitely many faces, we can just take the maximum of $c_F^{-1}$ over all of them to get the desired result.
\end{proof}

\begin{corollary}\label{corollary:aff-vs-norm}
    For any parameterised polyhedral complex $\Sigma$, there exist constants $K>0$ and $\delta>0$ such that we have $K^{-1}\cdot d_\aff(x,\Sigma) \leq d(x,\Sigma) \leq K\cdot d_\aff(x,\Sigma)$ for all $x \in V$ satisfying $d(x,\Sigma)<\delta$.
\end{corollary}

\begin{proof}
    Analogously to the Lemma, the upper bound is an immediate consequence of the defining functions being Lipschitz and does not require any restrictions on $\delta$.
    
    For the lower bound: pick a constant $K>0$ such that the conclusion of Lemma \ref{lemma:aff-vs-norm} holds for all the cells $C \in \Sigma$. Suppose that $\delta>0$ is small enough so that any two disjoint cells of $\Sigma$ have distance at least $2\delta$. Then for $d=d(x,\Sigma)<\delta$, there will exist a unique minimal cell $C_0$ of $\Sigma$ that intersects the closed ball $B_{d}(x)$, which will tautologically satisfy $d=d(x,C_0)$. This means that $d=d(x,C')$ for all cells $C'$ containing $C_0$, so by our choice of $K$, such cells will satisfy $d \leq K\cdot d_\aff(x,C')$. At the same time, we must have $d_\aff(x,C'')>K^{-1}\cdot d(x,C'')\geq K^{-1}\cdot d(x,\Sigma)= K^{-1}\cdot d$ for any cell $C''$ that does not contain $C_0$ by the choice of $\delta$. Therefore, if we take the minimum over all cells $C\in \Sigma$, we still get $d \leq K \cdot d_\aff(x,\Sigma)$, as desired. 
\end{proof}

In particular, this allows us to study two kinds of open neighbourhoods of polyhedral complexes: for $\delta>0$ and a parameterised polyhedral complex $\Sigma$, we can consider its ordinary $\delta$-neighbourhood $U(\Sigma,\delta)$ with respect to the norm distance and also its affine $\delta$-neighbourhood $U_\aff(\Sigma,\delta)$. Since we have shown that the two notions of distance are Lipschitz equivalent, we can use one type of neighbourhood to study the other, which is illustrated in the following results:

\begin{lemma}\label{lemma:nbhd-intersections}
    For any polyhedra $P_1$, $P_2$, there exists a constant $K>0$ such that $U(P_1,\delta) \cap U(P_2,\delta) \subseteq U(P_1 \cap P_2,K\delta)$ holds for all $\delta>0$.
\end{lemma}
\begin{proof}
    First, we observe that if $S_1$, $S_2$ are parameterisations for $P_1$ and $P_2$, there exists a parameterisation $S \subseteq S_1\cup S_2$ for $P_1 \cap P_2$. With respect to these, the inequality $d_\aff(x,P_1\cap P_2) \leq \min\{ d_\aff(x,P_1),d_\aff(x,P_2)\}$ clearly holds, so if we consider such parameterisations (which always exist by the discussion above) and replace $U$ by $U_\aff$ in the statement of the Lemma, it will hold for $K=1$. Therefore, it suffices to apply Lemma \ref{lemma:aff-vs-norm} to switch from one Lipschitz equivalent distance to another at the cost of a constant. 
\end{proof}

\begin{corollary}\label{corollary:nbhd-intersections}
    For any polyhedral complexes $\Sigma_1$, $\Sigma_2$, there exists a constant $K>0$ such that $U(\Sigma_1,\delta) \cap U(\Sigma_2,\delta) \subseteq U(\Sigma_1 \cap \Sigma_2,K\delta)$ holds for all sufficiently small $\delta>0$
\end{corollary}
\begin{proof}
    Same as Lemma \ref{lemma:nbhd-intersections}, but with Lemma \ref{lemma:aff-vs-norm} replaced by Corollary \ref{corollary:aff-vs-norm}. 
\end{proof}

Finally, we can also use these results to describe nice stratifications of small neighbourhoods of polyhedral complexes, along with projections defined over the neighbourhood (analogous to the data associated to a general Thom--Mather stratified space, but satisfying some nice explicit bounds that we shall make use of in the body of the paper):

\begin{lemma}\label{lemma:aff-nbhd-structure}
    For a parameterised polyhedral complex $\Sigma$ and all sufficiently small $\delta>0$, the closed neighbourhood $\overline{U_\aff(\Sigma,\delta)}$ admits a decomposition into $n$-dimensional polyhedral cells $C(\delta)$ with disjoint interiors indexed by $C \in \Sigma$, where $x \in C(\delta)$ if and only if $C$ is the minimal cell of $\Sigma$ satisfying $d_\aff(x,C)\leq \delta$.
\end{lemma}

\begin{proof}
    For a point $x \in \overline{U_\aff(\Sigma,\delta)}$, we can determine which $C(\delta)$ it belongs to for all $\delta$ small enough through the process that we used in the proof of Corollary \ref{corollary:aff-vs-norm}. Since $C(\delta)$ is cut out by linear inequalities, it is, by definition, a polyhedron. To see that it has non-empty interior, we can note that a small ball around a point in $\relint(C)$ will be contained inside $C(\delta)$ for sufficiently small $\delta$. 
\end{proof}

\begin{corollary}\label{corollary:aff-nbhd-structure}
    In the setting of Lemma \ref{lemma:aff-nbhd-structure}, there exists a constant $K>0$ such that whenever $x \in C(\delta)$ for a cell $C \in \Sigma$ and a sufficiently small $\delta>0$, there exists a point $x' \in C$ that satisfies $d_\aff(x',\partial C) \geq \delta$ and $d(x,x') \leq K\delta$.
\end{corollary}

\begin{proof}
    Let $C_\delta$ be the set of points in $C$ that have affine distance at least $\delta$ from $\partial C$, then it inherits a structure of a parameterised polyhedron from that of $C$. Moreover, it is clear that the points in $C(\delta)$ have affine distance at most $\delta$ from $C_\delta$. Therefore, it suffices to take $x'$ to be the orthogonal projection of $x$ onto $C_\delta$, this will have the first required property by design and the second one follows from Lemma \ref{lemma:aff-vs-norm}. 
\end{proof}
\section{Smoothing submanifolds of \texorpdfstring{$\R^n$}{TEXT}}
\label{section:appendix-smoothing}
In this appendix, we introduce some smoothing procedures for topological submanifolds of $\R^n$ that look like $\Sigma^\ba_\trans$ from Section \ref{section:combinatorics}, which is necessary as an intermediate step for constructing isotopies in Section \ref{section:potentials} from specific types of piecewise smooth isotopies considered there. We also use this framework to construct explicit smoothings of various piecewise linear objects that appear throughout the paper, such as $\partial C_{0,\trop}$ (in Section \ref{section:smoothing-applications}), or various identifications and scaling actions appearing in Section \ref{section:combinatorics} (Section \ref{section:ooga-booga-smoothing}).

Recall that if $Y_1$, \dots, $Y_r$ are manifolds with corners in $\R^n$, we say that they intersect transversely if for all open strata $S_j \subset Y_j$, the intersection $\bigcap_{j=1}^rS_j$ is a transverse intersection of smooth submanifolds of $\R^n$. We will be interested in subsets of $\R^n$ that are covered by the following definition: 

\begin{definition}\label{definition:int_strat_space}
    Suppose that $Y_1$, \dots, $Y_r$ are $n$-dimensional manifolds with corners in $\R^n$ that intersect transversely, then we call $X\coloneqq \bigcap_{j=1}^r\partial Y_j$ a \emph{transverse cornered complete intersection (TCCI)} of codimension $r$.  
\end{definition}

Observe that we can use the stratification of each $\partial Y_j$ into open lower dimensional manifolds to endow any TCCI with the structure of a smooth Whitney stratified space inside $\R^n$. 

The main cases of interest for us are going to be the following two examples:

\begin{example}
    For a full-dimensional manifold with corners $Y \subset \R^n$, its boundary $X=\partial Y$ is a codimension $1$ TCCI.
\end{example}
\begin{example}
    The complex of bounded cells in a codimension $r$ Batyrev--Borisov TTCI $\partial C_{0,\trop}=\bigcap_{j=1}^r \partial C_{0,\trop,j}$ from Section \ref{section:skeleta} is a compact codimension $r$ TCCI. 
\end{example}

We can observe that if $Y$ is a full-dimensional manifold with corners in $\R^n$, $X=\partial Y$ and $p \in X$ lies inside a codimension $m \geq 1$ stratum of $Y$, then picking standard boundary defining coordinates for $Y$ tells us that for a small neighbourhood $p\in U \subset \R^n$, the pair $(U,U\cap X)$ will be diffeomorphic to $(\R^n,C_m\times \R^{n-m})$, where $C_m\coloneqq\partial \R^m_{\geq0}$ and we are using the standard linear splitting $\R^n \cong \R^m \times \R^{n-m}$ to embed $\R^m$ as $\R^m \times \{0\}$. Using transversality, this approach generalises to give us the following local model for TCCI's:

\begin{lemma}\label{lemma:local_form_TCCI}
    Let $X=\bigcap_{j=1}^r \partial Y_j$ be a TCCI and let $p \in X$. Then there exists an open neighbourhood $p\in U \subset \R^n$ such that the pair $(U,U\cap X)$ is diffeomorphic to $(\R^n,C_{m_1}\times \dots \times C_{m_r}\times\R^{n-m})$, where $m_j$ are positive integers whose sum is $m$.
\end{lemma}

We introduce the following notation for these model spaces: for a vector $\mathbf{m}=(m_1, \dots, m_r)$ of positive integers that satisfy $m=\sum_j m_j \leq n$, we shall write $C_{n,\mathbf{m}}\coloneqq C_{m_1}\times C_{m_2} \times \dots \times C_{m_r} \times \R^{n-m}$. The existence of such a local model also immediately implies the following fact: 

\begin{corollary}
    Any codimension $r$ TCCI in $\R^n$ is a topological manifold of dimension $n-r$.
\end{corollary}

Therefore, the constructions that we explicitly describe fit into the classical framework of questions about smoothings of topological submanifolds of smooth manifolds. We follow the approach of \cite{Whitehead1961}; our new contribution is extending it to the setting where one needs to keep track of a polyhedral stratification of $\R^n$ and establishing certain existence and uniqueness results that hold for the class of TCCI's, but not in general.

\subsection{Smoothing TCCI's}\label{section:vanilla-smoothing} We begin by observation that $C_m=\partial\R^m_{\geq0} \subset \R^m$ can be identified with $\R^{m-1}$ by projecting onto $\R^m/\R v$ for any vector $v \in \inte(C_m)$. Therefore, a choice of such a vector $v$ naturally induces a smooth structure on $C_m$. Picking a normal vector also yields a tubular neighbourhood of $C_m$ inside $\R^m$, so we can construct an ambient smoothing of $C_m$ as the image of a sufficiently small (in $C^0$-sense) smooth section of this neighbourhood. Since the space of choices for $v$ is a contractible cone and the space of sections is convex, this construction turns out to be independent of choices. Clearly, one can also perform a similar procedure for a product of corners that serves as a local model for TCCI's by picking a tuple of vectors. Our strategy for a general TCCI will be patching such local constructions appropriately, which is possible because of convexity of certain choice spaces. 

First, we describe the normal geometry of TCCI's: let $X$ be a full-dimensional manifold with corners in $\R^n$ and suppose that $x \in \partial Y$ is in a codimension $m \geq 1$ stratum, then we can pick a chart $\phi$ identifying a pair $(U, U \cap Y)$ with $(\R^n,C_m \times \R^{n-m})$. Following  \cite[Definition 2.2]{joyce}, we can define\footnote{The minus sign in the definition is in order to match our convention from polyhedral geometry that we are considering outward pointing cones.} the \emph{normal cone of $Y$} at $x$ as $NC_xY\coloneqq -D\phi_0(C_m \times \R^{n-m}) \subset T_x\R^n$ (as remarked in the reference, this is independent of our choice of $\phi$). Analogously, define the \emph{conormal cone of $Y$} at $x$ as the dual cone $N^*C_xY \coloneqq (NC_xY)^\vee \subset T_x^*\R^n$. 

\begin{definition}\label{definition:positively-transverse-tuple}
    For a TCCI $X= \bigcap_{j=1}^r \partial Y_j$ and $x \in X$, we say that an $r$-tuple of vectors $\mathbf{v}=(v_1, \dots,v_r)$ in $T_x\R^n \cong \R^n$ is \emph{positively transverse to $X$ at $x$} if the matrix $M$ with entries $M_{ij}=\langle \alpha_i,v_j \rangle$ is positive definite for any choice of non-zero linear functionals $\alpha_i \in N^*C_x Y_i$. We denote the space of positively transverse tuples at $x$ as $\mG^+_x \subset (T_x\R^n)^r$.
\end{definition}

\begin{example}
    When $r=1$ and $X=\partial Y$, it follows from definition of the dual cone that $\mG^+_x$ is simply the interior of $NC_xY$.
\end{example}

Since the defining condition for $\mG^+_x$ is clearly preserved under taking convex combinations and one can construct a positively transverse tuple by working in local coordinates, it follows that $\mG^+_x$ is non-empty and convex. 

\begin{definition}\label{definition:neg-trans}
    For a TCCI $X=\bigcap_{j=1}^r \partial Y_j$, a \emph{positively transverse tuple} is a tuple of smooth functions $\mathbf{v}=(v_1,\dots,v_r)\colon X \rightarrow (\R^n)^r$ such that $\mathbf{v}(x) \in \mG^+_x$ for all $x \in X$ (we write the set of positively transverse tuples as $\Gamma(X,\mG^+)$). 
\end{definition}

\begin{lemma}\label{lemma:transverse-smooth-guy-exists}
    Let $X=\bigcap_{j=1}^r \partial Y_j$ be a compact TCCI, then $\Gamma(X,\mG^+)$ is non-empty and convex.
\end{lemma}

\begin{proof}
    For convexity, note that all the $\mathbf{v}$ are sections of a trivial rank-$nr$ vector bundle on $X$, so it makes sense to average them as such. Moreover, it was remarked above that $\mG_x^+$ are convex, so averaging two sections that have image inside $\mG^+$ will yield another element of $\Gamma(X,\mG^+)$. 

    Non-emptiness of $\Gamma(X,\mG^+)$ can be shown similarly, by appealing to non-emptiness of all the local spaces $\mG_x^+$ and a patching argument: by compactness of $X$, we can cover it by finitely many standard coordinate charts $U_1$, \dots, $U_N$ (with $U_i$ an open subset of $\R^n$ and $U_i \cap X$ a finite open cover of $X$) and pick a partition of unity $\{ \lambda_i\}_{i=1}^N$ subordinate to the cover. We can construct a positively transverse tuple $\mathbf{v}_i$ over $U_i$ as follows: for a local model $C_{n,\mathbf{m}}$, we pick a tuple $\mathbf{\widetilde{v}}$ that is positively transverse at the origin, it is then straightforward to check that the constant function with value $\mathbf{\widetilde{v}}$ is positively transverse to $C_{n,\mathbf{m}}$ at all points. After that, we can push this forward through a chart $\phi_i$ to obtain a positively transverse plane field $\mathbf{v}_i$ over $U_i$. Moreover, by using the convexity of all the individual spaces, one can make sense of the expression $\mathbf{v}(x)\coloneqq\sum_{i=1}^N \lambda_i(x)\mathbf{v}_i(x) \in (T_x\R^n)^r$ for all $x \in X$ and, moreover, it defines a positively transverse plane field.
\end{proof}

The idea behind this definition is that for any $\mathbf{v} \in \mG^+_x$, we will consider the space $V \coloneqq \spann(v_1,\dots,v_r) \subset T_x \R^n \cong \R^n$ and show that one can use the quotient map $q_V \colon \R^n \rightarrow \R^n/V$ as a local chart for $X$ near $x$ (see Lemma \ref{lemma:neg-transverse-means-transvers} for a precise statement). Before doing that, we introduce some more notation to relate the construction to the smoothing technique from \cite{Whitehead1961}. 

For a local model $C_{n,\mathbf{m}}$ from Lemma \ref{lemma:local_form_TCCI}, we can consider the following subset of the Grassmannian:
\begin{equation*}
    \mG_{n,\mathbf{m}}\coloneqq\{ V \in \Gr(n,r) \colon \textnormal{the restriction of } q_V\colon \R^n \rightarrow \R^n/V \textnormal{ to } C_{n,\mathbf{m}}\textnormal{ is a homeomorphism}\}.
\end{equation*}

Analogously, for a TCCI $X=\bigcap_{j=1}^r \partial Y_j$ and a point $x \in X$, we can pick a local chart $\phi$ as in Lemma \ref{lemma:local_form_TCCI} that identifies a neighbourhood of $x$ with the neighbourhood of the origin in $C_{n,\mathbf{m}}$. After identifying $T_0\R^n \cong \R^n$, this allows us to define the \emph{space of transverse planes at $x$} as
\begin{equation*}
    \mG_x\coloneqq(D_x\phi)^{+1}(\mG_{n,\mathbf{m}}) \subset \Gr(T_x\R^n,r).
\end{equation*}
Note that this is independent of the choice of local coordinates: if $\phi'$ were another chart identifying a neighbourhood of $x$ with the origin in $C_{n,\mathbf{m}'}$, then the vectors $\mathbf{m}$, $\mathbf{m}$ have to be the same up to rearranging their components, so by postcomposing $\phi'$ with a linear map permuting coordinates, we can assume that $\mathbf{m}=\mathbf{m'}$. The transition function $\psi=\phi' \circ \phi^{-1}$ is a diffeomorphism that maps the germ of $C_{n,\mathbf{m}}$ at the origin to itself, which means that $D_0\psi$ is a linear isomorphism that sends $C_{n,\mathbf{m}}$ to itself. In particular, this means that $D_0\psi$ induces an automorphism of $\mG_{n,\mathbf{m}}$, which implies that that the space $\mG_x$ is indeed well-defined.

\begin{remark}\label{remark:def-whitehead-equivalent}
    The space $\mG_x$ introduced above agrees with Whitehead's notion of \emph{transverse $r$-planes at $x$} introduced in \cite{Whitehead1961}. Unlike the general case of topological submanifolds of $\R^n$ that is discussed there, we can explicitly understand the local geometry for TCCI's, which makes the definition through local coordinates easier to work with. 
\end{remark}

With this in mind, we can make sense of what precisely we meant earlier by an appropriate choice of directions in which we can project (which is a slight variation on Whitehead's original definition):

\begin{definition}\label{definition:transverse-field}
    Let $X=\bigcap_{j=1}^r \partial Y_j$ be a TCCI, then we say that a smooth function $\nu\colon X \rightarrow \Gr(n,r)$ is a \emph{pre-smooth field of transverse planes on $X$} if $\nu(x) \in \mG_x$ for all $x \in X$ (we write $\nu \in \Gamma(X,\mG)$ for such plane fields).
\end{definition}

We now explain how our notion of positively transverse tuples fits into the picture:

\begin{lemma}\label{lemma:neg-transverse-means-transvers}
    For any TCCI $X$, there exists a canonical smooth map $\Gamma(X,\mG^+) \rightarrow \Gamma(X,\mG)$ given by $\mathbf{v} \mapsto \spann_\R(v_1,\dots,v_r)$.
\end{lemma}

\begin{proof}
    The statement is clearly local in nature, so it suffices to prove it in the case $X=C_{n,\mathbf{m}}$ at the point $x=0$. To explicitly specify the presentation of $X$ as a TCCI, we pick a dual basis for $(\R^n)^\vee$ consisting of vectors $\eta^0_1$, \dots, $\eta^0_{n-m}$ and $\eta^j_i$ for $1 \leq j \leq k$, $1 \leq i \leq m_j$ such that $Y_j \coloneqq \{x \in \R^n \colon \eta^i_j(x) \leq 0 \textnormal { for } 1 \leq i \leq m_j \} \cong \R_{\geq 0}^{m_j} \times \R^{n-m_j}$ are manifolds with corners satisfying $C_{n,\mathbf{m}}= \bigcap_{j=1}^r \partial Y_j$. Their conormal cones at the origin are then given by $N^*C_0Y_j=\cone(\eta_1^j,\dots,\eta^j_{m_j})$ for all $1 \leq j \leq r$. 
    
    Suppose that $\mathbf{v}$ is a positively transverse tuple at the origin, then $V=\spann_\R(v_1,\dots,v_r)$ must have dimension $r$, since the matrix $M$ from the definition being invertible tells us that there exists a surjective linear map from $V$ to $\R^r$. Therefore, we have a well-defined map $\mG^+_0\rightarrow \Gr(n,r)$. Note that if we prove that $q=q_V|_{C_{n,\mathbf{m}}}$ is injective then it has to be a homeomorphism: $C_{n,\mathbf{m}}$ and $\R^n/V$ are topological manifolds of the same dimension, so if $q$ is injective, it must be a homeomorphism onto some neighbourhood $0 \in U \subset \R^n/V$ by invariance of domain, and homogeneity of $q$ (i.e. $q(\lambda x)=\lambda q(x)$ for $x \in C_{n,\mathbf{m}}$, $\lambda \geq 0$) forces $U=\R^n/V$.

    Therefore, we just need to show that the map $q\colon C_{n,\mathbf{m}} \rightarrow V$ is injective. For the sake of contradiction, suppose that we have $q(x)=q(x')$ for $x \neq x'$, then we can write $x-x'=\sum_{j=1}^r t_j v_j$ for some real numbers $t_j$ that are not all zero. Since $\mathbf{v}$ is positively transverse, we must have $(\sum_{j=1}^r t_j\alpha_j)(x-x')>0$ for all non-zero linear functionals $\alpha_j \in N^*C_0Y_j=\cone(\eta^j_1,\dots,\eta^j_{m_j})$. Since both $x$ and $x'$ lie in $C_{n,\mathbf{m}}$, for any index $1 \leq j \leq r$, there must exist some $1 \leq i_j, i'_j \leq m_j$ such that $\eta^j_{i_j}(x)=0$ and $\eta^j_{i'_j}(x')=0$. In particular, this means that we have $\eta^j_{i_j}(x-x')=-\eta^j_{i_j}(x') \geq 0$ and $0 \geq \eta^j_{i'_j}(x)=\eta^j_{i'_j}(x-x')$. Therefore, by picking $\alpha_j=\eta^j_{i_j}$ when $t_j\geq 0$ and $\alpha_j=\eta^j_{i'_j}$ when $t_j<0$, we can achieve that $t_j \alpha_j(x-x') \leq 0$ holds for all $j$, which contradicts the inequality coming from positive transversality, so the map $q$ is injective, as desired. 
\end{proof}

Now that we understand how to use the data of a positively transverse tuple to produce a transverse plane field, we review the role of the latter in smoothing constructions. Recall that the Grassmannian carries a tautological vector bundle $\gamma_{n-r,r}\rightarrow\Gr(n,r)$, so given a pre-smooth field on $X$ we can define a smooth vector bundle $\mV_\nu\coloneqq\nu^*\gamma_{n-r,r}$ on $X$ (to be fully precise, there is a smooth vector bundle on some open neighbourhood where $\nu$ is defined and we look at its germ over $X$). 

Denote the inclusion of $X$ into $\R^n$ as $\iota\colon X \xhookrightarrow{} \R^n$, then we can use the fact that $\gamma_{n-r,r}$ is a subbundle of the trivial bundle $\underline{\R}^n \rightarrow \Gr(n,r)$ to get a map $\widetilde{\iota} \colon \Tot(\mV_\nu) \rightarrow \R^n$ by viewing $\Tot(\mV_\nu)$ as a subset of $X \times \Gr(n,r) \times \R^n$ and taking $\widetilde{\iota}\colon(x,V,v)\mapsto\iota(x)+v$. Note that a choice of an inner product on $\R^n$ gives us a metric on $\mV_\nu$, so we can consider $\mV^\varepsilon_\nu$, the disc bundle of radius $\varepsilon>0$ inside $\mV_\nu$. 

\begin{lemma}\label{lemma:TCCI_tub_nbhd}
    Let $X=\bigcap_{j=1}^r \partial Y_j$ be a compact codimension $r$ TCCI, then for $\varepsilon>0$ small enough, the restriction of $\widetilde{\iota}$ to $\Tot(\mV^{\varepsilon}_\nu)$ will be a homeomorphism onto its image. 
\end{lemma}

\begin{proof}
    Apply \cite[Theorem 1.5]{Whitehead1961} along with compactness of $X$. 
\end{proof}

Following the standard terminology, we call the image of this map $\mN^\varepsilon_\nu\coloneqq\widetilde{\iota}(\Tot(\mV^\varepsilon_\nu))$ for $\varepsilon>0$ small a \emph{$\nu$-neighbourhood of $X$}. Pushing forward the projection $\mV \rightarrow X$ gives us the \emph{$\nu$-projection} $\pi\colon\mN^\varepsilon_\nu \rightarrow X$. With this, we can upgrade the notion of pre-smooth transverse plane field (following \cite{Whitehead1961}; note that $\mN^\varepsilon_\nu$ naturally inherits a smooth structure from $\R^n$ as its open subset):

\begin{definition}\label{definition:sm-transverse-field}
    We say that a pre-smooth transverse plane field $\nu \in \Gamma(X,\mG)$ is \emph{smooth} if there exists a $\nu$-neighbourhood $\mN^\varepsilon_\nu$ such that the composition $\nu \circ \pi \colon \mN^\varepsilon_\nu \rightarrow \Gr(n,r)$ is smooth (or, equivalently, if $\nu$ defines a smooth foliation of $\mN^\varepsilon_\nu$ by affine planes after shrinking $\varepsilon$). 
\end{definition}

It is explained in the proof of \cite[Theorem 1.10]{Whitehead1961} how to construct smoothings of transverse plane fields that are Lipschitz, which is always the case for what we call pre-smooth plane fields. In particular, we get the following consequence of the Theorem that is natural to state in our setting:

\begin{corollary}\label{corollary:smoothing-plane-fields}
    Given a pre-smooth transverse plane field $\nu \in \Gamma(X,\mG)$ along a compact TCCI $X$ and $\varepsilon>0$, there exists a smooth transverse plane field $\nu' \in \Gamma(X,\mG)$ that is $\varepsilon$-close to $\nu$. 
\end{corollary}

Suppose that $\nu$ is a smooth transverse plane field along $X$, then for each $x \in X$, we can pick a small neighbourhood $x \in U_x \subset X$ such that the map $q_{\nu(x)} \colon \R^n \rightarrow \R^n/\nu(x)$ restricts to a homeomorphism from $U_x$ to some neighbourhood $V_x$ of the origin (alternatively, we could pick an inner product on $\R^n$ and consider orthogonal projections onto $\nu(x)^\perp \cong \R^n/\nu(x)$ like in the original reference). Smoothness of $\nu$ then guarantees that such charts assemble into a smooth structure on $X$:

\begin{theorem}[{\cite[Theorem 1.7]{Whitehead1961}}]
    The collection $\{(U_x,q_{\nu(x)})\}_{x \in X}$ described above defines a smooth atlas on $X$ such that $\pi \colon \mN^\varepsilon_\nu \rightarrow X$ becomes a smooth fibre bundle. 
\end{theorem}

Finally, this allows us to define what we mean by smoothing a TCCI:

\begin{definition}\label{definition:TCCI-smoothing}
    Let $\nu$ be a smooth transverse plane field along $X$, then a \emph{$\nu$-smoothing of $X$} is the image $X_\nu\coloneqq\im(s)$ of a smooth section $s\colon X \rightarrow \mN^\varepsilon_\nu$ of the fibre bundle $\mN^\varepsilon_\nu \xrightarrow{\pi} X$ for a sufficiently small $\varepsilon>0$. 
\end{definition}

\begin{theorem}[{\cite[Theorem 1.9]{Whitehead1961}}]\label{theorem:f-smoothing-exists}
    Given a smooth transverse plane field $\nu$ along $X$, a $\nu$-smoothing $X_\nu$ exists. 
\end{theorem}

The discussion preceding \cite[Theorem 1.9]{Whitehead1961} implies that the image of such a map is a smooth, embedded submanifold of $\R^n$. By construction of $\mN^\varepsilon_\nu$ from a vector bundle $\mV^\varepsilon_\nu$, it makes sense to take convex combinations of two sections, which tells us that the space of sections $\Gamma(X,\mN^\varepsilon_\nu)$ from Definition \ref{definition:TCCI-smoothing} is convex, hence it is contractible and any two $\nu$-smoothings are ambient isotopic in $\R^n$, showing that this notion does not depend on our choice of $s$ or the parameter $\varepsilon$, so it makes sense to refer to $X_\nu$ as \emph{the} $\nu$-smoothing of $X$.

\begin{corollary}\label{corollary:TCCI-smoothing-exists}
    Any compact TCCI in $\R^n$ admits a smoothing.
\end{corollary}
\begin{proof}
    By Lemma \ref{lemma:transverse-smooth-guy-exists}, there exists a positively transverse tuple $\mathbf{v} \in \Gamma(X,\mG^+)$, which induces a pre-smooth transverse plane field $\nu \in \Gamma(X,\mG)$ by Lemma \ref{lemma:neg-transverse-means-transvers}. This can be approximated by a smooth transverse plane field $\nu'$ by Corollary \ref{corollary:smoothing-plane-fields}, which provides the data sufficient to construct a smoothing of our TCCI by Theorem \ref{theorem:f-smoothing-exists}.
\end{proof}

We will now show that the smoothings defined this way are unique in a suitable sense. Recall that two (pre)-smooth plane fields $\nu_0$ and $\nu_1$ are \emph{transversally homotopic} if there exists a smooth one-parameter family $(\nu_t)_{0 \leq t \leq 1}$ of (pre)-smooth plane fields along $X$ interpolating between them. The geometric meaning of the notion is captured by the following result:

\begin{lemma}\label{lemma:htpic-implies-isotopic}
    If $\nu_0$ and $\nu_1$ are smoothly transversally homotopic transverse plane fields along $X$, then the $\nu_0$-smoothing $X_0$ is ambient isotopic to the $\nu_1$-smoothing $X_1$. 
\end{lemma}

\begin{proof}
    This follows from the proof of \cite[Theorem 1.11]{Whitehead1961}. 
\end{proof}

First, we can observe that our construction does not depend on the choice of smoothing in Corollary \ref{corollary:smoothing-plane-fields}, so the only input of our construction is the positively transverse tuple $\mathbf{v}$.

\begin{lemma}\label{lemma:TCCI-smoothing-unique}
    Let $\mathbf{v} \in \Gamma(X,\mG^+)$ be a positively transverse tuple along a compact TCCI $X$ with an associated transverse plane field $\nu \in \Gamma(X,\mG)$. Then for any two $\varepsilon$-smoothings $\nu'$ and $\nu''$ of $\nu$, the associated smoothings $X_{\nu'}$ and $X_{\nu''}$ are ambient isotopic. 
\end{lemma}

\begin{proof}
    By \cite[Theorem 1.3]{Whitehead1961}, $\nu'$ and $\nu''$ are transversally homotopic as pre-smooth plane fields. It is explained in Section 12 of loc. cit. that pre-smooth homotopies between smooth plane fields can be approximated by smooth homotopies relative to the endpoints (essentially, this is an application of the relative version of Theorem 1.10 from loc. cit.), so they are also transversally homotopic as smooth plane fields, and we are done by Lemma \ref{lemma:htpic-implies-isotopic}. 
\end{proof}

Therefore, it makes sense to speak of a smoothing $X_{\mathbf{v}}$ associated to a choice of a positive transverse tuple $\mathbf{v}$. In fact, we can use the same technique along with convexity of $\Gamma(X,\mG^+)$ proved in Lemma \ref{lemma:transverse-smooth-guy-exists} to show that this choice also does not matter.

\begin{proposition}\label{proposition:TCCI-smoothing-unique}
    For any positively transverse tuples $\mathbf{v}_0, \mathbf{v}_1 \in \Gamma(X,\mG^+)$, the associated smoothings $X_{\mathbf{v}_0}$ and $X_{\mathbf{v}_1}$ are ambient isotopic.
\end{proposition}
\begin{proof}
    By Lemma \ref{lemma:transverse-smooth-guy-exists}, $\mathbf{v}_t\coloneqq t\mathbf{v}_0+(1-t)\mathbf{v}_0$ is a positively transverse tuple for all $t \in [0,1]$, hence we get an associated homotopy of pre-smooth transverse plane fields $\nu_t$.  Pick two $\varepsilon$-smoothings $\nu_0'$ and $\nu_1'$ of $\nu_0$ and $\nu_1$, respectively (in the sense of Corollary \ref{corollary:smoothing-plane-fields}) so that we can take $X_{\mathbf{v}_i}=X_{\nu_i'}$ for $i=0,1$, then for all $\varepsilon>0$ sufficiently small, \cite[Theorem 1.3]{Whitehead1961} guarantees that there exist pre-smooth homotopies $\nu_0' \rightsquigarrow \nu_0$ and $\nu_1 \rightsquigarrow \nu_1'$. Therefore, concatenating these homotopies with the convex interpolation $\nu_0 \rightsquigarrow \nu_1$ yields a pre-smooth homotopy $\nu_0'\rightsquigarrow\nu_1'$. Since both the endpoints are smooth plane fields, we can use the same strategy as in Lemma \ref{lemma:TCCI-smoothing-unique} to approximate this homotopy by a smooth one, which implies the desired conclusion by Lemma \ref{lemma:htpic-implies-isotopic}.
\end{proof}

This tells us that we can drop the $\mathbf{v}$ from the the notation for smoothing for compact TCCI's, since the notion does not depend on the choice up to an ambient isotopy of $\R^n$. Therefore, in this setting, everything is canonical and it makes sense to speak of \emph{the} smoothing of a topological manifold $X \subset \R^n$. Remarkably, neither existence nor uniqueness hold in general, and a nice set of (counter)examples is presented at the end of \cite[Section 1]{Whitehead1961}.

\subsection{Smoothing in the presence of a stratification}\label{section:stratified-smoothing} In this section, we explain how to generalise results like Corollary \ref{corollary:TCCI-smoothing-exists} and Proposition \ref{proposition:TCCI-smoothing-unique} to the setting when one also wants to preserve the stratification of $\R^n$ into polyhedral domains, such as a stratification by a simplicial fan $\Sigma$ (which play a crucial role, for example, in Section \ref{section:potentials}). Even though we will be dealing exclusively with cases where all the strata are convex polyhedra, we make the following more general definition that works when $\R^n$ is replaced by a different ambient smooth manifold: 

\begin{definition}\label{definition:plain-stratification}
    We call a Whitney stratification $\mS$ of a manifold with corners $M$ into finitely many pieces \emph{plain} if the closure of each stratum is a submanifold with corners\footnote{Using the definitions of various types of smooth maps from \cite{joyce}, we get a notion of \emph{embedded submanifolds-with-corners} of a manifold-with-corners.} of $M$. 
\end{definition}

\begin{example}
    The trivial stratification of $M$ with a single stratum equal to $M$ is plain. 
\end{example}

\begin{example}
    A simplicial fan $\Sigma$ defines a plain stratification $\mS_\Sigma$ of $\R^n$, which we call the \emph{stratification induced by $\Sigma$}. More generally, a decomposition of $\R^n$ into finitely many simple polyhedra gives a plain stratification. 
\end{example}

The stratification $\mS_\Sigma$ can be refined in the following way: suppose that we are also given a simple polytope $P$ for which $\Sigma$ is the normal fan, so that its barycentric subdivision is well-defined. For cones $\tau \subset \sigma$, denote $\sigma^{\ba}_\tau\coloneqq\Sigma^\ba_\tau \cap \sigma$ (recall that we are using the convention $\Sigma^\ba_0=\Sigma^\ba$). More explicitly, $\sigma^\ba_{\tau} \cong [0,1]^{\dim(\sigma)-\dim(\tau)}$ is the face of $\sigma^\ba$ consisting of points where the expression $\sum_{\rho \in \tau(1)}v_\rho$ appears. For a positive number $t>0$, we define a stratification $\mS_{\Sigma}^t$ whose strata are $S^t(\sigma,\tau,\tau')\coloneqq\relint(t\cdot\sigma^\ba_\tau+\tau')$ for cones $\tau' \subseteq \tau \subseteq \sigma$ (see Figure \ref{fig:refined-strat} for an illustration). It is straightforward to check that the stratification defined this way is plain, that the closure of $S^t(\tau',\tau,\sigma)$ is a convex polyhedron lying inside $\sigma$ diffeomorphic to $[0,1]^{\dim(\sigma)-\dim(\tau)}\times\R_{\geq0}^{\dim(\tau')}$, and that the indexing is compatible with the poset structure on the cones of $\Sigma$ through the following relation:
\begin{equation*}
    \overline{S^t(\sigma,\tau,\tau')}=\bigcup_{\widetilde{\sigma} \subseteq \sigma, \tau \subseteq \widetilde{\tau}, \widetilde{\tau}' \subseteq \tau' } S^t(\widetilde{\sigma},\widetilde{\tau},\widetilde{\tau}').
\end{equation*}
We call $\mS^t_\Sigma$ the \emph{refined stratification induced by $\Sigma$ with parameter $t$}. 

\begin{figure}[ht]
    \centering
\begin{tikzpicture}[scale=1, every node/.style={font=\small}]
    
    \draw[->] (0,0) -- (-3,0);
    \draw[->] (0,0) -- (0,-3);
    \draw[->] (0,0) -- (3,3);
    \draw[->] (0,0) -- (0,3);
    
    \draw[dashed] (0,0) -- (-1,0) -- (-1,-1) -- (0,-1) -- cycle;
    
    \draw[dashed] (0,0) -- (0,-1) -- (1,0) -- (1,1) -- cycle;
    
    \draw[dashed] (0,0) -- (1,1) -- (1,2) -- (0,1) -- cycle;
    
    \draw[dashed] (0,0) -- (0,1) -- (-1,1) -- (-1,0) -- cycle;
    
    \draw[dashed] (-1,-1) -- (-3,-1);
    \draw[dashed] (-1,-1) -- (-1,-3);
    \draw[dashed] (1,0) -- (1,-3);
    \draw[dashed] (1,0) -- (3,2);
    \draw[dashed] (1,2) -- (2,3);
    \draw[dashed] (1,2) -- (1,3);
    \draw[dashed] (-1,1) -- (-3,1);
    \draw[dashed] (-1,1) -- (-1,3);
    
    \fill[blue, opacity=0.3] (-1,-1) -- (-3,-1) -- (-3,-3) -- (-1,-3) -- cycle;
    \node at (-4,-2){$\color{blue}S(\sigma_3,\sigma_3,\sigma_3)$};

    \fill[green, opacity=0.3] (-1,-1) -- (-3,-1) -- (-3,0) -- (-1,0) -- cycle;
    \node at (-4,-0.5){$\color{green}S(\sigma_3,\rho_4,\rho_4)$}; 

    \fill[orange, opacity=0.3] (-1,-1) -- (-1,0) -- (0,0) -- (0,-1) -- cycle;
    \node at (0.8,-0.55){$\color{orange}S(\sigma_3,0,0)$}; 
    
    \draw[red, very thick] (-1,-1) -- (-3,-1);
    \node at (-4,-1.25){$\color{red}S(\sigma_3,\sigma_3,\rho_4)$};

    \draw[violet, very thick] (-1,-1) -- (-1,0);
    \node at (-1.85,0.25) {$\color{violet}S(\sigma_3,\rho_4,0)$};

    \node[fill,circle,inner sep=1.5pt, color=black] at (-1,-1){};
    \node at (-1.85,-1.25) {$S(\sigma_3,\sigma_3,0)$};
  
\end{tikzpicture}

    \caption{The refined stratification $\mS^t_\Sigma$ associated to the fan of the first Hirzerbruch surface (as considered in Figure \ref{fig:bar-subdiv}) with a few regions labelled explicitly}
    \label{fig:refined-strat}
\end{figure}

\begin{definition}\label{definition:substratification}
    For two plain stratifications $\mS$, $\mS'$, we say that \emph{$\mS'$ is a refinement of $\mS$} if, for each $S \in \mS$, $\mS'(S)\coloneqq\{S' \in \mS' \colon S' \subset \overline{S}\}$ defines a plain stratification of $\overline{S}$.
\end{definition}

\begin{example}
    Any plain stratification is a refinement of the trivial stratification. The stratification $\mS^t_\Sigma$ is a refinement of $\mS_\Sigma$ for all $t>0$.
\end{example}
\begin{definition}\label{definition:s-smooth-manifolds}
    Suppose that $\mS$ is a plain stratification of $M$. We call a $X \subset M$ an \emph{$\mS$-smooth manifold with corners} if it is a topological manifold with boundary and $X_S\coloneqq X \cap \overline{S}$ is a smooth submanifold (with corners) of $\overline{S}$ for all $S \in \mS$ and, moreover, every closed stratum $\overline{V} \subset X_{S}$ that is not contained in $\partial\overline{S}$ intersects every boundary stratum boundary $\overline{S}'\subset\partial \overline{S}$ transversely inside $\overline{S}$. 
\end{definition}

Note that for a $\mS$-smooth manifold with corners $X$, we obtain a stratification of its (topological) boundary $\partial X$ by sets $\partial^{in}X_S \coloneqq \partial X \cap S$ for $S \in \mS$, which can be further decomposed into manifolds with corners following the boundary structure of $X_S$.  Also, transversality of intersections implies that we have $\codim_S(X_S)=\codim_M(X)$ whenever $X_S \neq \emptyset$. 

\begin{example}
    $\Sigma^\ba$ is a $\Sigma$-smooth manifold with corners in $\R^n$. A closed subset $X \subset M$ is a smooth submanifold with corners if and only if it is smooth with respect to the trivial stratification and $X$ intersects every boundary stratum of $M$ transversely. 
\end{example}

\begin{example}
    A typical non-example would the diagonal $\Delta \subset[0,1] \times [0,1]$, which is not smooth with respect to the trivial stratification, since it is codimension $1$ and hence is not allowed to intersect the $0$-dimensional boundary strata. 
\end{example}

Intuitively, this somewhat ad hoc class of manifolds should allow us to consider $X$ that look smooth along the tangent directions of the strata of $\mS$, but are allowed to bend in the normal direction. The notion also interacts well with refinements, as recorded by the following Lemma: 

\begin{lemma}\label{lemma:refined-stratification}
    Let $\mS$ be a plain stratification of $M$ with a refinement $\mS'$, then an $\mS$-smooth submanifold $Y$ is $\mS'$-smooth if and only if it intersects all the strata $S' \in \mS'\backslash \mS$ transversely.   
\end{lemma}

\begin{proof}
    This follows by \cite[Theorem 6.4]{joyce}, telling us that transverse intersections of manifolds with corners are manifolds with corners.
\end{proof}

\begin{corollary}\label{corollary:refined-stratification}
    If $\Sigma\subset\R^n$ is a simplicial fan and $X \subset \R^n$ is $\Sigma$-smooth and compact, then $X$ is also $\mS^t_\Sigma$-smooth for all $t>0$ small enough.
\end{corollary}

\begin{proof}
    The previous Lemma tells us that it suffices to check transversality of intersection of $X$ with the strata of $\mS^t_\Sigma$. But it follows from stability of transverse intersections under small perturbations that if $X \cap \tau'$ is transverse, then so is $X \cap S^t(\sigma,\tau,\tau')$ for sufficiently small $t$. 
\end{proof}

With this in mind, we can generalise the definition of a TCCI to this setting by replacing the notion of smoothness with $\mS$-smoothness: 

\begin{definition}\label{definition:sigma_TCCI}
    Suppose that $Y_1$, \dots, $Y_r$ are $n$-dimensional $\mS$-smooth manifolds with corners in $\R^n$ that intersect transversely\footnote{In the sense that the manifolds with corners $Y_{S,1}, \dots, Y_{S,r}$ intersect transversely inside $\overline{S}$ for all $S \in \mS$.}, then we call $X\coloneqq \bigcap_{j=1}^r \partial Y_j$ an \emph{$\mS$-transverse cornered complete intersection ($\mS$-TCCI)} of codimension $r$.
\end{definition}

\begin{example}
    $\partial\Sigma^\ba$ is a compact codimension $1$ $\Sigma$-TCCI; its subcomplex $\Sigma^\ba_{\trans}$ induced by a nef-partition into $r$ polytopes is a codimension $r$ $\Sigma$-TCCI. 
\end{example}

Now, we explain how all the basic definitions from the previous section extend to the setting of $\mS$-TCCI's by working one stratum at a time: first, suppose that $Y$ is $\mS$-smooth and suppose that $x \in \partial Y \cap S \subset Y_S$ for some $S \in \mS$ of dimension $l$ and $x$ is in a codimension $m$ stratum for $m \geq 1$. Then the \emph{$\mS$-normal cone} $NC^\mS_x Y \subset T_xS$ can be defined by picking a chart $\phi$ identifying a pair $(U,U \cap Y_S)$ for a neighbourhood $x \in U \subset S$ with a standard neighbourhood and then taking $NC^\mS_x Y \coloneqq -D\phi_0(C_m \times \R^{l-m}) \subset T_x S$. The \emph{$\mS$-conormal cone} $N^*C_x^\mS Y \subset T_x^*S$ is then defined as the dual cone. Given a $\mS$-TCCI $X=\bigcap_{j=1}^r \partial Y_j$ and $x \in X$, a tuple $\mathbf{v}(x)=(v_1(x),\dots,v_r(x)) \in T_x S$ is then called \emph{$\mS$-positively transverse to $X$} if the matrix from Definition \ref{definition:positively-transverse-tuple} is positive definite for all $\alpha_i \in N^*C_x^\mS Y_i$ (write this as $\mathbf{v}(x) \in \mG^{\mS,+}_x$). This naturally leads to a notion of global \emph{positively transverse $\mS$-tuples} $\mathbf{v} \in \Gamma(X,\mG^{\mS,+})$. We also obtain the following natural generalisation of Lemma \ref{lemma:transverse-smooth-guy-exists} to the stratified setting:

\begin{lemma}\label{lemma:sigma-transverse-guys-exist}
    Let $X=\bigcap_{j=1}^r \partial Y_j \subset \R^n$ be a compact $\mS$-TCCI, then the space of positively transverse $\mS$-tuples $\Gamma(X,\mG^{\mS,+})$ is non-empty and contractible. 
\end{lemma}

\begin{proof}
    Note that the condition on being tangent to $S \in \mS$ at all $x \in S$ is preserved by taking convex combinations, so the previous argument still shows convexity of the space. 

    The construction is also essentially identical, except we need the following observation: if $x \in X_S$ and $\mathbf{v}(x) \in \mG^{\mS,+}_x$, then there exists a neighbourhood of $x\in U \subset X$ such that $\mathbf{v}(x) \in \mG^{\mS,+}_y$ for all $y\in U$ (which follows from the fact that the strata of $X_{S'}$ not contained in $X_S$ with $S \subsetneq S'$ are transverse to $S$). This enables us to pick standard neighbourhoods inside $X_S$, take constant functions over them and observe that it they will also satisfy the condition over some thickenings of these neighbourhoods inside $X$, which means that the local construction of Lemma \ref{lemma:transverse-smooth-guy-exists} also works in this setting. 
\end{proof}

Suppose that $X$ is an $\mS$-TCCI and let $x \in \relint(X_S)$, then we have a local model for the TCCI $\relint(X_S)$ by Lemma \ref{lemma:local_form_TCCI}, hence we can pick a chart $\phi_{S}$ for it to identify the neighbourhood of $x$ in $X_S$ with a product of standard corners $C_{l,\mathbf{m}}$ for $l=\dim(S)$ and define the \emph{space of transverse $\mS$-planes to $X$ at $x$} as 
\begin{equation*}
    \mG^{\mS}_x\coloneqq(D_x\phi_S)^{-1}(\mG_{l,\mathbf{m}}) \subset \Gr(T_xS,r) \subset \Gr(T_x\R^n,r).
\end{equation*}
The notions of a (pre)-smooth field of transverse planes on $X$ then naturally extend to this stratified setting to a \emph{(pre)-smooth field of transverse $\mS$-planes on $X$} by replacing the condition $\nu(x) \in \mG_x$ with $\nu(x) \in\mG^{\mS}_x$. It is immediate from the non-stratified case that the map considered in Lemma \ref{lemma:neg-transverse-means-transvers} gives a canonical map $\Gamma(X,\mG^{\mS,+}) \rightarrow \Gamma(X,\mG^{\mS})$. We also obtain all associated notions, such as $\nu$-neighbourhoods and transverse homotopies \emph{relative to $\mS$}, in this new context, resulting in the following definition:

\begin{definition}\label{definition:sigma-smoothin}
    Let $X$ be an $\mS$-TCCI and $\nu \in \Gamma(X,\mG^{\mS})$ a smooth transverse field of $\mS$-planes along $X$, then a \emph{$\nu$-smoothing of $X$ relative to $\mS$} is the image of a smooth section of $\mN^\varepsilon_\nu \xrightarrow{\pi} X$ for a sufficiently small $\varepsilon>0$.
\end{definition}

Analogously to the case of TCCI's, this definition is independent of the choice of the section or $\varepsilon$ up to an ambient isotopy preserving $\mS$. 

\begin{corollary}\label{corollary:sigma-smoothing-exists}
    Any compact $\mS$-TCCI in $\R^n$ admits a smoothing. 
\end{corollary}

\begin{proof}
    The proof is exactly analogous to Corollary \ref{corollary:TCCI-smoothing-exists}, except one needs to perform the smoothing of the pre-smooth plane fields inductively from smaller strata to larger ones by employing the relative version of \cite[ Theorem 1.10]{Whitehead1961}. 
\end{proof}

The constructions of transverse homotopies also extend to this setting with the small alteration that we now consider isotopies that preserve $\mS$ (which becomes possible thanks to the condition $\nu(x) \subset T_xS$ for $x\in S$), so the discussion about uniqueness of smoothings also works with minor alterations and culminates in the following result:

\begin{proposition}\label{proposition:sigma-smoothing-unique}
    For any positively transverse $\mS$-tuples $\mathbf{v}_0,\mathbf{v}_1 \in \Gamma(X,\mG^{\mS,+})$ along a compact $\mS$-TCCI $X$, the two associated smoothings $X_{\mathbf{v}_0}$ and $X_{\mathbf{v}_1}$ relative to $\mS$ are ambient isotopic relative to $\mS$.
\end{proposition}

Therefore, analogously to the setting of TCCI's, it makes sense to speak of \emph{the} smoothing of an $\mS$-TCCI relative to $\mS$.

\subsection{Smoothings of tropical Batyrev--Borisov complete intersections}\label{section:smoothing-applications} 

In this section, we apply our general theory of smoothings to the spaces arising as skeleta in Section \ref{section:skeleta} (to make precise the intuition that the objects obtained from tailored complete intersections are indeed smoothings of the tropical complete intersections) and problems in Section \ref{section:potentials} (to explain how to construct genuine isotopies from stratified pre-isotopies). 

\begin{proposition}\label{proposition:smoothing-bbci}
    The manifold $\partial\widetilde{C}_{0,\beta}= \bigcap_{j=1}^r{\partial \widetilde{C}_{0,\beta,j}}$ is a smoothing of the complex of bounded cells $\partial C_{0,\trop} \subset \mA_{\trop}$ of a tropical smooth BBCI for all $\beta>0$ large enough. 
\end{proposition}
\begin{proof}
    Recall that, by definition, we have 
    \begin{equation*}
        \partial\widetilde{C}_{0,\beta}=\{ u \in M_\R \colon \widetilde{F}_{\beta,j}(u)=|c_{0,j}| \textnormal{ for all } 1 \leq j \leq r\},
    \end{equation*}
    where the functions $\widetilde{F}_{\beta,j}(u)$ are smoothings of the functions $\max_{\alpha \in \nabla_j\backslash \{0\}}\{e^{\beta l_\alpha(u)}\}$. Also, for any Riemannian metric $g$ on $M_\R$, one can consider a tuple $\mathbf{v}=(v_1,\dots,v_r)$ with $v_j(u)=\nabla_g\widetilde{F}_{\beta_0,j}(u)$ for all $j=1,\dots, r$, where $\beta_0>0$ is a large constant. Since the defining functions for $C_{0,\trop}$ are linear and the differentials of $\widetilde{F}_{\beta,j}$ are their positive linear combinations, they lie in the conormal cones of the defining tropical hypersurfaces.
    
    We can use an averaging construction analogous to the one in Lemma \ref{lemma:pre_isotopy_exists} to get a Riemannian metric $g$ defined on some neighbourhood of $C_{0,\trop}$ such that near each stratum, the defining linear functionals $\alpha$ are orthonormal with respect to it, so fix such a metric for the rest of the proof. In addition to that, we take $\beta_0$ large enough so that at $u$ in $C_{0,\trop}$, the functionals $\alpha$ such that $\widetilde{F}_{\alpha,\beta_0,j}(u) \neq 0$ for some $j$ are orthonormal. 
    This means that the $r$ subspaces $\spann_\R(\alpha \in \nabla_j  \backslash\{0\}\colon \widetilde{F}_{\alpha,\beta_0,j}(u) \neq 0)$ will be mutually $g$-orthogonal at $u$ and so the matrix $M$ from definition of positive transversality will always be diagonal with positive entries, hence the tuple $\mathbf{v}(u)$ associated to such a metric is positively transverse at $u$, so $\mathbf{v}$ defines a positively transverse tuple along $\partial C_{0,\trop}$ with a corresponding transverse pre-smooth plane field $\nu(u)=\spann_\R(v_1(u),\dots,v_r(u))$. 

    Therefore, we get a $\nu$-tubular neighbourhood $\mN^{4\varepsilon}_\nu$ for some $\varepsilon>0$. In fact, via estimates similar to the ones in the proof of \cite[Theorem 1.5]{Whitehead1961}, it follows that there exists an $\varepsilon>0$ for which there also exists some $\delta_0>0$ such that whenever $\nu'$ is $\delta_0$-close to $\nu$, $\mN^{2\varepsilon}_{\nu'}$ is a $\nu'$-tubular neighbourhood and $\mN^{\varepsilon}_\nu \subset \mN^{2\varepsilon}_{\nu'} \subset \mN^{4\varepsilon}_\nu$. 

    By Theorem \ref{theorem:trop_ci}, the distance of $\partial \widetilde{C}_{0,\beta}$ from $\partial C_{0,\trop}$ is $O(\beta^{-1})$, so for all sufficiently large $\beta$, we can guarantee that $\partial \widetilde{C}_{0,\beta} \subset \mN^{\varepsilon}_\nu$. By definition of the projection fro the tubular neighbourhood, any $u \in \partial \widetilde{C}_{0,\beta}$ can be written as $u=\pi(u)+v$ for a unique $\pi(u)\in \partial C_{0,\trop}$ and some $v \in f(\pi(u))$ of size $O(\varepsilon)$. Since $\pi(u)$ is on $\partial C_{0,\trop}$, $u$ is close to it and both $\beta, \beta_0$ are sufficiently large, the vectors $d\widetilde{F}_{\beta}(u)=(d\widetilde{F}_{\beta,1}(u),\dots,d\widetilde{F}_{\beta,r}(u))$ and $d\widetilde{F}_{\beta_0}(\pi(u))$ lie inside a product $\sigma_1 \times \dots \times \sigma_r$ for some transversal cone $\sigma=\sigma_1+\dots+\sigma_r$. By construction of $g$, the splitting $\R\sigma=\bigoplus_j \R \sigma_j$ is $g$-orthogonal, so the $(r \times r)$-matrix $\langle d\widetilde{F}_{\beta}(u),d\widetilde{F}_{\beta_0}(\pi(u)) \rangle_{g(u)}$ is diagonal with positive entries. Therefore, $T_u \partial \widetilde{C}_{0,\beta}=\ker(d\widetilde{F}_{\beta}(u))$ is transverse to $\nu(\pi(u))$ for all sufficiently large $\beta$ and so the manifold $\partial \widetilde{C}_{0,\beta}$ intersects the fibres of $\pi=\pi_\nu\colon\mN^\varepsilon_\nu \rightarrow \partial C_{0,\trop}$ transversely. Since the intersection number is locally constant thanks to the transversality statement, it suffices to check that the intersection is a point at a single point at a single point in every connected component of $C_{0,\trop}$. To do that, note that at points in the top-dimensional strata of $C_{0,\trop}$ and sufficiently far from the lower-dimensional strata, all the tailored polynomials $\widetilde{F}_{\beta},j$ and $\widetilde{F}_{\beta_0,j}$ are just monomials, so it is straightforward to check that the intersection with fibres is just a single point there, and that every connected component of $C_{0,\trop}$ will contain a top-dimensional stratum.

    Finally, suppose we are given a $\beta$ large enough so that the above discussion applies to it. Then, we can observe that the transversality of the intersection with fibres of $\pi_\nu$ is also an open condition, so there exists some $\delta(\beta)>0$ such that if $\nu'$ is $\delta(\beta)$-close to $\nu$, then $\partial \widetilde{C}_{0,\beta}$ will still intersect the fibres of $\pi_{\nu'}$ inside $\mN^{2\varepsilon}_{\nu'}$ transversely at a single point. By \cite[Theorem 1.10]{Whitehead1961}, there exists some smooth transverse plane field $\nu_\beta$ that is $\min\{\delta_0,\delta(\beta)\}$-close to $\nu$. By construction, $\partial \widetilde{C}_{0,\beta}$ is then the image of a smooth section of a $\nu_\beta$-tubular neighbourhood, so we're done.  
\end{proof}

The main application is passing from ambient isotopies that are piecewise smooth to smooth ones: suppose that $M$ is a manifold with corners with a plain stratification $\mS$, we shall call a map $H^t\colon M \times [0,1] \rightarrow M$ an \emph{$\mS$-isotopy} if it is a continuous isotopy and restricts to a smooth isotopy $H^t_S \colon \overline{S} \times [0,1] \rightarrow \overline{S}$ for each stratum $S \in \mS$. 

\begin{proposition}\label{proposition:smoothing-isotopies}
    Suppose that $Y \subset M$ is a compact $\mS$-TCCI and $H^t$ an $\mS$-isotopy with $H^0=\textnormal{id}$, $Y_0=Y$ and $Y_1=H^1(Y)$. Then there exists a smooth isotopy $\widetilde{H}^t \colon M \times [0,1] \rightarrow M$ preserving $\mS$ that takes the $\mS$-smoothing $\widetilde{Y}_0$ of $Y_0$ to the $\mS$-smoothing $\widetilde{Y}_1$ of $Y_1$. 
\end{proposition}

\begin{proof}
    Denote $Y_t\coloneqq H^t(Y)$, then each $Y_t$ is also a compact $\mS$-TCCI, so the smoothing of $\widetilde{Y}_t$ is well-defined by above. We also pick a smooth positively transverse $\mS$-tuple $\mathbf{v}_0$ along $Y_0$, then $\mathbf{v}_t\coloneqq H^t_*(\mathbf{v}_0)$ defines a positively transverse $\mS$-tuple along $Y_t$ for all $t$. 

    By composing with a non-decreasing smooth function from $[0,1]$ to itself that is identically equal to $0$ near $0$ and to $1$ near $1$, without loss of generality suppose that $H^t$ does not depend on $t$ near $0$ and $1$. Then we can extend it to an isotopy with domain $M\times\R$ by taking the constant extension along $(-\infty,0]$ and $[1,+\infty)$, note that it will depend on $t$ smoothly (we have not defined smoothings of topological manifolds with boundary, so this extension is necessary to replace a manifold with boundary $Y\times[0,1]$ with $Y\times \R$). Consider $Y_{\R}\coloneqq \bigcup_{t\in \R} Y_t \times \{t\} \subset M\times\R$ and note that the stratification $\mS$ induces a plain stratification $\mS_{\R}$ of $M\times\R$. Moreover, by techniques similar to the ones used in Corollary \ref{corollary:cornery_stuff1}, we can prove that $Y_{\R}$ is an $\mS$-TCCI that also fibres over $\R$. By setting $\mathbf{v}_{\R}(x,t)=\mathbf{v}_t(x)$ (with zero $\partial_t$-component), we get a positively transverse $\mS_{\R}$-tuple along $Y_{\R}$. This has an associated transverse $\mS$-plane field $\nu_{\R}$, which can be smoothed to $\nu_{\R}'$ as in Corollary \ref{corollary:smoothing-plane-fields} (note that we have assumed compactness of the $\mS$-TCCI in all the cases above, but we do not run into any problems here, since $\nu_\R$ is independent of the $t$-component for $t<0$ and $t>1$ and $Y_{\R}$ is compact in the other direction).

    Therefore, the choice of $\mathbf{v}_0$ also naturally gives us a smoothing $\widetilde{Y}_{\R}$ whose fibre over $t \in \R$ is a $\mathbf{v}_t$-smoothing $\widetilde{Y}_t$ of $Y_t$. This puts us in a setting analogous to the final step of the proof Lemma \ref{lemma:family_almost_ad}, where we have a sub-bundle $\widetilde{Y}_{\R}$ of the trivial fibre bundle $M\times\R$ over $\R$, so by integrating the lift of $\partial_t$ along the submersive projection like there, we get an ambient isotopy that takes $\widetilde{Y}_0$ to $\widetilde{Y}_1$ (one can make sure that the isotopy preserves $\mS$ thanks to the various transversality assumptions).
\end{proof}

The most important consequence of this statement is completing the proof of Proposition \ref{proposition:adapted_existence} by explaining how to use what we call \emph{stratified pre-isotopies} to construct genuine smooth isotopies preserving $\Sigma$: suppose that $\varphi$ is an adapted potential to a simple polytope $P$ with normal fan $\Sigma$. By adaptedness of $\varphi$, $\Phi(P)$ is $\Sigma$-smooth, so we can consider smoothings of both $\Phi(\partial P)$ and $\partial\Sigma^\ba$ as $\Sigma$-TCCI's. 

\begin{corollary}\label{corollary:smoothing-isotopies1}
    Suppose that $\varphi$ is an adapted potential such that there exists a pre-stratified isotopy $\{H^t_\sigma\}_{\sigma \in \Sigma}$ such that $H^t_\sigma$ deforms $\Phi(P) \cap \sigma$ to $\sigma^\ba$ for all $\sigma \in \Sigma$. Then $\varphi$ is strongly adapted.  
\end{corollary}

\begin{proof}
    By the previous Proposition, it is enough to produce a $\Sigma$-isotopy between $\partial Q$ for $Q\coloneqq\Phi(P)$ and $\partial\Sigma^\ba$, so we shall focus on that. 

    First, we describe how to deform $Q$ so that it looks like a product near the boundary of each cone (and identifications with products satisfy the right compatibility relations). In order to do that, let $\delta>0$ be small positive number and pick some smoothing $r^\delta$ of $\max\{x-\delta,0\}$ that is non-decreasing, agrees with the function for $x<\delta$ and $x>2\delta$, has the same support $(\delta,+\infty)$ and is $\delta^2$-close to it over the interval $(\delta,2\delta)$. Over each cone $\sigma\in\Sigma$, we then obtain a function $R_{\sigma}^\delta \colon \sigma \rightarrow \sigma$ given by 
        \begin{equation*}
            R_{\sigma}^\delta\left(\sum_{\rho \in \sigma(1)}{x_\rho v_\rho}\right)=\sum_{\rho \in \sigma(1)}r^\delta(x_\rho)\cdot v_\rho.
        \end{equation*}
    Observe that, by design, these functions glue to a continuous function $R^\delta\colon N_{\R} \rightarrow N_{\R}$, whose restriction to every cone is a weakly smooth map of manifolds with corners. We can check that $Q^\delta\coloneqq(R^\delta)^{-1}(Q)$ is going to be $\Sigma$-smooth for all sufficiently small $\delta>0$: let $\sigma$ be a cone of $\Sigma$ with rays $\sigma(1)=\{\rho_1,\dots,\rho_k\}$, then the region $Q_\sigma$ is given by $Q_\sigma =\{x \in \sigma \colon f_j(x)\leq1 \textnormal{ for all } j=1,\dots,k)\}$ for defining functions $f_j(x)\coloneqq \langle v_{\rho_j},\Phi^{-1}(x)\rangle$. Adaptedness of $\varphi$ guarantees that no other rays of $\Sigma$ will contribute to this equation and that $Q_\sigma$ is a manifold with corners that intersects the boundary $\partial \sigma$ transversely in the sense of Definition \ref{definition:s-smooth-manifolds}. In fact, both these facts follow from the hypersurfaces $H_j\coloneqq \{ x \in \sigma \colon f_j(x)=1\}$ intersecting transversely inside $\sigma$ and intersecting the boundary transversely (in the strong sense, i.e. any intersection $H_{i_1} \cap \dots \cap H_{i_l} \cap \relint(\tau)$ is transverse inside $\R\sigma$ for all $\{i_1,\dots,i_l\} \subseteq \{1,\dots,k\}$ and $\tau \subset \partial\sigma$). By construction, the map $r^\delta$ is $O(\delta)$-close to the identity (in $C^\infty$) over $[2\delta,+\infty)$, so it follows that the hypersurfaces $H^\delta_j$ with defining functions $f^\delta_j\coloneqq f_j \circ R^\delta_\sigma$ will still have this property for all sufficiently small $\delta>0$, so $Q^\delta_\sigma\coloneqq Q^\delta \cap \sigma$ will also be a manifold with corners inside $\sigma$ that intersects the boundary $\partial\sigma$ transversely and $Q^\delta$ is indeed $\Sigma$-smooth. Moreover, working inductively cone by cone, we can construct a $\Sigma$-isotopy that takes $Q$ to $Q^\delta$: it suffices to smoothly deform $r^\delta$ to the identity function, which induces deformations of each $R^\delta_\sigma$ to $\textnormal{id}_\sigma$ that glue on overlaps. Therefore, by Lemma \ref{proposition:smoothing-isotopies}, it suffices to construct a smooth isotopy that preserves $\Sigma$ and takes the smoothing of $\partial Q^\delta$ to the smoothing of $\partial \Sigma^\ba$.
    
    By construction of $R^\delta_\sigma$, we also have $(R^\delta_\sigma)^{-1}(\tau)=\overline{S^\delta(\sigma,\tau,\tau)}$ for any pair of cones $\tau \subseteq \sigma$, so Lemma \ref{lemma:refined-stratification} tells us that $Q^\delta$ is also $\mS^\delta_\Sigma$-smooth and covered by closed strata $\overline{Q^\delta(\sigma,\tau)}\coloneqq Q^\delta \cap \overline{S^\delta(\sigma,\tau,\tau)}$. Note that there is a projection $\pi^\ba\colon N_\R \rightarrow \delta\cdot\Sigma^\ba$ defined by mapping $\delta\cdot\sigma^\ba_\tau+\tau'$ onto $\delta\cdot\sigma^\ba_\tau$ on $S(\sigma,\tau,\tau')$, which is continuous and piecewise smooth with respect to the polyhedral stratification $\mS^\delta_\Sigma$. By construction, $\pi^\ba\colon\overline{Q^\delta(\sigma,\tau)}\rightarrow \delta\cdot\sigma^\ba_\tau$ is a submersive surjection, so a fibre bundle by Ehresmann's theorem (Theorem \ref{theorem:ehresmann_with_corners}). By construction of $R^\delta_\sigma$, this fibre bundle is trivial and $\pi^\ba$ agrees with a standard projection, so we get identifications $\overline{Q^\delta(\sigma,\tau)} \cong Q_\tau \times \sigma^\ba_\tau$ that are continuously compatible on overlaps.

    Pick a subdivision of $P$ into cubes as in Lemma \ref{lemma:bar_subdiv_homeo}, with the cube corresponding to a cone $\sigma$ being denoted as $P(\sigma)$ and its outward faces as $P(\sigma,\tau)\coloneqq P(\sigma)\cap F_\tau$ for $\tau \subsetneq \sigma$. This choice also induces a homeomorphism $\psi \colon \Sigma^\ba \rightarrow P$, which restricts to diffeomorphisms from $\sigma^\ba_\tau$ to $P(\sigma,\tau)$. This allows us to define a global isotopy $H^t \colon N_\R \times [0,1] \rightarrow N_\R$ by taking $H^t(x)\coloneqq H^t_\tau(\psi(\delta^{-1}\cdot\pi^\ba(x)),x-\pi^\ba(x))$ whenever $\pi^\ba(x) \in \Sigma^\ba_\tau$. This is well-defined, since $x-\pi^\ba(x) \in \tau$ under these conditions, induces a smooth isotopy of every closed stratum $\overline{S(\sigma,\tau,\tau)}$ and glues to a continuous isotopy thanks to the compatibility conditions on pre-stratified isotopies. Since $H^t_\tau$ takes $Q \cap \tau$ to $\tau^\ba$, the identifications of $\overline{Q^\delta(\sigma,\tau)}$ with products imply that the isotopy $H^t$ will take $Q^\delta$ to $\Sigma^\ba$ and it is, by construction, an $\mS^\delta_\Sigma$-isotopy. Therefore, by Proposition \ref{proposition:smoothing-isotopies}, we obtain a smooth isotopy that preserves $\mS^\delta_\Sigma$ and takes the smoothing of $\partial Q^\delta$ to the smoothing of $\partial \Sigma^\ba$. Because $\mS^\delta_\Sigma$ is a refinement of $\mS_\Sigma$, it will also preserve $\Sigma$, so we are done. 
\end{proof}

In the setting where the data comes from a tropical smooth BBCI associated to some nef-partition, we can repeat the same reasoning restricted to the subcomplex of transversal faces $\partial C_{0,\trop}$ of the polytope $P=C_{0,\trop,\textnormal{tot}}$ and combine it with Proposition \ref{proposition:smoothing-bbci} to get the following result:

\begin{corollary}\label{corollary:smoothing-isotopies2}
    For any potential $\varphi$ strongly adapted to $C_{0,\trop,\tot}$, there exists a smooth ambient isotopy of $N_\R$ that preserves $\Sigma$ and takes $\Phi(\partial\widetilde{C}_{0,\beta})$ to a smoothing of $\Sigma^\ba_\trans$. 
\end{corollary}

\subsection{Partial smoothings of barycentric subdivisions}\label{section:ooga-booga-smoothing} In this section, we explain how to manually construct a particular type of smoothing of $\Sigma^\ba$ and $\Sigma^\ba_\trans$, which play an important role in Section \ref{section:combinatorics}. The main goal is showing the existence of \emph{collections of coherent projections} and explaining how to smooth the \emph{nef scaling actions} introduced there. 

Recall that for a fan $\Sigma$ coming from a regular triangulation $\mT$ of $\partial \Delta^\vee$ and a cone $\sigma \in \Sigma$, we have a distinguished set of ray generators, which give us a basis $\{v_\rho \colon \rho \in \sigma(1)\}$ for $\R\sigma$ and hence also dual linear functionals $\eta^\rho_\sigma$ on the space. For an inclusion of cones $\sigma \subset \tau$, we clearly have $\eta^\rho_\tau|_\sigma=\eta^\rho_\sigma$ for all $\rho \in \sigma(1)$. Therefore, for any $\rho \in \Sigma(1)$, we can define a continuous piecewise linear function $h^\rho$ by letting $h^\rho(u)=0$ if $u \notin \st(\rho)$ and $h^\rho(u)=\eta^\rho_\sigma(u)$ if $u \in \sigma$ otherwise. The barycentric subdivision $\Sigma^\ba$ can then be alternatively described as 
\begin{equation*}
    \Sigma^\ba=\{ u \in N_\R \colon h^\rho(u) \leq 1 \textnormal{ for all } \rho \in \Sigma(1)\}.
\end{equation*}
We can also think of the function $h^\rho$ as a projection onto the ray $\rho$ (more precisely, the projection is obtained as $u \mapsto h^\rho(u)\cdot v_\rho$), and putting together the functions $h^\rho$ for all $\rho \in \sigma(1)$ gives a projection onto the cone $\sigma$ which agrees with the projection onto the second factor in the identification $\st(\sigma) \cong (N/\sigma)_\R \times \relint(\sigma)$ from Section \ref{section:std-covers-skeleton}. 

The construction relies on a standard combination of smoothing by convolution over $\st(\rho)$ and employing some cut-off functions near $\partial \overline{\st(\rho)}$ to make the function smooth on $N_\R$. We begin by picking an auxiliary inner product on $N_\R$ and let $\chi \colon N_\R \rightarrow \R_{\geq0}$ be a radially symmetric bump function that is supported on the unit ball and integrates to $1$; denote $\chi_\varepsilon(x)\coloneqq\varepsilon^{-n}\chi(x\cdot \varepsilon^{-1})$ for $\varepsilon>0$. We then define a smoothing of $h^\rho$ as the convolution $h_\varepsilon^\rho\coloneqq h^\rho *\chi_\varepsilon$. We can then observe that for sufficiently small $\varepsilon$, the level set $S^\rho_\varepsilon=\{u \in N_\R \colon h^\rho_\varepsilon(u)=1\}$ will be a subset of $\st(\rho)$ that is $O(\varepsilon)$-close to the corresponding level set of $h^\rho$. We can also check that it will be transverse to the radial vector field via the following calculation:

\begin{lemma}\label{lemma:convolution-est-1}
    There exists a constant $C>0$ such that for all $\rho \in \Sigma(1)$, $\varepsilon>0$ and $u \in S^\rho_\varepsilon$, we have $|D_uh^\rho_\varepsilon(u)-1| \leq C\varepsilon$. 
\end{lemma}

\begin{proof}
    By differentiating under the integral sign (we can do that since $h^\rho$ is $C^1$ on a set whose complement is null), we get:
    \begin{equation*}
        D_uh^\rho_\varepsilon(u)=\int_{N_\R} \chi_\varepsilon(v) D_{u+v}h^\rho(u)dv=\sum_{\sigma \in \Sigma(n)} \int_{\inte(\sigma)-u}\chi_\varepsilon(v)D_{u+v}h^\rho(u),
    \end{equation*}
    where the integrals in the second sum range over the points $v$ satisfying $u+v \in \inte(\sigma)$ for full-dimensional cones. But since the restriction of $h^\rho$ to these cones is given by a linear functional $\eta^\rho_\sigma$, we can explicitly compute that $D_{u+v}h^\rho(u)=\eta^\rho_\sigma(u)$ holds for the term indexed by $\sigma$. If we take $C$ to be equal to the maximum of operator norms of all the functionals $\eta^\sigma_\rho$, we then get a bound $|D_{u+v}h^\rho(u)-h^\rho(u+v)|=|D_{u+v}h^\rho(u)-\eta^\rho_\sigma(u+v)|=|\eta^\rho_\sigma(v)| \leq C \varepsilon$ that holds for all $u$ and for all $v \in B_\varepsilon(0)$, which gives us the required bound as follows:
    \begin{equation*}
        |D_uh^\rho_\varepsilon(u)-1|=|D_uh^\rho_\varepsilon(u)-h^\rho_\varepsilon(u)| \leq \int_{B_{\varepsilon}(0)}\chi_\varepsilon(v) \lvert D_{u+v}h^\rho(u+v)- h^\rho(u+v)\rvert dv \leq C\varepsilon.
    \end{equation*}
\end{proof}

 Therefore, each ray contained in $\st(\rho)$ intersects $S^\rho_\varepsilon$ transversely at a single point for all $\varepsilon>0$ sufficiently small, so we can define a smooth function $\widetilde{h}_\varepsilon^\rho \colon \st(\rho) \rightarrow \R_{> 0}$ homogeneous of degree $1$ by declaring it to be equal to $1$ along $S^\rho_\varepsilon$. The lemma can also be used to compare the differential of the homogenisation with that of the initial smoothing: 

\begin{corollary}\label{corollary:convolution-est-1}
    There exists a constant $C'>0$ such that for all $\rho \in \Sigma(1)$, sufficiently small $\varepsilon>0$ and $u \in S^\rho_\varepsilon$, we have $\lVert D_u\widetilde{h}^\rho_\varepsilon-D_uh^\rho_\varepsilon \rVert \leq C'\varepsilon$. 
\end{corollary}

\begin{proof}
    By definition of the homogenisation, if $v \in T_uS^\rho_\varepsilon$, then $D_u\widetilde{h}^\rho_\varepsilon(v)=D_uh^\rho_\varepsilon(v)=0$. Lemma \ref{lemma:convolution-est-1} tells us that we have a splitting $N_\R \cong T_uS^\rho_\varepsilon \oplus\R u$ and also that $\lvert 1-D_uh^\rho_\varepsilon(u)\rvert=\lvert D_u\widetilde{h}^\rho_\varepsilon(u)-D_uh^\rho_\varepsilon(u) \rvert \leq C\varepsilon$ holds for some constant. If we denote the orthogonal projection of $u$ onto $T_uS^\rho_\varepsilon=\ker(D_uh^\rho_\varepsilon)$ as $u'$, then the same estimate tells us that $\lVert D_uh^\rho_\varepsilon \rVert\cdot\lVert u' \rVert \geq 1-C\varepsilon$ and using the fact that $h^\rho$ is piecewise linear and hence Lipschitz, we get a uniform upper bound on $\lVert D_uh^\rho_\varepsilon \rVert$, which yields a uniform lower bound $\lVert u' \rVert \geq c >0$ for some constant for all sufficiently small $\varepsilon$. Therefore, since the splitting $T_uS^\rho_\varepsilon\oplus \R u'$ is orthogonal and $\lvert D_u\widetilde{h}^\rho_\varepsilon(u')-D_uh^\rho_\varepsilon(u') \rvert \leq C\varepsilon$, we can conclude that the operator norm of $D_u\widetilde{h}^\rho_\varepsilon-D_uh^\rho_\varepsilon$ is at most $\frac{C}{c}\varepsilon$, as desired. 
\end{proof}

Now, pick a non-zero cone $\sigma \in \Sigma$ with ray generators $\{\rho_1,\dots,\rho_r\}$, denote the associated ray generators $v_j\coloneqq v_{\rho_j}$ and the functions as $\widetilde{h}_j\coloneqq\widetilde{h}^{\rho_j}_\varepsilon$. They assemble into a map $\widetilde{h}_\sigma \colon \relint(\sigma) \rightarrow\R^r_{>0}$, which can be extended to a smooth function $\sigma \rightarrow \R_{\geq 0}^r$. It is also homogeneous of degree $1$, so it descends to a map onto the standard $(k-1)$-simplex $\widetilde{h}^\infty_\sigma \colon \partial_\infty \sigma \rightarrow \Delta^{k-1}$.

\begin{lemma}\label{lemma:convolutions-nice}
    The map $\widetilde{h}_\sigma \colon \sigma \rightarrow \R^r_{\geq 0}$ is a diffeomorphism of manifolds with corners for all non-zero $\sigma \in \Sigma$ and for all $\varepsilon>0$ small enough. 
\end{lemma}

\begin{proof}
    We note that it suffices to show that the map is a submersion over $\relint(\sigma)$ for each $\sigma$: it then follows that the extension is also a submersion by the same argument over $\relint(\tau)$ for $\tau \subset \partial \sigma$, so the extension is also a submersion, and since the domain and the codomain are of the same dimension, it follows that it is a local diffeomorphism. It is therefore a covering map, so the size of $\widetilde{h}_\sigma^{-1}(p)$ is locally constant, by checking that $\widetilde{h}_{\sigma}^{-1}(1,1,\dots,1)=\{v_1+\dots+v_r\}$ and using connectedness, we see that the map is injective. The map $\widetilde{h}^\infty_\sigma \colon \partial_\infty \sigma \rightarrow \Delta^{k-1}$ is also a local diffeomorphism, hence an open map, so its image is both open and closed, so it must be surjective by connectedness of $\Delta^{k-1}$, which also gives us surjectivity of $\widetilde{h}_\sigma$, thus showing that the map is indeed a diffeomorphism. 

    Let $u \in \relint(\sigma)$ and consider the $(k \times k)$-matrix $H$ with entries $H_{ij}=D_u\widetilde{h}_i(v_j)$, we will prove that it is invertible for all sufficiently small $\varepsilon>0$. Without loss of generality, we can permute the indices and assume that $0<\widetilde{h}_1(u) \leq \widetilde{h}_2(u) \leq \dots \leq \widetilde{h}_r(u)$. Denote the original piecewise linear functions as $h_i$ and the smoothings by convolution as $h_{i,\varepsilon}$.

    Note that for $x \in S^{\rho_i}_{\varepsilon}$, differentiating under the integral sign tells us that $D_xh_{i,\varepsilon}$ is a weighted average of values of $h_i(v)$ over a ball $B_\varepsilon(x)$. Thanks to the properties of the homogenisation and Corollary \ref{corollary:convolution-est-1}, it follows that there exists a constant $C_2>0$ such that for any $x \in \st(\rho_i)$ and any $\varepsilon$ small enough, we can get an estimate for $D_x\widetilde{h}_i$ in terms of an average:
    \begin{equation}\label{equation:conv-formula}
        \left\lVert D_x\widetilde{h}_i-\int_{B_{\varepsilon_i(x)}(x)}\chi_{\varepsilon_i(x)}(x-t)D_{t}h_idt \right\rVert \leq C_2 \varepsilon,
    \end{equation}
    where $\varepsilon_i(x)=\varepsilon\cdot\widetilde{h}_i(x)$ (note that we are using homogeneity of $h_i$ itself to scale and translate the integral from $B_{\varepsilon}(x\cdot \widetilde{h}_i(x)^{-1})$ to $B_{\varepsilon_i(x)}(x)$). 

    Denote the facet of $\sigma$ not containing $\rho_j$ as $\tau_j$. The functions $\widetilde{h}_i$ are all Lipschitz as remarked in the proof of Corollary \ref{corollary:convolution-est-1}, so pick a common Lipschitz constant $C_3$ that works for all of them. Therefore, considering the orthogonal projection of $u$ onto the facet $\tau_j$ yields an inequality $\widetilde{h}_j(u) \leq C_3d(u,\tau_j)$. In particular, this tells us that for all $0<\varepsilon<C_3^{-1}$ and $i \leq j$, one has $d(u,\tau_j) \geq C_3^{-1}\widetilde{h}_j(u) \geq C_3^{-1}\widetilde{h}_i(u)>\varepsilon_i(u)$, so the ball $B_{\varepsilon_i(u)}(u)$ does not intersect $\tau_j$. 

    This means that for sufficiently small $\varepsilon$, the ball $B_{\varepsilon_i(u)}(u)$ must be contained inside $\st(\cone(v_i,v_{i+1},\dots,v_r))$. By our definition of $h_i$ as a piecewise linear function, it is immediate that $D_th_i(v_j)=\delta_{ij}$ for all $t \in \st(\cone(v_i,\dots,v_r))$ and all $i\geq j$, so combining this with the formula \ref{equation:conv-formula}, we get a bound $|H_{ij}-\delta_{ij}| \leq C_4 \varepsilon$ for some constant $C_4>0$, all $\varepsilon \in (0,C_3^{-1})$ and $i\geq j$. In particular, we have shown that the matrix $H$ has distance $O(\varepsilon)$ from an upper triangular matrix with ones along the diagonal, from which the desired conclusion follows.
\end{proof}

One consequence of this Lemma is that $\{u \in N_\R \colon \widetilde{h}^\rho(u) \leq 1 \textnormal{ for all } \rho\in\Sigma(1)\}$ will be a manifold with corners that defines a \emph{partial smoothing} of the barycentric subdivision $\Sigma^\ba$. We can upgrade these to genuine smoothing as follows: for all $N>0$, denote the family smoothings of $\max\colon\R^N \rightarrow \R$ obtained via an appropriate convolution with a bump function $\chi_\varepsilon$ supported on a radius $\varepsilon$ ball as $\max_\varepsilon$. Let $\widetilde{h}_\varepsilon(u)\coloneqq\max_\varepsilon(\widetilde{h}^\rho_\varepsilon(u)\colon \rho \in \Sigma(1))$ and, in the setting of complete intersections defined via a nef partition, let $\widetilde{h}_{j,\varepsilon}(u)\coloneqq\max_\varepsilon(\widetilde{h}^\rho_\varepsilon(u)\colon \rho \in \Sigma_j(1))$, where $\Sigma_j$ is the fan on the triangulation $\mT_j$ of $\nabla_j$. Following the same strategy as in Proposition \ref{proposition:smoothing-bbci}, we can prove that these give us explicit defining functions for smoothings relative to $\Sigma$:

\begin{corollary}\label{corollary:trop-smoothings}
    For $\varepsilon>0$ small enough, the sets $\partial\widetilde{\Sigma}^\ba\coloneqq\{ u \in N_\R \colon \widetilde{h}_\varepsilon(u)=1\}$ and $\widetilde{\Sigma}^\ba_\trans\coloneqq\{ u \in N_\R \colon \widetilde{h}_{j,\varepsilon}(u)=1 \textnormal{ for all } j=1,\dots,r\}$ define smoothings of $\partial\Sigma^\ba$ and $\Sigma^\ba_\trans$ relative to $\Sigma$.
\end{corollary}
With this set-up in mind, we are also now ready to prove the first important result of this section, the existence of \emph{coherent projections} (Definition \ref{definition:coherent-projections}):

\begin{corollary}\label{corollary:coherent-projections}
    There exists a collection of coherent projections $\{ \pi_\sigma\}_{\sigma \in \Sigma\backslash 0}$ homogeneous of degree $1$.
\end{corollary}

\begin{proof}
    Pick a total order on rays of $\Sigma$, then for each $\sigma \in \Sigma\backslash\{0\}$, we can use the maps $\widetilde{h}^\rho$ for $\rho \in \sigma(1)$ sorted according to the ordering to get maps $\widetilde{h}_\sigma \colon \st(\sigma) \rightarrow \R^r_{>0}$ like in Lemma \ref{lemma:convolutions-nice}. The projection $\pi_\sigma$ can then be defined by composing this map with the inverse to the diffeomorphism $\relint(\sigma) \xrightarrow{\sim}\R^r_{>0}$. We claim that these are going to work for all $\varepsilon>0$ small enough. Note that the second property is satisfied by construction.

    As for the first one, if we denote $\sigma=\cone(v_1,\dots,v_r)$ as before, then the calculation for $H_{ij}$ from Lemma \ref{lemma:convolutions-nice} goes through for any $u\in\st(\sigma)$. In particular, it implies not just that the fibres of $\pi_\sigma$ (which are the same as fibres of $\widetilde{h}_\sigma$) are smooth submanifolds of $N_\R$, but also that they are transverse to $\R\sigma$. Therefore, the map $q_\sigma$ is a local diffeomorphism when restricted to a fibre $F_b\coloneqq\pi_\sigma^{-1}(b)$ for $b \in \relint(\sigma)$, so it suffices to prove that $q_\sigma|_{F_b}$ is a bijection. 

    Since the map is a local diffeomorphism, its image is open. We will show that it is also closed: suppose that we have a sequence of points $(x_m) \in F_b$ such that $q_\sigma(x_m)\rightarrow y$. This means that for $m$ large enough, we will have $\lVert q_\sigma(x)-y\rVert \leq 1$. By construction, $\widetilde{h}^\rho$ is also close to the piecewise linear function it is approximating, in particular, we have $|\widetilde{h}^\rho(x)-h^\rho(x)| \leq C \cdot \varepsilon \cdot \widetilde{h}^\rho(x)$ for all $\varepsilon>0$ small enough, $x\in \st(\sigma)$, $\rho \in \sigma(1)$ and some constant $C>0$. Therefore, $ c_\rho (1-C\varepsilon)\leq h^\rho(x_m) \leq c_\rho(1+C\varepsilon)$ for $\rho \in \sigma(1)$, meaning that for $\varepsilon>0$ small enough, the projections of $(x_m)$ onto $\sigma$ via $h_\sigma$ are confined to a compact neighbourhood of $b$. Therefore, all the terms $(x_m)$ must be eventually contained in a compact neighbourhood of a point $b \times y\in \sigma \times (N/\sigma)_\R \cong \st(\sigma)$ (we are using the piecewise smooth homeomorphism from Section \ref{section:combinatorics}), so it has a convergent subsequence that tends to some point $x \in \st(\sigma)$. By continuity, we get that $x \in F_b$ and also $q_\sigma(x)=y$, so the image is indeed closed. Therefore, by connectedness of the codomain $(N/\sigma)_\R$, the map $q_\sigma$ must be surjective. 

    Finally, observe that our argument also shows that every connected component also maps surjectively onto the base. By construction, we know that $F_b \cap \sigma=\{b\}$ so the only option is that $F_b$ is connected and $q_\sigma|_{F_b}$ is injective (since the number of preimages under a local diffeomorphism is locally constant). 
\end{proof}

We can also use a similar kind of smoothings to smooth certain scaling actions on $N_\R$: start from the observation that one can define an $\R^{\Sigma(1)}$ action on $N_\R$ by saying that the factor corresponding to a ray $\rho$ scales the $v_\rho$-component of all points $u \in \st(\rho)$ and does not affect the points $u \notin \st(\rho)$. To rigorously describe it, we shall present it in terms of flows. For $\rho \in \Sigma(1)$, let $V_\rho(x)\coloneqq h^\rho(x)\cdot v_\rho$, where the function $h^\rho$ is continuously extended by zero away from $\st(\rho)$. Denote the time $t$ flow of $V_\rho$ as $\Phi_\rho^t$, which will exist for all $t \in \R$ (this can be checked manually). Since $V_{\rho}(x) \in T_x\sigma$ whenever $x \in \relint(\sigma)$ and the restriction of $V_{\rho}|_\sigma$ to $\sigma$ is smooth, $\Phi^t_\rho$ defines a $\Sigma$-isotopy. Moreover, it is straightforward to check that $[V_\rho,V_{\rho'}]=0$ for all $\rho$, $\rho' \in \Sigma(1)$ (note that this expression makes sense despite the vector fields not being smooth, since can compute the Lie bracket separately along $\relint(\sigma)$ for all $\sigma \in \Sigma$). Therefore, the flows of the vector fields commute, so we have $\Phi^t_\rho \circ \Phi^{t'}_{\rho'}=\Phi^{t'}_{\rho'} \circ \Phi^{t}_{\rho}$ and the following definition makes sense: 

\begin{definition}\label{definition:total-scaling-action}
    The \emph{total scaling action} $A_\tot$ on $N_\R$ is defined as the action of $\R^{\Sigma(1)}$, where the function $t \colon \Sigma(1) \rightarrow \R$ acts by the composition of $\Phi^{t(\rho)}_\rho$ for all $\rho \in \Sigma(1)$.
\end{definition}

\begin{example}
    If we include $\R \rightarrow \R^{\Sigma(1)}$ diagonally, we get an action of $\R$ on $N_\R$, which is nothing else but the standard scaling action with $t \in \R$ acting by dilation with factor $e^t$.  
\end{example}

\begin{example}
    When our fan $\Sigma$ comes from a nef partition $\nabla=\nabla_1+\dots+\nabla_r$ with regular triangulations $\mT_j$, the rays in $\Sigma(1)$ are decomposed as $\Sigma(1)=\coprod_{j=1}^r \Sigma_j(1)$, hence we have $\R^{\Sigma(1)} \cong \R^{\Sigma_1(1)} \times \dots \R^{\Sigma_r(1)}$. Taking the diagonal inclusion of $\R$ into each factor gives us an injective group homomorphism $\R^r \rightarrow \R^{\Sigma(1)}$. The resulting $\R^r$-action $A_{nef}$ agrees with the one defined explicitly in Section \ref{section:combi-ci}, which we call the \emph{nef partition scaling action}. 
\end{example}

By construction, the restriction of the total scaling action onto each cone is smooth, but it might not be smooth as an action on $N_\R$. We can replace the functions $h^\rho$ by their appropriate smoothings to get rid of this issue: consider the non-homogeneous convolution smoothings $h^\rho_{\delta}$ defined above and let $r_\varepsilon(x)$ be a smoothing of $\max\{x-\varepsilon,0\}$ that is non-decreasing, has the same support, agrees with the function for $x>2\varepsilon$ and is $\varepsilon^2$-close to it over $(\varepsilon,2\varepsilon)$ (it is the same kind of smoothing as the one occurring in Corollary \ref{corollary:smoothing-isotopies2}). We can then note that $\widetilde{V}_{\rho}(x)\coloneqq r_\varepsilon(\widetilde{h}^\rho_\delta(x))\cdot v_\rho$ is a smooth vector field defined on $N_\R$. We record the properties of these smoothings in a Lemma:  

\begin{lemma}\label{lemma:vector-field-smoothings}
    Let $\rho, \rho' \in \Sigma(1)$ and suppose that $\varepsilon>0$ is a constant. Then for all sufficiently small $\delta>0$, we have the following:
    \begin{enumerate}
        \item The associated vector field $\widetilde{V}_\rho$ is supported in $\st(\rho)$.
        \item There exists a constant $C>0$ such that the inequality $\lVert \widetilde{V}_{\rho}(x) - V_\rho(x) \rVert \leq C\varepsilon$ holds.
        \item The flow $\widetilde{\Phi}^t_\rho$ of $\widetilde{V}_\rho$ is defined for all $t \in \R$ and gives us a smooth ambient isotopy of $N_\R$ preserving $\Sigma$.
        \item We have $[\widetilde{V}_{\rho},\widetilde{V}_{\rho'}]=0$, so the flows $\widetilde{\Phi}_{\rho}$ and $\widetilde{\Phi}_{\rho'}$ commute.
    \end{enumerate}
\end{lemma}

\begin{proof}
    Note that the support of $\widetilde{V}_\rho$ is, by construction, the set $\{x \in N_\R \colon h^\rho_\delta(x)>\varepsilon\}$. In particular, if $V_\rho(x) \neq 0$, then there must exist some $y \in B_\delta(x)$ such that $h^\rho(y)>\varepsilon$. As was noted above, the piecewise linear functions are Lipschitz with some common constant $L>0$, so we necessarily have $h^\rho(x) \geq h^\rho(y)-L\delta>\varepsilon-L\delta$. Therefore, for $\delta \leq L^{-1}\varepsilon$, the first conclusion holds. We shall prove that the remaining statements also hold under the same assumption.

    For the second property, observe that we have $|h^\rho_\delta(x)-h^\rho(x)| \leq L \delta$ for all $x\in N_\R$ and $|r_\varepsilon(s)-s|<\varepsilon$ for all $s \in \R$, which yields $\lVert \widetilde{V}_\rho(x)-V_\rho(x) \rVert = |r_\varepsilon(h^\rho_\delta(x))-h^\rho(x)|\cdot\lVert v_\rho \rVert \leq (L\delta +\varepsilon) \cdot \lVert v_\rho \rVert$. Therefore, if we take a constant $C>0$ to be larger than $2\max_\rho \lVert v_\rho \rVert$, the desired inequality will hold for the specified range of $\delta$. 
    
    The flowlines of $\widetilde{V}_{\rho}$ are rays parallel to $v_\rho$ or points and the bound $\lVert\widetilde{V}_{\rho}(x)\rVert <C'\lVert x \rVert$ for some constant $C'>0$ that follows from (2) tells us that it can not escape to infinity in finite time, which means that the flow is globally defined. Moreover, by property (1), $\widetilde{V}_{\rho}(x) \neq 0$ implies $x\in \st(\rho)$, so the vector filed is going to be tangent to $T_x\sigma$ whenever $x \in \relint(\sigma)$, which means that the flow preserves $\Sigma$. 
    
    Now, suppose that $\rho$ and $\rho'$ are distinct rays and suppose that $x \in \relint(\sigma)$; as with the non-smooth fields earlier, it suffices to evaluate the Lie bracket along $\sigma$. If at least one of the rays is not contained in $\sigma(1)$, then the corresponding vector field vanishes on a neighbourhood of $x$, so the Lie bracket vanishes at $x$. Otherwise, let $\sigma=\cone(v_1,\dots,v_r)$ with $v_1=v_{\rho}$, $v_2=v_{\rho'}$; then the dual functionals $x_j\coloneqq\eta^{\rho_j}_\sigma(x)$ give us coordinates on $\sigma$ with $v_{\rho_j}=\partial_j$, so that the coordinate expression for Lie derivatives yields 
    \begin{equation*}
[\widetilde{V}_\rho,\widetilde{V}_{\rho'}]=r_\varepsilon(h_\delta^\rho(x)) D_x(r_\varepsilon\circ h_\delta^{\rho'})(v_{\rho})\cdot v_{\rho'}-r_\varepsilon(h_\delta^{\rho'}(x)) D_x(r_\varepsilon\circ h_\delta^{\rho})(v_{\rho'})\cdot v_{\rho},
    \end{equation*}
    therefore we just need to check that $D_x(r_\varepsilon\circ h_\delta^{\rho})(v_{\rho'})=0$ whenever both of the functions $r_\varepsilon\circ h_\delta^{\rho}$ and $r_\varepsilon\circ h_\delta^{\rho'}$ do not vanish identically near $x$. 

    By the chain rule, it suffices to show that $D_x h_\delta^\rho(v_{\rho'})=0$ whenever $h_\delta^\rho(x)>\varepsilon$, $h_\delta^{\rho'}(x)>\varepsilon$. By differentiating under the integral sign, we see that $D_x h_\delta^\rho(v_{\rho'})$ is a weighted average of $D_y h^\rho(v_{\rho'})$ over $y \in B_\delta(x)$. By the condition on $\delta$, we obtain $h^\rho(y)>0$ and $h^{\rho'}(y)>0$ for all such $y$, therefore the ball $B_\delta(x)$ is contained in $\st(\cone(\rho,\rho'))$. In particular, all the linear functionals involved in the definition of $h^\rho$ over this ball vanish along $v_{\rho'}$ by design, which implies that $D_y h^\rho(v_{\rho'})=0$ for all such $y$ and therefore $D_x h_\delta^\rho(v_{\rho'})=0$, as desired. 
\end{proof}

Therefore, we can proceed analogously to what we did in our definition of total scaling action:

\begin{definition}\label{definition:total-scaling-action-sm}
    The \emph{$(\varepsilon,\delta)$-smoothing of the total scaling action} on $N_\R$ is the action $\widetilde{A}_\tot$ by $\R^{\Sigma(1)}$ where $t \colon \Sigma(1) \rightarrow \R$ acts via the composition of $\widetilde{\Phi}^{t(\rho)}_\rho$ for all $\rho \in \Sigma(1)$. 
\end{definition}

In particular, this gives us $\widetilde{A}_{nef}$, a smoothing of the nef partition scaling action $\R^r \curvearrowright N_\R$, which is precisely what we are after in Section \ref{section:combi-ci}. 

Finally, we can adjust the discussion of coherent projections to a setting that is equivariant with respect to the total scaling action. For a choice of smoothing parameters satisfying the above results, consider the sets $\sigma^{fr}\coloneqq\{ x \in \relint(\sigma) \colon h^\rho_{\delta}(x)>\varepsilon \textnormal{ for all } \rho \in \sigma(1)\}$ and $\st^{fr}(\sigma)\coloneqq \{ x \in N_\R \colon h^\rho_{\delta}(x)>\varepsilon \textnormal{ for all } \rho \in \sigma(1)\}$. By the above discussion, $\st^{fr}(\sigma)$ is precisely the set of points in $N_\R$ where the factors $\R^{\sigma(1)} \subset \R^{\Sigma(1)}$ act freely. By repeating the calculation from Lemma \ref{lemma:convolutions-nice} for the non-homogenised setting, we can show that level sets of $h^\rho_\delta$ for $\rho \in \sigma(1)$ that are contained in $\st^{fr}(\sigma)$ will be transverse to the fibres of $q_\sigma$ (in fact, everything simplifies, since the conditions on $\varepsilon$ and $\delta$ a priori rule out the possibility of balls spilling over into the wrong cones), and therefore to the orbits of the total scaling action. Combining this with observations like the ones from the proof of Corollary \ref{corollary:coherent-projections}, we get the following statement:

\begin{lemma}\label{lemma:orbit-space-total-scaling}
    The level set $\{x \in \st^{fr}(\sigma) \colon h^\rho_\delta(x)=c_\rho \textnormal{ for all } \rho \in \sigma(1)\}$ is a slice of the $\R^{\sigma(1)}$-action on $\st^{fr}(\sigma)$ whenever $c_\rho>\varepsilon$ for all $\rho$. 
\end{lemma}

This can be re-phrased in a slightly different language of the existence of coherent projections satisfying certain properties. Since quotienting out by the effect of the $\R^{\sigma(1)}$-action is the same as quotienting by the vector space $\R\sigma$, the following definition makes sense:

\begin{definition}\label{definition:coherent-projections-equiv}
    We say that a collection of maps $\{\pi_\sigma\}_{\sigma \in \Sigma \backslash 0}$ with $\pi_\sigma \colon \st^{fr}(\sigma) \rightarrow \sigma^{fr}$ is an \emph{equivariant coherent collection of projections} if
    \begin{enumerate}
        \item The fibres of $\pi_\sigma$ project diffeomorphically onto $(N/\sigma)_{\R}$ via the quotient map $q_\sigma\colon\st^{fr}(\sigma)\rightarrow (N/\sigma)_\R$.
        \item For any pair of cones $\sigma \subset \tau$, we have $\pi_\sigma \circ \pi_\tau=\pi_\sigma|_{\st^{fr}(\tau)}$.
        \item The map $\pi_\sigma$ is equivariant with respect to the smoothed total $\R^{\sigma(1)}$-action. 
    \end{enumerate}
\end{definition}

\begin{corollary}\label{corollary:coherent-projections-equiv}
    For any smoothing $\widetilde{A}_{\tot}$ of the total scaling action constructed as above, there exists a collection of equivariant coherent projections.
\end{corollary}

\begin{proof}
    The projection $\pi_\sigma$ will send a point $p \in \st^{fr}(\sigma)$ to the unique intersection of the level set $\{x \in N_\R \colon h^\rho_\delta(x)=h^\rho_\delta(p) \textnormal{ for all } \rho \in \sigma(1)\}$ with $\sigma^{fr}$; this will be well-defined and smooth since the maps $h^\rho_\delta$ give a diffeomorphism from $\sigma^{fr}$ to $(\varepsilon,\infty)^r$. Lemma \ref{lemma:orbit-space-total-scaling} shows that these projections satisfy the required conditions. 

    Alternatively, this can be seen as looking at the map $h^\tot_\delta \colon N_\R \rightarrow \R^{\Sigma(1)}$ whose $\rho$-component is given by $h^\rho_\delta$. Then the projections $\pi_\sigma$ can be defined via a relationship $\pi^\Sigma_\sigma \circ h^\tot_\delta=h_\delta^\tot \circ \pi_\sigma$, where $\pi^\Sigma_\sigma \colon \R^{\Sigma(1)} \rightarrow \R^{\sigma(1)}$ just picks out the rays in $\sigma$ (this is well-defined, since $h^\tot_\delta$ is an embedding over $\st^{fr}(\sigma)$ for all $\sigma \neq 0$). The vector fields $\widetilde{V}_\rho$ are pushed forward to $r_\varepsilon(x_\rho) \cdot \partial_{x_\rho}$, so the projections $\pi^\Sigma_\sigma$ are equivariant with respect to the action of $\R^{\sigma(1)}$, which means that $\pi_\sigma$ will also be equivariant. 
\end{proof}

\newpage
\printbibliography
\end{document}